\documentclass[a4paper, 12pt]{article}
\usepackage[T2A]{fontenc}
\usepackage[utf8x]{inputenc}
\usepackage[hidelinks]{hyperref}
\usepackage[english,russian]{babel}
\usepackage{amsmath}
\usepackage{amsfonts}
\renewcommand{\Im}{\operatorname{Im}}
\newcommand{\Ker}{\operatorname{Ker}}
\newcommand{\codim}{\operatorname{codim}}
\newcommand{\Rank}{\operatorname{rank}}
\newcommand{\sgn}{\operatorname{sign}}
\newcommand{\inter}{\operatorname{int}}
\newcommand{\card}{\operatorname{card}}
\newcommand{\diam}{\operatorname{diam}}
\newcommand{\supp}{\operatorname{supp}}
\newcommand{\supvrai}{\operatornamewithlimits{sup\,vrai}}
\newcommand{\e}{\mathfrak e}
\renewcommand{\l}{\mathfrak l}
\renewcommand{\k}{\mathfrak k}
\newcommand{\m}{\mathfrak m}
\newcommand{\n}{\mathfrak n}
\newcommand{\M}{\mathfrak M}
\newcommand{\s}{\mathcal s}
\newcommand{\Iota}{\mathfrak J}
\newcommand{\N}{\mathbb N}
\newcommand{\Z}{\mathbb Z}
\newcommand{\R}{\mathbb R}
\newcommand{\Nu}{\mathcal N}
\newcommand{\D}{\mathcal D}
\newcommand{\J}{\mathcal J}
\newcommand{\I}{\mathcal I}
\newcommand{\mes}{\operatorname{mes}}
\allowdisplaybreaks

\begin{document}

\author{ С. Н. Кудрявцев }
\title{Продолжение функций из неизотропных пространств Никольского --
Бесова и приближение их производных}
\date{}
\maketitle
\begin{abstract}
В статье рассмотрены неизотропные пространства Никольского и Бесова с
нормами, в определении которых вместо модулей непрерывности известных порядков
производных функций по координатным направлениям используются
"$L_p$-усредненные" модули непрерывности функций соответствующих порядков по
тем же направлениям. Для таких пространств функций, заданных в областях
определенного типа, построены непрерывные линейные отображения их в
обычные неизотропные пространства Никольского и Бесова в $ \R^d, $
являющиеся операторами продолжения функций, что влечет совпадение тех и
других пространств в упомянутых областях. В работе также найдена слабая
асимптотика аппроксимационных характеристик, относящихся к задаче
восстановления производной по значениям функций в заданном числе
точек, задаче С.Б. Стечкина для оператора дифференцирования,
задаче описания асимптотики поперечников для неизотропных классов
Никольского и Бесова в этих областях.
\end{abstract}

Ключевые слова: неизотропные пространства Никольского -- Бесова,
продолжение функций, эквивалентные нормы, восстановление производной,
приближение оператора, поперечник
\bigskip

\centerline{Введение}
\bigskip

В работе рассмотрен ряд задач теории функциональных пространств и теории
приближения для неизотропных пространств Никольского и Бесова функций, заданных
в областях из определенных классов. Эти задачи объединяет общая техника вывода
метрических соотношений, используемых для получения верхних оценок величин,
являющихся предметом изучения в этих задачах. Остановимся подробнее на
содержании работы.

При $ d \in \N, \alpha \in \R_+^d, 1 \le p,\theta < \infty $ для области $ D
\subset \R^d $ вводятся в рассмотрение
пространства $ (B_{p,\theta}^\alpha)^\prime(D) ((H_p^\alpha)^\prime(D)) $
с нормами
$$
\| f \|_{(B_{p,\theta}^\alpha)^\prime(D)} = \max\biggl(\| f \|_{L_p(D)},
\max_{j =1,\ldots,d} \left(\int_0^\infty
t^{-1 -\theta \alpha_j} (\Omega_j^{\prime l_j}(f, t)_{L_p(D)})^{\theta}\,dt
\right)^{1/\theta}\biggr),
$$
$$
\| f \|_{(H_p^\alpha)^\prime(D)} = \max(\| f \|_{L_p(D)}, \max_{j =1, \ldots, d}
\sup_{t \in \R_+} t^{-\alpha_j} \Omega_j^{\prime l_j}(f,t)_{L_p(D)}),
$$
где
\begin{multline*}
\Omega_j^{\prime l_j}(f,t)_{L_p(D)} =
\biggl((2 t)^{-1} \int_{ |\xi| \le t}
\| \Delta_{\xi e_j}^{l_j} f\|_{L_p(D_{l_j \xi e_j})}^p d\xi\biggr)^{1 /p},
t \in \R_+, l_j = \\
\min \{m \in \N: \alpha_j < m \}, j =1,\ldots,d.
\end{multline*}
Для ограниченных областей $D $ так называемого $ \alpha$-типа (см. п. 1.3.)
построены непрерывные линейные отображения пространств $ (B_{p,\theta}^\alpha)^\prime(D)
((H_p^\alpha)^\prime(D)) $ в пространства Бесова $ B_{p,\theta}^\alpha(\R^d) $
(Никольского $ H_p^\alpha(\R^d)), $ являющиеся операторами продолжения
функций, что влечет совпадение соответствующих пространств (см. п. 2.1.).
Публикации, в которых рассматривается такая задача, автору не известны. Из
близких работ по этой тематике приведем [1] -- [5] (см. также имеющуюся
там литературу), в которых изучается вопрос о продолжении за пределы
области определения гладких функций из неизотропных пространств с
сохранением класса. Отметим, что средства построения операторов
продолжения функций и схемы доказательства непрерывности таких
операторов, применяемые ниже, отличаются от тех, что использовались в упомянутых 
работах.

В случае ограниченной области $ D \alpha$-типа в статье при соответствующих
условиях на $ \lambda \in \Z_+^d $
установлена слабая асимптотика, т.е. найден порядок величины
наилучшей точности восстановления в пространстве $ L_q(D) $
производной $ \D^\lambda f $ по значениям в $ n $ точках функций
$ f $ из классов Никольского $ (\mathcal H_p^\alpha)^\prime(D) $ и Бесова
$ (\mathcal B_{p,\theta}^\alpha)^\prime(D), $ кроме того,
получена слабая асимптотика поведения в зависимости от $ \rho $
наилучшей точности приближения в $ L_q(D) $ оператора $ \D^\lambda $
на классах $ (\mathcal H_p^\alpha)^\prime(D) $ и
$ (\mathcal B_{p,\theta}^\alpha)^\prime(D) $ операторами, действующими из
$ L_s(D) $ в $ L_q(D), $ норма которых не превосходит $ \rho, $
и, наконец, описана слабая асимптотика колмогоровского, гельфандовского,
линейного, александровского и энтропийного $n$-поперечников в
$ L_q(D) $ классов, состоящих из  производных $ \D^\lambda f $ функций
$ f $ из единичного шара пространств $ (H_p^\alpha)^\prime(D) $ и
$ (B_{p,\theta}^\alpha)^\prime(D). $ В этой части настоящая работа продолжает
исследования, проводившиеся автором в [6] для упомянутых задач в
отношении неизотропных классов Никольского и Бесова функций, заданных
на кубе $ I^d. $ Проблематика, к которой относятся рассматриваемые задачи,
занимает значительное место в теории приближений, и обзор результатов по ней
в рамках этой статьи невозможен.

Отметим, что общие подходы при решении упомянутых задач сохранились такими же,
как в [6]. Однако для их реализации потребовались другие
средства приближения и существенные изменения схем вывода вспомогательных
соотношений, используемых для получения верхних оценок изучаемых величин,
в сравнении с теми, что применялись в [6]. Нижние оценки
рассматриваемых ниже величин устанавливаются, опираясь на нижние
оценки соответствующих величин, полученные в [6]. Как и
следовало ожидать, слабая асимптотика изучаемых величин в задачах о
восстановлении производной, приближении оператора дифференцирования и нахождении
порядков поперечников совпадает со слабой асимптотикой соответствующих им
величин, установленной в [6].

Работа состоит из введения и пяти параграфов, в первом из которых
рассматриваются некоторые средства для решения объявленных задач, во втором --
задача о продолжении, в третьем -- задача восстановления, в четв\"eртом --
задача приближения оператора дифференцирования, в пятом -- поперечники.
\bigskip

\centerline{\S 1. Предварительные сведения и вспомогательные
утверждения}
\bigskip

1.1. В этом пункте вводятся обозначения, относящиеся к
функциональным пространствам и классам, рассматриваемым в настоящей работе,
а также приводятся некоторые факты, необходимые в дальнейшем.

Для $ d \in \N $ через $ \Z_+^d $ обозначим множество
$$
\Z_+^d = \{\lambda = (\lambda_1, \ldots, \lambda_d) \in \Z^d:
\lambda_j \ge0, j =1, \ldots, d\}.
$$
Обозначим также при $ d \in \N $ для $ l \in \Z_+^d $ через
$ \Z_+^d(l) $ множество
$$
\Z_+^d(l) = \{ \lambda \in \Z_+^d: \lambda_j \le l_j, j =1, \ldots,d\}.
$$

Nапомним, что при $ n \in \N $ и $ 1 \le p  \le \infty $ через $ l_p^n $
обозначается пространство $ \R^n $ с фиксированной в н\"eм нормой
$$
\|x\|_{l_p^n} = \begin{cases} (\sum_{j =1}^n |x_j|^p)^{1/p} \text{ при } p < \infty; \\
\max_{j =1}^n |x_j| \text{ при } p = \infty, \end{cases} x \in \R^n.
$$

Для $ d \in \N, l \in \Z_+^d $ через $ \mathcal P^{d,l} $ будем
обозначать пространство вещественных  полиномов, состоящее из всех
функций $ f: \R^d \mapsto \R $ вида
$$
f(x) = \sum_{\lambda \in \Z_+^d(l)} a_{\lambda} \cdot x^{\lambda},
x \in \R^d,
$$
где $ a_{\lambda} \in \R, x^{\lambda} = x_1^{\lambda_1} \ldots x_d^{\lambda_d},
\lambda \in \Z_+^d(l). $

При $ d \in \N, l \in \Z_+^d $ для области $ D \subset \R^d $ через
$ \mathcal P^{d,l}(D) $ обозначим пространство функций $ f, $ определ\"eнных
в $ D, $ для каждой из которых существует полином $ g \in
\mathcal P^{d,l} $ такой, что сужение $ g \mid_D = f.$

Для топологического пространства $ T $ через $ C(T) $ будем
обозначать пространство непрерывных вещественных функций на $ T. $
Для множества $ A \subset T $ через $ \overline A $ обозначается замыкание
множества $ A, $ а через $ \inter A $  -- его внутренность.

Для измеримого по Лебегу множества $ D \subset \R^d $ при $ 1 \le p \le \infty $
через $ L_p(D),$ как обычно, обозначается
пространство всех вещественных измеримых на $ D $ функций $f,$
для которых определена норма
$$
\|f\|_{L_p(D)} = \begin{cases} (\int_D |f(x)|^p dx)^{1/p}, 1 \le p < \infty;\\
\supvrai_{x \in D}|f(x)|, p = \infty. \end{cases}
$$

При $ d \in \N $ для $ \lambda \in \Z_+^d $ через $ \D^\lambda $
будем обозначать оператор дифференцирования $ \D^\lambda =
\frac{\D^{|\lambda|}} {\D x_1^{\lambda_1} \ldots \D x_d^{\lambda_d}}, $
где $ |\lambda| = \sum_{j =1}^d \lambda_j. $

При обозначении известных пространств дифференцируемых функций будем
ориентироваться на [7].
Для области $ D \subset \R^d $ при $ 1 \le p \le \infty, l \in \Z_+^d $
через $ W_p^l(D) $ обозначается пространство всех функций $ f \in
L_p(D), $ для которых для каждого $ j =1,\ldots,d $  обобщенная
производная $ \D_j^{l_j} f = \frac{\D^{l_j}f }{ \D x_j^{l_j} } $
принадлежит $ L_p(D), $ с нормой
$$
\|f\|_{W_p^l(D)} = \max(\|f\|_{L_p(D)}, \max_{j =1,\ldots,d}
\|\D_j^{l_j}f\|_{L_p(D)}).
$$

Для $ x,y \in \R^d $ положим $ xy = x \cdot y = (x_1 y_1, \ldots,x_d y_d), $
а для $ x \in \R^d $ и $ A \subset \R^d $ определим
$$
x A = x \cdot A = \{xy: y \in A\}.
$$

Для $ x \in \R^d: x_j \ne 0, $ при $ j=1,\ldots,d,$ положим
$ x^{-1} = (x_1^{-1},\ldots,x_d^{-1}). $

Для $ d \in \N, x,y \in \R^d $ будем писать $ x \le y (x < y), $
если для каждого $ j =1,\ldots,d $ выполняется неравенство $ x_j
\le y_j (x_j < y_j). $

При $ d \in \N $ для $ x \in \R^d $ положим
$$
x_+ = ((x_1)_+, \ldots, (x_d)_+),
$$
где $ t_+ = \frac{1} {2} (t +|t|), t \in \R. $

Обозначим через $ \R_+^d $ множество $ x \in \R^d, $ для которых
$ x_j >0 $ при $ j =1,\ldots,d,$ и для $ a \in \R_+^d, x \in \R^d $
положим $ a^x = a_1^{x_1} \ldots a_d^{x_d}.$

При $ d \in \N $ для $ x \in \R^d $ через $\s(x) $ обозначим
множество $ \s(x) = \{j =1,\ldots,d: x_j \ne 0\}. $

Будем также обозначать через $ \chi_A $ характеристическую функцию
множества $ A \subset \R^d. $

При $ d \in \N $ определим множества
$$
I^d = \{x \in \R^d: 0 < x_j < 1, j =1,\ldots,d\},
$$
$$
\overline I^d = \{x \in \R^d: 0 \le x_j \le 1, j =1,\ldots,d\},
$$
$$
B^d = \{x \in \R^d: -1 \le x_j \le 1, j =1,\ldots,d\}.
$$

Через $ \e $ будем обозначать вектор в $ \R^d, $ задаваемый
равенством $ \e = (1,\ldots,1). $

Теперь привед\"eм некоторые факты, относящиеся к полиномам, которыми
мы будем пользоваться ниже.

В [6] содержится

Лемма 1.1.1

Пусть $ d \in \N, l \in \Z_+^d, \lambda \in \Z_+^d, 1 \le p,q \le \infty,
\rho, \sigma \in \R_+^d. $ Тогда существует константа
$ c_1(d,l,\lambda,\rho, \sigma) >0 $ такая, что для любых
измеримых по Лебегу множеств $ D,Q \subset \R^d, $  для которых
можно найти $ \delta \in \R_+^d $ и $ x^0 \in \R^d $ такие, что
$ D \subset (x^0 +\rho \delta B^d) $ и $ (x^0 +\sigma \delta I^d)
\subset Q, $ для любого полинома $ f \in \mathcal P^{d,l} $
выполняется неравенство
    \begin{equation*} \tag{1.1.1}
\| \D^\lambda f\|_{L_q(D)} \le c_1 \delta^{-\lambda -p^{-1} \e +q^{-1} \e}
\|f\|_{L_p(Q)}.
   \end{equation*}

Далее, напомним, что для области $ D \subset \R^d $ и вектора $ h \in \R^d $
через $ D_h $ обозначается множество
$$
D_h = \{x \in D: x +th \in D \forall t \in \overline I\}.
$$

Для функции $ f, $ заданной в области $ D \subset \R^d, $ и
вектора $ h \in \R^d $ определим в $ D_h $ е\"e разность $ \Delta_h f $ с шагом
$ h, $ полагая
$$
(\Delta_h f)(x) = f(x+h) -f(x), x \in D_h,
$$
а для $ l \in \N: l \ge 2, $ в $ D_{lh} $ определим $l$-ую
разность $ \Delta_h^l f $ функции $ f $ с шагом $ h $ равенством
$$
(\Delta_h^l f)(x) = (\Delta_h (\Delta_h^{l-1} f))(x), x \in
D_{lh},
$$
положим также $ \Delta_h^0 f = f. $

Как известно, справедливо равенство
$$
(\Delta_h^l f)(\cdot) = \sum_{k=0}^l C_l^k (-1)^{l-k} f(\cdot +kh),
C_l^k = \frac{l!} {k! (l-k)!}.
$$

При $ d \in \N $ для $ j =1,\ldots,d$ через $ e_j $ будем
обозначать вектор $ e_j = (0,\ldots,0,1_j,0,\ldots,0).$

Как показано в [8], справедлива

Лемма 1.1.2

Пусть $ d \in \N, l \in \Z_+^d. $ Тогда для любых $ \delta \in \R_+^d $ и
$ x^0 \in \R^d $ для $ Q = x^0 +\delta I^d $
существует единственный линейный оператор
$ P_{\delta, x^0}^{d,l}: L_1(Q) \mapsto \mathcal P^{d,l}, $
обладающий следующими свойствами:

1) для $ f \in \mathcal P^{d,l} $ имеет место равенство
\begin{equation*} \tag{1.1.2}
P_{\delta, x^0}^{d,l}(f \mid_Q) = f,
  \end{equation*}

2)
\begin{equation*}
\Ker P_{\delta,x^0}^{d,l} = \biggl\{\,f \in L_1(Q):
\int \limits_{Q} f(x) g(x) \,dx =0\ \forall g \in \mathcal P^{d,l}\,\biggr\},
\end{equation*}

прич\"eм существуют константы $ c_2(d,l) >0 $ и $ c_3(d,l) >0 $
такие, что

   3) при $ 1 \le p \le \infty $ для $ f \in L_p(Q) $ справедливо неравенство
  \begin{equation*} \tag{1.1.3}
\|P_{\delta, x^0}^{d,l}f \|_{L_p(Q)} \le c_2 \|f\|_{L_p(Q)},
  \end{equation*}

4) при $ 1 \le p < \infty $ для $ f \in L_p(Q) $ выполняется неравенство
\begin{equation*} \tag{1.1.4}
  \|f -P_{\delta, x^0}^{d,l}f \|_{L_p(Q)} \le c_3 \sum_{j =1}^d
\delta_j^{-1/p} \biggl(\int_{\delta_j B^1} \int_{Q_{(l_j +1) \xi e_j}}
|\Delta_{\xi e_j}^{l_j +1} f(x)|^p dx d\xi\biggr)^{1/p}.
\end{equation*}

Определим пространства и классы функций, изучаемые в настоящей работе (ср. с
[7]). Но прежде введ\"eм некоторые обозначения.

Пусть $ D$ --- область в $ \R^d$ и $ 1 \le p < \infty.$
Для $ f \in L_p(D) $ при $ j =1,\ldots,d,\ l \in \Z_+ $ обозначим модуль
непрерывности в $ L_p(D) $ порядка $ l $ по $ j $-му координатному направлению через
$$
\Omega_j^l(f,t)_{L_p(D)} =
\supvrai_{h \in t B^1} \| \Delta_{h e_j}^l f\|_{L_p(D_{l h e_j})}, t \in \R_+,
$$

а также введем в рассмотрение "усредненный" модуль непрерывности
в $ L_p(D) $ порядка $ l $ по $ j $-му координатному направлению, полагая
\begin{multline*}
\Omega_j^{\prime l}(f,t)_{L_p(D)} =
\biggl((2 t)^{-1} \int_{ t B^1} \| \Delta_{\xi e_j}^l f\|_{L_p(D_{l \xi e_j})}^p d\xi\biggr)^{1 /p} = \\
\biggl((2 t)^{-1} \int_{ t B^1} \int_{D_{l \xi e_j}} |\Delta_{\xi e_j}^l f(x)|^p dx
d\xi\biggr)^{1 /p}, t \in \R_+.
\end{multline*}

Из приведенных определений видно, что
\begin{multline*} \tag{1.1.5}
\Omega_j^{\prime l} (f, t)_{L_p(D)} \le
\Omega_j^l (f, t)_{L_p(D)}, t \in \R_+, 
f \in L_p(D), 1 \le p < \infty, \\
l \in \Z_+, D \text{ -- произвольная область в } \R^d.
\end{multline*}

Пусть теперь $ d \in \N, \alpha \in \R_+^d, 1 \le p < \infty $ и $ D $ --
область в $ \R^d. $ Тогда зададим вектор $ l = l(\alpha) \in \N^d, $ полагая
$ l_j = \min \{m \in \N: \alpha_j < m \}, j =1,\ldots,d, $ а также выберем
вектор $ \l \in \Z_+^d $ такой, что $ \l < \alpha, $ и через
$ (H_p^{\alpha})^{\l}(D) $ обозначим пространство всех функций
$ f \in W_p^{\l}(D), $ для которых при $ j =1,\ldots,d $ конечна величина
$$
\sup_{t \in \R_+} t^{-(\alpha_j -\l_j)} \cdot \Omega_j^{l_j -\l_j}(\D_j^{\l_j} f, t)_{L_p(D)} < \infty,
$$
а через $ (\mathcal H_p^{\alpha})^{\l}(D) $ -- множество функций $ f \in
W_p^{\l}(D), $ для которых при $ j =1,\ldots,d $ выполняется неравенство
$$
\sup_{t \in \R_+} t^{-(\alpha_j -\l_j)} \Omega_j^{l_j -\l_j}(\D_j^{\l_j} f,
t)_{L_p(D)} \le 1.
$$

В пространстве $ (H_p^\alpha)^{\l}(D) $ вводится норма
$$
\| f \|_{(H_p^\alpha)^{\l}(D)} = \max(\| f \|_{W_p^{\l}(D)}, \max_{j =1, \ldots, d}
\sup_{t \in \R_+} t^{-(\alpha_j -\l_j)} \Omega_j^{l_j -\l_j}(\D_j^{\l_j} f, t)_{L_p(D)}).
$$

При тех же условиях на $ \alpha, p, D $ обозначим через
$ (H_p^\alpha)^\prime(D) $ пространство всех функций
$ f \in L_p(D), $ обладающих тем свойством, что для любого
$ j =1,\ldots,d $ соблюдается условие
$$
\sup_{t \in \R_+} t^{-\alpha_j} \cdot \Omega_j^{\prime l_j}(f, t)_{L_p(D)} < \infty,
$$
а через $ (\mathcal H_p^\alpha)^\prime(D) $ -- множество функций
$ f \in L_p(D), $ обладающих тем свойством, что для любого
$ j =1,\ldots,d $ соблюдается неравенство
$$
\sup_{t \in \R_+} t^{-\alpha_j} \cdot \Omega_j^{\prime l_j}(f, t)_{L_p(D)} \le 1,
\text{где} l = l(\alpha).
$$
В пространстве $ (H_p^\alpha)^\prime(D) $ задается норма
$$
\| f \|_{(H_p^\alpha)^\prime(D)} = \max(\| f \|_{L_p(D)}, \max_{j =1, \ldots, d}
\sup_{t \in \R_+} t^{-\alpha_j} \Omega_j^{\prime l_j}(f,t)_{L_p(D)}).
$$

Для области $ D \subset \R^d $ при $ \alpha \in \R_+^d,\ 1 \le p < \infty,
\theta \in \R: 1 \le \theta < \infty, $
полагая, как и выше, $ l = l(\alpha) \in \N^d, $ и выбирая $ \l \in \Z_+^d:
\l < \alpha, $ через $ (B_{p,\theta}^\alpha)^{\l}(D) $ обозначим пространство
всех функций $ f \in W_p^{\l}(D), $ которые для любого $ j =1,\ldots,d $
удовлетворяют условию
$$
\int_0^\infty t^{-1 -\theta (\alpha_j -\l_j)}
(\Omega_j^{l_j -\l_j}(\D_j^{\l_j}f, t)_{L_p(D)})^{\theta}\,dt < \infty,
$$
а через $ (\mathcal B_{p,\theta}^\alpha)^{\l}(D) $ обозначим множество
всех функций $ f \in W_p^{\l}(D), $ для которых при любом $ j =1,\ldots,d $
соблюдается неравенство
$$
\int_0^\infty t^{-1 -\theta (\alpha_j -\l_j)}
(\Omega_j^{l_j -\l_j}(\D_j^{\l_j}f, t)_{L_p(D)})^{\theta}\,dt \le 1.
$$

В пространстве $ (B_{p,\theta}^\alpha)^{\l}(D) $ фиксируется норма
$$
\| f \|_{(B_{p,\theta}^\alpha)^{\l}(D)} = \max\biggl(\| f \|_{W_p^{\l}(D)},
\max_{j =1,\ldots,d} \left(\int_0^\infty
t^{-1 -\theta (\alpha_j -\l_j)} (\Omega_j^{l_j -\l_j}(\D_j^{\l_j}f,
t)_{L_p(D)})^{\theta}\,dt \right)^{1/\theta}\biggr).
$$

При $ \theta = \infty $ положим $ (B_{p,\infty}^\alpha)^{\l}(D) =
(H_p^\alpha)^{\l}(D). $

При тех же условиях на параметры  обозначим через
$ (B_{p,\theta}^\alpha)^\prime(D) $ пространство всех функций
$ f \in L_p(D), $ которые при $ l = l(\alpha) $ для каждого $ j = 1,\ldots,d $
подчинены условию
$$
\int_0^\infty t^{-1 -\theta \alpha_j}
(\Omega_j^{\prime l_j}(f, t)_{L_p(D)})^{\theta}\,dt < \infty,
$$
а через $ (\mathcal B_{p,\theta}^\alpha)^\prime(D) $ -- множество функций
$ f \in L_p(D), $ которые при $ j = 1,\ldots,d $ удовлетворяют неравенству
$$
\int_0^\infty t^{-1 -\theta \alpha_j}
(\Omega_j^{\prime l_j}(f, t)_{L_p(D)})^{\theta}\,dt \le 1.
$$

В пространстве $ (B_{p,\theta}^\alpha)^\prime(D) $ определяется норма
$$
\| f \|_{(B_{p,\theta}^\alpha)^\prime(D)} = \max\biggl(\| f \|_{L_p(D)},
\max_{j =1,\ldots,d} \left(\int_0^\infty
t^{-1 -\theta \alpha_j} (\Omega_j^{\prime l_j}(f, t)_{L_p(D)})^{\theta}\,dt
\right)^{1/\theta}\biggr).
$$
При $ \theta = \infty $ положим $ (B_{p,\infty}^\alpha)^\prime(D) =
(H_p^\alpha)^\prime(D). $

В случае, когда вектор $ \l = \l(\alpha) \in \Z_+^d $ имеет компоненты
$ (\l(\alpha))_j = \max\{m \in \Z_+: m < \alpha_j\}, j =1,\ldots,d,$
пространство $ (B_{p,\theta}^\alpha)^{\l}(D) ((H_p^\alpha)^{\l}(D))$ обычно
обозначается $ B_{p,\theta}^\alpha(D) (H_p^\alpha(D)).$

Неоднократно будет использоваться взятая из [7]

Лемма 1.1.3

Пусть $ d \in \N, l,m \in \Z_+^d, 1 \le p < \infty $ и $ D $ ---
область в $ \R^d. $ Тогда для $ f \in W_p^l(D) $ при
$ j =1,\lдdots,d $ и $ \xi \in \R $ справедлива оценка
\begin{equation*} \tag{1.1.6}
\| \Delta_{\xi e_j}^{l_j +m_j} f\|_{L_p(D_{(l_j +m_j) \xi e_j})} \le
|\xi|^{l_j} \| \Delta_{\xi e_j}^{m_j} \D_j^{l_j}f\|_{L_p(D_{m_j \xi e_j})}.
\end{equation*}

В тех же обозначениях, что и выше, принимая во внимание то обстоятельство,
что для $ f \in (B_{p,\theta}^\alpha)^{\l}(D) $ при $ j =1,\ldots,d $ и $ t >0 $
справедлива оценка (см., например, [6])
\begin{multline*}
t^{-(\alpha_j -\l_j)} \Omega_j^{l_j -\l_j}(\D_j^{\l_j}f, t)_{L_p(D)} \\
\le 2^{1/\theta +\alpha_j -\l_j}
\left(\int_0^\infty \tau^{-1 -\theta (\alpha_j -\l_j)}
(\Omega_j^{l_j -\l_j}(\D_j^{\l_j}f, \tau)_{L_p(D)})^{\theta}\, d\tau \right)^{1/\theta},
\end{multline*}
заключаем, что
\begin{multline*}
(B_{p,\theta}^\alpha)^{\l}(D) \subset (H_p^\alpha)^{\l}(D), \\
D \text{ -- область в} \R^d, \alpha \in \R_+^d, 1 \le p < \infty,
1 \le \theta < \infty, \l \in \Z_+^d: \l < \alpha.
\end{multline*}

Учитывая, что для $ f \in (B_{p, \theta}^\alpha)^\prime(D), \alpha \in
\R_+^d, 1 \le p,\theta < \infty, D $ -- область в $ \R^d, $ при
$ t \in \R_+, l = l(\alpha), j =1,\ldots,d $ выполняется неравенство
\begin{equation*}
t^{-\alpha_j} \Omega_j^{\prime l_j}(f, t)_{L_p(D)} \le 
2^{\alpha_j +1/\theta +1 /p} \biggl(\int_{\R_+}
\tau^{-1 -\theta \alpha_j} (\Omega_j^{\prime l_j}(f,
\tau)_{L_p(D)})^\theta d\tau \biggr)^{1/\theta},
\end{equation*}
заключаем, что
\begin{multline*} \tag{1.1.7}
(B_{p, \theta}^\alpha)^\prime(D) \subset
(H_p^\alpha)^\prime(D); 
(\mathcal B_{p, \theta}^\alpha)^\prime(D) \subset
c_4(\alpha) (\mathcal H_p^\alpha)^\prime(D); \\
\| f\|_{(H_p^\alpha)^\prime(D)} \le c_4(\alpha)
\| f\|_{(B_{p, \theta}^\alpha)^\prime(D)},
\end{multline*}
где $ c_4(\alpha) = \max_{j =1,\ldots,d} 2^{2+\alpha_j}. $
Из (1.1.5) и (1.1.6) следует, что для $ f \in (B_{p, \theta}^\alpha)^{\l}(D),
\l \in \Z_+^d: \l < \alpha, l = l(\alpha) $ при $ t \in \R_+, j =1,\ldots,d $
выполняется неравенство

\begin{multline*}
\Omega_j^{\prime l_j}(f, t)_{L_p(D)} \le
\Omega_j^{l_j}(f, t)_{L_p(D)} = \supvrai_{ \xi \in t B^1}
\| \Delta_{\xi e_j}^{\l_j +(l_j -\l_j)} f\|_{L_p(D_{(\l_j +(l_j -\l_j)) \xi e_j})} \le \\
\supvrai_{ \xi \in t B^1} |\xi|^{\l_j}
\| \Delta_{\xi e_j}^{l_j -\l_j} \D_j^{\l_j}f\|_{L_p(D_{(l_j -\l_j) \xi e_j})} \le \\
t^{\l_j} \cdot \supvrai_{ \xi \in t B^1}
\| \Delta_{\xi e_j}^{l_j -\l_j} \D_j^{\l_j}f\|_{L_p(D_{(l_j -\l_j) \xi e_j})} = 
t^{\l_j} \Omega_j^{l_j -\l_j}(\D_j^{\l_j} f, t)_{L_p(D)},
\end{multline*}
и, значит,
\begin{equation*} \tag{1.1.8}
(B_{p, \theta}^\alpha)^{\l}(D) \subset (B_{p, \theta}^\alpha)^\prime(D), \\
(\mathcal B_{p, \theta}^\alpha)^{\l}(D) \subset (\mathcal B_{p, \theta}^\alpha)^\prime(D)
\end{equation*}
и
\begin{multline*} \tag{1.1.9}
\| f\|_{(B_{p, \theta}^\alpha)^\prime(D)} \le
\| f\|_{(B_{p, \theta}^\alpha)^{\l}(D)},
\alpha \in \R_+^d, 1 \le p < \infty, 1 \le \theta \le \infty, \\
D \text{ -- произвольная область в } \R^d.
\end{multline*}

Обозначим через $ C^\infty(D) $ пространство бесконечно дифференцируемых
функций в области $ D \subset \R^d, $ а через $ C_0^\infty(D) $ -- пространство
функций $ f \in C^\infty(\R^d), $ у которых носитель $ \supp f \subset D. $

В заключение этого пункта введ\"eм ещ\"e несколько обозначений.
Для банахова пространства $ X $ (над $ \R$) обозначим $ B(X) = \{x \in X:
\|x\|_X \le 1\}. $

Для банаховых пространств $ X,Y $ через $ \mathcal B(X,Y) $ обозначим банахово
пространство, состоящее из непрерывных линейных операторов $ T: X \mapsto Y, $
с нормой
$$
\|T\|_{\mathcal B(X,Y)} = \sup_{x \in B(X)} \|Tx\|_Y.
$$
Отметим, что если $ X=Y,$ то $ \mathcal B(X,Y) $ является банаховой алгеброй.
\bigskip

1.2. В этом пункте вводятся в рассмотрение пространства кусочно-
полиномиальных функций и операторов в них, которые используются для построения
средств приближения для функций из изучаемых нами пространств.
Но сначала привед\"eм некоторые вспомогательные сведения.

Введ\"eм в рассмотрение систему разбиений единицы, используемую для построения
кусочно-полиномиальных функций из интересующих нас пространств.
Для этого обозначим через $ \psi^{1,0} $ характеристическую функцию
интервала $ I, $ т.е. функцию, определяемую равенством
$$
\psi^{1,0}(x) = \begin{cases} 1, & \text{ для } x \in I; \\
0, & \text{ для } x \in \R \setminus I.
\end{cases}
$$
При $ m \in \N $ положим
$$
\psi^{1,m}(x) = \int_I \psi^{1, m-1}(x-y) dy (\text{см.} [9]),
$$
а для $ d \in \N, m \in \Z_+^d $ определим
$$
\psi^{d,m}(x) = \prod_{j=1}^d \psi^{1,m_j}(x_j), x =
(x_1,\ldots,x_d) \in \R^d.
$$

Для $ d \in \N, m,n \in \Z^d: m \le n, $ обозначим
$$
\Nu_{m,n}^d = \{ \nu \in \Z^d: m \le \nu \le n \} = \prod_{j=1}^d
\Nu_{m_j,n_j}^1.
$$

Опираясь на определения, используя индукцию, нетрудно проверить
следующие свойства функций $ \psi^{d,m}, d \in \N, m \in \Z_+^d. $

1) При $ d \in \N, m \in \Z_+^d $
$$
\sgn \psi^{d,m}(x) = \begin{cases} 1, \text{ для } x \in ((m +\e) I^d); \\
0, \text{ для } x \in \R^d \setminus ((m +\e) I^d),
\end{cases}
$$

2) при $ d \in \N, m \in \Z_+^d $ для каждого $ \lambda \in
\Z_+^d(m) $ (обобщ\"eнная) производная $ \D^\lambda \psi^{d,m} \in
L_\infty(\R^d), $

3) при $ d \in \N, m \in \Z_+^d $ почти для всех $ x \in \R^d $
справедливо равенство
$$
\sum_{\nu \in \Z^d} \psi^{d,m}(x -\nu) =1,
$$

4) при $ m \in \N $ для всех $ x \in \R $ (при $ m =0 $ почти для всех $ x \in \R $)
имеет место равенство
\begin{equation*} \tag{1.2.1}
\psi^{1,m}(x) = \sum_{\mu \in \Nu_{0, m+1}^1} a_{\mu}^m
\psi^{1,m}(2x -\mu),
\end{equation*}
где $ a_\mu^m = 2^{-m} C_{m+1}^\mu. $
Используя разложение Ньютона для $ (1+1)^{m+1} $  и $ (-1+1)^{m+1}, $ легко
проверить, что при $ m \in \Z_+ $ выполняются равенства
\begin{equation*} \tag{1.2.2}
\sum_{\mu \in \Nu_{0,m +1}^1 \cap (2 \Z)} a_\mu^m =1,
\sum_{\mu \in \Nu_{0,m +1}^1 \cap (2 \Z +1)} a_\mu^m =1.
\end{equation*}

При $ d \in \N $ для $ t \in \R^d $ через $ 2^t $ будем обозначать
вектор $ 2^t = (2^{t_1}, \ldots, 2^{t_d}). $

Для $ d \in \N, m,\kappa \in \Z_+^d, \nu \in \Z^d $ обозначим
$$
g_{\kappa, \nu}^{d,m}(x) = \psi^{d,m}(2^\kappa x -\nu) =
\prod_{j =1}^d \psi^{1,m_j}( 2^{\kappa_j} x_j -\nu_j), x \in \R^d.
$$
Из первого среди приведенных выше свойств функций $ \psi^{d,m} $
следует, что при $ d \in \N, m,\kappa \in \Z_+^d, \nu \in \Z^d $
носитель $ \supp g_{\kappa,\nu}^{d,m} =
2^{-\kappa} \nu +2^{-\kappa} (m +\e) \overline I^d. $
При $ d \in \N, \kappa \in \Z_+^d, \nu \in \Z^d $ обозначим
$ Q_{\kappa, \nu}^d = 2^{-\kappa} \nu +2^{-\kappa} I^d, $
$ \overline Q_{\kappa, \nu}^d = 2^{-\kappa} \nu +2^{-\kappa} \overline I^d. $

Отметим некоторые полезные для нас свойства носителей функций
$ g_{\kappa,\nu}^{d,m}. $

При $ d \in \N, m,\kappa \in \Z_+^d $ для каждого $ \nu^\prime \in \Z^d $
имеет место равенство
\begin{equation*} \tag{1.2.3}
\{ \nu \in \Z^d: Q_{\kappa, \nu^\prime}^d \cap
\supp g_{\kappa, \nu}^{d,m} \ne \emptyset\} = \nu^\prime +\Nu_{-m,0}^d.
\end{equation*}

Из свойства 3) функций $ \psi^{d,m} $ вытекает, что при $ d \in \N,
m, \kappa \in \Z_+^d $ для любой области $ D \subset \R^d $ почти для всех
$ x \in D $ соблюдается равенство
\begin{equation*} \tag{1.2.4}
\sum_{ \nu \in \Z^d: \supp g_{\kappa, \nu}^{d,m} \cap D \ne \emptyset}
g_{\kappa, \nu}^{d,m}(x) =1.
\end{equation*}

При $ d \in \N $ для $ x,y \in \R^d $ будем обозначать
$$
(x,y) = \sum_{j =1}^d x_j y_j.
$$

Имея в виду свойство 2) функций $ \psi^{d,m}, $ отметим, что при
$ d \in \N, m,\kappa \in \Z_+^d, \nu \in \Z^d, \lambda \in \Z_+^d(m)$
выполняется равенство
\begin{multline*} \tag{1.2.5}
\| \D^\lambda g_{\kappa, \nu}^{d,m} \|_{L_\infty (\R^d)} =
2^{(\kappa, \lambda)} \| \D^\lambda \psi^{d,m} \|_{L_\infty(\R^d)} =
c_1(d,m,\lambda) 2^{(\kappa, \lambda)}.
\end{multline*}

Введем в рассмотрение следующие пространства кусочно-полиномиальных
функций.
При $ d \in \N, l \in \Z_+^d, m \in \N^d, \kappa \in \Z_+^d $ и области
$ D \subset \R^d, $ полагая
\begin{equation*}
N_\kappa^{d,m,D} = \{\nu \in \Z^d: \supp g_{\kappa, \nu}^{d,m}
\cap D \ne \emptyset\},
\end{equation*}
через $ \mathcal P_\kappa^{d,l,m,D} $ обозначим линейное пространство,
состоящее из функций $ f: \R^d \mapsto \R, $ для каждой из которых существует
набор полиномов
$ \{f_\nu \in \mathcal P^{d,l}, \nu \in N_\kappa^{d,m,D}\} $ такой, что
для $ x \in \R^d $ выполняется равенство
\begin{equation*} \tag{1.2.6}
f(x) = \sum_{\nu \in N_\kappa^{d,m,D}} f_\nu(x) g_{\kappa,\nu}^{d,m}(x).
\end{equation*}

Нетрудно проверить, что при $ d \in \N, l \in \Z_+^d, m \in \N^d,
\kappa \in \Z_+^d $ и ограниченной области $ D \subset \R^d $ отображение,
которое каждому набору полиномов $ \{f_\nu \in \mathcal P^{d, l},
\nu \in N_\kappa^{d,m,D} \} $ ставит в соответствие функцию $ f, $ задаваемую
равенством (1.2.6), является изоморфизмом прямого произведения
$ \card N_\kappa^{d,m,D} $ экземпляров пространства $ \mathcal P^{d,l} $
на пространство $ \mathcal P_\kappa^{d,l,m,D}. $

Предложение 1.2.1

Пусть $ d \in \N, l \in \Z_+^d, m \in \N^d, \kappa \in \Z_+^d, D $ --
ограниченная область в $ \R^d. $ Тогда при $ j =1,\ldots,d $ линейный
оператор $ H_\kappa^{j,d,l,m,D}:
\mathcal P_\kappa^{d,l,m,D} \mapsto \mathcal P_{\kappa +e_j}^{d,l,m,D}, $
значение которого на функции $ f \in \mathcal P_\kappa^{d,l,m,D}, $ задаваемой
равенством (1.2.6), определяется соотношением
\begin{multline*} \tag{1.2.7}
(H_\kappa^{j,d,l,m,D} f)(x) = \\
\sum_{\nu \in N_{\kappa +e_j}^{d,m,D}}
\biggl(\sum_{\substack{\nu^\prime \in N_\kappa^{d,m,D}, \mu_j \in \Nu_{0, m_j +1}^1: \\
 2 \nu^\prime_j +\mu_j = \nu_j, \nu^\prime_i = \nu_i, i = 1,\ldots,d, i \ne j }} a_{\mu_j}^{m_j}
f_{\nu^\prime}(x)\biggr) g_{\kappa +e_j,\nu}^{d,m}(x), x \in \R^d,
\end{multline*}
обладает тем свойством, что для $ f \in \mathcal P_\kappa^{d,l,m,D} $ выполняется равенство
\begin{equation*} \tag{1.2.8}
(H_\kappa^{j,d,l,m,D} f) \mid_{D} = f \mid_{D}.
\end{equation*}

Доказательство.

Прежде всего, заметим, что в условиях предложения для функции
$ f \in \mathcal P_\kappa^{d,l,m,D}, $ задаваемой равенством (1.2.6), ввиду
(1.2.1) при $ x \in D $ имеет место равенство
\begin{multline*}
f(x) = \sum_{\nu^\prime \in N_\kappa^{d,m,D}} f_{\nu^\prime}(x)
\biggl(\prod_{i =1}^d \psi^{1,m_i}(2^{\kappa_i} x_i -\nu^\prime_i)\biggr) = \\
\sum_{\nu^\prime \in N_\kappa^{d,m,D}} f_{\nu^\prime}(x)
\biggl(\sum_{\mu_j \in \Nu_{0, m_j +1}^1}
a_{\mu_j}^{m_j} \psi^{1,m_j}(2^{\kappa_j +1} x_j -2 \nu^\prime_j -\mu_j)\biggr)
\prod_{i =1,\ldots,d: i \ne j} \ psi^{1,m_i}(2^{\kappa_i} x_i -\nu^\prime_i) = \\
\sum_{\nu^\prime \in N_\kappa^{d,m,D}} \sum_{\mu_j \in \Nu_{0, m_j +1}^1}
a_{\mu_j}^{m_j} f_{\nu^\prime}(x) \psi^{1,m_j}(2^{\kappa_j +1} x_j -2 \nu^\prime_j -\mu_j)
\biggl(\prod_{i =1,\ldots,d: i \ne j} \ psi^{1,m_i}(2^{\kappa_i} x_i -\nu^\prime_i)\biggr) = \\
\sum_{\nu^\prime \in N_\kappa^{d,m,D}, \mu_j \in \Nu_{0, m_j +1}^1}
a_{\mu_j}^{m_j} f_{\nu^\prime}(x) \psi^{1,m_j}(2^{\kappa_j +1} x_j -2 \nu^\prime_j -\mu_j)
\biggl(\prod_{i =1,\ldots,d: i \ne j} \ psi^{1,m_i}(2^{\kappa_i} x_i -\nu^\prime_i)\biggr) = \\
\sum_{\nu \in \Z^d} \biggl(\sum_{\substack{\nu^\prime \in N_\kappa^{d,m,D}, \mu_j \in
\Nu_{0, m_j +1}^1: \\ 2 \nu^\prime_j +\mu_j =  \nu_j, \\ \nu^\prime_i = \nu_i,
i = 1,\ldots,d, i \ne j }} a_{\mu_j}^{m_j} f_{\nu^\prime}(x)
\psi^{1,m_j}(2^{\kappa_j +1} x_j -2 \nu^\prime_j -\mu_j)
\biggl(\prod_{i =1,\ldots,d: i \ne j} \ psi^{1,m_i}(2^{\kappa_i} x_i -\nu^\prime_i)\biggr)\biggr) = \\
\sum_{\nu \in \Z^d} \biggl(\sum_{\substack{\nu^\prime \in N_\kappa^{d,m,D}, \mu_j \in
\Nu_{0, m_j +1}^1: \\ 2 \nu^\prime_j +\mu_j = \nu_j, \\ \nu^\prime_i = \nu_i,
i = 1,\ldots,d, i \ne j }} a_{\mu_j}^{m_j} f_{\nu^\prime}(x)
\psi^{1,m_j}(2^{\kappa_j +1} x_j -\nu_j)
\biggl(\prod_{i =1,\ldots,d: i \ne j} \ psi^{1,m_i}(2^{\kappa_i} x_i -\nu_i)\biggr)\biggr) = \\
\sum_{\nu \in \Z^d} \biggl(\sum_{\substack{\nu^\prime \in N_\kappa^{d,m,D}, \mu_j \in
\Nu_{0, m_j +1}^1: \\ 2 \nu^\prime_j +\mu_j = \nu_j, \nu^\prime_i = \nu_i,
i = 1,\ldots,d, i \ne j }} a_{\mu_j}^{m_j} f_{\nu^\prime}(x)\biggr)
g_{\kappa +e_j,\nu}^{d,m}(x) = \\
\sum_{\nu \in \Z^d: \supp g_{\kappa +e_j,\nu}^{d,m} \cap D \ne \emptyset}
\biggl(\sum_{\substack{\nu^\prime \in N_\kappa^{d,m,D}, \mu_j \in \Nu_{0, m_j +1}^1: \\
2 \nu^\prime_j +\mu_j = \nu_j, \nu^\prime_i = \nu_i, i = 1,\ldots,d, i \ne j }}
a_{\mu_j}^{m_j} f_{\nu^\prime}(x)\biggr) g_{\kappa +e_j,\nu}^{d,m}(x) = \\
\sum_{\nu \in N_{\kappa +e_j}^{d,m,D}}
\biggl(\sum_{\substack{\nu^\prime \in N_\kappa^{d,m,D}, \mu_j \in \Nu_{0, m_j +1}^1: \\
2 \nu^\prime_j +\mu_j = \nu_j, \nu^\prime_i = \nu_i, i = 1,\ldots,d, i \ne j }}
a_{\mu_j}^{m_j} f_{\nu^\prime}(x)\biggr) g_{\kappa +e_j,\nu}^{d,m}(x).
\end{multline*}

Сопоставляя последнее равенство с (1.2.7), приходим к (1.2.8). $ \square $

Замечание.

В условиях предложения 1.2.1 при $ j =1,\ldots,d, $ если для $ \nu \in
N_{\kappa +e_j}^{d,m,D}, \mu_j \in \Nu_{0, m_j +1}^1 $ существует
$ \nu^\prime \in \Z^d, $ для которого $ 2 \nu^\prime_j +\mu_j = \nu_j,
\nu^\prime_i = \nu_i, i = 1,\ldots,d: i \ne j, $ то для такого $ \nu^\prime $
справедливо включение $ \nu^\prime \in N_{\kappa}^{d,m,D}. $

В самом деле, при соблюдении условий замечания, выбирая $ x \in D \cap
\supp g_{\kappa +e_j,\nu}^{d,m}, $ получаем, что
\begin{multline*}
2^{-\kappa_i} \nu_i \le x_i \le 2^{-\kappa_i} \nu_i +2^{-\kappa_i}
(m_i +1); \text{ при } i = 1,\ldots,d: i \ne j, \\
2^{-(\kappa_j +1)} \nu_j \le x_j \le 2^{-(\kappa_j +1)} \nu_j +
2^{-(\kappa_j +1)} (m_j +1),
\end{multline*}
откуда
\begin{multline*}
2^{-\kappa_i} \nu^\prime_i \le x_i \le 2^{-\kappa_i} \nu^\prime_i +
2^{-\kappa_i} (m_i +1); \text{ при } i = 1,\ldots,d: i \ne j, \\
2^{-\kappa_j} \nu^\prime_j = 2^{-\kappa_j} (\nu_j -\mu_j) /2 =
2^{-(\kappa_j +1)} \nu_j -2^{-(\kappa_j +1)} \mu_j \le \\ 2^{-(\kappa_j +1)} \nu_j
\le x_j \le 2^{-(\kappa_j +1)} \nu_j + 2^{-(\kappa_j +1)} (m_j +1) = \\
 2^{-\kappa_j} (\nu_j -\mu_j) /2 +
2^{-(\kappa_j +1)} \mu_j +2^{-(\kappa_j +1)} (m_j +1) = \\
2^{-\kappa_j} \nu^\prime_j +2^{-\kappa_j} (\mu_j +m_j +1) /2 \le
2^{-\kappa_j} \nu^\prime_j +2^{-\kappa_j} (m_j +1),
\end{multline*}
т.е. $ x \in D \cap \supp g_{\kappa,\nu^\prime}^{d,m}, $ а, значит,
$ \nu^\prime \in N_{\kappa}^{d,m,D}. $

При $ m \in \N^d, \k \in \Z_+^d, \nu \in \Z^d $ обозначим через
$ \M_{\k}^m(\nu) $ множество наборов чисел
\begin{multline*}
\M_{\k}^m(\nu) = \{ \m^{\k} = \{ \m_{j,s_j} \in \Nu_{0, m_j +1}^1,
s_j =1,\ldots,\k_j, j \in \s(\k)\}: \\
\nu_j 2^{-\iota_j} -\sum_{s_j =1}^{\iota_j} \m_{j, \k_j -s_j +1}
2^{-(\iota_j -s_j +1)} \in \Z \forall \iota_j = 1,\ldots, \k_j, j \in \s(\k)\}
\end{multline*}
и каждой паре $ \nu \in \Z^d, \m^{\k} \in \M_{\k}^m(\nu) $ сопоставим
$ \n_{\k}(\nu,\m^{\k}) \in \Z^d, $ полагая
\begin{equation*}
(\n_{\k}(\nu,\m^{\k}))_j = \begin{cases} \nu_j 2^{-\k_j} -\sum_{s_j =1}^{\k_j}
\m_{j, \k_j -s_j +1} 2^{-(\k_j -s_j +1)}, j \in \s(\k); \\
\nu_j, j \in \{1,\ldots,d\} \setminus \s(\k). \end{cases}
\end{equation*}

Для формулировки следующего утверждения при $ d \in \N $ для $ j \in
\{1,\ldots,d\} $ обозначим через $ \eta^j: \R^d \times \R^d \mapsto \R^d $
отображение, определяемое соотношением
$$
(\eta^j(\xi,x))_i = \begin{cases} \xi_i, i =1, \ldots, j; \\
x_i, i =j +1, \ldots, d, \end{cases}  \xi,x \in \R^d.
$$

Предложение 1.2.2

Пусть $ d \in \N, l \in \Z_+^d, m \in \N^d, \kappa \in \Z_+^d, \kappa^\prime
\in \Z_+^d: \kappa^\prime \le \kappa, D $ -- ограниченная область в $ \R^d. $
Тогда линейный оператор $ H_{\kappa, \kappa^\prime}^{d,l,m,D}:
\mathcal P_{\kappa^\prime}^{d,l,m,D} \mapsto \mathcal P_\kappa^{d,l,m,D}, $
значение которого для $ f \in \mathcal P_{\kappa^\prime}^{d,l,m,D} $
определяется равенством
\begin{equation*} \tag{1.2.9}
H_{\kappa, \kappa^\prime}^{d,l,m,D} f = \begin{cases} f, \text{ при } \kappa^\prime = \kappa; \\
(\prod_{j \in \s(\k)} (\prod_{\iota_j =1}^{\k_j}
H_{\eta^j(\kappa^\prime, \kappa) +(\k_j -\iota_j) e_j}^{j,d,l,m,D})) f,
\text{ при } \kappa^\prime \ne \kappa,
\end{cases}
\text{ где } \k = \kappa -\kappa^\prime,
\end{equation*}

обладает следующими свойствами:

1) для $ f \in \mathcal P_{\kappa^\prime}^{d,l,m,D} $ выполняется равенство
\begin{equation*} \tag{1.2.10}
(H_{\kappa, \kappa^\prime}^{d,l,m,D} f) \mid_{D} = f \mid_{D};
\end{equation*}

2) для $ f \in \mathcal P_{\kappa^\prime}^{d,l,m,D} $ вида
\begin{equation*} \tag{1.2.11}
f = \sum_{\nu^\prime \in N_{\kappa^\prime}^{d,m,D}}
f_{\kappa^\prime,\nu^\prime} g_{\kappa^\prime,\nu^\prime}^{d,m},
\{f_{\kappa^\prime,\nu^\prime} \in \mathcal P^{d,l},
\nu^\prime \in N_{\kappa^\prime}^{d,m,D}\},
\end{equation*}
имеет место представление
\begin{equation*} \tag{1.2.12}
H_{\kappa, \kappa^\prime}^{d,l,m,D} f = \sum_{\nu \in N_{\kappa}^{d,m,D}}
f_{\kappa,\nu} g_{\kappa,\nu}^{d,m},
\end{equation*}
где
\begin{multline*} \tag{1.2.13}
f_{\kappa,\nu} = \sum_{\m^{\k} \in \M_{\k}^m(\nu)} \biggl(\prod_{i \in \s(\k)}
(\prod_{s_i =1}^{\k_i} a_{\m_{i,s_i}}^{m_i})\biggr)
f_{\kappa^\prime, \n_{\k}(\nu,\m^{\k})} = \\
\sum_{\m^{\k} \in \M_{\k}^m(\nu)} A_{\m^{\k}}^m
f_{\kappa^\prime, \n_{\k}(\nu,\m^{\k})}, \nu \in N_{\kappa}^{d,m,D},
\end{multline*}
а
$$
A_{\m^{\k}}^m = \biggl(\prod_{i \in \s(\k)} (\prod_{s_i =1}^{\k_i}
a_{\m_{i,s_i}}^{m_i})\biggr), \m^{\k} \in \M_{\k}^m(\nu).
$$

Доказательство.

Для получения (1.2.10) достаточно, принимая во внимание (1.2.9),
последовательно применить (1.2.8).

Далее, покажем, что при любых $ \kappa^\prime \in \Z_+^d, \kappa =
\kappa^\prime +e_j, j =1,\ldots,d, $ выполняются равенства (1.2.12), (1.2.13).
Действительно, в этих условиях для $ f \in \mathcal P_{\kappa^\prime}^{d,l,m,D} $
вида (1.2.11), благодаря (1.2.9), (1.2.7), имеем
\begin{multline*} \tag{1.2.14}
H_{\kappa^\prime +e_j, \kappa^\prime}^{d,l,m,D} f = 
H_{\kappa^\prime}^{j,d,l,m,D} f = \\
\sum_{\nu \in N_{\kappa^\prime +e_j}^{d,m,D}} 
\biggl(\sum_{\substack{\nu^\prime \in N_{\kappa^\prime}^{d,m,D},
\m_j \in \Nu_{0, m_j +1}^1: \\ 2 \nu^\prime_j +\m_j = \nu_j,
\nu^\prime_i = \nu_i, i = 1,\ldots,d, i \ne j }} a_{\m_j}^{m_j}
f_{\kappa^\prime,\nu^\prime}\biggr) g_{\kappa^\prime +e_j,\nu}^{d,m} = 
\sum_{\nu \in N_{\kappa^\prime +e_j}^{d,m,D}} f_{\kappa^\prime +e_j,\nu}
g_{\kappa^\prime +e_j,\nu}^{d,m},
\end{multline*}
причем для $ \nu \in N_{\kappa^\prime +e_j}^{d,m,D} $ в силу замечания после
предложения 1.2.1 соблюдается равенство
\begin{multline*} \tag{1.2.15}
f_{\kappa^\prime +e_j,\nu} = \sum_{\substack{\nu^\prime \in N_{\kappa^\prime}^{d,m,D},
\m_j \in \Nu_{0, m_j +1}^1: \\ 2 \nu^\prime_j +\m_j = \nu_j,
\nu^\prime_i = \nu_i, i = 1,\ldots,d, i \ne j }} a_{\m_j}^{m_j}
f_{\kappa^\prime,\nu^\prime} = \\
\sum_{\substack{\nu^\prime \in \Z^d, \m_j \in \Nu_{0, m_j +1}^1: \\ (\nu_j -\m_j) /2 \in \Z,
\nu^\prime_j = (\nu_j -\m_j) /2, \nu^\prime_i = \nu_i, i = 1,\ldots,d, i \ne j }}
a_{\m_j}^{m_j} f_{\kappa^\prime,\nu^\prime} = \\
\sum_{ \m_j \in \Nu_{0, m_j +1}^1: (\nu_j -\m_j) /2 \in \Z}
\sum_{\substack{\nu^\prime \in \Z^d: \\ \nu^\prime_j = (\nu_j -\m_j) /2, \nu^\prime_i =
\nu_i, i = 1,\ldots,d, i \ne j }} a_{\m_j}^{m_j} f_{\kappa^\prime,\nu^\prime} = \\
\sum_{ \m^{e_j} \in \M_{e_j}^m(\nu)}
a_{\m_{j,1}}^{m_j} f_{\kappa^\prime,\n_{e_j}(\nu,\m^{e_j})},
\end{multline*}
т.е. имеют место (1.2.12), (1.2.13) при $ \kappa = \kappa^\prime +e_j,
j =1,\ldots,d. $

Теперь установим, что если для $ \kappa, \kappa^\prime \in \Z_+^d $ таких, что
$ \k = \kappa -\kappa^\prime \in \Z_+^d, $ соблюдаются равенства (1.2.12),
(1.2.13), то при $ j =1,\ldots,d $ для $ f \in \mathcal P_{\kappa^\prime}^{d,l,m,D} $
вида (1.2.11) имеет место представление
\begin{equation*} \tag{1.2.16}
H_\kappa^{j,d,l,m,D} (H_{\kappa, \kappa^\prime}^{d,l,m,D} f) =
\sum_{\nu \in N_{\kappa +e_j}^{d,m,D}}
f_{\kappa +e_j,\nu} g_{\kappa +e_j,\nu}^{d,m},
\end{equation*}
где
\begin{equation*} \tag{1.2.17}
f_{\kappa +e_j,\nu} = \sum_{\m^{\k +e_j} \in \M_{\k +e_j}^m(\nu)}
\biggl(\prod_{i \in \s(\k +e_j)} (\prod_{s_i =1}^{(\k +e_j)_i} a_{\m_{i,s_i}}^{m_i})\biggr)
f_{\kappa^\prime, \n_{\k +e_j}(\nu,\m^{\k +e_j})}, \nu \in N_{\kappa +e_j}^{d,m,D}.
\end{equation*}
В самом деле, для $ \kappa, \kappa^\prime \in \Z_+^d: \k = \kappa -
\kappa^\prime \in \Z_+^d, $ для которых выполняются (1.2.12), (1.2.13), при
$ j =1,\ldots,d $ для $ f $ вида (1.2.11) с учетом (1.2.12), (1.2.7)
имеем
\begin{multline*}
H_\kappa^{j,d,l,m,D} (H_{\kappa, \kappa^\prime}^{d,l,m,D} f) = \\
\sum_{\nu \in N_{\kappa +e_j}^{d,m,D}}
\biggl(\sum_{\substack{\rho \in N_{\kappa}^{d,m,D}, \m_j \in \Nu_{0, m_j +1}^1: \\ 2 \rho_j +\m_j =
\nu_j, \rho_i = \nu_i, i = 1,\ldots,d, i \ne j }} a_{\m_j}^{m_j}
f_{\kappa,\rho}\biggr) g_{\kappa +e_j,\nu}^{d,m} = 
\sum_{\nu \in N_{\kappa +e_j}^{d,m,D}} f_{\kappa +e_j,\nu} g_{\kappa +e_j,\nu}^{d,m},
\end{multline*}
причем при $ j =1,\ldots,d $ для $ \nu \in N_{\kappa +e_j}^{d,m,D} $ вследствие
замечания после предложения 1.2.1, а также в силу (1.2.13) выполняется
равенство
\begin{multline*} \tag{1.2.18}
f_{\kappa +e_j,\nu} =
\sum_{\substack{\rho \in N_{\kappa}^{d,m,D}, \m_j \in \Nu_{0, m_j +1}^1: \\ 2 \rho_j +\m_j =
\nu_j, \rho_i = \nu_i, i = 1,\ldots,d, i \ne j }} a_{\m_j}^{m_j}
f_{\kappa,\rho} = \\
\sum_{\substack{\rho \in \Z^d, \m_j \in \Nu_{0, m_j +1}^1: \\
(\nu_j -\m_j) /2 \in \Z, \rho_j = (\nu_j -\m_j) /2, \rho_i = \nu_i,
i = 1,\ldots,d, i \ne j }} a_{\m_j}^{m_j} f_{\kappa,\rho} = \\
\sum_{ \m_j \in \Nu_{0, m_j +1}^1: (\nu_j -\m_j) /2 \in \Z}
\sum_{\substack{\rho \in \Z^d: \\ \rho_j = (\nu_j -\m_j) /2, \rho_i = \nu_i,
i = 1,\ldots,d, i \ne j }} a_{\m_j}^{m_j} f_{\kappa,\rho} = \\
\sum_{ \m_j \in \Nu_{0, m_j +1}^1: (\nu_j -\m_j) /2 \in \Z}
\sum_{\substack{\rho \in \Z^d: \\ \rho_j = (\nu_j -\m_j) /2, \rho_i = \nu_i,
i = 1,\ldots,d, i \ne j }} a_{\m_j}^{m_j}
\biggl(\sum_{\m^{\k} \in \M_{\k}^m(\rho)} \biggl(\prod_{i \in \s(\k)}
(\prod_{s_i =1}^{\k_i} a_{\m_{i,s_i}}^{m_i})\biggr)
f_{\kappa^\prime, \n_{\k}(\rho,\m^{\k})}\biggr).
\end{multline*}
Пользуясь тем, что при $ j =1,\ldots,d $ для $ \nu \in N_{\kappa +e_j}^{d,m,D},
\m_j \in \Nu_{0, m_j +1}^1: (\nu_j -\m_j) /2 \in \Z, \rho \in \Z^d:
\rho_j = (\nu_j -\m_j) /2, \rho_i = \nu_i, i = 1,\ldots,d, i \ne j, $ 
множество
\begin{multline*}
\M_{\k}^m(\rho) = \{ \m^{\k} = \{ \m_{i,s_i} \in \Nu_{0, m_i +1}^1,
s_i =1,\ldots,\k_i, i \in \s(\k)\}:\\ \rho_i 2^{-\iota_i} -
\sum_{s_i =1}^{\iota_i} \m_{i, \k_i -s_i +1} 2^{-(\iota_i -s_i +1)} \in \Z\
\forall \iota_i = 1,\ldots, \k_i, i \in \s(\k)\} = \\
\{ \m^{\k} = \{ \m_{i,s_i} \in \Nu_{0, m_i +1}^1,
s_i =1,\ldots,\k_i, i \in \s(\k)\}: \nu_i 2^{-\iota_i} -
\sum_{s_i =1}^{\iota_i} \m_{i, \k_i -s_i +1} 2^{-(\iota_i -s_i +1)} \in \Z\ \\
\forall \iota_i = 1,\ldots, \k_i, i \in \s(\k) \setminus \{j\},
\nu_j 2^{-\iota_j -1} -\m_j 2^{-\iota_j -1} -
\sum_{s_j =1}^{\iota_j} \m_{j, \k_j -s_j +1} 2^{-(\iota_j -s_j +1)} \in \Z\ \\
\forall \iota_j = 1,\ldots, \k_j, j \in \s(\k) \} = \\
\{ \m^{\k} = \{ \m_{i,s_i} \in \Nu_{0, m_i +1}^1,
s_i =1,\ldots,\k_i, i \in \s(\k)\}: \nu_i 2^{-\iota_i} -
\sum_{s_i =1}^{\iota_i} \m_{i, \k_i -s_i +1} 2^{-(\iota_i -s_i +1)} \in \Z\ \\
\forall \iota_i = 1,\ldots, \k_i, i \in \s(\k) \setminus \{j\},
\nu_j 2^{-\iota^\prime_j} -\m_j 2^{-\iota^\prime_j} -
\sum_{s^\prime_j =2}^{\iota^\prime_j} \m_{j, \k_j +1 -s^\prime_j +1} 2^{-(\iota^\prime_j -s^\prime_j +1)} \in \Z\ \\
\forall \iota^\prime_j = 2,\ldots, \k_j +1, j \in \s(\k) \},
\end{multline*}
а для $ \m^{\k} \in \M_{\k}^m(\rho) $ мультииндекс $ \n_{\k}(\rho,\m^{\k}) $ 
имеет компоненты
\begin{multline*}
(\n_{\k}(\rho,\m^{\k}))_i = \begin{cases} \rho_i 2^{-\k_i} -
\sum_{s_i =1}^{\k_i} \m_{i, \k_i -s_i +1} 2^{-(\k_i -s_i +1)}, i \in \s(\k); \\
\rho_i, i \in \{1,\ldots,d\} \setminus \s(\k) \end{cases}
= \\
\begin{cases} \begin{cases} \nu_i 2^{-\k_i} -
\sum_{s_i =1}^{\k_i} \m_{i, \k_i -s_i +1} 2^{-(\k_i -s_i +1)}, i \in
\s(\k) \setminus \{j\}; \\
\nu_i, i \in \{1,\ldots,d\} \setminus \s(\k); \\
\nu_i 2^{-\k_i -1} -\m_i 2^{-\k_i -1} -
\sum_{s_i =1}^{\k_i} \m_{i, \k_i -s_i +1} 2^{-(\k_i -s_i +1)}, i = j, \\
\end{cases}
\text{ при } j \in \s(\k); \\
\begin{cases} \nu_i 2^{-\k_i} -
\sum_{s_i =1}^{\k_i} \m_{i, \k_i -s_i +1} 2^{-(\k_i -s_i +1)}, i \in
\s(\k); \\
\nu_i, i \in (\{1,\ldots,d\} \setminus \s(\k)) \setminus \{j\}; \\
(\nu_i -\m_i) /2, i =j, \\
\end{cases}
\text{ при } j \in (\{1,\ldots,d\} \setminus \s(\k)), \\
\end{cases} = \\
(\n_{\k +e_j}(\nu,\m^{\k +e_j}))_i, i = 1,\ldots,d, \text{ при } \\
\m^{\k +e_j} = (\m^{\k}, \m_j) = \\
\begin{cases} \{\m_{i,s_i}, s_i =1,\ldots,\k_i,
i \in \s(\k), \m_{j,\k_j +1} = \m_j\}, j \in \s(\k); \\
\{\m_{i,s_i}, s_i =1,\ldots,\k_i,
i \in \s(\k), \m_{j,1} = \m_j\}, j \in \{1,\ldots,d\} \setminus \s(\k),
\end{cases}
\end{multline*}
причем $ \m^{\k +e_j} \in \M_{\k +e_j}^m(\nu), $
из (1.2.18) получаем, что
\begin{multline*}
f_{\kappa +e_j,\nu} = \\ \sum_{ \m_j \in \Nu_{0, m_j +1}^1: (\nu_j -\m_j) /2 \in \Z}
a_{\m_j}^{m_j} \biggl(\sum_{\m^{\k} \in \M_{\k}^m(\rho) \mid_{\rho = \nu -((\nu_j +\m_j) /2) e_j}}
(\prod_{i \in \s(\k)} (\prod_{s_i =1}^{\k_i} a_{\m_{i,s_i}}^{m_i}))
f_{\kappa^\prime, \n_{\k +e_j}(\nu,\m^{\k +e_j})}\biggr) = \\
\sum_{ \m_j \in \Nu_{0, m_j +1}^1: (\nu_j -\m_j) /2 \in \Z}
\sum_{\m^{\k} \in \M_{\k}^m(\rho) \mid_{\rho = \nu -((\nu_j +\m_j) /2) e_j}}
a_{\m_j}^{m_j} (\prod_{i \in \s(\k)} (\prod_{s_i =1}^{\k_i} a_{\m_{i,s_i}}^{m_i}))
f_{\kappa^\prime, \n_{\k +e_j}(\nu,\m^{\k +e_j})} = \\
\sum_{\substack{ (\m^{\k}, \m_j) | \m_j \in \Nu_{0, m_j +1}^1:\\ (\nu_j -\m_j) /2 \in \Z,
\m^{\k} \in \M_{\k}^m(\rho) \mid_{\rho = \nu -((\nu_j +\m_j) /2) e_j}}}
a_{\m_j}^{m_j} (\prod_{i \in \s(\k)} (\prod_{s_i =1}^{\k_i} a_{\m_{i,s_i}}^{m_i}))
f_{\kappa^\prime, \n_{\k +e_j}(\nu,\m^{\k +e_j})} = \\
\sum_{ \m^{\k +e_j} \in \M_{\k +e_j}^m(\nu)}
(\prod_{i \in \s(\k +e_j)} (\prod_{s_i =1}^{(\k +e_j)_i} a_{\m_{i,s_i}}^{m_i}))
f_{\kappa^\prime, \n_{\k +e_j}(\nu,\m^{\k +e_j})},
\nu \in N_{\kappa +e_j}^{d,m,D}. .
\end{multline*}
Тем самым установлено соблюдение (1.2.16), (1.2.17) при выполнении
указанных выше условий.

Для завершения доказательства справедливости (1.2.12), (1.2.13) в
общей ситуации достаточно для $ \kappa^\prime, \kappa \in \Z_+^d: \kappa^\prime
\le \kappa, $ исходя из (1.2.9), выбрать последовательности
$ \kappa^i \in \Z_+^d, j_i \in \Nu_{1,d}^1, i =0,\ldots,K $ при $ K =
\sum_{j =1}^d \kappa_j -\kappa^\prime_j, $ для которых выполняются условия
$$
\kappa^0 = \kappa^\prime, \kappa^K = \kappa,
$$
\begin{equation*} \tag{1.2.19}
H_{\kappa^{i +1}, \kappa^\prime}^{d,l,m,D} =
H_{\kappa^{i}}^{j_i,d,l,m,D} H_{\kappa^{i}, \kappa^\prime}^{d,l,m,D}
(\kappa^{i +1} = \kappa^i +e_{j_i}), i = 0,\ldots,K -1,
\end{equation*}
и применить (1.2.14), (1.2.15), а затем соответствующее число раз -- (1.2.19),
(1.2.16), (1.2.17). $ \square $

Лемма 1.2.3

При $ d \in \N, \nu \in \Z^d, \k \in \Z_+^d, m \in \N^d $ имеет место
равенство
\begin{equation*} \tag{1.2.20}
\sum_{\m^{\k} \in \M_{\k}^m(\nu)} \prod_{i \in \s(\k)}
(\prod_{s_i =1}^{\k_i} a_{\m_{i,s_i}}^{m_i}) =1.
\end{equation*}

Доказательство.

Поскольку в условиях леммы множество $ \M_{\k}^m(\nu) $ представляется в виде
\begin{multline*}
\M_{\k}^m(\nu) = \{ \m^{\k} = \{ \m_{j,s_j} \in \Nu_{0, m_j +1}^1,
s_j =1,\ldots,\k_j, j \in \s(\k)\}: \\ \nu_j 2^{-\iota_j} -
\sum_{s_j =1}^{\iota_j} \m_{j, \k_j -s_j +1} 2^{-(\iota_j -s_j +1)} \in \Z \
\forall \iota_j = 1,\ldots, \k_j, j \in \s(\k)\} = \\
\prod_{j \in \S(\k)} \{ \{\m_{j,s_j} \in \Nu_{0, m_j +1}^1, s_j =1,\ldots,\k_j\}:
\nu_j 2^{-\iota_j} -\sum_{s_j =1}^{\iota_j} \m_{j, \k_j -s_j +1}
2^{-(\iota_j -s_j +1)} \in \Z \\
\forall \iota_j = 1,\ldots, \k_j\} = \prod_{j \in \s(\k)} \M_{\k_j}^{m_j}(\nu_j),
\end{multline*}
то
\begin{multline*}
\sum_{\m^{\k} \in \M_{\k}^m(\nu)} \prod_{i \in \s(\k)}
(\prod_{s_i =1}^{\k_i} a_{\m_{i,s_i}}^{m_i}) = 
\sum_{\{ \{\m_{i,s_i}, s_i =1,\ldots,\k_i\}, i \in \s(\k)\} \in
\prod_{i \in \s(\k)} \M_{\k_i}^{m_i}(\nu_i)} \prod_{i \in \s(\k)}
(\prod_{s_i =1}^{\k_i} a_{\m_{i,s_i}}^{m_i}) = \\
\prod_{i \in \s(\k)} \biggl(\sum_{ \{\m_{i,s_i}, s_i =1,\ldots,\k_i\} \in
\M_{\k_i}^{m_i}(\nu_i)} ( \prod_{s_i =1}^{\k_i} a_{\m_{i,s_i}}^{m_i})\biggr).
\end{multline*}
Отсюда видим, что (1.2.20) достаточно установить в случае $ d =1. $ В этом
случае при $ \nu \in \Z, \k \in \N, m \in \N, $ меняя индексацию, имеем
\begin{multline*} \tag{1.2.21}
\mathcal A_{\k}^m(\nu) = \sum_{ \{\m_{s}, s =1,\ldots,\k\} \in \M_{\k}^{m}(\nu)}
( \prod_{s =1}^{\k} a_{\m_{s}}^{m}) = \\
\sum_{ \{\m_{s} \in \Nu_{0,m +1}^1, s =1,\ldots,\k\}:
\nu 2^{-\iota} -\sum_{s =1}^{\iota} \m_{\k -s +1} 2^{-(\iota -s +1)} \in \Z
\forall \iota = 1,\ldots, \k} ( \prod_{s =1}^{\k} a_{\m_{s}}^{m}) = \\
\sum_{\substack{ \{\m_{s} \in \Nu_{0,m +1}^1, s =1,\ldots,\k -1\},
\m_{\k} \in \Nu_{0,m +1}^1: \\ (\nu -\m_{\k}) /2 \in \Z, \nu 2^{-\iota} -
\m_{\k} 2^{-\iota} -\sum_{s =2}^{\iota} \m_{\k -s +1} 2^{-(\iota -s +1)} \in \Z
\forall \iota = 2,\ldots, \k}} ( \prod_{s =1}^{\k} a_{\m_{s}}^{m}) = \\
\sum_{ \substack{\{\m_{s} \in \Nu_{0,m +1}^1, s =1,\ldots,\k -1\},
\m_{\k} \in \Nu_{0,m +1}^1: \\ (\nu -\m_{\k}) /2 \in \Z,
((\nu -\m_{\k}) /2) 2^{-\iota^\prime} -\sum_{s^\prime =1}^{\iota^\prime}
\m_{\k -1 -s^\prime +1} 2^{-(\iota^\prime -s^\prime +1)} \in \Z
\forall \iota^\prime = 1,\ldots, \k -1}} ( \prod_{s =1}^{\k} a_{\m_{s}}^{m}) = \\
\sum_{\m_{\k} \in \Nu_{0,m +1}^1: (\nu -\m_{\k}) /2 \in \Z}
\sum_{\substack{ \{\m_{s} \in \Nu_{0,m +1}^1, s =1,\ldots,\k -1\}: \\
((\nu -\m_{\k}) /2) 2^{-\iota^\prime} -\sum_{s^\prime =1}^{\iota^\prime}
\m_{\k -1 -s^\prime +1} 2^{-(\iota^\prime -s^\prime +1)} \in \Z
\forall \iota^\prime = 1,\ldots, \k -1}} a_{\m_{\k}}^m ( \prod_{s =1}^{\k -1}
a_{\m_{s}}^{m}) = \\
\sum_{\m_{\k} \in \Nu_{0,m +1}^1: (\nu -\m_{\k}) /2 \in \Z} a_{\m_{\k}}^m
\sum_{\substack{ \{\m_{s} \in \Nu_{0,m +1}^1, s =1,\ldots,\k -1\}: \\
((\nu -\m_{\k}) /2) 2^{-\iota^\prime} -\sum_{s^\prime =1}^{\iota^\prime}
\m_{\k -1 -s^\prime +1} 2^{-(\iota^\prime -s^\prime +1)} \in \Z
\forall \iota^\prime = 1,\ldots, \k -1}} ( \prod_{s =1}^{\k -1} a_{\m_{s}}^{m}) = \\
\sum_{\m_{\k} \in \Nu_{0,m +1}^1: (\nu -\m_{\k}) /2 \in \Z} a_{\m_{\k}}^m
\mathcal A_{\k -1}^m((\nu -\m_{\k}) /2).
\end{multline*}
Учитывая, что в силу (1.2.2) при $ m \in \N $ для любого $ \nu \in \Z $
справедливо равенство
\begin{equation*}
\mathcal A_1^m(\nu) = \sum_{\m \in \Nu_{0,m +1}^1: (\nu -\m) /2 \in \Z} a_{\m}^m =
\sum_{\m \in \Nu_{0,m +1}^1: \m \in (\nu +2\Z)} a_{\m}^m =1,
\end{equation*}
используя (1.2.21) и (1.2.2), по индукции относительно $ \k $ получаем, что
$ \mathcal A_{\k}^m(\nu) =1, m,\k \in \N, \nu \in \Z, $ что завершает вывод (1.2.20). $ \square $
\bigskip

1.3. В этом пункте определяются средства приближения функций из рассматриваемых
нами пространств.

Для формулировки предложения 1.3.1 потребуются следующие обозначения.
При $ d \in \N, l, \kappa \in \Z_+^d, \nu \in \Z^d $ определим линейный
оператор $ S_{\kappa,\nu}^{d,l}: L_1(Q_{\kappa,\nu}^d) \mapsto \mathcal P^{d,l}, $
полагая $ S_{\kappa, \nu}^{d,l} = P_{\delta, x^0}^{d,l} $ при
$ \delta = 2^{-\kappa}, x^0 = 2^{-\kappa} \nu $ (см. лемму 1.1.2).
Отметим, что в ситуации, когда $ Q_{\kappa,\nu}^d \subset D, $ где $ D $ --
область в $ \R^d, $ для $ f \in L_1(D) $ вместо $ S_{\kappa, \nu}^{d,l}( f \mid_{Q_{\kappa, \nu}^d}) $
будем писать $ S_{\kappa, \nu}^{d,l} f.$

Для области $ D \subset \R^d $ и $ \kappa \in \Z_+^d $ таких, что множество
\begin{equation*} \tag{1.3.1}
\{\nu^\prime \in \Z^d: Q_{\kappa,\nu^\prime}^d \subset D\} \ne \emptyset,
\end{equation*}
для каждого $ \nu \in \Z^d $ фиксируем некоторый $ \nu_\kappa^D(\nu) \in \Z^d, $
для которого $ Q_{\kappa,\nu_\kappa^D(\nu)}^d \subset D, $ а
$$
\|\nu -\nu_\kappa^D(\nu)\|_{l_\infty^d} = \min_{\nu^\prime \in \Z^d:
Q_{\kappa,\nu^\prime}^d \subset D} \|\nu -\nu^\prime\|_{l_\infty^d},
$$
и при $ d \in \N, l \in \Z_+^d, m \in \N^d $ для ограниченной области
$ D \subset \R^d $ и $ \kappa \in \Z_+^d, $ удовлетворяющих (1.3.1),
определим линейный непрерывный оператор
$ E_\kappa^{d,l,m,D}: L_1(D) \mapsto \mathcal P_\kappa^{d,l,m,D} \cap
L_\infty(\R^d) $ равенством
\begin{equation*}
E_\kappa^{d,l,m,D} f = \sum_{\nu \in N_\kappa^{d,m,D}}
(S_{\kappa, \nu_\kappa^D(\nu)}^{d,l} f ) g_{\kappa, \nu}^{d,m}, f \in L_1(D).
\end{equation*}

При $ k \in \Z_+, \alpha \in \R_+^d $ определим $ \kappa(k,\alpha) $
как вектор, имеющий компоненты
$ (\kappa(k,\alpha))_j = [k /\alpha_j], j =1,\ldots,d, ([a] \text{ -- целая часть} a),$
и при $ d \in \N, l \in \Z_+^d, m \in \N^d, \alpha \in \R_+^d, k \in \Z_+, $
ограниченной области $ D \subset \R^d, $ для которых при $ \kappa = \kappa(k,\alpha) $
выполняется (1.3.1), обозначим
$$
E_k^{d,l,m,D,\alpha} = E_\kappa^{d,l,m,D}
$$
при $ \kappa = \kappa(k,\alpha). $

Обозначим еще через $ \mathcal I^D $ линейное отображение, которое каждой
функции $ f, $ заданной на области $ D \subset \R^d, $ сопоставляет функцию
$ \mathcal I^D f, $ определяемую
на $ \R^d $ равенством
\begin{equation*}
(\mathcal I^D f)(x) = \begin{cases} f(x), \text{ при } x \in D; \\
0, \text{ при } x \in \R^d \setminus D.
\end{cases}
\end{equation*}

Предложениее 1.3.1

Пусть $ d \in \N, \alpha \in \R_+^d, l \in \N^d, m \in \N^d, 1 \le p < \infty $
и ограниченная область $ D \subset \R^d $ такова, что существуют $ k^0(D,\alpha)
\in \Z_+, \gamma^0(D,\alpha) \in \R_+^d, $ для которых при любом $ k \in \Z_+:
k \ge k^0, $ при $ \kappa = \kappa(k,\alpha) $ для каждого $ \nu \in \Z^d:
Q_{\kappa,\nu}^d \cap D \ne \emptyset, $ существует $ \nu^\prime \in \Z^d $
такой, что
\begin{equation*} \tag{1.3.2}
Q_{\kappa,\nu^\prime}^d \subset D \cap (2^{-\kappa} \nu +\gamma^0 2^{-\kappa} B^d).
\end{equation*}
Тогда для любой функции $ f \in L_p(D) $ в $ L_p(D) $ имеет место равенство
\begin{equation*} \tag{1.3.3}
f = (E_{k^0}^{d,l -\e,m,D,\alpha} f) \mid_D +\sum_{k = k^0 +1}^\infty
(E_k^{d,l -\e,m,D,\alpha} f -E_{k -1}^{d,l -\e,m,D,\alpha} f) \mid_D.
\end{equation*}

Доказательство.

Сначала отметим свойства некоторых вспомогательных множеств и других
объектов, которые понадобятся для доказательства предложения.

Для $ d \in \N, \kappa \in \Z_+^d, m \in \N^d $ и области $ D \subset \R^d $
обозначим
$$
G_\kappa^{d,m,D} = \cup_{\nu \in N_\kappa^{d,m,D}} \supp g_{\kappa,\nu}^{d,m},
$$
и рассмотрим представление
\begin{multline*} \tag{1.3.4}
G_\kappa^{d,m,D} = G_\kappa^{d,m,D} \cap \R^d = G_\kappa^{d,m,D} \cap
((\cup_{n \in \Z^d} Q_{\kappa,n}^d) \cup A_\kappa^d) = \\
(\cup_{n \in \Z^d} (G_\kappa^{d,m,D} \cap Q_{\kappa,n}^d)) \cup
(G_\kappa^{d,m,D} \cap A_\kappa^d) = \\
(\cup_{n \in \Z^d: Q_{\kappa,n}^d \cap G_\kappa^{d,m,D} \ne \emptyset}
(G_\kappa^{d,m,D} \cap Q_{\kappa,n}^d)) \cup (G_\kappa^{d,m,D} \cap A_\kappa^d), \\
\text{где} \mes A_\kappa^d =0, A_\kappa^d \cap Q_{\kappa,n}^d = \emptyset,
Q_{\kappa,n}^d \cap Q_{\kappa,n^\prime}^d = \emptyset, n,n^\prime \in
\Z^d: n \ne n^\prime.
\end{multline*}
Заметим, что если при $ n \in \Z^d $ пересечение
\begin{equation*}
Q_{\kappa,n}^d \cap G_\kappa^{d,m,D} \ne \emptyset,
\end{equation*}
т.е.
\begin{equation*}
Q_{\kappa,n}^d \cap (\cup_{\nu \in N_\kappa^{d,m,D}} \supp g_{\kappa,\nu}^{d,m}) =
\cup_{\nu \in N_\kappa^{d,m,D}} (Q_{\kappa,n}^d \cap \supp g_{\kappa,\nu}^{d,m})
\ne \emptyset,
\end{equation*}
то существует $ \nu \in N_\kappa^{d,m,D}, $ для которого
\begin{equation*}
Q_{\kappa,n}^d \cap \supp g_{\kappa,\nu}^{d,m} \ne \emptyset.
\end{equation*}
Для $ n \in \Z^d, \nu \in N_\kappa^{d,m,D}: Q_{\kappa,n}^d \cap
\supp g_{\kappa,\nu}^{d,m} \ne \emptyset, $ выбирая $ x \in
Q_{\kappa,n}^d \cap \supp g_{\kappa,\nu}^{d,m}, $ при $ j =1,\ldots,d $ имеем
$$
2^{-\kappa_j} n_j < x_j \le 2^{-\kappa_j} \nu_j +2^{-\kappa_j} (m_j +1);
2^{-\kappa_j} \nu_j \le x_j < 2^{-\kappa_j} n_j +2^{-\kappa_j},
$$
или
$$
n_j < \nu_j +m_j +1;
\nu_j < n_j +1,
$$
откуда
$$
n_j +1 \le \nu_j +m_j +1;
\nu_j \le n_j,
$$
следовательно,
$$
2^{-\kappa_j} n_j +2^{-\kappa_j} \le 2^{-\kappa_j} \nu_j +2^{-\kappa_j} (m_j +1);
2^{-\kappa_j} \nu_j \le 2^{-\kappa_j} n_j,
$$
а это значит, что
\begin{equation*} \tag{1.3.5}
Q_{\kappa,n}^d \subset \overline Q_{\kappa,n}^d \subset
\supp g_{\kappa,\nu}^{d,m} \subset G_\kappa^{d,m,D},
n \in \Z^d, \nu \in N_\kappa^{d,m,D}: Q_{\kappa,n}^d \cap
\supp g_{\kappa,\nu}^{d,m} \ne \emptyset.
\end{equation*}
Из (1.3.4) и (1.3.5) получаем
\begin{multline*} \tag{1.3.6}
G_\kappa^{d,m,D} = (\cup_{n \in \Z^d: Q_{\kappa,n}^d \cap G_\kappa^{d,m,D} \ne
\emptyset} Q_{\kappa,n}^d) \cup (G_\kappa^{d,m,D} \cap A_\kappa^d), \\
d \in \N, m \in \N^d, \kappa \in \Z_+^d, D  \text{ -- область в} \R^d.
\end{multline*}

В условиях предложения при $ k \ge k^0, \kappa = \kappa(k,\alpha) $
для $ \nu \in N_\kappa^{d,m,D}, $ выбирая точку $ z \in (2^{-\kappa} \nu +
2^{-\kappa} (m +\e) \overline I^d) \cap D, $ и учитывая открытость множества
$ D, $ найдем точку $ y \in (2^{-\kappa} \nu +2^{-\kappa} (m +\e) I^d) \cap D. $
Затем возьмем клетку $ \overline Q_{\kappa,\rho}^d, \rho \in \Z^d, $
содержащую точку $ y. $ Тогда поскольку
$ (2^{-\kappa} \nu +2^{-\kappa} (m +\e) I^d) \cap D $ является
окрестностью точки $ y, $ принадлежащей замыканию $ Q_{\kappa,\rho}^d, $ то
множество $ (2^{-\kappa} \nu +2^{-\kappa} (m +\e) I^d) \cap D \cap
Q_{\kappa,\rho}^d \ne \emptyset. $
Принимая во внимание это обстоятельство и (1.3.5), получаем, что
\begin{equation*} \tag{1.3.7}
Q_{\kappa,\rho}^d \subset \overline Q_{\kappa,\rho}^d \subset
\supp g_{\kappa,\nu}^{d,m},
\end{equation*}
и $ Q_{\kappa,\rho}^d \cap D \ne \emptyset. $
Исходя из условий предложения, согласно  (1.3.2) выберем $ \nu^\prime \in \Z^d, $
для которого
\begin{equation*} \tag{1.3.8}
Q_{\kappa,\nu^\prime}^d \subset D \cap
(2^{-\kappa} \rho +\gamma^0 2^{-\kappa} B^d).
\end{equation*}
Из (1.3.8) и определения $ \nu_\kappa^D(\nu) $ имеем
\begin{equation*} \tag{1.3.9}
\|\nu -\nu_\kappa^D(\nu)\|_{l_\infty^d} \le \|\nu -\nu^\prime\|_{l_\infty^d} \le
\|\nu -\rho\|_{l_\infty^d} +\|\\rho -\nu^\prime\|_{l_\infty^d}.
\end{equation*}
Ввиду (1.3.7) получаем
$$
| 2^{-\kappa_J} \nu_j -2^{-\kappa_J} \rho_j | \le 2^{-\kappa_J} (m_j +1),
$$
или
$$
| \nu_j -\rho_j | \le (m_j +1), j =1,\ldots,d,
$$
т.е.
\begin{equation*} \tag{1.3.10}
\|\nu -\rho\|_{l_\infty^d} \le c_1(d,m).
\end{equation*}
Кроме того, из (1.3.8) следует, что
$$
| 2^{-\kappa_J} \nu_j^\prime -2^{-\kappa_J} \rho_j | \le 2^{-\kappa_J} \gamma^0_j,
$$
или
$$
| \nu_j^\prime -\rho_j | \le \gamma^0_j, j =1,\ldots,d,
$$
т.е.
\begin{equation*} \tag{1.3.11}
\| \nu^\prime -\rho\|_{l_\infty^d} \le c_2(d,\gamma^0).
\end{equation*}
Объединяя (1.3.9), (1.3.10), (1.3.11), приходим к выводу, что
\begin{multline*} \tag{1.3.12}
\|\nu -\nu_\kappa^D(\nu)\|_{l_\infty^d} \le c_1(d,m) +c_2(d,\gamma^0) =
c_3(d,m,D,\alpha),\\ \nu \in N_\kappa^{d,m,D},
\kappa = \kappa(k, \alpha), k \ge k^0.
\end{multline*}

Отметим еще, что в условиях предложения при $ k \ge k^0, \kappa = \kappa(k,\alpha), $
для $ n \in \Z^d: Q_{\kappa,n}^d \cap G_\kappa^{d,m,D} \ne \emptyset,
\nu \in N_\kappa^{d,m,D}: \supp g_{\kappa, \nu}^{d,m} \cap Q_{\kappa,n}^d
\ne \emptyset, $ для $ x \in \overline Q_{\kappa,\nu_\kappa^D(\nu)}^d $ в силу
(1.3.5) и (1.3.12) при $ j =1,\ldots, d $ выполняется неравенство
\begin{multline*}
| x_j -2^{-\kappa_j} n_j | \le \\
 | x_j -2^{-\kappa_j} (\nu_\kappa^D(\nu))_j |
+| 2^{-\kappa_j} (\nu_\kappa^D(\nu))_j -2^{-\kappa_j} \nu_j |
+| 2^{-\kappa_j} \nu_j -2^{-\kappa_j} n_j |\le \\
 2^{-\kappa_j} +2^{-\kappa_j}
\|\nu -\nu_\kappa^D(\nu)\|_{l_\infty^d} +2^{-\kappa_j} (m_j +1) \le
2^{-\kappa_j} (1 +c_3 +m_j +1) = \\
\gamma_j^1(d,m,D,\alpha) 2^{-\kappa_j},
\end{multline*}
т.е.
\begin{equation*} \tag{1.3.13}
Q_{\kappa,\nu_\kappa^D(\nu)}^d \subset \overline Q_{\kappa,\nu_\kappa^D(\nu)}^d
\subset (2^{-\kappa} n +\gamma^1 2^{-\kappa} B^d).
\end{equation*}

В условиях предложения при $ k \ge k^0,  \kappa = \kappa(k,\alpha),
n \in \Z^d: Q_{\kappa,n}^d \cap G_\kappa^{d,m,D} \ne \emptyset, $ зададим
$$
x_{\kappa,n}^{d,m,D,\alpha} = 2^{-\kappa} n -\gamma^1 2^{-\kappa};\\
\delta_{\kappa,n}^{d,m,D,\alpha} = 2 \gamma^1 2^{-\kappa}
$$
и определим клетку $ D_{\kappa,n}^{d,m,D,\alpha} $ равенством
$$
D_{\kappa,n}^{d,m,D,\alpha} = x_{\kappa,n}^{d,m,D,\alpha} +
\delta_{\kappa,n}^{d,m,D,\alpha} I^d = \inter (2^{-\kappa} n +\gamma^1 2^{-\kappa} B^d).
$$

Из приведенных определений с учетом того, что $ \gamma^1 > \e, $ видно, что
справедливо включение
\begin{equation*} \tag{1.3.14}
Q_{\kappa, n}^d \subset D_{\kappa,n}^{d,m,D,\alpha}, n \in \Z^d:
Q_{\kappa,n}^d \cap G_\kappa^{d,m,D} \ne \emptyset, \kappa = \kappa(k,\alpha),
k \ge k^0.
\end{equation*}
При $ k \ge k^0, \kappa = \kappa(k,\alpha), n \in \Z^d: Q_{\kappa,n}^d \cap
G_\kappa^{d,m,D} \ne \emptyset, $ введем в рассморение линейный оператор
$ \mathcal S_{\kappa,n}^{d,l -\e,m,D,\alpha}, $ полагая
\begin{equation*}
\mathcal S_{\kappa,n}^{d,l -\e,m,D,\alpha} = P_{\\delta,x^0}^{d,l -\e}
\end{equation*}
при $ x^0 = x_{\kappa,n}^{d,m,D,\alpha}, \delta = \delta_{\kappa,n}^{d,m,D,\alpha}. $

Учитывая (1.3.14), нетрудно видеть, что в условиях предложения существует
константа $ c_4(d,m,D,\alpha) >0 $ такая, что для любого $ k \in \Z_+:
k \ge k^0, $ при $ \kappa = \kappa(k,\alpha) $ для каждого $ x \in \R^d $
число
\begin{equation*} \tag{1.3.15}
\card \{ n \in \Z^d: Q_{\kappa,n}^d \cap G_\kappa^{d,m,D} \ne \emptyset,
x \in D_{\kappa,n}^{d,m,D,\alpha} \} \le c_4.
\end{equation*}

Из (1.3.13) и определения $ D_{\kappa,n}^{d,m,D,\alpha} $ следует, что
при $ k \ge k^0, \kappa = \kappa(k,\alpha) $ для $ n \in \Z^d:
Q_{\kappa,n}^d \cap G_\kappa^{d,m,D} \ne \emptyset, \nu \in N_\kappa^{d,m,D}:
\supp g_{\kappa, \nu}^{d,m} \cap Q_{\kappa,n}^d \ne \emptyset, $ имеет место
включение
\begin{equation*} \tag{1.3.16}
Q_{\kappa,\nu_\kappa^D(\nu)}^d \subset D_{\kappa,n}^{d,m,D,\alpha}.
\end{equation*}

Из (1.2.3) вытекает, что при $ \kappa \in \Z_+^d $ для $ n \in \Z^d:
Q_{\kappa,n}^d \cap G_\kappa^{d,m,D} \ne \emptyset, $ верно неравенство
\begin{equation*} \tag{1.3.17}
\card \{ \nu \in N_\kappa^{d,m,D}: \supp g_{\kappa, \nu}^{d,m} \cap
Q_{\kappa,n}^d \ne \emptyset \} \le c_5(d,m).
\end{equation*}

Пусть теперь в условиях предложения $ f \in L_p(D) $ и $ k \in \Z_+: k \ge k^0,
\kappa = \kappa(k,\alpha). $
Тогда, полагая $ F = \mathcal I^D f, $ ввиду (1.2.4) имеем
\begin{multline*} \tag{1.3.18}
\biggl\|f -((E_{k^0}^{d,l -\e,m,D,\alpha} f) \mid_D
+\sum_{\k = k^0 +1}^k (E_{\k}^{d,l -\e,m,D,\alpha} f
-E_{\k -1}^{d,l -\e,m,D,\alpha} f) \mid_D) \biggr\|_{L_p(D)}^p = \\
\| F \mid_D -(E_{k}^{d,l -\e,m,D,\alpha} f) \mid_D \|_{L_p(D)}^p =
\| F \mid_D -(E_\kappa^{d,l -\e,m,D} f) \mid_D \|_{L_p(D)}^p = \\
\biggl\| (F (\sum_{ \nu \in \Z^d: \supp g_{\kappa, \nu}^{d,m} \cap D \ne \emptyset}
g_{\kappa, \nu}^{d,m})) \mid_D -(E_\kappa^{d,l -\e,m,D} f) \mid_D \biggr\|_{L_p(D)}^p \le \\
\biggl\| F (\sum_{ \nu \in N_\kappa^{d,m,D}}
g_{\kappa, \nu}^{d,m}) -(E_\kappa^{d,l -\e,m,D} f) \biggr\|_{L_p(\R^d)}^p = \\
\biggl\| F (\sum_{ \nu \in N_\kappa^{d,m,D}}
g_{\kappa, \nu}^{d,m}) -\sum_{\nu \in N_\kappa^{d,m,D}}
(S_{\kappa, \nu_\kappa^D(\nu)}^{d,l -\e} f) g_{\kappa, \nu}^{d,m} \biggr\|_{L_p(\R^d)}^p = \\
\biggl\| F (\sum_{ \nu \in N_\kappa^{d,m,D}}
g_{\kappa, \nu}^{d,m}) -\sum_{\nu \in N_\kappa^{d,m,D}}
(S_{\kappa, \nu_\kappa^D(\nu)}^{d,l -\e} F) g_{\kappa, \nu}^{d,m} \biggr\|_{L_p(\R^d)}^p = \\
\int_{\R^d}\biggl| \sum_{\nu \in N_\kappa^{d,m,D}} (F(x) -
(S_{\kappa, \nu_\kappa^D(\nu)}^{d,l -\e} F)(x))
g_{\kappa, \nu}^{d,m}(x)\biggr|^p dx = \\
\int_{G_\kappa^{d,m,D}}\biggl| \sum_{\nu \in N_\kappa^{d,m,D}} (F(x) -
(S_{\kappa, \nu_\kappa^D(\nu)}^{d,l -\e} F)(x))
g_{\kappa, \nu}^{d,m}(x)\biggr|^p dx.
\end{multline*}

В силу (1.3.6) (см. также (1.3.4)) выводим
\begin{multline*} \tag{1.3.19}
\int_{G_\kappa^{d,m,D}}\biggl| \sum_{\nu \in N_\kappa^{d,m,D}} (F(x) -
(S_{\kappa, \nu_\kappa^D(\nu)}^{d,l -\e} F)(x))
g_{\kappa, \nu}^{d,m}(x)\biggr|^p dx = \\
\int_{\cup_{n \in \Z^d: Q_{\kappa,n}^d \cap G_\kappa^{d,m,D} \ne \emptyset}
Q_{\kappa,n}^d}\biggl| \sum_{\nu \in N_\kappa^{d,m,D}} (F(x) -
(S_{\kappa, \nu_\kappa^D(\nu)}^{d,l -\e} F)(x))
g_{\kappa, \nu}^{d,m}(x)\biggr|^p dx = \\
\sum_{n \in \Z^d: Q_{\kappa,n}^d \cap G_\kappa^{d,m,D} \ne \emptyset}
\int_{Q_{\kappa,n}^d} \biggl| \sum_{\nu \in N_\kappa^{d,m,D}} (F(x) -
(S_{\kappa, \nu_\kappa^D(\nu)}^{d,l -\e} F)(x))
g_{\kappa, \nu}^{d,m}(x)\biggr|^p dx = \\
\sum_{n \in \Z^d: Q_{\kappa,n}^d \cap G_\kappa^{d,m,D} \ne \emptyset}
\int_{Q_{\kappa,n}^d} \biggl| \sum_{\nu \in N_\kappa^{d,m,D}:
\supp g_{\kappa, \nu}^{d,m} \cap Q_{\kappa,n}^d \ne \emptyset} (F(x) -
(S_{\kappa, \nu_\kappa^D(\nu)}^{d,l -\e} F)(x))
g_{\kappa, \nu}^{d,m}(x)\biggr|^p dx = \\
\sum_{n \in \Z^d: Q_{\kappa,n}^d \cap G_\kappa^{d,m,D} \ne \emptyset}
\biggl\| \sum_{\nu \in N_\kappa^{d,m,D}:
\supp g_{\kappa, \nu}^{d,m} \cap Q_{\kappa,n}^d \ne \emptyset} (F -
(S_{\kappa, \nu_\kappa^D(\nu)}^{d,l -\e} F))
g_{\kappa, \nu}^{d,m} \biggr\|_{L_p(Q_{\kappa,n}^d)}^p \le \\
\sum_{n \in \Z^d: Q_{\kappa,n}^d \cap G_\kappa^{d,m,D} \ne \emptyset}
\biggl( \sum_{\nu \in N_\kappa^{d,m,D}:
\supp g_{\kappa, \nu}^{d,m} \cap Q_{\kappa,n}^d \ne \emptyset}
\| (F -(S_{\kappa, \nu_\kappa^D(\nu)}^{d,l -\e} F))
g_{\kappa, \nu}^{d,m} \|_{L_p(Q_{\kappa,n}^d)}\biggr)^p.
\end{multline*}

Далее, для $ n \in \Z^d: Q_{\kappa,n}^d \cap G_\kappa^{d,m,D} \ne \emptyset,
\nu \in N_\kappa^{d,m,D}: \supp g_{\kappa, \nu}^{d,m} \cap Q_{\kappa,n}^d
\ne \emptyset, $ с учетом (1.2.5) получаем
\begin{multline*} \tag{1.3.20}
\| (F -(S_{\kappa, \nu_\kappa^D(\nu)}^{d,l -\e} F))
g_{\kappa, \nu}^{d,m} \|_{L_p(Q_{\kappa,n}^d)} \le 
\| F -S_{\kappa, \nu_\kappa^D(\nu)}^{d,l -\e} F \|_{L_p(Q_{\kappa,n}^d)} \le \\
\|F -\mathcal S_{\kappa,n}^{d,l -\e,m,D,\alpha} F \|_{L_p(Q_{\kappa,n}^d)} +
\|\mathcal S_{\kappa,n}^{d,l -\e,m,D,\alpha} F
-S_{\kappa, \nu_\kappa^D(\nu)}^{d, l-\e} F \|_{L_p(Q_{\kappa,n}^d)}.
\end{multline*}

Учитывая (1.3.14), на основании (1.1.4) для $ n \in \Z^d: Q_{\kappa,n}^d
\cap G_\kappa^{d,m,D} \ne \emptyset, $ заключаем, что
\begin{multline*} \tag{1.3.21}
\|F -\mathcal S_{\kappa,n}^{d,l -\e,m,D,\alpha} F \|_{L_p(Q_{\kappa,n}^d)} \le
\|F -\mathcal S_{\kappa,n}^{d,l -\e,m,D,\alpha} F \|_{L_p(D_{\kappa,n}^{d,m,D,\alpha})} \\
\le c_6 \sum_{j =1}^d 2^{\kappa_j /p} \biggl(\int_{ c_7 2^{-\kappa_j} B^1}
\int_{ (D_{\kappa,n}^{d,m,D,\alpha})_{l_j \xi e_j}}
|\Delta_{\xi e_j}^{l_j} F(x)|^p dx d\xi\biggr)^{1/p}.
\end{multline*}

Принимая во внимание (1.3.14), (1.3.16), (1.1.1), (1.1.2), (1.1.3) и снова
(1.3.16), а также, опираясь на (1.1.4), для $ n \in \Z^d: Q_{\kappa,n}^d \cap
G_\kappa^{d,m,D} \ne \emptyset, \nu \in N_\kappa^{d,m,D}:
\supp g_{\kappa, \nu}^{d,m} \cap Q_{\kappa,n}^d \ne \emptyset, $
находим, что
\begin{multline*} \tag{1.3.22}
\|\mathcal S_{\kappa,n}^{d,l -\e,m,D,\alpha} F
-S_{\kappa,\nu_\kappa^D(\nu)}^{d, l-\e} F \|_{L_p(Q_{\kappa,n}^d)} \le \\
c_8 \|S_{\kappa, \nu_\kappa^D(\nu)}^{d, l -\e}
(\mathcal S_{\kappa,n}^{d,l -\e,m,D,\alpha} F -F) \|_{L_p(Q_{\kappa, \nu_\kappa^D(\nu)}^d)} \le 
c_9 \|F -\mathcal S_{\kappa,n}^{d,l -\e,m,D,\alpha} F \|_{L_p(D_{\kappa,n}^{d,m,D,\alpha})} \le \\
c_{10} \sum_{j =1}^d 2^{\kappa_j /p} \biggl(\int_{ c_7 2^{-\kappa_j} B^1}
\int_{ (D_{\kappa,n}^{d,m,D,\alpha})_{l_j \xi e_j}}
|\Delta_{\xi e_j}^{l_j} F(x)|^p dx d\xi\biggr)^{1/p}.
\end{multline*}

Соединяя (1.3.20), (1.3.21) и (1.3.22), для $ n \in \Z^d: Q_{\kappa,n}^d
\cap G_\kappa^{d,m,D} \ne \emptyset, \nu \in N_\kappa^{d,m,D}:
\supp g_{\kappa, \nu}^{d,m} \cap Q_{\kappa,n}^d \ne \emptyset, $ получаем
\begin{equation*}
\|(F -S_{\kappa, \nu_\kappa^D(\nu)}^{d, l-\e} F)
g_{\kappa,\nu}^{d,m} \|_{L_p(Q_{\kappa,n}^d)} \le \\
c_{11} \sum_{j =1}^d 2^{\kappa_j /p} \biggl(\int_{ c_7 2^{-\kappa_j} B^1}
\int_{ (D_{\kappa,n}^{d,m,D,\alpha})_{l_j \xi e_j}}
|\Delta_{\xi e_j}^{l_j} F(x)|^p dx d\xi\biggr)^{1/p}.
\end{equation*}

Подставляя эту оценку в (1.3.19) и применяя (1.3.17) и неравенство
Г\"eльдера, а затем используя (1.3.15) и неравенство Г\"eльдера, выводим
\begin{multline*}
\int_{G_\kappa^{d,m,D}}\biggl| \sum_{\nu \in N_\kappa^{d,m,D}} (F(x) -
(S_{\kappa, \nu_\kappa^D(\nu)}^{d,l -\e} F)(x))
g_{\kappa, \nu}^{d,m}(x)\biggr|^p dx \le \\
\sum_{\substack{n \in \Z^d: \\ Q_{\kappa,n}^d \cap G_\kappa^{d,m,D} \ne \emptyset}}
\biggl( \sum_{\substack{\nu \in N_\kappa^{d,m,D}: \\
\supp g_{\kappa, \nu}^{d,m} \cap Q_{\kappa,n}^d \ne \emptyset}}
c_{11} \sum_{j =1}^d 2^{\kappa_j /p} \biggl(\int_{ c_7 2^{-\kappa_j} B^1}
\int_{ (D_{\kappa,n}^{d,m,D,\alpha})_{l_j \xi e_j}}
|\Delta_{\xi e_j}^{l_j} F(x)|^p dx d\xi\biggr)^{1/p}\biggr)^p \le \\
\sum_{n \in \Z^d: Q_{\kappa,n}^d \cap G_\kappa^{d,m,D} \ne \emptyset}
c_{12}^p \biggl(\sum_{j =1}^d 2^{\kappa_j} \int_{ c_7 2^{-\kappa_j} B^1}
\int_{ (D_{\kappa,n}^{d,m,D,\alpha})_{l_j \xi e_j}}
|\Delta_{\xi e_j}^{l_j} F(x)|^p dx d\xi\biggr) \le \\
c_{12}^p \sum_{j =1}^d 2^{\kappa_j } \int_{ c_7 2^{-\kappa_j} B^1}
\biggl(\sum_{n \in \Z^d: Q_{\kappa,n}^d \cap G_\kappa^{d,m,D} \ne \emptyset}
\int_{D_{\kappa,n}^{d,m,D,\alpha}} |\Delta_{\xi e_j}^{l_j} F(x)|^p dx \biggr) d\xi = \\
c_{12}^p \sum_{j =1}^d 2^{\kappa_j } \int_{ c_7 2^{-\kappa_j} B^1}
\int_{ \R^d} \biggl(\sum_{n \in \Z^d: Q_{\kappa,n}^d \cap G_\kappa^{d,m,D} \ne
\emptyset} \chi_{D_{\kappa,n}^{d,m,D,\alpha}}(x)\biggr)
|\Delta_{\xi e_j}^{l_j} F(x)|^p dx d\xi \\
\le \biggl(c_{13} \sum_{j =1}^d 2^{\kappa_j /p} \biggl(\int_{c_7 2^{-\kappa_j} B^1}
\int_{ \R^d} |\Delta_{\xi e_j}^{l_j} F(x)|^p dx d\xi\biggr)^{1/p}\biggr)^p.
\end{multline*}
Объединяя последнее неравенство с (1.3.18), для $ f \in L_p(D) $ при
$ k \in \Z_+: k \ge k^0, $ приходим к соотношению
\begin{multline*}
\| f -((E_{k^0}^{d,l -\e,m,D,\alpha} f) \mid_D
+\sum_{\k = k^0 +1}^k (E_{\k}^{d,l -\e,m,D,\alpha} f
-E_{\k -1}^{d,l -\e,m,D,\alpha} f) \mid_D) \|_{L_p(D)} \le \\
c_{14} \sum_{j =1}^d 2^{k / (\alpha_j p)} \biggl(\int_{c_{15} 2^{-k / \alpha_j} B^1}
\int_{ \R^d} |\Delta_{\xi e_j}^{l_j} F(x)|^p dx d\xi\biggr)^{1/p} \le \\
c_{16} \sum_{j =1}^d \Omega_j^{l_j}(\mathcal I^D f,
c_{15} 2^{-k / \alpha_j})_{L_p(\R^d)} \to 0 \text{ при } k \to \infty,
\end{multline*}
что влеч\"eт (1.3.3). $ \square $

Как известно, имеет место

Лемма 1.3.2

Пусть $ d \in \N, \lambda \in \Z_+^d, D $ --- область в $ \R^d $ и
функция $ f \in C^\infty(D), $ а $ g \in L_1(D), $ прич\"eм
для каждого $ \mu \in \Z_+^d(\lambda) $ обобщ\"eнная производная
$ \D^\mu g \in L_1(D). $ Тогда в пространстве обобщ\"eнных функций в
области $ D $ имеет место равенство
\begin{equation*} \tag{1.3.23}
\D^\lambda (fg) = \sum_{ \mu \in \Z_+^d(\lambda)} C_\lambda^\mu
\D^{\lambda -\mu} f \D^\mu g.
\end{equation*}

Лемма 1.3.3

Пусть $ d \in \N, \alpha \in \R_+^d, l \in \Z_+^d, m \in \N^d, \lambda \in
\Z_+^d(m), 1 \le s, q \le \infty $ и область $ D \subset \R^d $ удовлетворяет
условиям предложения 1.3.1. Тогда существует константа
$ c_{17}(d,l,m,D,\alpha,\lambda,s,q) > 0 $ такая, что для любой функции
$ f \in L_s(D) $ при $ k \in \Z_+: k \ge k^0(D,\alpha), $ справедливо
неравенство
\begin{equation*} \tag{1.3.24}
\| \D^\lambda E_k^{d,l,m,D,\alpha} f \|_{L_q(\R^d)} \le
c_{17} 2^{k (\alpha^{-1}, \lambda +(s^{-1} -q^{-1})_+ \e)} \| f\|_{L_s(D)}.
\end{equation*}

Доказательство.

При доказательсве леммы будем использовать объекты и связанные с ними
факты из доказательства предложения 1.3.1. В условиях леммы для $ f \in
L_s(D) $ при $ s \le q, k \in \Z_+: k \ge k^0, \kappa = \kappa(k,\alpha), $ ,
принимая во внимание (1.3.23), имеем
\begin{multline*} \tag{1.3.25}
\| \D^\lambda (E_k^{d,l,m,D,\alpha} f) \|_{L_q(\R^d)} =
\| \D^\lambda E_\kappa^{d,l,m,D} f \|_{L_q(\R^d)} = \\
\biggl\| \D^\lambda (\sum_{\nu \in N_\kappa^{d,m,D}}
(S_{\kappa, \nu_\kappa^D(\nu)}^{d,l} f) g_{\kappa, \nu}^{d,m}) \biggr\|_{L_q(\R^d)} = 
\biggl\| \sum_{\nu \in N_\kappa^{d,m,D}} \D^\lambda
((S_{\kappa, \nu_\kappa^D(\nu)}^{d,l} f) g_{\kappa, \nu}^{d,m}) \biggr\|_{L_q(\R^d)} = \\
\biggl\| \sum_{\nu \in N_\kappa^{d,m,D}} \sum_{ \mu \in \Z_+^d(\lambda)}
C_\lambda^\mu (\D^\mu S_{\kappa, \nu_\kappa^D(\nu)}^{d,l} f)
\D^{\lambda -\mu} g_{\kappa, \nu}^{d,m} \biggr\|_{L_q(\R^d)} = \\
\biggl\| \sum_{ \mu \in \Z_+^d(\lambda)} C_\lambda^\mu
\sum_{\nu \in N_\kappa^{d,m,D}}
(\D^\mu S_{\kappa, \nu_\kappa^D(\nu)}^{d,l} f)
\D^{\lambda -\mu} g_{\kappa, \nu}^{d,m} \biggr\|_{L_q(\R^d)} \le \\
\sum_{ \mu \in \Z_+^d(\lambda)} C_\lambda^\mu
\biggl\| \sum_{\nu \in N_\kappa^{d,m,D}}
(\D^\mu S_{\kappa, \nu_\kappa^D(\nu)}^{d,l} f)
\D^{\lambda -\mu} g_{\kappa, \nu}^{d,m} \biggr\|_{L_q(\R^d)}.
\end{multline*}

Оценивая правую часть (1.3.25), при $ \mu \in \Z_+^d(\lambda) $ с учетом
(1.3.6), (1.3.4) получаем (ср. с (1.3.19))
\begin{multline*} \tag{1.3.26}
\biggl\| \sum_{\nu \in N_\kappa^{d,m,D}}
(\D^\mu S_{\kappa, \nu_\kappa^D(\nu)}^{d,l} f)
\D^{\lambda -\mu} g_{\kappa, \nu}^{d,m} \biggr\|_{L_q(\R^d)}^q = \\
\int_{\R^d} \biggl| \sum_{\nu \in N_\kappa^{d,m,D}}
(\D^\mu S_{\kappa, \nu_\kappa^D(\nu)}^{d,l} f)
\D^{\lambda -\mu} g_{\kappa, \nu}^{d,m} \biggr|^q dx = \\
\int_{G_\kappa^{d,m,D}} \biggl| \sum_{\nu \in N_\kappa^{d,m,D}}
(\D^\mu S_{\kappa, \nu_\kappa^D(\nu)}^{d,l} f)
\D^{\lambda -\mu} g_{\kappa, \nu}^{d,m} \biggr|^q dx = \\
\sum_{n \in \Z^d: Q_{\kappa,n}^d \cap G_\kappa^{d,m,D} \ne \emptyset}
\int_{Q_{\kappa,n}^d} \biggl| \sum_{\nu \in N_\kappa^{d,m,D}}
(\D^\mu S_{\kappa, \nu_\kappa^D(\nu)}^{d,l} f)
\D^{\lambda -\mu} g_{\kappa, \nu}^{d,m} \biggr|^q dx = \\
\sum_{n \in \Z^d: Q_{\kappa,n}^d \cap G_\kappa^{d,m,D} \ne \emptyset}
\int_{Q_{\kappa,n}^d} \biggl| \sum_{\nu \in N_\kappa^{d,m,D}: Q_{\kappa,n}^d \cap
\supp g_{\kappa, \nu}^{d,m} \ne \emptyset}
(\D^\mu S_{\kappa, \nu_\kappa^D(\nu)}^{d,l} f)
\D^{\lambda -\mu} g_{\kappa, \nu}^{d,m} \biggr|^q dx = \\
\sum_{n \in \Z^d: Q_{\kappa,n}^d \cap G_\kappa^{d,m,D} \ne \emptyset}
\biggl\| \sum_{\nu \in N_\kappa^{d,m,D}: Q_{\kappa,n}^d \cap
\supp g_{\kappa, \nu}^{d,m} \ne \emptyset}
(\D^\mu S_{\kappa, \nu_\kappa^D(\nu)}^{d,l} f)
\D^{\lambda -\mu} g_{\kappa, \nu}^{d,m} \biggr\|_{L_q(Q_{\kappa,n}^d)}^q \le \\
\sum_{n \in \Z^d: Q_{\kappa,n}^d \cap G_\kappa^{d,m,D} \ne \emptyset}
\biggl(\sum_{\nu \in N_\kappa^{d,m,D}: Q_{\kappa,n}^d \cap
\supp g_{\kappa, \nu}^{d,m} \ne \emptyset}
\| (\D^\mu S_{\kappa, \nu_\kappa^D(\nu)}^{d,l} f)
\D^{\lambda -\mu} g_{\kappa, \nu}^{d,m} \|_{L_q(Q_{\kappa,n}^d)}\biggr)^q.
\end{multline*}

Для оценки правой части (1.3.26), используя сначала (1.2.5), а, затем
применяя (1.1.1), при $ n \in \Z^d: Q_{\kappa,n}^d \cap G_\kappa^{d,m,D} \ne
\emptyset, \nu \in N_\kappa^{d,m,D}: Q_{\kappa,n}^d \cap
\supp g_{\kappa, \nu}^{d,m} \ne \emptyset, $ выводим
\begin{multline*} \tag{1.3.27}
\| (\D^\mu S_{\kappa, \nu_\kappa^D(\nu)}^{d,l} f)
\D^{\lambda -\mu} g_{\kappa, \nu}^{d,m} \|_{L_q(Q_{\kappa,n}^d)} \le 
\| \D^{\lambda -\mu} g_{\kappa, \nu}^{d,m} \|_{L_\infty(\R^d)}
\| \D^\mu S_{\kappa, \nu_\kappa^D(\nu)}^{d,l} f \|_{L_q(Q_{\kappa,n}^d)} = \\
c_{18} 2^{(\kappa, \lambda -\mu)}
\| \D^\mu S_{\kappa, \nu_\kappa^D(\nu)}^{d,l} f \|_{L_q(Q_{\kappa,n}^d)} \le 
c_{19} 2^{(\kappa, \lambda -\mu)} 2^{(\kappa, \mu +s^{-1} \e -q^{-1} \e)}
\| S_{\kappa, \nu_\kappa^D(\nu)}^{d,l} f\|_{L_s(Q_{\kappa,n}^d)} = \\
c_{19} 2^{(\kappa, \lambda +s^{-1} \e -q^{-1} \e)}
\| S_{\kappa, \nu_\kappa^D(\nu)}^{d,l} f\|_{L_s(Q_{\kappa,n}^d)}.
\end{multline*}

Далее, принимая во внимание, что в силу (1.3.13) имеет место включение
\begin{multline*} 
Q_{\kappa,n}^d \subset (2^{-\kappa} \nu_\kappa^D(\nu) +(\gamma^1 +\e) 2^{-\kappa} B^d), \\
n \in \Z^d: Q_{\kappa,n}^d \cap G_\kappa^{d,m,D} \ne \emptyset,
\nu \in N_\kappa^{d,m,D}: \supp g_{\kappa, \nu}^{d,m} \cap Q_{\kappa,n}^d
\ne \emptyset,
\end{multline*} 
на основании (1.1.1), (1.1.3), (1.3.16) заключаем, что для
$ n \in \Z^d: Q_{\kappa,n}^d \cap G_\kappa^{d,m,D} \ne \emptyset,
\nu \in N_\kappa^{d,m,D}: \supp g_{\kappa, \nu}^{d,m} \cap Q_{\kappa,n}^d
\ne \emptyset, $ выполняется неравенство
\begin{multline*} \tag{1.3.28}
\| S_{\kappa, \nu_\kappa^D(\nu)}^{d,l} f\|_{L_s(Q_{\kappa,n}^d)} \le
c_{20} \| S_{\kappa, \nu_\kappa^D(\nu)}^{d,l} f
\|_{L_s(Q_{\kappa, \nu_\kappa^D(\nu)}^d)} \le \\
c_{21} \| f \|_{L_s(Q_{\kappa, \nu_\kappa^D(\nu)}^d)} \le c_{21}
\| f\|_{L_s(D \cap D_{\kappa, n}^{d,m,D,\alpha})}.
\end{multline*}

Объединяя (1.3.27) и (1.3.28), находим, что при $ n \in \Z^d: Q_{\kappa,n}^d
\cap G_\kappa^{d,m,D} \ne \emptyset,
\nu \in N_\kappa^{d,m,D}: Q_{\kappa,n}^d \cap \supp g_{\kappa, \nu}^{d,m} \ne
\emptyset, $ имеет место неравенство
$$
\| (\D^\mu S_{\kappa, \nu_\kappa^D(\nu)}^{d,l} f)
\D^{\lambda -\mu} g_{\kappa, \nu}^{d,m} \|_{L_q(Q_{\kappa,n}^d)} \le 
c_{22} 2^{(\kappa, \lambda +s^{-1} \e -q^{-1} \e)}
\| f\|_{L_s(D \cap D_{\kappa,n}^{d,m,D,\alpha})}.
$$

Подставляя эту оценку в (1.3.26) и применяя (1.3.17), а, затем используя
неравенство Г\"eльдера с показателем $ s /q \le 1 $ и оценку (1.3.15), приходим к
неравенству
\begin{multline*} \tag{1.3.29}
\biggl\| \sum_{\nu \in N_\kappa^{d,m,D}}
(\D^\mu S_{\kappa, \nu_\kappa^D(\nu)}^{d,l} f)
\D^{\lambda -\mu} g_{\kappa, \nu}^{d,m} \biggr\|_{L_q(\R^d)}^q \le \\
\sum_{n \in \Z^d: Q_{\kappa,n}^d \cap G_\kappa^{d,m,D} \ne \emptyset}
\biggl(\sum_{\substack{\nu \in N_\kappa^{d,m,D}: \\ Q_{\kappa,n}^d \cap
\supp g_{\kappa, \nu}^{d,m} \ne \emptyset}}
c_{22} 2^{(\kappa, \lambda +s^{-1} \e -q^{-1} \e)}
\| f\|_{L_s(D \cap D_{\kappa,n}^{d,m,D,\alpha})}\biggr)^q \le \\
(c_{22} 2^{(\kappa, \lambda +s^{-1} \e -q^{-1} \e)})^q
\sum_{n \in \Z^d: Q_{\kappa,n}^d \cap G_\kappa^{d,m,D} \ne \emptyset}
\biggl( c_5 \| f\|_{L_s(D \cap D_{\kappa,n}^{d,m,D,\alpha})}\biggr)^q = \\
(c_{23} 2^{(\kappa, \lambda +s^{-1} \e -q^{-1} \e)})^q
\sum_{n \in \Z^d: Q_{\kappa,n}^d \cap G_\kappa^{d,m,D} \ne \emptyset}
\biggl( \int_{D \cap D_{\kappa,n}^{d,m,D,\alpha}} | f(x)|^s dx\biggr)^{q /s} \le \\
(c_{23} 2^{(\kappa, \lambda +s^{-1} \e -q^{-1} \e)})^q
\biggl(\sum_{n \in \Z^d: Q_{\kappa,n}^d \cap G_\kappa^{d,m,D} \ne \emptyset}
\int_{D \cap D_{\kappa,n}^{d,m,D,\alpha}} | f(x)|^s dx\biggr)^{q /s} = \\
(c_{23} 2^{(\kappa, \lambda +s^{-1} \e -q^{-1} \e)})^q
\biggl(\sum_{n \in \Z^d: Q_{\kappa,n}^d \cap G_\kappa^{d,m,D} \ne \emptyset}
\int_D \chi_{ D_{\kappa,n}^{d,m,D,\alpha}}(x) | f(x)|^s dx\biggr)^{q /s} = \\
(c_{23} 2^{(\kappa, \lambda +s^{-1} \e -q^{-1} \e)})^q
\biggl(\int_D \biggl(\sum_{n \in \Z^d: Q_{\kappa,n}^d \cap G_\kappa^{d,m,D} \ne \emptyset}
\chi_{ D_{\kappa,n}^{d,m,D,\alpha}}(x)\biggr)  | f(x)|^s dx\biggr)^{q /s} \le \\
(c_{23} 2^{(\kappa, \lambda +s^{-1} \e -q^{-1} \e)})^q
\biggl(\int_D c_4 | f(x)|^s dx\biggr)^{q /s} = 
(c_{24} 2^{(\kappa, \lambda +s^{-1} \e -q^{-1} \e)} \| f\|_{L_s(D)})^q.
\end{multline*}

Соединяя (1.3.25) с (1.3.29). получаем
\begin{multline*}
\| \D^\lambda E_k^{d,l,m,D,\alpha} f \|_{L_q(\R^d)} \le
\sum_{ \mu \in \Z_+^d(\lambda)} C_\lambda^\mu
c_{24} 2^{(\kappa, \lambda +s^{-1} \e -q^{-1} \e)} \| f\|_{L_s(D)} = \\
c_{25} 2^{(\kappa, \lambda +s^{-1} \e -q^{-1} \e)} \| f\|_{L_s(D)} \le 
c_{17} 2^{k (\alpha^{-1}, \lambda +(s^{-1} -q^{-1}) \e)} \| f\|_{L_s(D)},
\end{multline*}
что совпадает с (1.3.24) при $ s \le q. $ Справедливость (1.3.24) при $ q < s $
вытекает из соблюдения (1.3.24) при $ q = s $ и неравенства
\begin{multline*}
\| \D^\lambda E_k^{d,l,m,D,\alpha} f \|_{L_q(\R^d)} =
\| \D^\lambda E_k^{d,l,m,D,\alpha} f \|_{L_q(G)} \le \\
(\mes G)^{1 /q -1 /s} \| \D^\lambda E_k^{d,l,m,D,\alpha} f \|_{L_s(G)} \le
(\mes G)^{1 /q -1 /s} \| \D^\lambda E_k^{d,l,m,D,\alpha} f \|_{L_s(\R^d)},
\end{multline*}
где $ G $ -- ограниченная область в $ \R^d $ такая, что $ G_\kappa^{d,m,D}
\subset G $ при $ \kappa = \kappa(k,\alpha), k \in \Z_+: k \ge k^0. \square $
Отметим, что в лемме 1.3.3 при $ s \le q $ и в предложении 1.3.1 условие
ограниченности области $ D $ можно опустить.

Введем следующее обозначение. Пусть $ d \in \N, \alpha \in \R_+^d, D $ --
ограниченная область в $ \R^d $ и число $ k \in \N $ таковы, что при
$ \kappa^\prime = \kappa(k -1,\alpha) $ существует $ \nu^\prime \in \Z^d, $ для
которого $ Q_{\kappa^\prime,\nu^\prime}^d \subset D, $ а также пусть
$ l \in \Z_+^d, m \in \N^d. $ Тогда положим
\begin{equation*} \tag{1.3.30}
\mathcal E_k^{d,l,m,D,\alpha} = E_\kappa^{d,l,m,D} -
H_{\kappa, \kappa^\prime}^{d,l,m,D} E_{\kappa^\prime}^{d,l,m,D},
\text{ при } \kappa = \kappa(k,\alpha), \kappa^\prime = \kappa(k -1,\alpha).
\end{equation*}

При $ d \in \n, \alpha \in \R_+^d $ будем говорить, что область
$ D \subset \R^d $ является областью $ \alpha $-типа, если существуют константы
$ K^0(D,\alpha) \in \Z_+, \Gamma^0(D,\alpha) \in \R_+^d,
c_{26}(D,\alpha) >0 $ такие, что для $ k \in \Z_+: k \ge K^0, $ при $ \kappa =
\kappa(k,\alpha) $ выполняются следующие условия:

1) для каждого $ \nu \in \Z^d: Q_{\kappa,\nu}^d \cap D \ne \emptyset, $
существует $ \nu^\prime \in \Z^d $ такой, что соблюдается включение
\begin{equation*}
Q_{\kappa,\nu^\prime}^d \subset D \cap (2^{-\kappa} \nu +\Gamma^0 2^{-\kappa} B^d);
\end{equation*}

2) для любых $ \nu, \nu^\prime \in \Z^d $ таких, что $ Q_{\kappa,\nu}^d \subset D,
Q_{\kappa,\nu^\prime}^d \subset D, $ существуют последовательнности
$ \nu^\iota \in \Z^d, j^\iota \in \Nu_{1,d}^1, \epsilon^\iota \in \{-1,1\},
\iota =0,\ldots,\Iota, $ со следующими свойствами:

\begin{equation*} \tag{1.3.31}
\Iota \le c_{26} \|\nu -\nu^\prime\|_{l_\infty^d};
\end{equation*}

\begin{equation*} \tag{1.3.32}
\nu^0 = \nu, \nu^{\Iota} = \nu^\prime, \nu^{\iota +1} = \nu^\iota +\epsilon^\iota
e_{j^\iota}, \iota =0,\ldots,\Iota -1;
\end{equation*}

при $ \iota =0,\ldots,\Iota -1 $ справедливо включение
\begin{equation*} \tag{1.3.33}
Q^\iota = \inter (\overline Q_{\kappa,\nu^\iota}^d \cup
\overline Q_{\kappa,\nu^{\iota +1}}^d) \subset D.
\end{equation*}

Предложение 1.3.4

Пусть $ d \in \N, l \in \N^d, m \in \N^d, \lambda \in \Z_+^d(m),
\alpha \in \R_+^d, D \subset \R^d $ -- ограниченная область $ \alpha $-типа,
$ 1 \le p < \infty, 1 \le q \le \infty. $
Тогда существуют константы $ c_{27}(d,l,m,D,\alpha,\lambda,p,q) > 0,
c_{28}(d,m,D,\alpha) > 0 $ такие, что при $ k \in \Z_+: k > K^0(D,\alpha) $
(см. определение области $ \alpha $-типа), для $ f \in L_p(D) $ соблюдается
неравенство
\begin{multline*} \tag{1.3.34}
\| \D^\lambda \mathcal E_k^{d,l -\e,m,D,\alpha} f \|_{L_q(\R^d)}
\le c_{27} 2^{k (\alpha^{-1}, \lambda  +(p^{-1} -q^{-1})_+ \e)}
\sum_{j =1}^d \Omega_j^{\prime l_j}(f, c_{28} 2^{-k /\alpha_j})_{L_p(D)}.
\end{multline*}

Доказательство.

В условиях предложения пусть $ f \in L_p(D), k \in \Z_+: k > K^0, p \le q. $
Тогда, полагая $ \kappa = \kappa(k,\alpha), \kappa^\prime = \kappa(k -1,\alpha),
\k = \kappa -\kappa^\prime, $ и принимая во внимание (1.2.11), (1.2.12),
(1.2.13), (1.2.20), (1.3.23), имеем
\begin{multline*} \tag{1.3.35}
\| \D^\lambda \mathcal E_k^{d,l -\e,m,D,\alpha} f \|_{L_q(\R^d)} =
\| \D^\lambda (E_\kappa^{d,l -\e,m,D} f -
H_{\kappa, \kappa^\prime}^{d,l -\e,m,D} (E_{\kappa^\prime}^{d,l -\e,m,D} f)) \|_{L_q(\R^d)} = \\
\biggl\| \D^\lambda \biggl(\sum_{\nu \in N_\kappa^{d,m,D}}
(S_{\kappa, \nu_\kappa^D(\nu)}^{d,l -\e} f) g_{\kappa, \nu}^{d,m} -
H_{\kappa, \kappa^\prime}^{d,l -\e,m,D} \biggl(\sum_{\nu^\prime \in N_{\kappa^\prime}^{d,m,D}}
(S_{\kappa^\prime, \nu_{\kappa^\prime}^D(\nu^\prime)}^{d,l -\e} f)
g_{\kappa^\prime, \nu^\prime}^{d,m}\biggr)\biggr)  \biggr\|_{L_q(\R^d)} = \\
\biggl\| \D^\lambda \biggl(\sum_{\nu \in N_\kappa^{d,m,D}}
(S_{\kappa, \nu_\kappa^D(\nu)}^{d,l -\e} f) g_{\kappa, \nu}^{d,m} -
\sum_{\nu \in N_{\kappa}^{d,m,D}}
\biggl(\sum_{\m^{\k} \in \M_{\k}^m(\nu)} A_{\m^{\k}}^m
S_{\kappa^\prime, \nu_{\kappa^\prime}^D(\n_{\k}(\nu,\m^{\k}))}^{d,l -\e} f\biggr)
g_{\kappa, \nu}^{d,m}\biggr) \biggr\|_{L_q(\R^d)} = \\
\biggl\| \D^\lambda \biggl(\sum_{\nu \in N_\kappa^{d,m,D}}
(\sum_{\m^{\k} \in \M_{\k}^m(\nu)} A_{\m^{\k}}^m)
(S_{\kappa, \nu_\kappa^D(\nu)}^{d,l -\e} f) g_{\kappa, \nu}^{d,m} - \\
\sum_{\nu \in N_{\kappa}^{d,m,D}}
(\sum_{\m^{\k} \in \M_{\k}^m(\nu)} A_{\m^{\k}}^m
S_{\kappa^\prime, \nu_{\kappa^\prime}^D(\n_{\k}(\nu,\m^{\k}))}^{d,l -\e} f)
g_{\kappa, \nu}^{d,m}\biggr) \biggr\|_{L_q(\R^d)} = \\
\biggl\| \D^\lambda \biggl(\sum_{\nu \in N_\kappa^{d,m,D}}
((\sum_{\m^{\k} \in \M_{\k}^m(\nu)} A_{\m^{\k}}^m
S_{\kappa, \nu_\kappa^D(\nu)}^{d,l -\e} f) -
(\sum_{\m^{\k} \in \M_{\k}^m(\nu)} A_{\m^{\k}}^m
S_{\kappa^\prime, \nu_{\kappa^\prime}^D(\n_{\k}(\nu,\m^{\k}))}^{d,l -\e} f))
g_{\kappa, \nu}^{d,m}\biggr) \biggr\|_{L_q(\R^d)} = \\
\biggl\| \D^\lambda \biggl(\sum_{\nu \in N_\kappa^{d,m,D}}
\biggl(\sum_{\m^{\k} \in \M_{\k}^m(\nu)} A_{\m^{\k}}^m
((S_{\kappa, \nu_\kappa^D(\nu)}^{d,l -\e} f) -
(S_{\kappa^\prime, \nu_{\kappa^\prime}^D(\n_{\k}(\nu,\m^{\k}))}^{d,l -\e} f))\biggr)
g_{\kappa, \nu}^{d,m}\biggr) \biggr\|_{L_q(\R^d)} = \\
\biggl\| \sum_{\nu \in N_\kappa^{d,m,D}} \D^\lambda
\biggl(\biggl(\sum_{\m^{\k} \in \M_{\k}^m(\nu)} A_{\m^{\k}}^m
((S_{\kappa, \nu_\kappa^D(\nu)}^{d,l -\e} f) -
(S_{\kappa^\prime, \nu_{\kappa^\prime}^D(\n_{\k}(\nu,\m^{\k}))}^{d,l -\e} f))\biggr)
g_{\kappa, \nu}^{d,m}\biggr) \biggr\|_{L_q(\R^d)} = \\
\biggl\| \sum_{\nu \in N_\kappa^{d,m,D}} \sum_{ \mu \in \Z_+^d(\lambda)}
C_\lambda^\mu \D^\mu \biggl(\sum_{\m^{\k} \in \M_{\k}^m(\nu)} A_{\m^{\k}}^m
((S_{\kappa, \nu_\kappa^D(\nu)}^{d,l -\e} f) -
(S_{\kappa^\prime, \nu_{\kappa^\prime}^D(\n_{\k}(\nu,\m^{\k}))}^{d,l -\e} f))\biggr)
\D^{\lambda -\mu} g_{\kappa, \nu}^{d,m} \biggr\|_{L_q(\R^d)} = \\
\biggl\| \sum_{ \mu \in \Z_+^d(\lambda)} C_\lambda^\mu
\sum_{\nu \in N_\kappa^{d,m,D}}
\D^\mu \biggl(\sum_{\m^{\k} \in \M_{\k}^m(\nu)} A_{\m^{\k}}^m
((S_{\kappa, \nu_\kappa^D(\nu)}^{d,l -\e} f) -
(S_{\kappa^\prime, \nu_{\kappa^\prime}^D(\n_{\k}(\nu,\m^{\k}))}^{d,l -\e} f))\biggr)
\D^{\lambda -\mu} g_{\kappa, \nu}^{d,m} \biggr\|_{L_q(\R^d)} = \\
\biggl\| \sum_{ \mu \in \Z_+^d(\lambda)} C_\lambda^\mu
\sum_{\nu \in N_\kappa^{d,m,D}}
\biggl(\sum_{\m^{\k} \in \M_{\k}^m(\nu)} A_{\m^{\k}}^m
\D^\mu (S_{\kappa, \nu_\kappa^D(\nu)}^{d,l -\e} f -
S_{\kappa^\prime, \nu_{\kappa^\prime}^D(\n_{\k}(\nu,\m^{\k}))}^{d,l -\e} f)\biggr)
\D^{\lambda -\mu} g_{\kappa, \nu}^{d,m} \biggr\|_{L_q(\R^d)} \le \\
\sum_{ \mu \in \Z_+^d(\lambda)} C_\lambda^\mu
\biggl\| \sum_{\nu \in N_\kappa^{d,m,D}}
\biggl(\sum_{\m^{\k} \in \M_{\k}^m(\nu)} A_{\m^{\k}}^m
\D^\mu (S_{\kappa, \nu_\kappa^D(\nu)}^{d,l -\e} f -
S_{\kappa^\prime, \nu_{\kappa^\prime}^D(\n_{\k}(\nu,\m^{\k}))}^{d,l -\e} f)\biggr)
\D^{\lambda -\mu} g_{\kappa, \nu}^{d,m} \biggr\|_{L_q(\R^d)}.
\end{multline*}

Оценивая правую часть (1.3.35), при $ \mu \in \Z_+^d(\lambda) $ с учетом
(1.3.6), (1.3.4) получаем (ср. с (1.3.19))
\begin{multline*} \tag{1.3.36}
\biggl\| \sum_{\nu \in N_\kappa^{d,m,D}}
\biggl(\sum_{\m^{\k} \in \M_{\k}^m(\nu)} A_{\m^{\k}}^m
\D^\mu (S_{\kappa, \nu_\kappa^D(\nu)}^{d,l -\e} f -
S_{\kappa^\prime, \nu_{\kappa^\prime}^D(\n_{\k}(\nu,\m^{\k}))}^{d,l -\e} f)\biggr)
\D^{\lambda -\mu} g_{\kappa, \nu}^{d,m} \biggr\|_{L_q(\R^d)}^q = \\
\int_{\R^d} \biggl| \sum_{\nu \in N_\kappa^{d,m,D}}
\biggl(\sum_{\m^{\k} \in \M_{\k}^m(\nu)} A_{\m^{\k}}^m
\D^\mu (S_{\kappa, \nu_\kappa^D(\nu)}^{d,l -\e} f -
S_{\kappa^\prime, \nu_{\kappa^\prime}^D(\n_{\k}(\nu,\m^{\k}))}^{d,l -\e} f)\biggr)
\D^{\lambda -\mu} g_{\kappa, \nu}^{d,m} \biggr|^q dx = \\
\int_{G_\kappa^{d,m,D}} \biggl| \sum_{\nu \in N_\kappa^{d,m,D}}
\biggl(\sum_{\m^{\k} \in \M_{\k}^m(\nu)} A_{\m^{\k}}^m
\D^\mu (S_{\kappa, \nu_\kappa^D(\nu)}^{d,l -\e} f -
S_{\kappa^\prime, \nu_{\kappa^\prime}^D(\n_{\k}(\nu,\m^{\k}))}^{d,l -\e} f)\biggr)
\D^{\lambda -\mu} g_{\kappa, \nu}^{d,m} \biggr|^q dx = \\
\sum_{n \in \Z^d: Q_{\kappa,n}^d \cap G_\kappa^{d,m,D} \ne \emptyset}
\int_{Q_{\kappa,n}^d} \biggl| \sum_{\nu \in N_\kappa^{d,m,D}}
\biggl(\sum_{\m^{\k} \in \M_{\k}^m(\nu)} A_{\m^{\k}}^m
\D^\mu (S_{\kappa, \nu_\kappa^D(\nu)}^{d,l -\e} f -\\
S_{\kappa^\prime, \nu_{\kappa^\prime}^D(\n_{\k}(\nu,\m^{\k}))}^{d,l -\e} f)\biggr)
\D^{\lambda -\mu} g_{\kappa, \nu}^{d,m} \biggr|^q dx = \\
\sum_{n \in \Z^d: Q_{\kappa,n}^d \cap G_\kappa^{d,m,D} \ne \emptyset}
\int_{Q_{\kappa,n}^d} \biggl| \sum_{\substack{\nu \in N_\kappa^{d,m,D}:\\ Q_{\kappa,n}^d \cap
\supp g_{\kappa, \nu}^{d,m} \ne \emptyset}}
\biggl(\sum_{\m^{\k} \in \M_{\k}^m(\nu)} A_{\m^{\k}}^m
\D^\mu (S_{\kappa, \nu_\kappa^D(\nu)}^{d,l -\e} f -\\
S_{\kappa^\prime, \nu_{\kappa^\prime}^D(\n_{\k}(\nu,\m^{\k}))}^{d,l -\e} f)\biggr)
\D^{\lambda -\mu} g_{\kappa, \nu}^{d,m} \biggr|^q dx = \\
\sum_{n \in \Z^d: Q_{\kappa,n}^d \cap G_\kappa^{d,m,D} \ne \emptyset}
\biggl\| \sum_{\substack{\nu \in N_\kappa^{d,m,D}: \\ Q_{\kappa,n}^d \cap
\supp g_{\kappa, \nu}^{d,m} \ne \emptyset}}
\biggl(\sum_{\m^{\k} \in \M_{\k}^m(\nu)} A_{\m^{\k}}^m
\D^\mu (S_{\kappa, \nu_\kappa^D(\nu)}^{d,l -\e} f -\\
S_{\kappa^\prime, \nu_{\kappa^\prime}^D(\n_{\k}(\nu,\m^{\k}))}^{d,l -\e} f)\biggr)
\D^{\lambda -\mu} g_{\kappa, \nu}^{d,m} \biggr\|_{L_q(Q_{\kappa,n}^d)}^q \le \\
\sum_{n \in \Z^d: Q_{\kappa,n}^d \cap G_\kappa^{d,m,D} \ne \emptyset}
\biggl(\sum_{\substack{\nu \in N_\kappa^{d,m,D}: \\ Q_{\kappa,n}^d \cap
\supp g_{\kappa, \nu}^{d,m} \ne \emptyset}}
\biggl\|\biggl(\sum_{\m^{\k} \in \M_{\k}^m(\nu)} A_{\m^{\k}}^m
\D^\mu (S_{\kappa, \nu_\kappa^D(\nu)}^{d,l -\e} f -\\
S_{\kappa^\prime, \nu_{\kappa^\prime}^D(\n_{\k}(\nu,\m^{\k}))}^{d,l -\e} f)\biggr)
\D^{\lambda -\mu} g_{\kappa, \nu}^{d,m} \biggr\|_{L_q(Q_{\kappa,n}^d)}\biggr)^q \le \\
\sum_{n \in \Z^d: Q_{\kappa,n}^d \cap G_\kappa^{d,m,D} \ne \emptyset}
\biggl(\sum_{\substack{\nu \in N_\kappa^{d,m,D}: \\ Q_{\kappa,n}^d \cap
\supp g_{\kappa, \nu}^{d,m} \ne \emptyset}}
\sum_{\m^{\k} \in \M_{\k}^m(\nu)} A_{\m^{\k}}^m
\biggl\| \D^\mu (S_{\kappa, \nu_\kappa^D(\nu)}^{d,l -\e} f -\\
S_{\kappa^\prime, \nu_{\kappa^\prime}^D(\n_{\k}(\nu,\m^{\k}))}^{d,l -\e} f)
\D^{\lambda -\mu} g_{\kappa, \nu}^{d,m} \biggr\|_{L_q(Q_{\kappa,n}^d)}\biggr)^q.
\end{multline*}

Для оценки правой части (1.3.36), используя сначала (1.2.5), а,
затем применяя (1.1.1), при $ \mu \in \Z_+^d(\lambda), n \in \Z^d: Q_{\kappa,n}^d \cap
G_\kappa^{d,m,D} \ne \emptyset, \nu \in N_\kappa^{d,m,D}: Q_{\kappa,n}^d \cap
\supp g_{\kappa, \nu}^{d,m} \ne \emptyset, \m^{\k} \in \M_{\k}^m(\nu) $
выводим
\begin{multline*} \tag{1.3.37}
\| \D^\mu (S_{\kappa, \nu_\kappa^D(\nu)}^{d,l -\e} f -
S_{\kappa^\prime, \nu_{\kappa^\prime}^D(\n_{\k}(\nu,\m^{\k}))}^{d,l -\e} f)
\D^{\lambda -\mu} g_{\kappa, \nu}^{d,m} \|_{L_q(Q_{\kappa,n}^d)} \le \\
\| \D^{\lambda -\mu} g_{\kappa, \nu}^{d,m} \|_{L_\infty(\R^d)}
\| \D^\mu (S_{\kappa, \nu_\kappa^D(\nu)}^{d,l -\e} f -
S_{\kappa^\prime, \nu_{\kappa^\prime}^D(\n_{\k}(\nu,\m^{\k}))}^{d,l -\e} f)\|_{L_q(Q_{\kappa,n}^d)} = \\
c_{29} 2^{(\kappa, \lambda -\mu)}
\| \D^\mu (S_{\kappa, \nu_\kappa^D(\nu)}^{d,l -\e} f -
S_{\kappa^\prime, \nu_{\kappa^\prime}^D(\n_{\k}(\nu,\m^{\k}))}^{d,l -\e} f)\|_{L_q(Q_{\kappa,n}^d)} \le \\
c_{30} 2^{(\kappa, \lambda -\mu)} 2^{(\kappa, \mu +p^{-1} \e -q^{-1} \e)}
\| S_{\kappa, \nu_\kappa^D(\nu)}^{d,l -\e} f -
S_{\kappa^\prime, \nu_{\kappa^\prime}^D(\n_{\k}(\nu,\m^{\k}))}^{d,l -\e} f\|_{L_p(Q_{\kappa,n}^d)} = \\
c_{30} 2^{(\kappa, \lambda +p^{-1} \e -q^{-1} \e)}
\| S_{\kappa, \nu_\kappa^D(\nu)}^{d,l -\e} f -
S_{\kappa^\prime, \nu_{\kappa^\prime}^D(\n_{\k}(\nu,\m^{\k}))}^{d,l -\e} f\|_{L_p(Q_{\kappa,n}^d)}.
\end{multline*}

Оценим норму в правой части (1.3.37). Для этого, фиксировав $ n \in \Z^d:
Q_{\kappa,n}^d \cap G_\kappa^{d,m,D} \ne \emptyset, \nu \in N_\kappa^{d,m,D}:
Q_{\kappa,n}^d \cap \supp g_{\kappa, \nu}^{d,m} \ne \emptyset,
\m^{\k} \in \M_{\k}^m(\nu), $ заметим, что
\begin{multline*} \tag{1.3.38}
D \supset Q_{\kappa^\prime, \nu_{\kappa^\prime}^D(\n_{\k}(\nu,\m^{\k}))}^d =
2^{-\kappa^\prime} \nu_{\kappa^\prime}^D(\n_{\k}(\nu,\m^{\k})) +
2^{-\kappa^\prime} I^d = 
 2^{-\kappa} 2^{\k} \nu_{\kappa^\prime}^D(\n_{\k}(\nu,\m^{\k})) +\\
2^{-\kappa^\prime} I^d \supset 2^{-\kappa} 2^{\k} \nu_{\kappa^\prime}^D(\n_{\k}(\nu,\m^{\k})) +
2^{-\kappa} I^d = Q_{\kappa, 2^{\k} \nu_{\kappa^\prime}^D(\n_{\k}(\nu,\m^{\k}))}^d.
\end{multline*}
Учитывая, что и $ Q_{\kappa, \nu_\kappa^D(\nu)}^d \subset D, $ выберем
последовательности $ \nu^\iota \in \Z^d, j^\iota \in \Nu_{1,d}^1, \epsilon^\iota
\in \{-1, 1\}, \iota =0,\ldots,\Iota, $ для которых $ \nu^0 = \nu_\kappa^D(\nu),
\nu^{\Iota} = 2^{\k} \nu_{\kappa^\prime}^D(\n_{\k}(\nu,\m^{\k})) $ и
соблюдаются соотношения (1.3.31), (1.3.32), (1.3.33). Тогда
\begin{multline*} \tag{1.3.39}
\| S_{\kappa, \nu_\kappa^D(\nu)}^{d,l -\e} f -
S_{\kappa^\prime, \nu_{\kappa^\prime}^D(\n_{\k}(\nu,\m^{\k}))}^{d,l -\e} f\|_{L_p(Q_{\kappa,n}^d)} = \\
\| S_{\kappa, \nu^0}^{d,l -\e} f -
S_{\kappa^\prime, \nu_{\kappa^\prime}^D(\n_{\k}(\nu,\m^{\k}))}^{d,l -\e} f\|_{L_p(Q_{\kappa,n}^d)} = \\
\| S_{\kappa, \nu^0}^{d,l -\e} f -\sum_{\iota =1}^{\Iota}
(S_{\kappa, \nu^\iota}^{d,l -\e} f -S_{\kappa, \nu^\iota}^{d,l -\e} f) -
S_{\kappa^\prime, \nu_{\kappa^\prime}^D(\n_{\k}(\nu,\m^{\k}))}^{d,l -\e} f\|_{L_p(Q_{\kappa,n}^d)} = \\
\| \sum_{\iota =0}^{\Iota -1} (S_{\kappa, \nu^\iota}^{d,l -\e} f -
S_{\kappa, \nu^{\iota +1}}^{d,l -\e} f) +
S_{\kappa, \nu^{\Iota}}^{d,l -\e} f -
S_{\kappa^\prime, \nu_{\kappa^\prime}^D(\n_{\k}(\nu,\m^{\k}))}^{d,l -\e} f\|_{L_p(Q_{\kappa,n}^d)} \le \\
\sum_{\iota =0}^{\Iota -1} \| S_{\kappa, \nu^\iota}^{d,l -\e} f -
S_{\kappa, \nu^{\iota +1}}^{d,l -\e} f\|_{L_p(Q_{\kappa,n}^d)} +
\| S_{\kappa, \nu^{\Iota}}^{d,l -\e} f -
S_{\kappa^\prime, \nu_{\kappa^\prime}^D(\n_{\k}(\nu,\m^{\k}))}^{d,l -\e} f\|_{L_p(Q_{\kappa,n}^d)}.
\end{multline*}

Для проведения оценки слагаемых в правой части (1.3.39) отметим некоторые
полезные для нас факты.
При $ \iota =0,\ldots,\Iota, j =1,\ldots,d $ для $ x \in
\overline Q_{\kappa,\nu^\iota}^d $ ввиду (1.3.32), (1.3.13), (1.3.31)
справедливо неравенство
\begin{multline*} \tag{1.3.40}
| x_j -2^{-\kappa_j} n_j | \\
\le | x_j -2^{-\kappa_j} (\nu^\iota)_j | +
| 2^{-\kappa_j} (\nu^\iota)_j -2^{-\kappa_j} n_j | \le 2^{-\kappa_j} +
| 2^{-\kappa_j} (\nu^\iota)_j -2^{-\kappa_j} n_j | = \\
2^{-\kappa_j} +
| 2^{-\kappa_j} (\nu^\iota)_j -\sum_{i =0}^{\iota -1} (2^{-\kappa_j} (\nu^i)_j
-2^{-\kappa_j} (\nu^i)_j) -2^{-\kappa_j} n_j |=\\
 2^{-\kappa_j} +
| \sum_{i =0}^{\iota -1} (2^{-\kappa_j} (\nu^{i +1})_j
-2^{-\kappa_j} (\nu^i)_j) +2^{-\kappa_j} (\nu^0)_j -2^{-\kappa_j} n_j | \le \\
2^{-\kappa_j} +\sum_{i =0}^{\iota -1} | 2^{-\kappa_j} (\nu^{i +1})_j
-2^{-\kappa_j} (\nu^i)_j| +| 2^{-\kappa_j} (\nu^0)_j -2^{-\kappa_j} n_j | \le \\
2^{-\kappa_j} +\sum_{i =0}^{\iota -1} 2^{-\kappa_j} \| \nu^{i +1} -\nu^i \|_{l_\infty^d}
+| 2^{-\kappa_j} (\nu^0)_j -2^{-\kappa_j} n_j | = \\
2^{-\kappa_j} +\sum_{i =0}^{\iota -1} 2^{-\kappa_j} \| \epsilon^i e_{j^i} \|_{l_\infty^d}
+| 2^{-\kappa_j} (\nu_\kappa^D(\nu))_j  -2^{-\kappa_j} n_j | = \\
2^{-\kappa_j} +\sum_{i =0}^{\iota -1} 2^{-\kappa_j}
+| 2^{-\kappa_j} (\nu_\kappa^D(\nu))_j  -2^{-\kappa_j} n_j | = \\
2^{-\kappa_j} (\iota +1) +| 2^{-\kappa_j} (\nu_\kappa^D(\nu))_j -
2^{-\kappa_j} n_j | \le 
2^{-\kappa_j} (\Iota +1) +\Gamma_j^1 2^{-\kappa_j} \le \\
2^{-\kappa_j} (c_{26} \| \nu_\kappa^D(\nu) -
2^{\k} \nu_{\kappa^\prime}^D(\n_{\k}(\nu,\m^{\k})) \|_{l_\infty^d} +1) +
\Gamma_j^1 2^{-\kappa_j}.
\end{multline*}

Оценивая правую часть (1.3.40), имеем
\begin{multline*} \tag{1.3.41}
\| \nu_\kappa^D(\nu) -
2^{\k} \nu_{\kappa^\prime}^D(\n_{\k}(\nu,\m^{\k})) \|_{l_\infty^d} = \\
\| \nu_\kappa^D(\nu) -\nu +\nu -2^{\k} \n_{\k}(\nu,\m^{\k}) +2^{\k} \n_{\k}(\nu,\m^{\k})
-2^{\k} \nu_{\kappa^\prime}^D(\n_{\k}(\nu,\m^{\k})) \|_{l_\infty^d} \le \\
\| \nu_\kappa^D(\nu) -\nu \|_{l_\infty^d} +\| \nu -
2^{\k} \n_{\k}(\nu,\m^{\k}) \|_{l_\infty^d} +\| 2^{\k} \n_{\k}(\nu,\m^{\k})
-2^{\k} \nu_{\kappa^\prime}^D(\n_{\k}(\nu,\m^{\k})) \|_{l_\infty^d},
\end{multline*}
Замечая, что
\begin{equation*} \tag{1.3.42}
\k_j = [k / \alpha_j] -[(k -1) / \alpha_j] \le 1 +1 / \alpha_j, j = 1,\ldots,d,
\end{equation*}
получаем, что для $ \nu \in N_\kappa^{d,m,D}, \m^{\k} \in \M_{\k}^m(\nu) $
при $ j \in \Nu_{1,d}^1 \setminus \s(\k) $ имеет место равенство
\begin{equation*}
\nu_j -(2^{\k} \n_{\k}(\nu,\m^{\k}))_j = \nu_j -2^{\k_j} (\n_{\k}(\nu,\m^{\k}))_j =
\nu_j -1 \cdot \nu_j =0,
\end{equation*}
а при $ j \in \s(\k) $ справедливо соотношение
\begin{multline*}
| \nu_j -(2^{\k} \n_{\k}(\nu,\m^{\k}))_j | = | \nu_j -2^{\k_j} (\n_{\k}(\nu,\m^{\k}))_j | =\\
| \nu_j -2^{\k_j} (\nu_j 2^{-\k_j} - \sum_{s_j =1}^{\k_j} \m_{j, \k_j -s_j +1} 2^{-(\k_j -s_j +1)}) | = \\
| \nu_j -\nu_j +\sum_{s_j =1}^{\k_j} \m_{j, \k_j -s_j +1} 2^{s_j -1} | = 
\sum_{s_j =1}^{\k_j} \m_{j, \k_j -s_j +1} 2^{s_j -1} \le \\
\sum_{s_j =1}^{\k_j} (m_j +1) 2^{s_j -1} \le (m_j +1) 2^{\k_j} \le (m_j +1) 2^{1 +1 / \alpha_j},
\end{multline*}
и, значит,
\begin{equation*} \tag{1.3.43}
\| \nu -2^{\k} \n_{\k}(\nu,\m^{\k}) \|_{l_\infty^d} \le \max_{j \in \Nu_{1,d}^1}
(m_j +1) 2^{1 +1 / \alpha_j}.
\end{equation*}
Учитывая, что для $ \nu \in N_\kappa^{d,m,D}, \m^{\k} \in \M_{\k}^m(\nu) $
мультииндекс $ \n_{\k}(\nu,\m^{\k}) \in N_{\kappa^\prime}^{d,m,D}, $ согласно
(1.3.12), (1.3.42) находим, что
\begin{multline*} \tag{1.3.44}
\| 2^{\k} \n_{\k}(\nu,\m^{\k})
-2^{\k} \nu_{\kappa^\prime}^D(\n_{\k}(\nu,\m^{\k})) \|_{l_\infty^d} =
\| 2^{\k} (\n_{\k}(\nu,\m^{\k})
-\nu_{\kappa^\prime}^D(\n_{\k}(\nu,\m^{\k}))) \|_{l_\infty^d} \le \\
\| 2^{\k} \|_{l_\infty^d} \| \n_{\k}(\nu,\m^{\k})
-\nu_{\kappa^\prime}^D(\n_{\k}(\nu,\m^{\k})) \|_{l_\infty^d} \le 
c_3 \max_{j \in \Nu_{1,d}^1} 2^{1 +1 / \alpha_j}.
\end{multline*}
Используя для оценки правой части (1.3.41) неравенства (1.3.12), (1.3.43),
(1.3.44), выводим
\begin{multline*} \tag{1.3.45}
\| \nu_\kappa^D(\nu) -
2^{\k} \nu_{\kappa^\prime}^D(\n_{\k}(\nu,\m^{\k})) \|_{l_\infty^d} \le 
c_3 +(\max_{j \in \Nu_{1,d}^1} (m_j +1) 2^{1 +1 / /\alpha_j}) +c_3
\max_{j \in \Nu_{1,d}^1} 2^{1 +1 / \alpha_j} = \\
c_{31}(d,m,D,\alpha),
\nu \in N_\kappa^{d,m,D}, \m^{\k} \in \M_{\k}^m(\nu).
\end{multline*}
Объединяя (1.3.40), (1.3.45). видим, что для $ \iota =0,\ldots,\Iota $ имеет место включение
\begin{equation*}
\overline Q_{\kappa,\nu^\iota}^d \subset (2^{-\kappa} n +\Gamma^2
2^{-\kappa} B^d)
\end{equation*}
с $ \Gamma^2 = (c_{26} c_{31} +1) \e +\Gamma^1, $
из которого с учетом (1.3.33) следует, что
\begin{equation*} \tag{1.3.46}
Q_{\kappa,\nu^\iota}^d \cup Q_{\kappa,\nu^{\iota +1}}^d \subset Q^\iota \subset
D \cap (2^{-\kappa} n +\Gamma^2 2^{-\kappa} B^d), \iota =0,\ldots,\Iota -1.
\end{equation*}

Заметим еще, что для $ x \in
\overline Q_{\kappa^\prime, \nu_{\kappa^\prime}^D(\n_{\k}(\nu,\m^{\k}))}^d $
при $ j =1,\ldots,d, $ благодаря (1.3.42), (1.3.45), (1.3.13), соблюдается
неравенство
\begin{multline*}
| x_j -2^{-\kappa_j} n_j | = | x_j -2^{-\kappa^\prime_j}
(\nu_{\kappa^\prime}^D(\n_{\k}(\nu,\m^{\k})))_j +2^{-\kappa^\prime_j}
(\nu_{\kappa^\prime}^D(\n_{\k}(\nu,\m^{\k})))_j -2^{-\kappa_j} n_j | \le \\
| x_j -2^{-\kappa^\prime_j} (\nu_{\kappa^\prime}^D(\n_{\k}(\nu,\m^{\k})))_j |
+| 2^{-\kappa^\prime_j} (\nu_{\kappa^\prime}^D(\n_{\k}(\nu,\m^{\k})))_j -
2^{-\kappa_j} n_j | \le \\
2^{-\kappa^\prime_j} +| 2^{-\kappa^\prime_j}
(\nu_{\kappa^\prime}^D(\n_{\k}(\nu,\m^{\k})))_j -2^{-\kappa_j} n_j | = \\
2^{\k_j} 2^{-\kappa_j} +| 2^{-\kappa_j} 2^{\k_j}
(\nu_{\kappa^\prime}^D(\n_{\k}(\nu,\m^{\k})))_j -2^{-\kappa_j} n_j | \le \\
2^{1 +1 / \alpha_j} 2^{-\kappa_j} +2^{-\kappa_j} | 2^{\k_j}
(\nu_{\kappa^\prime}^D(\n_{\k}(\nu,\m^{\k})))_j -n_j | \le \\
2^{1 +1 / \alpha_j} 2^{-\kappa_j} +2^{-\kappa_j} | (2^{\k}
\nu_{\kappa^\prime}^D(\n_{\k}(\nu,\m^{\k})))_j -(\nu_\kappa^D(\nu))_j |
+2^{-\kappa_j} | (\nu_\kappa^D(\nu))_j -n_j | \le\\
 2^{1 +1 / \alpha_j} 2^{-\kappa_j} +2^{-\kappa_j} \| 2^{\k}
\nu_{\kappa^\prime}^D(\n_{\k}(\nu,\m^{\k})) -\nu_\kappa^D(\nu) \|_{l_\infty^d}
+| 2^{-\kappa_j} (\nu_\kappa^D(\nu))_j -2^{-\kappa_j} n_j | \le \\
2^{1 +1 / \alpha_j} 2^{-\kappa_j} +c_{31} 2^{-\kappa_j}
+2^{-\kappa_j} \Gamma^1_j = 2^{-\kappa_j} \Gamma^3_j
\end{multline*}
с $ \Gamma^3_j = 2^{1 +1 / \alpha_j} +c_{31} +\Gamma^1_j, $
т.е.
\begin{equation*} \tag{1.3.47}
Q_{\kappa^\prime, \nu_{\kappa^\prime}^D(\n_{\k}(\nu,\m^{\k}))}^d \subset
\overline Q_{\kappa^\prime, \nu_{\kappa^\prime}^D(\n_{\k}(\nu,\m^{\k}))}^d
\subset (2^{-\kappa} n +\Gamma^3 2^{-\kappa} B^d).
\end{equation*}

Итак, задавая координаты вектора $ \Gamma^4 \in \R_+^d $ соотношением
$ \Gamma^4_j = \max(\Gamma^2_j, \Gamma^3_j), j =1,\ldots,d, $ и обозначая
$$
D_{\kappa,n}^{\prime d,m,D,\alpha} = 2^{-\kappa} n +\Gamma^4 2^{-\kappa} B^d,
n \in \Z^d: Q_{\kappa,n}^d \cap G_\kappa^{d,m,D} \ne \emptyset,
$$
в силу (1.3.46), (1.3.47) имеем
\begin{equation*} \tag{1.3.48}
Q_{\kappa,\nu^\iota}^d \cup Q_{\kappa,\nu^{\iota +1}}^d \subset Q^\iota \subset
D \cap D_{\kappa,n}^{\prime d,m,D,\alpha}, \iota =0,\ldots,\Iota -1.
\end{equation*}
\begin{equation*} \tag{1.3.49}
Q_{\kappa^\prime, \nu_{\kappa^\prime}^D(\n_{\k}(\nu,\m^{\k}))}^d \subset
D \cap D_{\kappa,n}^{\prime d,m,D,\alpha},
\nu \in N_\kappa^{d,m,D}: Q_{\kappa,n}^d \cap \supp g_{\kappa,\nu}^{d,m} \ne
\emptyset, \m^{\k} \in \M_{\k}^m(\nu).
\end{equation*}

Из приведенных определений с учетом того, что $ \Gamma^4 > \e, $ видно, что
при $ n \in \Z^d: Q_{\kappa,n}^d \cap G_\kappa^{d,m,D} \ne \emptyset, $
справедливо включение
\begin{equation*} \tag{1.3.50}
Q_{\kappa, n}^d \subset D_{\kappa, n}^{\prime d,m,D,\alpha},
\end{equation*}
а из (1.3.48), (1.3.50) следует, что
\begin{equation*} \tag{1.3.51}
Q_{\kappa,n}^d \subset (2^{-\kappa} \nu^{\iota} +(\Gamma^4 +\e) 2^{-\kappa} B^d),
\iota =0,\ldots,\Iota.
\end{equation*}

Учитывая (1.3.50), нетрудно видеть, что существует константа
$ c_{32}(d,m,D,\alpha) >0 $ такая, что для каждого $ x \in \R^d $ число
\begin{equation*} \tag{1.3.52}
\card \{ n \in \Z^d: Q_{\kappa,n}^d \cap G_\kappa^{d,m,D} \ne \emptyset,
x \in D_{\kappa,n}^{\prime d,m,D,\alpha} \} \le c_{32}.
\end{equation*}

Отметим еще, что вследствие (1.3.31), (1.3.45) справедлива оценка
\begin{equation*} \tag{1.3.53}
\Iota \le c_{26} \| \nu^0 -\nu^{\Iota} \|_{l_\infty^d} = c_{26} \| \nu_\kappa^D(\nu) -
2^{\k} \nu_{\kappa^\prime}^D(\n_{\k}(\nu,\m^{\k})) \|_{l_\infty^d} \le c_{26}
c_{31} = c_{33}(d,m,D,\alpha).
\end{equation*}

Принимая во внимание, что при $ \iota =0,\ldots,\Iota -1 $ ввиду (1.3.33),
(1.3.32) выполняется равенство
\begin{equation*}
Q^\iota = \begin{cases} 2^{-\kappa} \nu^\iota +2^{-(\kappa -e_{j^\iota})} I^d,
\text{ при } \epsilon^\iota =1; \\
2^{-\kappa} \nu^{\iota +1} +2^{-(\kappa -e_{j^\iota})} I^d,

\text{ при } \epsilon^\iota =-1,
\end{cases}
\end{equation*}
определим линейный оператор $ S^\iota: L_1(Q^\iota) \mapsto \mathcal P^{d,l-\e}, $
полагая
$$
S^\iota = P_{\delta,x^0}^{d,l -\e}
$$
при $ \delta = 2^{-(\kappa -e_{j^\iota})}, $
\begin{equation*}
x^0 = \begin{cases} 2^{-\kappa} \nu^\iota, \text{ при } \epsilon^\iota =1; \\
2^{-\kappa} \nu^{\iota +1}, \text{ при } \epsilon^\iota =-1,
\end{cases}
\iota =0,\ldots,\Iota -1.
\end{equation*}
Теперь проведем оценку слагаемых в правой части (1.3.39).
При $ \iota =0,\ldots,\Iota -1, $ благодаря (1.3.51) применяя (1.1.1),
а, затем используя (1.1.2), (1.1.3), (1.3.48), и наконец, пользуясь (1.1.4)
и учитывая (1.3.48), приходим к неравенству
\begin{multline*} \tag{1.3.54}
\| S_{\kappa, \nu^\iota}^{d,l -\e} f -
S_{\kappa, \nu^{\iota +1}}^{d,l -\e} f\|_{L_p(Q_{\kappa,n}^d)} = \\
\| S_{\kappa, \nu^\iota}^{d,l -\e} f -S^\iota f +S^\iota f -
S_{\kappa, \nu^{\iota +1}}^{d,l -\e} f\|_{L_p(Q_{\kappa,n}^d)} \le \\
\| S_{\kappa, \nu^\iota}^{d,l -\e} f -S^\iota f \|_{L_p(Q_{\kappa,n}^d)}
+\| S^\iota f -S_{\kappa, \nu^{\iota +1}}^{d,l -\e} f\|_{L_p(Q_{\kappa,n}^d)} \le \\
c_{34} \| S_{\kappa, \nu^\iota}^{d,l -\e} f -S^\iota f \|_{L_p(Q_{\kappa,\nu^\iota}^d)}
+c_{34} \| S^\iota f -S_{\kappa, \nu^{\iota +1}}^{d,l -\e} f
\|_{L_p(Q_{\kappa,\nu^{\iota +1}}^d)} = \\
c_{34} \| S_{\kappa, \nu^\iota}^{d,l -\e} (f -S^\iota f)
\|_{L_p(Q_{\kappa,\nu^\iota}^d)}
+c_{34} \| S_{\kappa, \nu^{\iota +1}}^{d,l -\e} (f -S^\iota f)
\|_{L_p(Q_{\kappa,\nu^{\iota +1}}^d)} \le \\
c_{35} \| f -S^\iota f \|_{L_p(Q_{\kappa,\nu^\iota}^d)}
+c_{35} \| f -S^\iota f\|_{L_p(Q_{\kappa,\nu^{\iota +1}}^d)} \le 
c_{36} \| f -S^\iota f \|_{L_p(Q^\iota)} \le \\ c_{37} \sum_{j =1}^d
((2^{-(\kappa -e_{j^\iota})})_j)^{-1/p} \biggl(\int_{(2^{-(\kappa -e_{j^\iota})})_j B^1}
\int_{Q_{l_j \xi e_j}^\iota}| \Delta_{\xi e_j}^{l_j} f(x)|^p dx d\xi\biggr)^{1/p} \le \\
c_{38} \sum_{j =1}^d 2^{\kappa_j /p} \biggl(\int_{2^{-\kappa_j +1} B^1}
\int_{Q_{l_j \xi e_j}^\iota}| \Delta_{\xi e_j}^{l_j} f(x)|^p dx d\xi\biggr)^{1/p} \le \\
c_{38} \sum_{j =1}^d 2^{\kappa_j /p} \biggl(\int_{2^{-\kappa_j +1} B^1}
\int_{(D \cap D_{\kappa,n}^{\prime d,m,D,\alpha})_{l_j \xi e_j}}
| \Delta_{\xi e_j}^{l_j} f(x)|^p dx d\xi\biggr)^{1/p} \le \\
c_{38} \sum_{j =1}^d 2^{\kappa_j /p} \biggl(\int_{2^{-\kappa_j +1} B^1}
\int_{D_{l_j \xi e_j} \cap D_{\kappa,n}^{\prime d,m,D,\alpha}}
| \Delta_{\xi e_j}^{l_j} f(x)|^p dx d\xi\biggr)^{1/p}.\
\end{multline*}

Для оценки последнего слагаемого в правой части (1.3.39), ввиду (1.3.51)
используя (1.1.1), затем пользуясь (1.1.2), (1.1.3), учитывая (1.3.38),
применяя (1.1.4) и принимая во внимание (1.3.49), получаем
\begin{multline*} \tag{1.3.55}
\| S_{\kappa, \nu^{\Iota}}^{d,l -\e} f -
S_{\kappa^\prime, \nu_{\kappa^\prime}^D(\n_{\k}(\nu,\m^{\k}))}^{d,l -\e} f\|_{L_p(Q_{\kappa,n}^d)} \le \\
c_{39} \| S_{\kappa, \nu^{\Iota}}^{d,l -\e} f -
S_{\kappa^\prime, \nu_{\kappa^\prime}^D(\n_{\k}(\nu,\m^{\k}))}^{d,l -\e} f\|_{L_p(Q_{\kappa,\nu^{\Iota}}^d)} = \\
c_{39} \| S_{\kappa, \nu^{\Iota}}^{d,l -\e} (f -
S_{\kappa^\prime, \nu_{\kappa^\prime}^D(\n_{\k}(\nu,\m^{\k}))}^{d,l -\e} f)\|_{L_p(Q_{\kappa,\nu^{\Iota}}^d)} \le \\
c_{40} \| f -S_{\kappa^\prime, \nu_{\kappa^\prime}^D(\n_{\k}(\nu,\m^{\k}))}^{d,l -\e} f
\|_{L_p(Q_{\kappa,\nu^{\Iota}}^d)} \le \\
c_{40} \| f -S_{\kappa^\prime, \nu_{\kappa^\prime}^D(\n_{\k}(\nu,\m^{\k}))}^{d,l -\e} f
\|_{L_p(Q_{\kappa^\prime, \nu_{\kappa^\prime}^D(\n_{\k}(\nu,\m^{\k}))}^d)} \le \\
c_{41} \sum_{j =1}^d 2^{\kappa^\prime_j /p} \biggl(\int_{2^{-\kappa^\prime_j} B^1}
\int_{(Q_{\kappa^\prime, \nu_{\kappa^\prime}^D(\n_{\k}(\nu,\m^{\k}))}^d)_{l_j \xi e_j}}
| \Delta_{\xi e_j}^{l_j} f(x)|^p dx d\xi\biggr)^{1/p} \le \\
c_{41} \sum_{j =1}^d 2^{\kappa_j /p} \biggl(\int_{c_{42} 2^{-\kappa_j} B^1}
\int_{(D \cap D_{\kappa,n}^{\prime d,m,D,\alpha})_{l_j \xi e_j}}
| \Delta_{\xi e_j}^{l_j} f(x)|^p dx d\xi\biggr)^{1/p} \le \\
c_{41} \sum_{j =1}^d 2^{\kappa_j /p} \biggl(\int_{c_{42} 2^{-\kappa_j} B^1}
\int_{D_{l_j \xi e_j} \cap D_{\kappa,n}^{\prime d,m,D,\alpha}}
| \Delta_{\xi e_j}^{l_j} f(x)|^p dx d\xi\biggr)^{1/p},
\end{multline*}
где $ c_{42} = \max_{j =1,\ldots,d} 2^{1 +1 /\alpha_j} > 2. $
Объединяя (1.3.39), (1.3.54), (1.3.55) и учитывая (1.3.53), при $ n \in \Z^d:
Q_{\kappa,n}^d \cap G_\kappa^{d,m,D} \ne \emptyset, \nu \in N_\kappa^{d,m,D}:
Q_{\kappa,n}^d \cap \supp g_{\kappa, \nu}^{d,m} \ne \emptyset,
\m^{\k} \in \M_{\k}^m(\nu) $ имеем
\begin{multline*} \tag{1.3.56}
\| S_{\kappa, \nu_\kappa^D(\nu)}^{d,l -\e} f -
S_{\kappa^\prime, \nu_{\kappa^\prime}^D(\n_{\k}(\nu,\m^{\k}))}^{d,l -\e} f\|_{L_p(Q_{\kappa,n}^d)} \le \\
\Iota c_{38} \sum_{j =1}^d 2^{\kappa_j /p} \biggl(\int_{2^{-\kappa_j +1} B^1}
\int_{D_{l_j \xi e_j} \cap D_{\kappa,n}^{\prime d,m,D,\alpha}}
| \Delta_{\xi e_j}^{l_j} f(x)|^p dx d\xi\biggr)^{1/p} +\\
c_{41} \sum_{j =1}^d 2^{\kappa_j /p} \biggl(\int_{c_{42} 2^{-\kappa_j} B^1}
\int_{D_{l_j \xi e_j} \cap D_{\kappa,n}^{\prime d,m,D,\alpha}}
| \Delta_{\xi e_j}^{l_j} f(x)|^p dx d\xi\biggr)^{1/p} \le \\
c_{33} c_{38} \sum_{j =1}^d 2^{\kappa_j /p} \biggl(\int_{c_{42} 2^{-\kappa_j} B^1}
\int_{D_{l_j \xi e_j} \cap D_{\kappa,n}^{\prime d,m,D,\alpha}}
| \Delta_{\xi e_j}^{l_j} f(x)|^p dx d\xi\biggr)^{1/p} +\\
c_{41} \sum_{j =1}^d 2^{\kappa_j /p} \biggl(\int_{c_{42} 2^{-\kappa_j} B^1}
\int_{D_{l_j \xi e_j} \cap D_{\kappa,n}^{\prime d,m,D,\alpha}}
| \Delta_{\xi e_j}^{l_j} f(x)|^p dx d\xi\biggr)^{1/p} = \\
c_{43} \sum_{j =1}^d 2^{\kappa_j /p} \biggl(\int_{c_{42} 2^{-\kappa_j} B^1}
\int_{D_{l_j \xi e_j} \cap D_{\kappa,n}^{\prime d,m,D,\alpha}}
| \Delta_{\xi e_j}^{l_j} f(x)|^p dx d\xi\biggr)^{1/p}.
\end{multline*}
Из (1.3.37) и (1.3.56) находим, что при $ \mu \in \Z_+^d(\lambda), n \in \Z^d:
Q_{\kappa,n}^d \cap G_\kappa^{d,m,D} \ne \emptyset, \nu \in N_\kappa^{d,m,D}:
Q_{\kappa,n}^d \cap \supp g_{\kappa, \nu}^{d,m} \ne \emptyset,
\m^{\k} \in \M_{\k}^m(\nu) $ выполняется неравенство
\begin{multline*}
\| \D^\mu (S_{\kappa, \nu_\kappa^D(\nu)}^{d,l -\e} f -
S_{\kappa^\prime, \nu_{\kappa^\prime}^D(\n_{\k}(\nu,\m^{\k}))}^{d,l -\e} f)
\D^{\lambda -\mu} g_{\kappa, \nu}^{d,m} \|_{L_q(Q_{\kappa,n}^d)} \le \\
c_{44} 2^{(\kappa, \lambda +p^{-1} \e -q^{-1} \e)}
\sum_{j =1}^d 2^{\kappa_j /p} \biggl(\int_{c_{42} 2^{-\kappa_j} B^1}
\int_{D_{l_j \xi e_j} \cap D_{\kappa,n}^{\prime d,m,D,\alpha}}
| \Delta_{\xi e_j}^{l_j} f(x)|^p dx d\xi\biggr)^{1/p}.
\end{multline*}

Подставляя эту оценку в (1.3.36) и учитывая (1.2.20), (1.3.17), а
затем, применяя неравенство Г\"eльдера с показателем $ q $ и
неравенство Г\"eльдера с показателем $ p/q \le 1, $ и, наконец,
принимая во внимание (1.3.52) и используя неравенство Г\"eльдера с
показателем $ 1/q, $ при $ \mu \in \Z_+^d(\lambda) $ выводим

\begin{multline*} \tag{1.3.57}
\biggl\| \sum_{\nu \in N_\kappa^{d,m,D}}
\biggl(\sum_{\m^{\k} \in \M_{\k}^m(\nu)} A_{\m^{\k}}^m
\D^\mu (S_{\kappa, \nu_\kappa^D(\nu)}^{d,l -\e} f -
S_{\kappa^\prime, \nu_{\kappa^\prime}^D(\n_{\k}(\nu,\m^{\k}))}^{d,l -\e} f)\biggr)
\D^{\lambda -\mu} g_{\kappa, \nu}^{d,m} \biggr\|_{L_q(\R^d)}^q \le \\
\sum_{n \in \Z^d: Q_{\kappa,n}^d \cap G_\kappa^{d,m,D} \ne \emptyset}
\biggl(\sum_{\substack{\nu \in N_\kappa^{d,m,D}:\\ Q_{\kappa,n}^d \cap
\supp g_{\kappa, \nu}^{d,m} \ne \emptyset}}
\sum_{\m^{\k} \in \M_{\k}^m(\nu)} A_{\m^{\k}}^m
c_{44} 2^{(\kappa, \lambda +p^{-1} \e -q^{-1} \e)}
\sum_{j =1}^d 2^{\kappa_j /p}\times\\ 
\biggl(\int_{c_{42} 2^{-\kappa_j} B^1}
\int_{D_{l_j \xi e_j} \cap D_{\kappa,n}^{\prime d,m,D,\alpha}}
| \Delta_{\xi e_j}^{l_j} f(x)|^p dx d\xi\biggr)^{1/p}\biggr)^q \le \\
\sum_{n \in \Z^d: Q_{\kappa,n}^d \cap G_\kappa^{d,m,D} \ne \emptyset}
\biggl(c_5 c_{44} 2^{(\kappa, \lambda +p^{-1} \e -q^{-1} \e)}
\sum_{j =1}^d 2^{\kappa_j /p}\times\\ 
\biggl(\int_{c_{42} 2^{-\kappa_j} B^1}
\int_{D_{l_j \xi e_j} \cap D_{\kappa,n}^{\prime d,m,D,\alpha}}
| \Delta_{\xi e_j}^{l_j} f(x)|^p dx d\xi\biggr)^{1/p}\biggr)^q = \\
(c_{45} 2^{(\kappa, \lambda +p^{-1} \e -q^{-1} \e)})^q
\sum_{n \in \Z^d: Q_{\kappa,n}^d \cap G_\kappa^{d,m,D} \ne \emptyset}
\biggl(\sum_{j =1}^d 2^{\kappa_j /p}\times\\ 
\biggl(\int_{c_{42} 2^{-\kappa_j} B^1}
\int_{D_{l_j \xi e_j} \cap D_{\kappa,n}^{\prime d,m,D,\alpha}}
| \Delta_{\xi e_j}^{l_j} f(x)|^p dx d\xi\biggr)^{1/p}\biggr)^q \le \\
(c_{46} 2^{(\kappa, \lambda +p^{-1} \e -q^{-1} \e)})^q
\sum_{n \in \Z^d: Q_{\kappa,n}^d \cap G_\kappa^{d,m,D} \ne \emptyset}
\sum_{j =1}^d 2^{\kappa_j q /p}\times\\ 
\biggl(\int_{c_{42} 2^{-\kappa_j} B^1}
\int_{D_{l_j \xi e_j} \cap D_{\kappa,n}^{\prime d,m,D,\alpha}}
| \Delta_{\xi e_j}^{l_j} f(x)|^p dx d\xi\biggr)^{q/p} = \\
(c_{46} 2^{(\kappa, \lambda +p^{-1} \e -q^{-1} \e)})^q
\sum_{j =1}^d 2^{\kappa_j q /p}\times\\
\sum_{n \in \Z^d: Q_{\kappa,n}^d \cap G_\kappa^{d,m,D} \ne \emptyset}
\biggl(\int_{c_{42} 2^{-\kappa_j} B^1}
\int_{D_{l_j \xi e_j} \cap D_{\kappa,n}^{\prime d,m,D,\alpha}}
| \Delta_{\xi e_j}^{l_j} f(x)|^p dx d\xi\biggr)^{q/p} \le \\
(c_{46} 2^{(\kappa, \lambda +p^{-1} \e -q^{-1} \e)})^q
\sum_{j =1}^d 2^{\kappa_j q /p}\times\\
\biggl(\sum_{n \in \Z^d: Q_{\kappa,n}^d \cap G_\kappa^{d,m,D} \ne \emptyset}
\int_{c_{42} 2^{-\kappa_j} B^1}
\int_{D_{l_j \xi e_j} \cap D_{\kappa,n}^{\prime d,m,D,\alpha}}
| \Delta_{\xi e_j}^{l_j} f(x)|^p dx d\xi\biggr)^{q/p} = \\
(c_{46} 2^{(\kappa, \lambda +p^{-1} \e -q^{-1} \e)})^q
\sum_{j =1}^d 2^{\kappa_j q /p}\times\\
\biggl(\int_{c_{42} 2^{-\kappa_j} B^1}
\sum_{n \in \Z^d: Q_{\kappa,n}^d \cap G_\kappa^{d,m,D} \ne \emptyset}
\int_{D_{l_j \xi e_j}} \chi_{D_{\kappa,n}^{\prime d,m,D,\alpha}}(x)
| \Delta_{\xi e_j}^{l_j} f(x)|^p dx d\xi\biggr)^{q/p} = \\
(c_{46} 2^{(\kappa, \lambda +p^{-1} \e -q^{-1} \e)})^q
\sum_{j =1}^d 2^{\kappa_j q /p}\times\\
\biggl(\int_{c_{42} 2^{-\kappa_j} B^1}
\int_{D_{l_j \xi e_j}}
\biggl(\sum_{n \in \Z^d: Q_{\kappa,n}^d \cap G_\kappa^{d,m,D} \ne \emptyset}
\chi_{D_{\kappa,n}^{\prime d,m,D,\alpha}}(x)\biggr)
| \Delta_{\xi e_j}^{l_j} f(x)|^p dx d\xi\biggr)^{q/p} \le \\
(c_{46} 2^{(\kappa, \lambda +p^{-1} \e -q^{-1} \e)})^q
\sum_{j =1}^d 2^{\kappa_j q /p}
\biggl(\int_{c_{42} 2^{-\kappa_j} B^1}
\int_{D_{l_j \xi e_j}} c_{32}
| \Delta_{\xi e_j}^{l_j} f(x)|^p dx d\xi\biggr)^{q/p} \le \\
(c_{47} 2^{(\kappa, \lambda +p^{-1} \e -q^{-1} \e)})^q
\biggl(\sum_{j =1}^d 2^{\kappa_j /p}
(\int_{c_{42} 2^{-\kappa_j} B^1}
\int_{D_{l_j \xi e_j}} | \Delta_{\xi e_j}^{l_j} f(x)|^p dx d\xi)^{1/p}\biggr)^q = \\
(c_{47} 2^{(\kappa, \lambda +p^{-1} \e -q^{-1} \e)}
\sum_{j =1}^d 2^{\kappa_j /p}
\biggl(\int_{c_{42} 2^{-\kappa_j} B^1}
\int_{D_{l_j \xi e_j}} | \Delta_{\xi e_j}^{l_j} f(x)|^p dx d\xi\biggr)^{1/p})^q \le \\
( c_{47} 2^{k (\alpha^{-1}, \lambda  +(p^{-1} -q^{-1})_+ \e)}
\sum_{j =1}^d 2^{k /(\alpha_j p)}
\biggl(\int_{c_{28} 2^{-k /\alpha_j} B^1}
\int_{D_{l_j \xi e_j}} | \Delta_{\xi e_j}^{l_j} f(x)|^p dx d\xi\biggr)^{1/p})^q \le \\
(c_{48} 2^{k (\alpha^{-1}, \lambda  +(p^{-1} -q^{-1})_+ \e)}
\sum_{j =1}^d \Omega_j^{\prime l_j}(f, c_{28} 2^{-k /\alpha_j})_{L_p(D)})^q.
\end{multline*}

Соединяя (1.3.35) с (1.3.57), приходим к (1.3.34) при $ p \le q. $
Для получения неравенства (1.3.34) при $ q < p, $ фиксировав ограниченную
область $ D^\prime $ в $ \R^d, $ для которой $ \supp h \subset D^\prime $
для $ h \in \mathcal P_\kappa^{d,l,m,D} $ при $ \kappa \in \Z_+^d, $ благодаря
неравенству Г\"eльдера и соблюдению (1.3.34) при $ q = p, $ имеем
\begin{multline*}
\| \D^\lambda \mathcal E_k^{d,l -\e,m,D,\alpha} f \|_{L_q(\R^d)} =
\| \D^\lambda \mathcal E_k^{d,l -\e,m,D,\alpha} f \|_{L_q(D^\prime)} \le \\
c(D^\prime) \| \D^\lambda \mathcal E_k^{d,l -\e,m,D,\alpha} f \|_{L_p(D^\prime)} = 
c(D^\prime) \| \D^\lambda \mathcal E_k^{d,l -\e,m,D,\alpha} f \|_{L_p(\R^d)} \le \\
c_{27} 2^{k (\alpha^{-1}, \lambda +(p^{-1} -q^{-1})_+ \e)}
\sum_{j =1}^d \Omega_j^{\prime l_j}(f, c_{28} 2^{-k /\alpha_j})_{L_p(D)}, f \in L_p(D). \square
\end{multline*}

Следствие

В условиях предложения 1.3.4 для $ f \in L_p(D) $ при $ k \in \Z_+:  
k > K^0(D,\alpha), $ имеет место неравенство
\begin{multline*} \tag{1.3.58}
\| \D^\lambda ((E_k^{d,l -\e,m,D,\alpha} f) \mid_D)
-\D^\lambda ((E_{k -1}^{d,l -\e,m,D,\alpha} f) \mid_D) \|_{L_q(D)}\\
\le c_{27} 2^{k (\alpha^{-1}, \lambda  +(p^{-1} -q^{-1})_+ \e)}
\sum_{j =1}^d \Omega_j^{\prime l_j}(f, c_{28} 2^{-k /\alpha_j})_{L_p(D)}.
\end{multline*}

В самом деле, при соблюдении условий предложения 1.3.4 для $ f \in L_p(D) $
при $ k \in \Z_+: k > K^0, $ полагая $ \kappa = \kappa(k,\alpha), \kappa^\prime =
\kappa(k -1,\alpha), $ согласно (1.2.10), (1.3.30), (1.3.34) имеем
\begin{multline*}
\| \D^\lambda ((E_k^{d,l -\e,m,D,\alpha} f) \mid_D)
-\D^\lambda ((E_{k -1}^{d,l -\e,m,D,\alpha} f) \mid_D) \|_{L_q(D)} = \\
\| \D^\lambda ((E_\kappa^{d,l -\e,m,D} f) \mid_D)
-\D^\lambda ((E_{\kappa^\prime}^{d,l -\e,m,D} f) \mid_D) \|_{L_q(D)} = \\
\| \D^\lambda ((E_\kappa^{d,l -\e,m,D} f) \mid_D
-(E_{\kappa^\prime}^{d,l -\e,m,D} f) \mid_D) \|_{L_q(D)} = \\
\| \D^\lambda ((E_\kappa^{d,l -\e,m,D} f) \mid_D
-(H_{\kappa, \kappa^\prime}^{d,l -\e,m,D}(E_{\kappa^\prime}^{d,l -\e,m,D}
f)) \mid_D) \|_{L_q(D)} = \\
\| \D^\lambda ((E_\kappa^{d,l -\e,m,D} f
-H_{\kappa, \kappa^\prime}^{d,l -\e,m,D}(E_{\kappa^\prime}^{d,l -\e,m,D}
f)) \mid_D) \|_{L_q(D)} = \\
\| \D^\lambda ((\mathcal E_k^{d,l -\e,m,D,\alpha} f) \mid_D) \|_{L_q(D)} = 
\| (\D^\lambda (\mathcal E_k^{d,l -\e,m,D,\alpha} f)) \mid_D \|_{L_q(D)} \le \\
\| \D^\lambda \mathcal E_k^{d,l -\e,m,D,\alpha} f \|_{L_q(\R^d)}
\le c_{27} 2^{k (\alpha^{-1}, \lambda  +(p^{-1} -q^{-1})_+ \e)}
\sum_{j =1}^d \Omega_j^{\prime l_j}(f, c_{28} 2^{-k /\alpha_j})_{L_p(D)}.
\end{multline*}

Предложение 1.3.5

Пусть выполнены условия предложения 1.3.4. Тогда существуют константы
$ c_{49}(d,l,m,D,\alpha,\lambda,p,q) >0, c_{50}(d,m,D,\alpha) >0 $ такие, что
если для функции $ f \in L_p(D) $ функция
\begin{multline*} \tag{1.3.59}
K(f,t) = K_{p,q}^{d,l,m,D,\alpha,\lambda}(f,t) =\\
 t^{-(\alpha^{-1}, \lambda +(p^{-1} -q^{-1})_+ \e) -1} 
\sum_{j =1}^d \Omega_j^{\prime l_j}(f, c_{50} t^{1 /\alpha_j})_{L_p(D)}
\in L_1(I),
\end{multline*}
то для любого $ k \in \Z_+: k \ge K^0 $ справедливо неравенство
\begin{equation*} \tag{1.3.60}
\|\D^\lambda f -\D^\lambda ((E_k^{d,l -\e,m,D,\alpha} f) \mid_D)\|_{L_q(D)} \le
c_{49} \int_0^{2^{-k}} K(f,t) dt.
\end{equation*}

Доказательство.

Прежде всего заметим, что в условиях предложения 1.3.5 в силу предложения
1.3.1 для $ f \in L_p(D) $ имеет место равенство (1.3.3) с константой $ K^0 $
вместо $ k^0 $ (см. определение области $\alpha$-типа).
Далее, пусть $ f \in L_p(D) $ и соблюдается условие (1.3.59).
Тогда согласно (1.3.58) при $ k \in \Z_+: k > K^0, $ выполняется
неравенство
\begin{multline*} \tag{1.3.61}
\| \D^\lambda ((E_k^{d,l -\e,m,D,\alpha} f) \mid_D)
-\D^\lambda ((E_{k -1}^{d,l -\e,m,D,\alpha} f) \mid_D) \|_{L_q(D)}\\
\le c_{27} 2^{k (\alpha^{-1}, \lambda  +(p^{-1} -q^{-1})_+ \e)}
\sum_{j =1}^d \Omega_j^{\prime l_j}(f, c_{28} 2^{-k /\alpha_j})_{L_p(D)} \le \\
c_{49} \int_{2^{-k}}^{2^{-k +1}} t^{-(\alpha^{-1}, \lambda +(p^{-1} -q^{-1})_+ \e) -1}
\sum_{j =1}^d \Omega_j^{\prime l_j}(f, c_{50} t^{1 /\alpha_j})_{L_p(D)} dt.
\end{multline*}

Из (1.3.61) и (1.3.59) вытекает, что ряд
$$
\D^\lambda ((E_{K^0}^{d,l -\e,m,D,\alpha} f) \mid_D) +
\sum_{k = K^0 +1}^\infty (\D^\lambda ((E_k^{d,l -\e,m,D,\alpha} f) \mid_D)
-\D^\lambda ((E_{k -1}^{d,l -\e,m,D,\alpha} f) \mid_D))
$$
сходится в $ L_q(D). $ Принимая
во внимание это обстоятельство, а также равенство (1.3.3), для
любой функции $ \phi \in C_0^\infty(D) $ имеем
\begin{multline*}
\langle \D^\lambda f, \phi \rangle = 
(-1)^{|\lambda|} \int_{D} f \D^\lambda \phi dx  = 
(-1)^{|\lambda|} \int_{D} E_{K^0}^{d,l -\e,m,D,\alpha} f
\D^\lambda \phi dx +\\
\sum_{k = K^0 +1}^\infty (-1)^{|\lambda|} \int_{D}
(E_k^{d,l -\e,m,D,\alpha} f -E_{k -1}^{d,l -\e,m,D,\alpha} f)
\D^\lambda \phi dx=\\
\int_{D} \D^\lambda (E_{K^0}^{d,l -\e,m,D,\alpha} f) \phi dx \\
+\sum_{k = K^0 +1}^\infty \int_{D} (\D^\lambda (E_k^{d,l -\e,m,D,\alpha} f)
-\D^\lambda (E_{k -1}^{d,l -\e,m,D,\alpha} f)) \phi dx \\
= \int_{D} \biggl(\D^\lambda ((E_{K^0}^{d,l -\e,m,D,\alpha} f) \mid_D)+\\
\sum_{k = K^0 +1}^\infty (\D^\lambda ((E_k^{d,l -\e,m,D,\alpha} f) \mid_D) -
\D^\lambda ((E_{k -1}^{d,l -\e,m,D,\alpha} f) \mid_D)) \biggr) \phi dx.
\end{multline*}

А это значит, что для обобщ\"eнной производной $ \D^\lambda f $
функции $ f $ в $ L_q(D) $ имеет место равенство
\begin{equation*}
\D^\lambda f = \D^\lambda ((E_{K^0}^{d,l -\e,m,D,\alpha} f) \mid_D)
+\sum_{k = K^0 +1}^\infty \biggl(\D^\lambda ((E_k^{d,l -\e,m,D,\alpha} f) \mid_D) -
\D^\lambda ((E_{k -1}^{d,l -\e,m,D,\alpha} f) \mid_D)\biggr).
\end{equation*}

Исходя из этого равенства, используя  (1.3.61), при $ k \in \Z_+: k \ge K^0 $
выводим оценку
\begin{multline*}
\| \D^\lambda f -\D^\lambda ((E_k^{d,l -\e,m,D,\alpha} f) \mid_D) \|_{L_q(D)} \le \\
\sum_{j = k +1}^\infty \biggl\| \D^\lambda((E_j^{d,l -\e,m,D,\alpha} f) \mid_D)
-\D^\lambda((E_{j -1}^{d,l -\e,m,D,\alpha} f) \mid_D)\biggr\|_{L_q(D)} \\
\le \sum_{j = k +1}^\infty c_{49} \int_{2^{-j}}^{2^{-j +1}} K(f,t) dt =
c_{49} \int_0^{2^{-k}} K(f,t) dt,
\end{multline*}
что совпадает с (1.3.60). $ \square $

Предложение 1.3.6

Пусть $ d \in \N, \alpha \in \R_+^d, l = l(\alpha), m \in \N^d, \lambda \in
\Z_+^d(m), D \subset \R^d $ -- ограниченная область $ \alpha $-типа.
Пусть ещ\"e $ 1 \le p < \infty, 1 \le q \le \infty $ и соблюдается условие
\begin{equation*} \tag{1.3.62}
1 -(\alpha^{-1}, \lambda +(p^{-1} -q^{-1})_+ \e) >0.
\end{equation*}
Тогда существует константа $ c_{51}(d,\alpha,m,D,\lambda,p,q) >0 $ такая, что
для любой функции $ f \in (\mathcal H_p^\alpha)^\prime(D) $
при $ k \in \Z_+: k \ge K^0, $ имеет место неравенство
\begin{equation*} \tag{1.3.63}
\| \D^\lambda f -\D^\lambda E_k^{d,l -\e,m,D,\alpha} f \|_{L_q(D)}
\le c_{51} 2^{-k(1 -(\alpha^{-1}, \lambda +(p^{-1} -q^{-1})_+ \e))}.
\end{equation*}

Доказательство.

Прежде всего заметим, что для параметров $ d,\alpha,l,m,\lambda,p,q,D, $
указанных в условии предложения, при соблюдении неравенства (1.3.62)
для $ f \in (H_p^\alpha)^\prime(D) $ при $ t \in I $ справедливо соотношение
\begin{multline*} \tag{1.3.64}
K(f,t) = t^{-(\alpha^{-1}, \lambda +(p^{-1} -q^{-1})_+ \e) -1}
\sum_{j =1}^d \Omega_j^{\prime l_j}(f, c_{50} t^{1 /\alpha_j})_{L_p(D)} \\
\le c_{52} t^{-(\alpha^{-1}, \lambda +(p^{-1} -q^{-1})_+ \e)}
\max_{j =1,\ldots,d} \sup_{\tau >0} \tau^{-\alpha_j}
\Omega_j^{\prime l_j}(f, \tau)_{L_p(D)} \in L_1(I).
\end{multline*}
Принимая во внимание (1.3.59) и (1.3.64), в соответствии с (1.3.60)
заключаем, что для $ f \in (\mathcal H_p^\alpha)^\prime(D) $
при $ k \in \Z_+: k \ge K^0, $ имеет место оценка
\begin{multline*}
\| \D^\lambda f -\D^\lambda ((E_k^{d,l -\e,m,D,\alpha} f) \mid_D)\|_{L_q(D)} \le \\
c_{49} \int_0^{2^{-k}} c_{52} t^{-(\alpha^{-1}, \lambda +
(p^{-1} -q^{-1})_+ \e)} dt 
= c_{51} 2^{-k(1 -(\alpha^{-1}, \lambda +(p^{-1} -q^{-1})_+ \e))}. \square
\end{multline*}

Теорема 1.3.7

Пусть $ d \in \N, \alpha \in \R_+^d, l = l(\alpha), \lambda \in \Z_+^d. $
Пусть ещ\"e $ 1 \le p < \infty, 1 \le q \le \infty $ и соблюдается условие
(1.3.62). Тогда существуют константы $ c_{53}(d,\alpha,\lambda,p,q) >0 $ и
$ c_{54}(d,\alpha) >0 $ такие, что для любого $ \delta \in \R_+^d $ и любого
$ x^0 \in \R^d $ для $ Q = x^0 +\delta I^d $ для любой функции
$ f \in (H_p^\alpha)^\prime(Q) $ справедливо неравенство
\begin{multline*} \tag{1.3.65}
\| \D^\lambda f -\D^\lambda P_{\delta,  x^0}^{d, l -\e} f \|_{L_q(Q)} \le \\
c_{53} \delta^{-\lambda -p^{-1} \e +q^{-1} \e} \int_0^1
t^{-(\alpha^{-1}, \lambda +(p^{-1} -q^{-1})_+ \e) -1} \times \\
\sum_{j =1}^d \delta_j^{-1/p}  t^{-p^{-1} /\alpha_j}
\biggl(\int_{ c_{54} \delta_j t^{1 /\alpha_j} B^1}
\int_{ Q_{l_j \xi e_j}} | \Delta_{\xi e_j}^{l_j} f(x)|^p dx d\xi \biggr)^{1/p} dt.
\end{multline*}

Доказательство.

Нетрудно видеть, что при любом $ \alpha \in \R_+^d $ куб $ I^d $ является
ограниченной областью $ \alpha $-типа, причем константа $ K^0 = K^0(D,\alpha) $
из определения области $ \alpha $-типа равна $0.$
Фиксировав $ m \in \N^d, $ для которого $ \lambda \in \Z_+^d(m), $ учитывая, что
для $ f \in (H_p^\alpha)^\prime(I^d) $
в силу (1.3.62) имеет место (1.3.64) при $ D = I^d, $ а, следовательно,
соблюдается (1.3.59), согласно (1.3.60) при $ k =0, $ с уч\"eтом того, что ввиду
(1.2.4) для $ x \in I^d $ выполняется равенство
$ (E_0^{d,l -\e,,m,I^d,\alpha} f)(x) = (P_{\e,0}^{d,l -\e} f)(x),$
получаем (1.3.65) при $ \delta = \e, x^0 =0. $ А отсюда выводим неравенство
(1.3.65) при произвольных $ \delta \in \R_+^d $ и $ x^0 \in \R^d $
пут\"eм перехода от функции $ f \in (H_p^\alpha)^\prime(Q) $ к функции
$ (h_{\delta,x^0} f)(\cdot) = f(x^0 +\delta \cdot) \in
(H_p^\alpha)^\prime(I^d), $ и, учитывая тот факт, что
$ P_{\delta,x^0}^{d,l -\e} f = h_{\delta, x^0}^{-1} P_{\e, 0}^{d, l -\e}
h_{\delta, x^0} f, f \in L_1(Q). \square $

Теорема 1.3.8

Пусть $ d \in \N, \alpha \in \R_+^d, 1 \le p < \infty, D \subset \R^d $ --
ограниченная область $ \alpha $-типа и соблюдается условие
\begin{equation*} \tag{1.3.66}
1 -(\alpha^{-1}, p^{-1} \e) >0.
\end{equation*}
Тогда имеет место включение
\begin{equation*} \tag{1.3.67}
(H_p^\alpha)^\prime(D) \subset C(D).
\end{equation*}

Доказательство.

Для доказательства (1.3.67) достаточно заметить, что при выполнении (1.3.66)
и $ m \in \N^d $ в силу (1.3.63) каждая функция $ f \in (\mathcal H_p^\alpha)^\prime(D) $,
а, следовательно, и каждая функция $ f \in (H_p^\alpha)^\prime(D) $
является пределом в $ L_\infty(D) $ равномерно сходящейся на $ \overline D $
последовательности непрерывных в $ \R^d $ функций $ \{E_k^{d,l -\e,m,D,\alpha} f,
k \in \Z_+: k \ge K^0\}. \square $
\bigskip

\centerline{\S 2. Продолжение функций из пространств
$ (H_p^\alpha)^\prime(D) $ и $ (B_{p,\theta}^\alpha)^\prime(D) $}
\bigskip

2.1. В этом пункте доказывается один из важнейших результатов работы

Теорема 2.1.1

Пусть $ d \in \N, \alpha \in \R_+^d, 1 \le p < \infty, 1 \le \theta \le \infty $
и $ D \subset \R^d $ -- ограниченная область $ \alpha $-типа.
Тогда существует непрерывное линейное отображение
$ E^{d, \alpha, p, \theta,D}: (B_{p, \theta}^\alpha)^\prime(D) \mapsto
L_p(\R^d) $ такое, что для любой функции $ f \in (B_{p, \theta}^\alpha)^\prime(D) $
соблюдается равенство
\begin{equation*} \tag{2.1.1}
(E^{d, \alpha, p, \theta,D} f) \mid_{D} = f,
\end{equation*}
а при $ \l \in \Z_+^d: \l < \alpha, $ существует константа
$ C_0(d,\alpha,p,\theta,D,\l) > 0, $ для которой выполняется неравенство
\begin{equation*} \tag{2.1.2}
\| E^{d, \alpha, p, \theta,D} f\|_{(B_{p, \theta}^\alpha)^{\l}(\R^d)} \le
C_0 \| f\|_{(B_{p, \theta}^\alpha)^\prime(D)}.
\end{equation*}

Доказательство.

В условиях теоремы положим $ l = l(\alpha) $ и фиксируем $ m \in \N^d $ так,
чтобы $ l \in \Z_+^d(m). $ Отметим еще, что в качестве $ K^0 $ рассматривается
константа из определения облати $ \alpha $-типа. Сначала рассмотрим случай,
когда $ \theta \ne \infty. $

Прежде всего заметим, что для $ f \in (B_{p, \theta}^\alpha)^\prime(D) $ ряд
\begin{equation*}
\sum_{k = K^0 +1}^\infty \| \mathcal E_k^{d,l -\e,m,D,\alpha} f\|_{L_p(\R^d)}
\end{equation*}
сходится, ибо в силу (1.3.34) при $ q = p, \lambda =0, $ и неравенства Г\"eльдера
имеет место оценка
\begin{multline*} \tag{2.1.3}
\sum_{ k \in \N: n \le k \le n +r}
\| \mathcal E_k^{d,l-\e,m,D,\alpha} f \|_{L_p(\R^d)} \le
\sum_{ k \in \N: n \le k \le n +r} c_{1} \sum_{j =1}^d
\Omega_j^{\prime l_j}(f, c_{2} 2^{-k /\alpha_j})_{L_p(D)} = \\
c_{1} \sum_{j =1}^d \sum_{ k \in \n: n \le k \le n +r}
\Omega_j^{\prime l_j}(f, c_{2} 2^{-k /\alpha_j})_{L_p(D)} = \\
c_{1} \sum_{j =1}^d \sum_{ k \in \N: n \le k \le n +r}
2^{-k} 2^k \Omega_j^{\prime l_j}(f, c_2 2^{-k /\alpha_j})_{L_p(D)} \le \\ 
c_{1} \sum_{j =1}^d \biggl( \sum_{ k \in \N: n \le k \le n +r}
2^{-\theta^\prime k} \biggr)^{1 / \theta^\prime}
\biggl( \sum_{ k \in \N: n \le k \le n +r} (2^k \Omega_j^{\prime l_j}(f, c_2
2^{-k /\alpha_j})_{L_p(D)})^{\theta} \biggr)^{1 / \theta} \le \\ 
c_{3} \sum_{j =1}^d \biggl( \sum_{ k \in \N: n \le k \le n +r}
2^{-\theta^\prime k} \biggr)^{1 / \theta^\prime} \biggl( \sum_{k \in \N}
((2^{-(k -1) /\alpha_j})^{-\alpha_j}
\Omega_j^{\prime l_j}(f, c_2 2^{-k /\alpha_j})_{L_p(D)})^{\theta} \biggr)^{1 / \theta} \le \\
c_{4} \sum_{j =1}^d \biggl( \sum_{ k \in \N: n \le k \le n +r}
2^{-\theta^\prime k} \biggr)^{1 / \theta^\prime} \biggl( \sum_{k \in \N}
\int_{2^{-k /\alpha_j}}^{2^{-(k -1) /\alpha_j}} t^{-1 -\theta \alpha_j}
(\Omega_j^{\prime l_j}(f, c_2 t)_{L_p(D)})^{\theta} dt \biggr)^{1 / \theta} \le \\
c_{4} \sum_{j =1}^d \biggl( \sum_{ k \in \N: n \le k \le n +r}
2^{-\theta^\prime k} \biggr)^{1 / \theta^\prime}
\biggl(\int_{\R_+} t^{-1 -\theta \alpha_j}
(\Omega_j^{\prime l_j}(f, c_2 t)_{L_p(D)})^{\theta} dt \biggr)^{1 / \theta} \le \\
c_{5} \sum_{j =1}^d \biggl( \sum_{ k \in \N: n \le k \le n +r}
2^{-\theta^\prime k} \biggr)^{1 / \theta^\prime}
\biggl(\int_{\R_+} t^{-1 -\theta \alpha_j}
(\Omega_j^{\prime l_j}(f, t)_{L_p(D)})^{\theta} dt \biggr)^{1 / \theta},\\ 
n \in \N: n > K^0, r \in \Z_+,
\end{multline*}
где $ \theta^\prime = \theta/(\theta -1). $

Поэтому для $ f \in (B_{p,\theta}^\alpha)^\prime(D) $ ряд
$ \sum_{ k = K^0 +1}^\infty \mathcal E_k^{d,l-\e,m,D,\alpha} f $
сходится в $ L_p(\R^d). $ Определим отображение $ E^{d,\alpha,p,\theta,D}:
(B_{p,\theta}^\alpha)^\prime(D) \mapsto L_p(\R^d), $ полагая для
$ f \in (B_{p,\theta}^\alpha)^\prime(D) $ значение
\begin{equation*} \tag{2.1.4}
E^{d,\alpha,p,\theta,D} f = E_{K^0}^{d,l -\e,m,D,\alpha} f +
\sum_{ k = K^0 +1}^\infty \mathcal E_k^{d,l -\e,m,D,\alpha} f =
E_{K^0}^{d,l -\e,m,D,\alpha} f +E^{\prime d,\alpha,p,\theta,D} f,
\end{equation*}
где
\begin{equation*} \tag{2.1.5}
E^{\prime d,\alpha,p,\theta,D} f = \sum_{ k = K^0 +1}^\infty
\mathcal E_k^{d,l-\e,m,D,\alpha} f.
\end{equation*}
Ясно, что отображение $ E^{d,\alpha,p,\theta,D} $ линейно, а, вследствие
(2.1.4), (1.3.30), (1.2.10), (1.3.3), с учетом обозначений $ \kappa = \kappa(k,\alpha),
\kappa^\prime = \kappa(k -1,\alpha), k > K^0, $ в $ L_p(D) $ выполняется
равенство
\begin{multline*}
(E^{d,\alpha,p,\theta,D} f) \mid_{D} =
(E_{K^0}^{d,l -\e,m,D,\alpha} f) \mid_D +\sum_{ k = K^0 +1}^\infty
(\mathcal E_k^{d,l-\e,m,D,\alpha} f) \mid_{D} = \\
(E_{K^0}^{d,l -\e,m,D,\alpha} f) \mid_D +\sum_{ k = K^0 +1}^\infty
((E_\kappa^{d,l -\e,m,D} f) \mid_D -
(H_{\kappa, \kappa^\prime}^{d,l -\e,m,D} (E_{\kappa^\prime}^{d,l -\e,m,D}
f)) \mid_D) = \\
(E_{K^0}^{d,l -\e,m,D,\alpha} f) \mid_D +\sum_{ k = K^0 +1}^\infty
((E_\kappa^{d,l -\e,m,D} f) \mid_D -
(E_{\kappa^\prime}^{d,l -\e,m,D} f) \mid_D) = \\
(E_{K^0}^{d,l -\e,m,D,\alpha} f) \mid_D +\sum_{ k = K^0 +1}^\infty
((E_k^{d,l -\e,m,D,\alpha} f) \mid_D -
(E_{k -1}^{d,l -\e,m,D,\alpha} f) \mid_D) = \\
(E_{K^0}^{d,l -\e,m,D,\alpha} f) \mid_D +\sum_{ k = K^0 +1}^\infty
(E_k^{d,l -\e,m,D,\alpha} f -
E_{k -1}^{d,l -\e,m,D,\alpha} f) \mid_D = f,
\end{multline*}
т.е. соблюдается (2.1.1).

Покажем, что при $ \l \in \Z_+^d: \l < \alpha, $ отображение, которое каждой
функции $ f \in (B_{p,\theta}^\alpha)^\prime(D) $ ставит в соответствие функцию
$ E_{K^0}^{d,l -\e,m,D,\alpha} f $
является непрерывным отображением пространства $ (B_{p,\theta}^\alpha)^\prime(D) $
в $ (B_{p,\theta}^\alpha)^{\l}(\R^d). $ Действительно, для $ \mathcal F =
E_{K^0}^{d,l -\e,m,D,\alpha} f, f \in (B_{p,\theta}^\alpha)^\prime(D), $
при $ j =1,\ldots,d $ и $ h \in \R, $ согласно (1.1.6), (1.3.24), с одной
стороны имеет место неравенство
\begin{multline*}
\| \Delta_{h e_j}^{l_j -\l_j} (\D_j^{\l_j} \mathcal F)\|_{L_p(\R^d)} =
\| \Delta_{h e_j}^{l_j -\l_j} (\D_j^{\l_j} E_{K^0}^{d,l -\e,m,D,\alpha} f) \|_{L_p(\R^d)} \le \\
|h|^{l_j -\l_j} \| \D^{l_j e_j} E_{K^0}^{d, l -\e,m,D, \alpha} f \|_{L_p(\R^d)} \le \\
|h|^{l_j -\l_j} c_6 2^{K^0 (\alpha^{-1}, l_j e_j)} \| f\|_{L_p(D)} =
c_7 |h|^{l_j -\l_j} \| f \|_{L_p(D)},
\end{multline*}
а с другой стороны --
\begin{equation*}
\| \Delta_{h e_j}^{l_j -\l_j} (\D_j^{\l_j} \mathcal F) \|_{L_p(\R^d)} \le
c_8 \| \D_j^{\l_j} \mathcal F \|_{L_p(\R^d)} =
c_8 \| \D_j^{\l_j} E_{K^0}^{d, l -\e,m,D,\alpha} f \|_{L_p(\R^d)} \le
c_{9} \| f \|_{L_p(D)}.
\end{equation*}
Отсюда видим, что при $ j =1,\ldots,d, t \in \R_+ $ выполняются неравенства
\begin{equation*}
\Omega_j^{l_j -\l_j}(\D_j^{\l_j} \mathcal F, t)_{L_p(\R^d)} \le
c_7 t^{l_j -\l_j} \| f \|_{L_p(D)}; \\
\Omega_j^{l_j -\l_j}(\D_j^{\l_j} \mathcal F, t)_{L_p(\R^d)} \le
c_{9} \| f \|_{L_p(D)}.
\end{equation*}
Используя эти неравенства, получаем, что
\begin{multline*}
\biggl(\int_{\R_+} t^{-1 -\theta (\alpha_j -\l_j)}
(\Omega_j^{l_j -\l_j}(\D_j^{\l_j} \mathcal F, t)_{L_p(\R^d)})^{\theta} dt \biggr)^{1 / \theta} = \\
\biggl(\int_{I} t^{-1 -\theta (\alpha_j -\l_j)}
(\Omega_j^{l_j -\l_j}(\D_j^{\l_j} \mathcal F, t)_{L_p(\R^d)})^{\theta} dt +
\int_1^\infty t^{-1 -\theta (\alpha_j -\l_j)}
(\Omega_j^{l_j -\l_j}(\D_j^{\l_j} \mathcal F, t)_{L_p(\R^d)})^{\theta} dt \biggr)^{1 / \theta} \le \\
\biggl(\int_{I} t^{-1 -\theta (\alpha_j -\l_j)}
(c_7 t^{l_j -\l_j} \| f \|_{L_p(D)})^{\theta} dt +
\int_1^\infty t^{-1 -\theta (\alpha_j -\l_j)}
(c_{9} \| f \|_{L_p(D)})^{\theta} dt \biggr)^{1 / \theta} = \\
\biggl(c_7^\theta \int_{I} t^{-1 +\theta (l_j -\alpha_j)} dt +c_{9}^\theta
\int_1^\infty t^{-1 -\theta (\alpha_j -\l_j)} dt \biggr)^{1 / \theta} \| f \|_{L_p(D)} = \\
c_{10} \| f \|_{L_p(D)} \le c_{10} \| f \|_{(B_{p,\theta}^\alpha)^\prime(D)},
f \in (B_{p,\theta}^\alpha)^\prime(D), j =1,\ldots,d.
\end{multline*}
Из сказанного с учетом (1.3.24) вытекает, что для $ f \in (B_{p,\theta}^\alpha)^\prime(D) $
справедливо неравенство
\begin{equation*} \tag{2.1.6}
\| E_{K^0}^{d,l -\e,m,D,\alpha} f \|_{(B_{p,\theta}^\alpha)^{\l}(\R^d)}
\le c_{11} \| f \|_{(B_{p,\theta}^\alpha)^\prime(D)}.
\end{equation*}

Теперь проверим, что существует константа
$ c_{12}(d,\alpha,p,\theta,D) > 0 $ такая, что для
$ f \in (B_{p, \theta}^\alpha)^\prime(D) $ выполняется неравенство
\begin{equation*} \tag{2.1.7}
\| E^{\prime d,\alpha,p,\theta,D} f\|_{(B_{p, \theta}^\alpha)^{\l}(\R^d)} \le
c_{12} \| f\|_{(B_{p, \theta}^\alpha)^\prime(D)}.
\end{equation*}

В самом деле, в силу (2.1.5), (2.1.3) справедливо неравенство
\begin{multline*} \tag{2.1.8}
\| E^{\prime d,\alpha,p,\theta,D} f \|_{L_p(\R^d)} \le \sum_{ k = K^0 +1}^\infty
\| \mathcal E_k^{d,l-\e,m,D,\alpha} f\|_{L_p(\R^d)} \le \\
c_5 \sum_{j =1}^d \biggl(\sum_{ k = K^0 +1}^\infty 2^{-\theta^\prime k}\biggr)^{1 / \theta^\prime}
\biggl(\int_{\R_+} t^{-1 -\theta \alpha_j}
(\Omega_j^{\prime l_j}(f, t)_{L_p(D)})^{\theta} dt \biggr)^{1 / \theta} \le \\
c_{13} \| f \|_{(B_{p, \theta}^\alpha)^\prime(D)}.
\end{multline*}

Далее, пользуясь тем, что для $ f \in (B_{p,\theta}^\alpha)^\prime(D) $
при $ i =1, \ldots,d, k \in \Z_+: k > K^0, $ ввиду (1.3.34) справедлива оценка
\begin{multline*}
\| \D_i^{\l_i} \mathcal E_k^{d,l -\e,m,D,\alpha} f \|_{L_p(\R^d)}
\le c_{14} 2^{k (\alpha^{-1}, \l_i e_i)}
\sum_{j =1}^d \Omega_j^{\prime l_j}(f, c_2 2^{-k /\alpha_j})_{L_p(D)} \le \\
c_{15} 2^{k \l_i /\alpha_i} \sum_{j =1}^d 2^{-k} (c_2 2^{-k /\alpha_j})^{-\alpha_j}
\Omega_j^{\prime l_j}(f, c_2 2^{-k /\alpha_j})_{L_p(D)} \le \\
c_{16} 2^{-k(1 -\l_i /\alpha_i)} \max_{j =1,\ldots,d} \sup_{t \in \R_+}
t^{-\alpha_j} \Omega_j^{\prime l_j}(f, t)_{L_p(D)},
\end{multline*}
с учетом (1.1.7) заключаем, что ряд
\begin{equation*}
\sum_{k = K^0 +1}^\infty \D_i^{\l_i} \mathcal E_k^{d,l -\e,m,D,\alpha} f
\end{equation*}
сходится в $ L_p(\R^d). $ Отсюда и из (2.1.5) приходим к выводу, что в $ L_p(\R^d) $
имеет место равенство
\begin{equation*} \tag{2.1.9}
\D_i^{\l_i} E^{\prime d,\alpha,p,\theta,D} f = \sum_{ k = K^0 +1}^\infty
\D_i^{\l_i} \mathcal E_k^{d,l-\e,m,D,\alpha} f, f \in
(B_{p,\theta}^\alpha)^\prime(D), i =1, \ldots,d.
\end{equation*}
Принимая во внимание (2.1.9), (1.3.34), а также неравенство Г\"eльдера и (2.1.3),
для $ f \in (B_{p,\theta}^\alpha)^\prime(D) $ при $ i =1, \ldots,d $ имеем
\begin{multline*} \tag{2.1.10}
\| \D_i^{\l_i} E^{\prime d,\alpha,p,\theta,D} f \|_{L_p(\R^d)}
= \| \sum_{ k = K^0 +1}^\infty \D_i^{\l_i}
\mathcal E_k^{d,l-\e,m,D,\alpha} f \|_{L_p(\R^d)} \le \\
\sum_{ k = K^0 +1}^\infty \| \D_i^{\l_i}
\mathcal E_k^{d,l-\e,m,D,\alpha} f \|_{L_p(\R^d)} \le \\
\sum_{ k = K^0 +1}^\infty c_{14} 2^{k (\alpha^{-1}, \l_i e_i)}
\sum_{j =1}^d \Omega_j^{\prime l_j}(f, c_2 2^{-k /\alpha_j})_{L_p(D)} = \\
c_{14} \sum_{j =1}^d \sum_{ k = K^0 +1}^\infty 2^{k \l_i /\alpha_i}
\Omega_j^{\prime l_j}(f, c_2 2^{-k /\alpha_j})_{L_p(D)} = \\
c_{14} \sum_{j =1}^d \sum_{ k = K^0 +1}^\infty 2^{-k(1 -\l_i /\alpha_i)} 2^k
\Omega_j^{\prime l_j}(f, c_2 2^{-k /\alpha_j})_{L_p(D)} \le \\
c_{14} \sum_{j =1}^d \biggl( \sum_{ k = K^0 +1}^\infty
2^{-k \theta^\prime (1 -\l_i /\alpha_i)} \biggr)^{1 / \theta^\prime}
\biggl( \sum_{ k = K^0 +1}^\infty (2^k \Omega_j^{\prime l_j}(f, c_2
2^{-k /\alpha_j})_{L_p(D)})^{\theta} \biggr)^{1 / \theta} \le \\
c_{17} \sum_{j =1}^d \biggl( \sum_{ k = K^0 +1}^\infty
2^{-k \theta^\prime (1 -\l_i /\alpha_i)} \biggr)^{1 / \theta^\prime}
\biggl(\int_{\R_+} t^{-1 -\theta \alpha_j}
(\Omega_j^{\prime l_j}(f, t)_{L_p(D)})^{\theta} dt \biggr)^{1 / \theta} \le \\
c_{18} \| f\|_{(B_{p,\theta}^\alpha)^\prime(D)}.
\end{multline*}

Наконец, для функции $ F = E^{\prime d,\alpha,p,\theta,D} f, f \in (B_{p,\theta}^\alpha)^\prime(D), $
при $ i =1, \ldots,d $ оценим \\
\begin{multline*} \tag{2.1.11}
\biggl(\int_{\R_+} t^{-1 -\theta (\alpha_i -\l_i)}
(\Omega_i^{l_i -\l_i}(\D_i^{\l_i} F, t)_{L_p(\R^d)})^{\theta} dt\biggr)^{1 /\theta} =\\
\biggl(\int_{I} t^{-1 -\theta (\alpha_i -\l_i)} (\Omega_i^{l_i -\l_i}(\D_i^{\l_i} F, t)_{L_p(\R^d)})^{\theta} dt +
\int_1^\infty t^{-1 -\theta (\alpha_i -\l_i)}
(\Omega_i^{l_i -\l_i}(\D_i^{\l_i} F, t)_{L_p(\R^d)})^{\theta} dt\biggr)^{1 /\theta}.
\end{multline*}

Учитывая, что при $ i =1, \ldots,d $ для $ t \in \R_+ $
справедливо неравенство
\begin{multline*}
\Omega_i^{l_i -\l_i}(\D_i^{\l_i} F, t)_{L_p(\R^d)} =
\sup_{ \{ h \in \R: h \in t B^1 \}} \| \Delta_{h e_i}^{l_i -\l_i} \D_i^{\l_i} F \|_{L_p(\R^d)}\le \\
\sup_{ \{ h \in \R: h \in t B^1 \}} c_{19} \| \D_i^{\l_i} F \|_{L_p(\R^d)} =
c_{19} \| \D_i^{\l_i} F \|_{L_p(\R^d)},
\end{multline*}
находим, что
\begin{multline*}
\int_1^\infty t^{-1 -\theta (\alpha_i -\l_i)}
(\Omega_i^{l_i -\l_i}(\D_i^{\l_i} F, t)_{L_p(\R^d)})^{\theta} dt \le \\
\int_1^\infty t^{-1 -\theta (\alpha_i -\l_i)}
(c_{19} \| \D_i^{\l_i} F \|_{L_p(\R^d)})^{\theta} dt =\\
(c_{19} \| \D_i^{\l_i} F \|_{L_p(\R^d)})^{\theta}
\int_1^\infty t^{-1 -\theta (\alpha_i -\l_i)} dt =
(c_{20} \| \D_i^{\l_i} F \|_{L_p(\R^d)})^{\theta}.
\end{multline*}

Подставляя последнюю оценку в (2.1.11) и применяя неравенство Гёльдера,
получаем, что

\begin{multline*} \tag{2.1.12}
\biggl(\int_{\R_+} t^{-1 -\theta (\alpha_i -\l_i)}
(\Omega_i^{l_i -\l_i}(\D_i^{\l_i} F, t)_{L_p(\R^d)})^{\theta} dt\biggr)^{1 /\theta} \le \\
\biggl(\int_{I} t^{-1 -\theta (\alpha_i -\l_i)} (\Omega_i^{l_i -\l_i}(\D_i^{\l_i} F, t)_{L_p(\R^d)})^{\theta} dt +
(c_{20} \| \D_i^{\l_i} F \|_{L_p(\R^d)})^{\theta}\biggr)^{1 /\theta} \le \\
(\int_{I} t^{-1 -\theta (\alpha_i -\l_i)}
(\Omega_i^{l_i -\l_i}(\D_i^{\l_i} F, t)_{L_p(\R^d)})^{\theta} dt)^{1 /\theta} +
c_{20} \| \D_i^{\l_i} F \|_{L_p(\R^d)}.
\end{multline*}

Таким образом, приходим к необходимости оценки
$$ \biggl(\int_{I} t^{-1 -\theta (\alpha_i -\l_i)} 
(\Omega_i^{l_i -\l_i}(\D_i^{\l_i} F, t)_{L_p(\R^d)})^{\theta} dt\biggr)^{1 /\theta} $$
для $ i =1, \ldots,d.$ Провед\"eм эту оценку.

Пользуясь тем, что для любого $ i =1,\ldots, d $ имеет место представление
$ I = \cup_{k \in \N} [2^{-k /\alpha_i}, 2^{-(k -1) /\alpha_i}), $
получаем, что при $ i =1,\ldots,d $ соблюдается неравенство
\begin{multline*} \tag{2.1.13}
\biggl(\int_I t^{-1 -\theta (\alpha_i -\l_i)}
(\Omega_i^{l_i -\l_i}(\D_i^{\l_i} F, t)_{L_p(\R^d)})^{\theta} dt\biggr)^{1/\theta}  = \\
\biggl(\sum_{ k \in \N} \int_{[2^{-k /\alpha_i}, 2^{-(k -1) /\alpha_i})}
t^{-1 -\theta (\alpha_i -\l_i)}
(\Omega_i^{l_i -\l_i}(\D_i^{\l_i} F, t)_{L_p(\R^d)})^{\theta} dt\biggr)^{1 /\theta} \le \\
\biggl(\sum_{ k \in \N} \!\int_{[2^{-k /\alpha_i}, 2^{-(k -1) /\alpha_i})}
(2^{-k /\alpha_i}\!)^{-1 -\theta(\alpha_i -\l_i)}
(\Omega_i^{l_i -\l_i}(\D_i^{\l_i} F, 2^{-(k -1) /\alpha_i}\!)_{L_p(\R^d)})^{\theta} dt\biggr)^{1 /\theta} \le \\
\biggl(\sum_{ k \in \N} c_{21}^{\theta} 2^{k \theta(1 -\l_i /\alpha_i)}
(\Omega_i^{l_i -\l_i}(\D_i^{\l_i} F, c_{22} 2^{-k /\alpha_i})_{L_p(\R^d)})^{\theta}\biggr)^{1 /\theta}.
\end{multline*}

При $ i =1,\ldots,d, k \in \N, $ чтобы оценить
$ \Omega_i^{l_i -\l_i}(\D_i^{\l_i} F, c_{22} 2^{-k /\alpha_i})_{L_p(\R^d)}, $
для $ h \in (c_{22} 2^{-k /\alpha_i} B^1) $ с уч\"eтом (2.1.9), (2.1.10) имеем
\begin{multline*} \tag{2.1.14}
\| \Delta_{h e_i}^{l_i -\l_i} \D_i^{\l_i} F \|_{L_p(\R^d)} =
\biggl\| \Delta_{h e_i}^{l_i -\l_i} \biggl(\sum_{ \kappa = K^0 +1}^\infty
\D_i^{\l_i} \mathcal E_\kappa^{d, l -\e,m,D,\alpha} f\biggr) \biggr\|_{L_p(\R^d)} = \\
\biggl\| \sum_{ \kappa = K^0 +1}^\infty \Delta_{h e_i}^{l_i -\l_i}
(\D_i^{\l_i} \mathcal E_\kappa^{d, l -\e,m,D,\alpha} f) \biggr\|_{L_p(\R^d)} \le \\
\sum_{ \kappa = K^0 +1}^\infty \| \Delta_{h e_i}^{l_i -\l_i}
(\D_i^{\l_i} \mathcal E_\kappa^{d, l -\e,m,D,\alpha} f)\|_{L_p(\R^d)} = \\
 \sum_{ \kappa: K^0 +1 \le \kappa \le k} \| \Delta_{h e_i}^{l_i -\l_i}
(\D_i^{\l_i} \mathcal E_\kappa^{d, l -\e,m,D,\alpha} f)\|_{L_p(\R^d)} +\\
\sum_{ \kappa \ge \max(K^0 +1, k +1)} \| \Delta_{h e_i}^{l_i -\l_i}
(\D_i^{\l_i} \mathcal E_\kappa^{d, l -\e,m,D,\alpha} f)\|_{L_p(\R^d)}.
\end{multline*}

При оценке первой суммы в правой части (2.1.14), пользуясь (1.1.6), (1.3.34)
при $ q = p, \lambda = l_i e_i, $ выводим
\begin{multline*} \tag{2.1.15}
\sum_{ \kappa: K^0 +1 \le \kappa \le k} \| \Delta_{h e_i}^{l_i -\l_i}
(\D_i^{\l_i} \mathcal E_\kappa^{d, l -\e,m,D,\alpha} f)\|_{L_p(\R^d)} \le \\
\sum_{ \kappa: K^0 +1 \le \kappa \le k} |h|^{l_i -\l_i}
\| \D^{l_i e_i} (\mathcal E_\kappa^{d, l -\e,m,D,\alpha} f)\|_{L_p(\R^d)} \le \\
\sum_{ \kappa: K^0 +1 \le \kappa \le k} c_{23} 2^{-k(l_i -\l_i) /\alpha_i}
c_{24} 2^{\kappa l_i /\alpha_i}
\sum_{j =1}^d \Omega_j^{\prime l_j}(f, c_2 2^{-\kappa /\alpha_j})_{L_p(D)} \le \\
c_{25} 2^{-k(l_i -\l_i) /\alpha_i} \sum_{j =1}^d
\sum_{ \kappa =1}^k 2^{\kappa l_i /\alpha_i}
\Omega_j^{\prime l_j}(f, c_2 2^{-\kappa /\alpha_j})_{L_p(D)},\\
 h \in (c_{22} 2^{-k /\alpha_i} B^1), k \in \N, i =1,\ldots,d.
\end{multline*}

Оценивая вторую сумму в правой части (2.1.14), в силу (1.3.34) при $ q = p,
\lambda = \l_i e_i $ получаем
\begin{multline*} \tag{2.1.16}
\sum_{ \kappa \ge \max(K^0 +1, k +1)} \| \Delta_{h e_i}^{l_i -\l_i}
(\D_i^{\l_i} \mathcal E_\kappa^{d, l -\e,m,D,\alpha} f)\|_{L_p(\R^d)} \le \\
\sum_{ \kappa \ge \max(K^0 +1, k +1)} c_{26}
\| \D_i^{\l_i} \mathcal E_\kappa^{d, l -\e,m,D,\alpha} f\|_{L_p(\R^d)} \le \\
c_{26} \sum_{ \kappa \ge \max(K^0 +1, k +1)}
c_{14} 2^{\kappa \l_i /\alpha_i} \sum_{j =1}^d
\Omega_j^{\prime l_j}(f, c_2 2^{-\kappa /\alpha_j})_{L_p(D)} \le \\
c_{27} \sum_{j =1}^d \sum_{ \kappa = k +1}^\infty
2^{\kappa \l_i /\alpha_i} \Omega_j^{\prime l_j}(f, c_2 2^{-\kappa /\alpha_j})_{L_p(D)},\\
 h \in (c_{22} 2^{-k /\alpha_i} B^1), k \in \N, i =1,\ldots,d.
\end{multline*}

Соединяя (2.1.14), (2.1.15), (2.1.16), заключаем, что для $ k \in \N,
i =1,\ldots,d $ соблюдается неравенство
\begin{multline*}
\Omega_i^{l_i -\l_i}(\D_i^{\l_i} F, c_{22} 2^{-k /\alpha_i})_{L_p(\R^d)} =
\sup_{ h \in (c_{22} 2^{-k /\alpha_i} B^1)}
\| \Delta_{h e_i}^{l_i -\l_i} \D_i^{\l_i} F \|_{L_p(\R^d)} \le \\
c_{25} 2^{-k(l_i -\l_i) /\alpha_i} \sum_{j =1}^d
\sum_{ \kappa =1}^k 2^{\kappa l_i /\alpha_i}
\Omega_j^{\prime l_j}(f, c_2 2^{-\kappa /\alpha_j})_{L_p(D)} + \\
c_{27} \sum_{j =1}^d \sum_{ \kappa = k +1}^\infty 2^{\kappa \l_i /\alpha_i}
\Omega_j^{\prime l_j}(f, c_2 2^{-\kappa /\alpha_j})_{L_p(D)} \le \\
c_{28} \sum_{j =1}^d \biggl(2^{-k(l_i -\l_i) /\alpha_i}
\sum_{ \kappa =1}^k 2^{\kappa l_i /\alpha_i}
\Omega_j^{\prime l_j}(f, c_2 2^{-\kappa /\alpha_j})_{L_p(D)} + \\
\sum_{ \kappa = k +1}^\infty 2^{\kappa \l_i /\alpha_i}
\Omega_j^{\prime l_j}(f, c_2 2^{-\kappa /\alpha_j})_{L_p(D)}\biggr).
\end{multline*}

При $ i =1,\ldots,d,$ подставляя эту оценку в (2.1.13) и применяя неравенство
Г\"eльдера, приходим к соотношению
\begin{multline*} \tag{2.1.17}
\biggl(\int_I t^{-1 -\theta (\alpha_i -\l_i)}
(\Omega_i^{l_i -\l_i}(\D_i^{\l_i} F, t)_{L_p(\R^d)})^{\theta} dt\biggr)^{1/\theta} \le \\
c_{21} \biggl(\sum_{ k \in \N} 2^{k \theta(1 -\l_i /\alpha_i)}
(c_{28} \sum_{j =1}^d \biggl(2^{-k(l_i -\l_i) /\alpha_i}
\sum_{ \kappa =1}^k 2^{\kappa l_i /\alpha_i}
\Omega_j^{\prime l_j}(f, c_2 2^{-\kappa /\alpha_j})_{L_p(D)} + \\
\sum_{ \kappa = k +1}^\infty 2^{\kappa \l_i /\alpha_i}
\Omega_j^{\prime l_j}(f, c_2 2^{-\kappa /\alpha_j})_{L_p(D)})\biggr)^\theta\biggr)^{1 /\theta} \le \\
c_{29} \biggl(\sum_{ k \in \N} 2^{k \theta(1 -\l_i /\alpha_i)}
c_{30} \sum_{j =1}^d \biggl(2^{-k(l_i -\l_i) /\alpha_i}
\sum_{ \kappa =1}^k 2^{\kappa l_i /\alpha_i}
\Omega_j^{\prime l_j}(f, c_2 2^{-\kappa /\alpha_j})_{L_p(D)} + \\
\sum_{ \kappa = k +1}^\infty 2^{\kappa \l_i /\alpha_i}
\Omega_j^{\prime l_j}(f, c_2 2^{-\kappa /\alpha_j})_{L_p(D)}\biggr)^\theta\biggr)^{1 /\theta} = \\
c_{31} \biggl(\sum_{ k \in \N} \sum_{j =1}^d 2^{k \theta(1 -\l_i /\alpha_i)}
\biggl(2^{-k(l_i -\l_i) /\alpha_i} \sum_{ \kappa =1}^k 2^{\kappa l_i /\alpha_i}
\Omega_j^{\prime l_j}(f, c_2 2^{-\kappa /\alpha_j})_{L_p(D)} + \\
\sum_{ \kappa = k +1}^\infty 2^{\kappa \l_i /\alpha_i}
\Omega_j^{\prime l_j}(f, c_2 2^{-\kappa /\alpha_j})_{L_p(D)}\biggr)^\theta\biggr)^{1 /\theta} = \\
c_{31} \biggl(\sum_{j =1}^d \sum_{ k \in \N} 2^{k \theta(1 -\l_i /\alpha_i)}
\biggl(2^{-k(l_i -\l_i) /\alpha_i} \sum_{ \kappa =1}^k 2^{\kappa l_i /\alpha_i}
\Omega_j^{\prime l_j}(f, c_2 2^{-\kappa /\alpha_j})_{L_p(D)} + \\
\sum_{ \kappa = k +1}^\infty 2^{\kappa \l_i /\alpha_i}
\Omega_j^{\prime l_j}(f, c_2 2^{-\kappa /\alpha_j})_{L_p(D)}\biggr)^\theta\biggr)^{1 /\theta} \le \\
c_{31} \biggl(\sum_{j =1}^d \sum_{ k \in \N} 2^{k \theta(1 -\l_i /\alpha_i)}
c_{32} \biggl((2^{-k(l_i -\l_i) /\alpha_i} \sum_{ \kappa =1}^k 2^{\kappa l_i /\alpha_i}
\Omega_j^{\prime l_j}(f, c_2 2^{-\kappa /\alpha_j})_{L_p(D)})^\theta + \\
(\sum_{ \kappa = k +1}^\infty 2^{\kappa \l_i /\alpha_i}
\Omega_j^{\prime l_j}(f, c_2 2^{-\kappa /\alpha_j})_{L_p(D)})^\theta\biggr)\biggr)^{1 /\theta} = \\
c_{33} \biggl(\sum_{j =1}^d \sum_{ k \in \N}
\biggl(2^{-k(l_i /\alpha_i -1) \theta} (\sum_{ \kappa =1}^k 2^{\kappa l_i /\alpha_i}
\Omega_j^{\prime l_j}(f, c_2 2^{-\kappa /\alpha_j})_{L_p(D)})^\theta + \\
2^{k \theta(1 -\l_i /\alpha_i)} (\sum_{ \kappa = k +1}^\infty 2^{\kappa \l_i /\alpha_i}
\Omega_j^{\prime l_j}(f, c_2 2^{-\kappa /\alpha_j})_{L_p(D)})^\theta\biggr)\biggr)^{1 /\theta} = \\
c_{33} \biggl(\sum_{j =1}^d \sum_{ k \in \N}
2^{-k(l_i /\alpha_i -1) \theta} \biggl(\sum_{ \kappa =1}^k 2^{\kappa l_i /\alpha_i}
\Omega_j^{\prime l_j}(f, c_2 2^{-\kappa /\alpha_j})_{L_p(D)}\biggr)^\theta + \\
\sum_{j =1}^d \sum_{ k \in \N} 2^{k \theta(1 -\l_i /\alpha_i)}
\biggl(\sum_{ \kappa = k +1}^\infty 2^{\kappa \l_i /\alpha_i}
\Omega_j^{\prime l_j}(f, c_2 2^{-\kappa /\alpha_j})_{L_p(D)}\biggr)^\theta\biggr)^{1 /\theta}.
\end{multline*}

Далее, проводя оценку правой части (2.1.17), фиксируем $ \epsilon \in I, $
для которого при $ i =1, \ldots,d $ выполняются неравенства
$ l_i /\alpha_i -1 -\epsilon >0, 1 -\l_i /\alpha_i -\epsilon >0. $
Тогда при $ i,j =1, \ldots, d, k \in \N, $ применяя неравенство Г\"eльдера, имеем
\begin{multline*} \tag{2.1.18}
\biggl(\sum_{ \kappa =1}^k 2^{\kappa l_i /\alpha_i}
\Omega_j^{\prime l_j}(f, c_2 2^{-\kappa /\alpha_j})_{L_p(D)}\biggr)^\theta = \\
\biggl(\sum_{ \kappa =1}^k 2^{\kappa l_i /\alpha_i -\kappa -\epsilon \kappa}
2^{\epsilon \kappa} 2^\kappa
\Omega_j^{\prime l_j}(f, c_2 2^{-\kappa /\alpha_j})_{L_p(D)}\biggr)^\theta \le \\
\biggl(\biggl(\sum_{ \kappa =1}^k 2^{\kappa (l_i /\alpha_i -1 -\epsilon) \theta^\prime}\biggr)^{1 /\theta^\prime}
\biggl(\sum_{ \kappa =1}^k 2^{\kappa \epsilon \theta} 2^{\kappa \theta}
(\Omega_j^{\prime l_j}(f, c_2 2^{-\kappa /\alpha_j})_{L_p(D)})^\theta\biggr)^{1 /\theta}\biggr)^\theta \le \\
(c_{34} 2^{k(l_i /\alpha_i -1 -\epsilon) \theta^\prime})^{\theta /\theta^\prime}
\sum_{ \kappa =1}^k 2^{\kappa \epsilon \theta} 2^{\kappa \theta}
(\Omega_j^{\prime l_j}(f, c_2 2^{-\kappa /\alpha_j})_{L_p(D)})^\theta = \\
c_{35} 2^{k(l_i /\alpha_i -1 -\epsilon) \theta}
\sum_{ \kappa =1}^k 2^{\kappa \epsilon \theta} 2^{\kappa \theta}
(\Omega_j^{\prime l_j}(f, c_2 2^{-\kappa /\alpha_j})_{L_p(D)})^\theta.
\end{multline*}
Используя (2.1.18) и принимая во внимание вывод неравенства (2.1.3),
получаем
\begin{multline*} \tag{2.1.19}
\sum_{ k \in \N}
2^{-k(l_i /\alpha_i -1) \theta} \biggl(\sum_{ \kappa =1}^k 2^{\kappa l_i /\alpha_i}
\Omega_j^{\prime l_j}(f, c_2 2^{-\kappa /\alpha_j})_{L_p(D)}\biggr)^\theta \le \\
\sum_{ k \in \N}
2^{-k(l_i /\alpha_i -1) \theta} c_{35} 2^{k(l_i /\alpha_i -1 -\epsilon) \theta}
\sum_{ \kappa =1}^k 2^{\kappa \epsilon \theta} 2^{\kappa \theta}
(\Omega_j^{\prime l_j}(f, c_2 2^{-\kappa /\alpha_j})_{L_p(D)})^\theta = \\
c_{35} \sum_{ k \in \N} 2^{-k \epsilon \theta}
\sum_{ \kappa =1}^k 2^{\kappa \epsilon \theta} 2^{\kappa \theta}
(\Omega_j^{\prime l_j}(f, c_2 2^{-\kappa /\alpha_j})_{L_p(D)})^\theta = \\
c_{35} \sum_{ k \in \N} \sum_{ \kappa =1}^k
2^{-k \epsilon \theta} 2^{\kappa \epsilon \theta} 2^{\kappa \theta}
(\Omega_j^{\prime l_j}(f, c_2 2^{-\kappa /\alpha_j})_{L_p(D)})^\theta = \\
c_{35} \sum_{ k \in \N, \kappa \in \N: \kappa \le k}
2^{-k \epsilon \theta} 2^{\kappa \epsilon \theta} 2^{\kappa \theta}
(\Omega_j^{\prime l_j}(f, c_2 2^{-\kappa /\alpha_j})_{L_p(D)})^\theta = \\
c_{35} \sum_{\kappa \in \N} \biggl(\sum_{k = \kappa}^\infty
2^{-k \epsilon \theta}\biggr) 2^{\kappa \epsilon \theta} 2^{\kappa \theta}
(\Omega_j^{\prime l_j}(f, c_2 2^{-\kappa /\alpha_j})_{L_p(D)})^\theta \le \\
c_{35} \sum_{\kappa \in \N} c_{36} 2^{-\kappa \epsilon \theta}
2^{\kappa \epsilon \theta} 2^{\kappa \theta}
(\Omega_j^{\prime l_j}(f, c_2 2^{-\kappa /\alpha_j})_{L_p(D)})^\theta = \\
c_{37} \sum_{\kappa \in \N} (2^{\kappa }
\Omega_j^{\prime l_j}(f, c_2 2^{-\kappa /\alpha_j})_{L_p(D)})^\theta \le \\
c_{38} \int_{\R_+} t^{-1 -\theta \alpha_j}
(\Omega_j^{\prime l_j}(f, t)_{L_p(D)})^\theta dt, i,j =1,\ldots,d.
\end{multline*}

Продолжая оценку правой части (2.1.17), при $ i, j =1, \ldots,d, k \in \N, $
пользуясь с учетом (2.1.3) неравенством Г\"eльдера. выводим
\begin{multline*} \tag{2.1.20}
\biggl(\sum_{ \kappa = k +1}^\infty 2^{\kappa \l_i /\alpha_i}
\Omega_j^{\prime l_j}(f, c_2 2^{-\kappa /\alpha_j})_{L_p(D)}\biggr)^\theta = \\
\biggl(\sum_{ \kappa = k +1}^\infty 2^{\kappa \l_i /\alpha_i} 2^{\kappa \epsilon}
2^{-\kappa} 2^{-\kappa \epsilon} 2^{\kappa}
\Omega_j^{\prime l_j}(f, c_2 2^{-\kappa /\alpha_j})_{L_p(D)}\biggr)^\theta \le \\
\biggl(\biggl(\sum_{ \kappa = k +1}^\infty 2^{-\kappa(1 -\l_i /\alpha_i -\epsilon) \theta^\prime}\biggr)^{1 /\theta^\prime}
\biggl(\sum_{ \kappa = k +1}^\infty 2^{-\kappa \epsilon \theta} 2^{\kappa \theta}
(\Omega_j^{\prime l_j}(f, c_2 2^{-\kappa /\alpha_j})_{L_p(D)})^\theta\biggr)^{1 /\theta}\biggr)^\theta \le \\
(c_{39} 2^{-k (1 -\l_i /\alpha_i -\epsilon) \theta^\prime})^{\theta /\theta^\prime}
\sum_{ \kappa = k +1}^\infty 2^{-\kappa \epsilon \theta} 2^{\kappa \theta}
(\Omega_j^{\prime l_j}(f, c_2 2^{-\kappa /\alpha_j})_{L_p(D)})^\theta = \\
c_{40} 2^{-k(1 -\l_i /\alpha_i -\epsilon) \theta}
\sum_{ \kappa = k +1}^\infty 2^{-\kappa \epsilon \theta} 2^{\kappa \theta}
(\Omega_j^{\prime l_j}(f, c_2 2^{-\kappa /\alpha_j})_{L_p(D)})^\theta.
\end{multline*}

Пользуясь (2.1.20), находим, что
\begin{multline*} \tag{2.1.21}
\sum_{ k \in \N} 2^{k \theta(1 -\l_i /\alpha_i)}
\biggl(\sum_{ \kappa = k +1}^\infty 2^{\kappa \l_i /\alpha_i}
\Omega_j^{\prime l_j}(f, c_2 2^{-\kappa /\alpha_j})_{L_p(D)}\biggr)^\theta \le \\
\sum_{ k \in \N} 2^{k \theta(1 -\l_i /\alpha_i)}
c_{40} 2^{-k(1 -\l_i /\alpha_i -\epsilon) \theta}
\sum_{ \kappa = k +1}^\infty 2^{-\kappa \epsilon \theta} 2^{\kappa \theta}
(\Omega_j^{\prime l_j}(f, c_2 2^{-\kappa /\alpha_j})_{L_p(D)})^\theta = \\
c_{40} \sum_{ k \in \N} 2^{k \epsilon \theta}
\sum_{ \kappa = k +1}^\infty 2^{-\kappa \epsilon \theta} 2^{\kappa \theta}
(\Omega_j^{\prime l_j}(f, c_2 2^{-\kappa /\alpha_j})_{L_p(D)})^\theta \le \\
c_{40} \sum_{ k \in \N} 2^{k \epsilon \theta}
\sum_{ \kappa = k}^\infty 2^{-\kappa \epsilon \theta}
(2^\kappa \Omega_j^{\prime l_j}(f, c_2 2^{-\kappa /\alpha_j})_{L_p(D)})^\theta = \\
c_{40} \sum_{ k \in \N, \kappa \in \N: \kappa \ge k} 2^{k \epsilon \theta}
2^{-\kappa \epsilon \theta}
(2^\kappa \Omega_j^{\prime l_j}(f, c_2 2^{-\kappa /\alpha_j})_{L_p(D)})^\theta = \\
c_{40} \sum_{\kappa \in \N} \biggl(\sum_{k =1}^\kappa 2^{k \epsilon \theta}\biggr)
2^{-\kappa \epsilon \theta}
(2^\kappa \Omega_j^{\prime l_j}(f, c_2 2^{-\kappa /\alpha_j})_{L_p(D)})^\theta \le \\
c_{40} \sum_{\kappa \in \N} c_{41} 2^{\kappa \epsilon \theta}
2^{-\kappa \epsilon \theta}
(2^\kappa \Omega_j^{\prime l_j}(f, c_2 2^{-\kappa /\alpha_j})_{L_p(D)})^\theta = \\
c_{42} \sum_{\kappa \in \N} (2^\kappa
\Omega_j^{\prime l_j}(f, c_2 2^{-\kappa /\alpha_j})_{L_p(D)})^\theta \le \\
c_{43} \int_{\R_+} t^{-1 -\theta \alpha_j}
(\Omega_j^{\prime l_j}(f, t)_{L_p(D)})^\theta dt (\text{см. вывод} (2.1.3)), i,j =1,\ldots,d.
\end{multline*}

Соединяя (2.1.17), (2.1.19), (2.1.21) и применяя неравенство Г\"eльдера с
показателем $ 1 /\theta \le 1, $ приходим к неравенству
\begin{multline*}
\biggl(\int_I t^{-1 -\theta (\alpha_i -\l_i)}
(\Omega_i^{l_i -\l_i}(\D_i^{\l_i} F, t)_{L_p(\R^d)})^{\theta} dt\biggr)^{1/\theta} \le \\
c_{33} \biggl(\sum_{j =1}^d c_{44} \int_{\R_+} t^{-1 -\theta \alpha_j}
(\Omega_j^{\prime l_j}(f, t)_{L_p(D)})^\theta dt\biggr)^{1 /\theta} \le \\
c_{45} \sum_{j =1}^d \biggl(\int_{\R_+} t^{-1 -\theta \alpha_j}
(\Omega_j^{\prime l_j}(f, t)_{L_p(D)})^\theta dt\biggr)^{1 /\theta},
\end{multline*}
которое в сочетании с (2.1.12) влечет неравенство
\begin{multline*} \tag{2.1.22}
\biggl(\int_{\R_+} t^{-1 -\theta (\alpha_i -\l_i)}
(\Omega_i^{l_i -\l_i}(\D_i^{\l_i} F, t)_{L_p(\R^d)})^{\theta} dt\biggl)^{1 /\theta} \le \\
c_{46} \biggl(\sum_{j =1}^d \biggl(\int_{\R_+} t^{-1 -\theta \alpha_j}
(\Omega_j^{\prime l_j}(f, t)_{L_p(D)})^\theta dt\biggr)^{1 /\theta} +
\| \D_i^{\l_i} F \|_{L_p(\R^d)}\biggr) = \\
c_{46} \biggl(\sum_{j =1}^d \biggl(\int_{\R_+} t^{-1 -\theta \alpha_j}
(\Omega_j^{\prime l_j}(f, t)_{L_p(D)})^\theta dt\biggr)^{1 /\theta} +
\| \D_i^{\l_i} (E^{\prime d,\alpha,p,\theta,D} f) \|_{L_p(\R^d)}\biggr),\\
i =1,\ldots,d.
\end{multline*}
Из (2.1.8), (2.1.10), (2.1.22) следует (2.1.7). Соединяя (2.1.6), (2.1.7) и
учитывая (2.1.4), приходим к (2.1.2) при $ \theta \ne \infty. $
При $ \theta = \infty $ доаказательство теоремы проводится по той же схеме с
заменой в соответствующих местах операции суммирования на операцию
взятия супремума или максимума. $ \square $

Следствие

В условиях теоремы 2.1.1 имеет место включение
\begin{equation*} \tag{2.1.23}
(B_{p, \theta}^\alpha)^\prime(D) \subset (B_{p, \theta}^\alpha)^{\l}(D),
\end{equation*}
и для любой функции $ f \in (B_{p, \theta}^\alpha)^\prime(D) $
выполняется неравенство
\begin{equation*} \tag{2.1.24}
\| f\|_{(B_{p, \theta}^\alpha)^{\l}(D)} \le C_0
\| f\|_{(B_{p, \theta}^\alpha)^\prime(D)}.
\end{equation*}

Для получения (2.1.23), (2.1.24) достаточно применить теорему 2.1.1
и неравенство
\begin{equation*}
\| (E^{d,\alpha,p,\theta,D} f) \mid_{D} \|_{(B_{p, \theta}^\alpha)^{\l}(D)} \le
\| E^{d,\alpha,p,\theta,D} f\|_{(B_{p, \theta}^\alpha)^{\l}(\R^d)}.
\end{equation*}

Из (1.1.8), (1.1.9) и (2.1.23), (2.1.24) вытекает, что в условиях теоремы 2.1.1

имеет место равенство $ (B_{p, \theta}^\alpha)^\prime(D) =
(B_{p, \theta}^\alpha)^{\l}(D) $ и нормы
$ \| \cdot \|_{(B_{p, \theta}^\alpha)^{\l}(D)},
\| \cdot \|_{(B_{p, \theta}^\alpha)^\prime(D)} $ эквивалентны.
\bigskip

2.2. В этом пункте приведем используемые в следующих параграфах соотношения
между нормами образов и прообразов при некоторых отображениях рассматриваемых
пространств.

При $ d \in \N $ для $ \delta \in \R_+^d $ и $ x^0 \in \R^d $ обозначим через
$ h_{\delta, x^0} $ отображение, которое каждой функции $ f, $ заданной на
некотором множестве $ S \subset \R^d, $  ставит в соответствие функцию
$ h_{\delta, x^0} f, $ определяемую на множестве $ \{ x \in \R^d: x^0 +\delta
x \in S\} = \delta^{-1} (S -x^0) $ равенством $ (h_{\delta, x^0} f)(x) =
f(x^0 +\delta x). $ Так как для $ \delta \in \R_+^d, x^0 \in \R^d $
отображение $ \R^d \ni x \mapsto x^0 +\delta x \in \R^d $ ---
взаимно однозначно, то отображение $ h_{\delta, x^0} $ является
биекцией на себя  множества всех функций с областью определения
в $ \R^d. $ При этом обратное  отображение $ h_{\delta, x^0}^{-1} $
задается равенством
\begin{equation*} \tag{2.2.1}
(h_{\delta, x^0}^{-1} f)(x) = f(\delta^{-1} (x -x^0)) =
(h_{\delta^\prime, x^{\prime 0}} f)(x)  \text{с} \delta^\prime = \delta^{-1},
x^{\prime 0} =-\delta^{-1} x^0.
\end{equation*}

Отметим, что при $ 1 \le p \le \infty $ для $ f \in L_p(x^0 +\delta D), $ где
$ D $ -- область в $ \R^d, \delta \in \R_+^d, x^0 \in \R^d, $ имеет место
равенство
\begin{multline*} \tag{2.2.2}
\| h_{\delta,x^0} f\|_{L_p(D)} = \biggl(\int_D |f(x^0 +\delta y)|^p dy\biggr)^{1/p}
= \biggl(\int_{ x^0 +\delta D} |f(x)|^p \delta^{-\e} dx\biggr)^{1/p} = \\
\delta^{-p^{-1} \e} \biggl(\int_{x^0 +\delta D} |f(x)|^p dx\biggr)^{1/p} =
\delta^{-p^{-1} \e} \|f\|_{L_p(x^0 +\delta D)},
\end{multline*}
а, следовательно, для $ f \in L_p(D) $ выполняется равенство
\begin{equation*} \tag{2.2.3}
\| h_{\delta,x^0}^{-1} f\|_{L_p(x^0 +\delta D)} = \delta^{p^{-1} \e}
\| h_{\delta,x^0} h_{\delta,x^0}^{-1} f \|_{L_p(D)} =
\delta^{p^{-1} \e} \|f\|_{L_p(D)}.
\end{equation*}

Лемма  2.2.1

Пусть $ d \in \N, l \in \Z_+^d, D $ -- область в $ \R^d, 1 \le p < \infty,
\delta \in \R_+^d, x^0 \in \R^d. $ Тогда при $ j =1,\ldots,d, t \in \R_+ $
для $ f \in L_p(x^0 +\delta D) $ справедливо равенство
\begin{equation*} \tag{2.2.4} 
\Omega_j^{\prime l_j}((h_{\delta, x^0} f),t)_{L_p(D)} =
\delta^{-p^{-1} \e} \Omega_j^{\prime l_j}(f, \delta_j t)_{L_p(x^0 +\delta D)}.
\end{equation*}

Доказательство.

В условиях леммы при $ j =1,\ldots,d $ для $ f \in L_p(x^0 +\delta D), $
делая замену переменных $ \eta = \delta_j^{-1} \xi, y = \delta^{-1}(x -x^0), $
находим, что соблюдается равенство
\begin{multline*}
\int_{ t B^1} \int_{ D_{l_j \eta e_j}} | (\Delta_{\eta e_j}^{l_j}
(h_{\delta, x^0} f))(y)|^p dy d\eta = \\
\int_{\{(\eta,y): \eta \in t B^1, y \in D_{l_j \eta e_j}\}}
|\sum_{i =0}^{l_j} C_{l_j}^i (-1)^{l_j -i}
(h_{\delta, x^0} f)(y +i \eta e_j)|^p d\eta dy = \\
\int_{\{(\eta,y): \eta \in t B^1, y \in D_{l_j \eta e_j}\}}
|\sum_{i =0}^{l_j} C_{l_j}^i (-1)^{l_j -i}
f(x^0 +\delta y +i \delta_j \eta e_j)|^p d\eta dy = \\
\int_{\{ (\xi,x): \xi \in \delta_j t B^1, x \in (x^0 +\delta D)_{l_j \xi e_j} \}}
|\sum_{i =0}^{l_j} C_{l_j}^i (-1)^{l_j -i}
f(x +i \xi e_j)|^p \delta^{-\e} \delta_j^{-1} d\xi dx = \\
\delta^{-\e} \delta_j^{-1} \int_{\delta_j t B^1} \int_{ (x^0 +\delta D)_{l_j \xi e_j}}
|\Delta_{\xi e_j}^{l_j} f(x)|^p dx d\xi.
\end{multline*}
Откуда выводим
\begin{multline*}
\Omega_j^{\prime l_j}((h_{\delta, x^0} f),t)_{L_p(D)} =
((2 t)^{-1} \int_{ t B^1} \int_{D_{l_j \eta e_j}}
| (\Delta_{\eta e_j}^{l_j} (h_{\delta, x^0} f))(y)|^p dy d\eta)^{1 /p} = \\
((2 t)^{-1} \delta^{-\e} \delta_j^{-1} \int_{\delta_j t B^1}
\int_{ (x^0 +\delta D)_{l_j \xi e_j}}
|\Delta_{\xi e_j}^{l_j} f(x)|^p dx d\xi)^{1 /p} = \\
\delta^{-p^{-1} \e} ((2 \delta_j t)^{-1} \int_{\delta_j t B^1}
\int_{ (x^0 +\delta D)_{l_j \xi e_j}}
|\Delta_{\xi e_j}^{l_j} f(x)|^p dx d\xi)^{1 /p} = \\
\delta^{-p^{-1} \e} \Omega_j^{\prime l_j}(f, \delta_j t)_{L_p(x^0 +\delta D)},
t \in \R_+, f \in L_p(x^0 +\delta D), j =1,\ldots,d,
\end{multline*}
что совпадает с (2.2.4). $ \square $

Лемма 2.2.2

Пусть $ d \in \N, D $ -- область в $ \R^d, \alpha \in \R_+^d, 1 \le p < \infty,
1 \le \theta \le \infty, \delta \in \R_+^d, x^0 \in \R^d. $ Тогда существуют
константы $ c_1(d,\alpha,p,\delta) > 0, c_2(d,\alpha,p,\delta) > 0 $
такие, что для любой функции $ f \in (B_{p,\theta}^\alpha)^\prime(x^0 +\delta D) $
соблюдается неравенство
\begin{equation*} \tag{2.2.5}
\| h_{\delta, x^0} f \|_{(B_{p,\theta}^\alpha)^\prime(D)} \le
c_1 \| f \|_{(B_{p,\theta}^\alpha)^\prime(x^0 +\delta D)},
\end{equation*}
а для $ f \in (B_{p,\theta}^\alpha)^\prime(D) $ выполняется неравенство
\begin{equation*} \tag{2.2.6}
\| h_{\delta, x^0}^{-1} f \|_{(B_{p,\theta}^\alpha)^\prime(x^0 +\delta D)} \le
c_2 \| f \|_{(B_{p,\theta}^\alpha)^\prime(D)}.
\end{equation*}

Доказательство.

В условиях леммы, полагая $ l = l(\alpha), $ для
$ f \in (B_{p,\theta}^\alpha)^\prime(x^0 +\delta D) $ при $ j =1,\ldots,d $
в силу (2.2.4) функция
$ t^{-1 -\theta \alpha_j} (\Omega_j^{\prime l_j}((h_{\delta, x^0}
f),t)_{L_p(D)})^\theta $ суммируема на $ \R_+, $ и справедливо соотношение
\begin{multline*} \tag{2.2.7}
\biggl(\int_{\R_+} t^{-1 -\theta \alpha_j} (\Omega_j^{\prime l_j}((h_{\delta, x^0}
f),t)_{L_p(D)})^\theta dt\biggr)^{1 /\theta} = \\
\biggl(\int_{\R_+} t^{-1 -\theta \alpha_j} (\delta^{-p^{-1} \e}
\Omega_j^{\prime l_j}(f, \delta_j t)_{L_p(x^0 +\delta D)})^\theta dt\biggr)^{1 /\theta} = \\
\delta^{-p^{-1} \e} \biggl(\int_{\R_+} t^{-1 -\theta \alpha_j}
(\Omega_j^{\prime l_j}(f, \delta_j t)_{L_p(x^0 +\delta D)})^\theta dt\biggr)^{1 /\theta} = \\
\delta^{-p^{-1} \e} \biggl(\int_{\R_+} (\delta_j^{-1} \tau)^{-1 -\theta \alpha_j}
(\Omega_j^{\prime l_j}(f, \delta_j \delta_j^{-1} \tau)_{L_p(x^0 +\delta D)})^\theta
\delta_j^{-1} d\tau\biggr)^{1 /\theta} = \\
\delta^{-p^{-1} \e} \biggl(\int_{\R_+} \delta_j^{\theta \alpha_j} \tau^{-1 -\theta \alpha_j}
(\Omega_j^{\prime l_j}(f, \tau)_{L_p(x^0 +\delta D)})^\theta d\tau\biggr)^{1 /\theta} = \\
\delta^{-p^{-1} \e} \delta_j^{\alpha_j} \biggl(\int_{\R_+} \tau^{-1 -\theta \alpha_j}
(\Omega_j^{\prime l_j}(f, \tau)_{L_p(x^0 +\delta D)})^\theta d\tau\biggr)^{1 /\theta}.
\end{multline*}
Объединяя неравенства (2.2.7) и (2.2.2), приходим к (2.2.5).
Наконец, для $ f \in (B_{p,\theta}^\alpha)^\prime(D), $ применяя (2.2.1),
(2.2.5), получаем, что $ h_{\delta,x^0}^{-1} f \in
(B_{p,\theta}^\alpha)^\prime(x^0 +\delta D) $ и
соблюдается (2.2.6). $ \square $
\bigskip

\centerline{\S 3. Восстановление частных производных функций}
\centerline{из неизотропных классов Никольского и Бесова}
\centerline{по значениям функций в заданном числе точек}
\bigskip

3.1. В этом пункте будет установлена оценка сверху величины
наилучшей точности линейного восстановления в пространстве $ L_q(D) $
частной производной $ \D^\lambda f $ по значениям в $ n $ точках функций $ f $
из классов $ (\mathcal H_p^\alpha)^\prime(D) $ и
$ (\mathcal B_{p,\theta}^\alpha)^\prime(D), $ определенных в ограниченной
области $ D \alpha $-типа.

Но сначала опишем постановку задачи.

Пусть $ T $ -- топологическое пространство, $ C(T) $ ---
пространство непрерывных вещественных функций, заданных на $ T, X $ --
банахово пространство над $\R $ и $ U: D(U) \mapsto X $ -- линейный оператор с
областью определения $ D(U) \subset C(T). $
Пусть ещ\"e $ \mathcal K \subset D(U) $ -- некоторый класс функций.

Для $ n \in \N $ через $ \Phi_n(C(T)) $ обозначим совокупность всех отображений
$ \phi: C(T) \mapsto \R^n, $ для каждого из которых существует набор точек
$ \{t^j \in T, j =1,\ldots,n\} $ такой, что
$ \phi(f) = (f(t^1), \ldots, f(t^n)), f \in C(T), $ а также через
$ \mathcal A^n(X) (\overline{\mathcal A}^n(X)) $ обозначим
множество всех отображений (всех линейных отображений)
$ A: \R^n \mapsto X. $

Тогда при $ n \in \N $ положим
\begin{equation*}
\sigma_n(U,\mathcal K,X) = \inf_{A \in \mathcal  A^n(X), \phi \in \Phi_n(C(T))}
\sup_{f \in \mathcal K} \|Uf -A \circ \phi(f) \|_X,
\end{equation*}
а
\begin{equation*}
\overline \sigma_n(U,\mathcal K,X) = \inf_{  A \in \overline{\mathcal A}^n(X),
\phi \in \Phi_n(C(T))} \sup_{f \in \mathcal K} \|Uf -A \circ \phi(f) \|_X.
\end{equation*}

Для доказательства основного утверждения этого пункта, теоремы
3.1.2, потребуется взятая из [6]

Лемма 3.1.1

Пусть $ d \in \N $ и $ l \in \Z_+^d. $ Тогда справедливы следующие
утверждения:

1) отображение $ \J^{d,l}: \mathcal P^{d,l} \mapsto \R^{(l+\e)^{\e}}, $
которое каждому полиному $ f \in \mathcal P^{d,l} $ ставит в соответствие
набор вещественных чисел $ y = \{y_\lambda = f(\lambda), \lambda \in \Z_+^d(l)\}, $
является изоморфизмом;

2) для $ \rho, \sigma, \tau \in \R_+^d   $ существует константа
$ c_1(d,l,\rho,\sigma,\tau) >0 $ такая, что для любого $ \delta \in \R_+^d, $
для любых точек $ x^0, y^0 \in \R^d $ таких, что $ y^0 \in
(x^0 +\sigma \delta B^d),$ и любого множества $ Q \subset
(y^0 +\tau \delta B^d), $ для любого полинома $ f \in \mathcal P^{d,l} $
справедлива оценка
\begin{equation*} \tag{3.1.1}
\sup_{x \in Q} | f(x)| \le c_1 \max_{\lambda \in \Z_+^d(l)} | f(x^0
+\rho \delta \lambda)|.
\end{equation*}

Теорема 3.1.2

Пусть $ d \in \N, \alpha \in \R_+^d, D \subset \R^d $ -- ограниченная область
$ \alpha $-типа, $ \lambda \in \Z_+^d, 1 \le p < \infty, 1 \le q \le \infty $
таковы, что выполняются условия (1.3.66) и (1.3.62). Пусть ещ\"e $ T = D,
X = L_q(D), D(U) = \{f \in C(D): \D^\lambda f \in L_q(D)\}, U = \D^\lambda,
\mathcal K = (\mathcal H_p^\alpha)^\prime(D), (\mathcal B_{p,\theta}^\alpha)^\prime(D), $
где $ \theta \in \R: \theta \ge 1. $ Тогда существуют константы
$ c_2(d,\alpha,D,p,q,\lambda) >0, n_0 \in \N $ такие, что для любого
$ n \in \N: n  \ge n_0, $ соблюдается неравенство
\begin{equation*} \tag{3.1.2}
\overline \sigma_n(U,\mathcal K,X) \le c_2 n^{-(1 -(\alpha^{-1}, \lambda +
(p^{-1} -q^{-1})_+ \e)) /(\alpha^{-1},\e)}.
\end{equation*}

Доказательство.

Учитывая включение (1.1.7), доказательство достаточно провести
лишь в случае $ \mathcal K = (\mathcal H_p^\alpha)^\prime(D). $ Так что
рассмотрим этот случай.

Прежде всего отметим, что соблюдение условий (1.3.66) и (1.3.62)
обеспечивает включение $ \mathcal K \subset D(U), $ и, следовательно,
величина $ \overline \sigma_n(U,\mathcal K,X) $ определена.

В условиях теоремы фиксируем $ l = l(\alpha) \in \N^d, m \in \N^d: \lambda \in
\Z_+^d(m) $ и константу $ K^0 = K^0(D,\alpha) \in \Z_+ $ (см. определение
области $ \alpha $-типа).
При $ k \in \Z_+: k \ge K^0, $ положим $ n(k,\alpha,D) =
l^{\e} \cdot \card \Nu_\kappa(D) $ при $ \kappa = \kappa(k,\alpha), $
где $ \Nu_\kappa(D) = \{\nu \in \Z^d: Q_{\kappa,\nu}^d \subset D\}. $

Пусть теперь $ n \in \N: n \ge n(K^0,\alpha,D). $ Тогда выберем $ k \in
\Z_+: k \ge K^0, $ так, чтобы имело место соотношение
\begin{equation*} \tag{3.1.3}
n(k,\alpha,D) \le n \le n(k +1,\alpha,D).
\end{equation*}

Далее, для $ \kappa = \kappa(k,\alpha) $ построим систему точек
\begin{equation*}
x_{\kappa,\nu}^\mu = 2^{-\kappa} (\nu +\frac{1} {4} \e
+\frac {1} {2} l^{-1} \mu) \in Q_{\kappa,\nu}^d, \nu \in \Nu_\kappa(D),
\mu \in \Z_+^d(l -\e),
\end{equation*}
и определим отображение $ \phi \in \Phi_{n(k,\alpha,D)}(C(D)), $ полагая для $ f \in C(D) $ значение
$$
\phi(f) = \{ f(x_{\kappa,\nu}^\mu), \nu  \in \Nu_\kappa(D), \mu \in \Z_+^d(l -\e)\}.
$$

Пользуясь леммой 3.1.1, возьм\"eм систему полиномов $ \{\pi_\mu^{d,l-\e} \in
\mathcal P^{d,l-\e}, \mu \in \Z_+^d(l-\e)\} $ такую, что для $ \mu, \mu^\prime \in \Z_+^d(l-\e) $ значение
\begin{equation*}
\pi_\mu^{d,l-\e}(\mu^\prime) = \begin{cases} 1, \text{ при } \mu^\prime = \mu; \\
0, \text{ при } \mu^\prime \ne \mu, \end{cases}
\end{equation*}
и обозначим через $ R_{\kappa,\nu,\mu}^{d,l-\e} $ полином
из $ \mathcal P^{d,l-\e} $ равный $ R_{\kappa,\nu,\mu}^{d,l-\e}(x) =
\pi_\mu^{d,l-\e}(2l 2^\kappa (x -x_{\kappa,\nu}^0)), \nu \in \Nu_\kappa(D),
\mu \in \Z_+^d(l-\e). $

Нетрудно видеть, что при $ \nu \in \Nu_\kappa(D), \mu, \mu^\prime \in \Z_+^d(l-\e) $
значение $ R_{\kappa,\nu,\mu}^{d,l-\e} (x_{\kappa,\nu}^{\mu^\prime}) =
\pi_\mu^{d,l-\e}(\mu^\prime). $

Определим для $ \nu \in \Nu_\kappa(D) $ линейный оператор
$ R_{\kappa,\nu}^{d,l-\e}: \R^{n(k,\alpha,D)} \mapsto \mathcal P^{d,l-\e}, $
полагая для $ t = \{t_{\nu^\prime,\mu} \in \R, \nu^\prime \in \Nu_\kappa(D),
\mu \in \Z_+^d(l-\e)\} $ значение
\begin{equation*}
R_{\kappa,\nu}^{d,l-\e} t = \sum_{\mu \in \Z_+^d(l-\e)}
t_{\nu,\mu} R_{\kappa,\nu,\mu}^{d,l-\e}.
\end{equation*}

Принимая во внимание сказанное выше,  получаем, что при
$ \nu \in \Nu_\kappa(D) $ для $ f \in C(D) $ и $ \mu \in \Z_+^d(l-\e) $ имеет
место равенство
\begin{equation*} \tag{3.1.4}
(R_{\kappa,\nu}^{d,l-\e} \circ \phi(f))  (x_{\kappa,\nu}^\mu)
= f(x_{\kappa,\nu}^\mu).
\end{equation*}

Теперь определим оператор $ A \in \overline{\mathcal A}^{n(k,\alpha,D)}(X), $
полагая для $ t = \{t_{\nu,\mu} \in \R, \nu \in \Nu_\kappa(D), \mu \in \Z_+^d(l-\e)\} $
значение
\begin{equation*}
A t = \D^\lambda \biggl(\sum_{\nu \in N_\kappa^{d,m,D}}
(R_{\kappa, \nu_\kappa^D(\nu)}^{d,l-\e} t) g_{\kappa, \nu}^{d,m}\biggr)
\text{(см. п. 1.3.).}
\end{equation*}

Тогда для $ f \in \mathcal K $ имеем
\begin{multline*} \tag{3.1.5}
\| \D^\lambda f -A \circ \phi(f) \|_X \le
\|\D^\lambda f -\D^\lambda E_k^{d,l -\e,m,D,\alpha}(f) \|_{L_q(D)} \\
+\|\D^\lambda E_\kappa^{d,l-\e,m,D} f -A \circ \phi(f)\|_{L_q(D)}.
\end{multline*}

Оценку первого слагаемого в правой части (3.1.5) да\"eт (1.3.63).
При оценке второго слагаемого в правой части (3.1.5) будем действовать так же,
как при доказательстве предложения 1.3.4, опираясь на те же вспомогательные
объекты, что и там. Используя (1.3.23), имеем
\begin{multline*} \tag{3.1.6}
\|\D^\lambda E_\kappa^{d,l-\e,m,D} f -A \circ \phi(f)\|_{L_q(D)} = \\
\|\D^\lambda E_\kappa^{d,l-\e,m,D} f -\D^\lambda (\sum_{\nu \in N_\kappa^{d,m,D}}
(R_{\kappa, \nu_\kappa^D(\nu)}^{d,l-\e} \phi(f)) g_{\kappa, \nu}^{d,m}) \|_{L_q(D)} \le \\
\|\D^\lambda E_\kappa^{d,l-\e,m,D} f -\D^\lambda (\sum_{\nu \in N_\kappa^{d,m,D}}
(R_{\kappa, \nu_\kappa^D(\nu)}^{d,l-\e} \phi(f)) g_{\kappa, \nu}^{d,m}) \|_{L_q(\R^d)} = \\
\|\D^\lambda ((E_\kappa^{d,l-\e,m,D} f) -\sum_{\nu \in N_\kappa^{d,m,D}}
(R_{\kappa, \nu_\kappa^D(\nu)}^{d,l-\e} \phi(f)) g_{\kappa, \nu}^{d,m}) \|_{L_q(\R^d)} = \\
\|\D^\lambda (\sum_{\nu \in N_\kappa^{d,m,D}}
(S_{\kappa, \nu_\kappa^D(\nu)}^{d,l -\e} f) g_{\kappa, \nu}^{d,m}
-\sum_{\nu \in N_\kappa^{d,m,D}}
(R_{\kappa, \nu_\kappa^D(\nu)}^{d,l-\e} \phi(f)) g_{\kappa, \nu}^{d,m}) \|_{L_q(\R^d)} = \\
\|\D^\lambda (\sum_{\nu \in N_\kappa^{d,m,D}}
((S_{\kappa, \nu_\kappa^D(\nu)}^{d,l -\e} f)
-(R_{\kappa, \nu_\kappa^D(\nu)}^{d,l-\e} \phi(f))) g_{\kappa, \nu}^{d,m}) \|_{L_q(\R^d)} = \\
\| \sum_{\nu \in N_\kappa^{d,m,D}}
\D^\lambda (((S_{\kappa, \nu_\kappa^D(\nu)}^{d,l -\e} f)
-(R_{\kappa, \nu_\kappa^D(\nu)}^{d,l-\e} \phi(f))) g_{\kappa, \nu}^{d,m}) \|_{L_q(\R^d)} = \\
\| \sum_{\nu \in N_\kappa^{d,m,D}} \sum_{\mu \in \Z_+^d(\lambda)} C_\lambda^\mu
\D^\mu ((S_{\kappa, \nu_\kappa^D(\nu)}^{d,l -\e} f)
-(R_{\kappa, \nu_\kappa^D(\nu)}^{d,l-\e} \phi(f)))
\D^{\lambda -\mu} g_{\kappa, \nu}^{d,m} \|_{L_q(\R^d)} = \\
\| \sum_{\mu \in \Z_+^d(\lambda)} C_\lambda^\mu
\sum_{\nu \in N_\kappa^{d,m,D}}
\D^\mu ((S_{\kappa, \nu_\kappa^D(\nu)}^{d,l -\e} f)
-(R_{\kappa, \nu_\kappa^D(\nu)}^{d,l-\e} \phi(f)))
\D^{\lambda -\mu} g_{\kappa, \nu}^{d,m} \|_{L_q(\R^d)} \le \\
\sum_{\mu \in \Z_+^d(\lambda)} C_\lambda^\mu
\| \sum_{\nu \in N_\kappa^{d,m,D}}
\D^\mu ((S_{\kappa, \nu_\kappa^D(\nu)}^{d,l -\e} f)
-(R_{\kappa, \nu_\kappa^D(\nu)}^{d,l-\e} \phi(f)))
\D^{\lambda -\mu} g_{\kappa, \nu}^{d,m} \|_{L_q(\R^d)}.
\end{multline*}

Оценивая правую часть (3.1.6), при $ \mu \in \Z_+^d(\lambda) $ с учетом
(1.3.6), (1.3.4) получаем
\begin{multline*} \tag{3.1.7}
\biggl\| \sum_{\nu \in N_\kappa^{d,m,D}}
\D^\mu ((S_{\kappa, \nu_\kappa^D(\nu)}^{d,l -\e} f) -
(R_{\kappa, \nu_\kappa^D(\nu)}^{d,l -\e} \phi(f)))
\D^{\lambda -\mu} g_{\kappa, \nu}^{d,m} \biggr\|_{L_q(\R^d)}^q = \\
\int_{\R^d} \biggl| \sum_{\nu \in N_\kappa^{d,m,D}}
\D^\mu ((S_{\kappa, \nu_\kappa^D(\nu)}^{d,l -\e} f) -
(R_{\kappa, \nu_\kappa^D(\nu)}^{d,l -\e} \phi(f)))
\D^{\lambda -\mu} g_{\kappa, \nu}^{d,m} \biggr|^q dx = \\
\int_{G_\kappa^{d,m,D}} \biggl| \sum_{\nu \in N_\kappa^{d,m,D}}
\D^\mu ((S_{\kappa, \nu_\kappa^D(\nu)}^{d,l -\e} f) -
(R_{\kappa, \nu_\kappa^D(\nu)}^{d,l -\e} \phi(f)))
\D^{\lambda -\mu} g_{\kappa, \nu}^{d,m} \biggr|^q dx = \\
\sum_{\substack{n \in \Z^d:\\ Q_{\kappa,n}^d \cap G_\kappa^{d,m,D} \ne \emptyset}}
\int_{Q_{\kappa,n}^d} \biggl| \sum_{\nu \in N_\kappa^{d,m,D}}
\D^\mu ((S_{\kappa, \nu_\kappa^D(\nu)}^{d,l -\e} f) -
(R_{\kappa, \nu_\kappa^D(\nu)}^{d,l -\e} \phi(f)))
\D^{\lambda -\mu} g_{\kappa, \nu}^{d,m} \biggr|^q dx = \\
\sum_{\substack{n \in \Z^d:\\ Q_{\kappa,n}^d \cap G_\kappa^{d,m,D} \ne \emptyset}}
\int_{Q_{\kappa,n}^d} \biggl| \sum_{\nu \in N_\kappa^{d,m,D}: Q_{\kappa,n}^d \cap
\supp g_{\kappa, \nu}^{d,m} \ne \emptyset}
\D^\mu ((S_{\kappa, \nu_\kappa^D(\nu)}^{d,l -\e} f) -\\
(R_{\kappa, \nu_\kappa^D(\nu)}^{d,l -\e} \phi(f)))
\D^{\lambda -\mu} g_{\kappa, \nu}^{d,m} \biggr|^q dx = \\
\sum_{n \in \Z^d: Q_{\kappa,n}^d \cap G_\kappa^{d,m,D} \ne \emptyset}
\biggl\| \sum_{\nu \in N_\kappa^{d,m,D}: Q_{\kappa,n}^d \cap
\supp g_{\kappa, \nu}^{d,m} \ne \emptyset}
\D^\mu ((S_{\kappa, \nu_\kappa^D(\nu)}^{d,l -\e} f) -\\
(R_{\kappa, \nu_\kappa^D(\nu)}^{d,l -\e} \phi(f)))
\D^{\lambda -\mu} g_{\kappa, \nu}^{d,m} \biggr\|_{L_q(Q_{\kappa,n}^d)}^q \le \\
\sum_{n \in \Z^d: Q_{\kappa,n}^d \cap G_\kappa^{d,m,D} \ne \emptyset}
\biggl(\sum_{\nu \in N_\kappa^{d,m,D}: Q_{\kappa,n}^d \cap
\supp g_{\kappa, \nu}^{d,m} \ne \emptyset}
\| \D^\mu (S_{\kappa, \nu_\kappa^D(\nu)}^{d,l -\e} f -\\
R_{\kappa, \nu_\kappa^D(\nu)}^{d,l -\e} \phi(f))
\D^{\lambda -\mu} g_{\kappa, \nu}^{d,m} \|_{L_q(Q_{\kappa,n}^d)}\biggr)^q.
\end{multline*}

Для оценки правой части (3.1.7), используя сначала (1.2.5), а
затем, применяя (1.1.1), при $ \mu \in \Z_+^d(\lambda),
n \in \Z^d: Q_{\kappa,n}^d \cap G_\kappa^{d,m,D} \ne \emptyset,
\nu \in N_\kappa^{d,m,D}: Q_{\kappa,n}^d \cap \supp g_{\kappa, \nu}^{d,m} \ne \emptyset, $
выводим
\begin{multline*} \tag{3.1.8}
\| \D^\mu (S_{\kappa, \nu_\kappa^D(\nu)}^{d,l -\e} f -
R_{\kappa, \nu_\kappa^D(\nu)}^{d,l -\e} \phi(f))
\D^{\lambda -\mu} g_{\kappa, \nu}^{d,m} \|_{L_q(Q_{\kappa,n}^d)} \le \\
\| \D^{\lambda -\mu} g_{\kappa, \nu}^{d,m} \|_{L_\infty(\R^d)}
\| \D^\mu (S_{\kappa, \nu_\kappa^D(\nu)}^{d,l -\e} f -
R_{\kappa, \nu_\kappa^D(\nu)}^{d,l -\e} \phi(f))\|_{L_q(Q_{\kappa,n}^d)} = \\
c_3 2^{(\kappa, \lambda -\mu)}
\| \D^\mu (S_{\kappa, \nu_\kappa^D(\nu)}^{d,l -\e} f -
R_{\kappa, \nu_\kappa^D(\nu)}^{d,l -\e} \phi(f))\|_{L_q(Q_{\kappa,n}^d)} \le \\
c_4 2^{(\kappa, \lambda -\mu)} 2^{(\kappa, \mu -q^{-1} \e)}
\| S_{\kappa, \nu_\kappa^D(\nu)}^{d,l -\e} f -
R_{\kappa, \nu_\kappa^D(\nu)}^{d,l -\e} \phi(f)\|_{L_\infty(Q_{\kappa,n}^d)} = \\
c_4 2^{(\kappa, \lambda -q^{-1} \e)}
\| S_{\kappa, \nu_\kappa^D(\nu)}^{d,l -\e} f -
R_{\kappa, \nu_\kappa^D(\nu)}^{d,l -\e} \phi(f)\|_{L_\infty(Q_{\kappa,n}^d)}.
\end{multline*}

Далее, с учетом того, что в силу (1.3.13) имеет место включение
\begin{multline*}
Q_{\kappa,n}^d \subset (2^{-\kappa} \nu_\kappa^D(\nu) +
(\gamma^1 +\e) 2^{-\kappa} B^d),\\
n \in \Z^d: Q_{\kappa,n}^d \cap G_\kappa^{d,m,D} \ne \emptyset,
\nu \in N_\kappa^{d,m,D}: \supp g_{\kappa, \nu}^{d,m} \cap Q_{\kappa,n}^d
\ne \emptyset,
\end{multline*}
применяя (3.1.1), затем используя (3.1.4), (1.3.65), наконец, делая
замену переменной, получаем, что при $ n \in \Z^d: Q_{\kappa,n}^d \cap
G_\kappa^{d,m,D} \ne \emptyset,
\nu \in N_\kappa^{d,m,D}: \supp g_{\kappa, \nu}^{d,m} \cap Q_{\kappa,n}^d
\ne \emptyset,, $ выполняется неравенство
\begin{multline*} \tag{3.1.9}
\| S_{\kappa, \nu_\kappa^D(\nu)}^{d,l -\e} f -
R_{\kappa, \nu_\kappa^D(\nu)}^{d,l -\e} \phi(f)\|_{L_\infty(Q_{\kappa,n}^d)}.= \\
\sup_{x \in Q_{\kappa,n}^d}
| (S_{\kappa, \nu_\kappa^D(\nu)}^{d,l -\e} f)(x) -
(R_{\kappa, \nu_\kappa^D(\nu)}^{d,l -\e} \phi(f))(x)| \le \\
c_5 \max_{\mu \in \Z_+^d(l-\e)}
\biggl| (S_{\kappa, \nu_\kappa^D(\nu)}^{d,l -\e} f)
(x_{\kappa,\nu_\kappa^D(\nu)}^0 +\frac {1} {2} l^{-1} 2^{-\kappa} \mu) -
(R_{\kappa, \nu_\kappa^D(\nu)}^{d,l -\e} \phi(f))
(x_{\kappa,\nu_\kappa^D(\nu)}^0 +\frac {1} {2} l^{-1} 2^{-\kappa} \mu) \biggr| = \\
c_5 \max_{\mu \in \Z_+^d(l-\e)}
| (S_{\kappa, \nu_\kappa^D(\nu)}^{d,l -\e} f)
(x_{\kappa,\nu_\kappa^D(\nu)}^\mu) -
(R_{\kappa, \nu_\kappa^D(\nu)}^{d,l -\e} \phi(f))
(x_{\kappa,\nu_\kappa^D(\nu)}^\mu)| = \\
c_5 \max_{\mu \in \Z_+^d(l-\e)}
| (S_{\kappa, \nu_\kappa^D(\nu)}^{d,l -\e} f)
(x_{\kappa,\nu_\kappa^D(\nu)}^\mu) -f(x_{\kappa,\nu_\kappa^D(\nu)}^\mu)| \le \\
c_5 \sup_{x \in Q_{\kappa,\nu_\kappa^D(\nu)}^d}
| (S_{\kappa, \nu_\kappa^D(\nu)}^{d,l -\e} f)(x) -f(x)| = \\
c_5 \| f -S_{\kappa, \nu_\kappa^D(\nu)}^{d,l -\e} f\|_{L_\infty
(Q_{\kappa,\nu_\kappa^D(\nu)}^d)} \le 
c_6 2^{(\kappa, p^{-1} \e)} \int_0^1
t^{-(\alpha^{-1}, p^{-1} \e) -1}  \\
\times \sum_{j =1}^d 2^{\kappa_j /p}  t^{-p^{-1} /\alpha_j}
\biggl(\int_{ c_7 2^{-\kappa_j} t^{1 /\alpha_j} B^1}
\int_{ (Q_{\kappa,\nu_\kappa^D(\nu)}^d)_{l_j \xi e_j}}
| \Delta_{\xi e_j}^{l_j} f(x)|^p dx d\xi \biggr)^{1/p} dt \le \\
c_6 2^{k(\alpha^{-1}, p^{-1} \e)} \int_0^1 t^{-(\alpha^{-1}, p^{-1} \e) -1} \\
\times \sum_{j =1}^d (2^{-k} t)^{-p^{-1} /\alpha_j}
\biggl(\int_{ c_8 (2^{-k} t)^{1/\alpha_j} B^1}
\int_{(Q_{\kappa,\nu_\kappa^D(\nu)}^d)_{l_j \xi e_j}}
| \Delta_{\xi e_j}^{l_j} f(x)|^p dx d\xi\biggr)^{1/p} dt = \\
c_6 \int_0^{2^{-k}} t^{-(\alpha^{-1}, p^{-1} \e) -1} \\
\times \sum_{j =1}^d t^{-p^{-1} /\alpha_j} \biggl(\int_{c_8 t^{1/\alpha_j} B^1}
\int_{(Q_{\kappa,\nu_\kappa^D(\nu)}^d)_{l_j \xi e_j}}
| \Delta_{\xi e_j}^{l_j} f(x)|^p dx d\xi\biggr)^{1/p} dt.
\end{multline*}

Объединяя (3.1.8), (3.1.9) и учитывая (1.3.48) при $ \iota =0 $ и то
обстоятельство, что $ \nu^0 = \nu_\kappa^D(\nu), $ выводим
\begin{multline*} \tag{3.1.10}
\| \D^\mu (S_{\kappa, \nu_\kappa^D(\nu)}^{d,l -\e} f -
R_{\kappa, \nu_\kappa^D(\nu)}^{d,l -\e} \phi(f))
\D^{\lambda -\mu} g_{\kappa, \nu}^{d,m} \|_{L_q(Q_{\kappa,n}^d)} \le \\
c_4 2^{(\kappa, \lambda -q^{-1} \e)}
c_6 \int_0^{2^{-k}} t^{-(\alpha^{-1}, p^{-1} \e) -1} \\
\times \sum_{j =1}^d t^{-p^{-1} /\alpha_j} \biggl(\int_{c_8 t^{1/\alpha_j} B^1}
\int_{(Q_{\kappa,\nu_\kappa^D(\nu)}^d)_{l_j \xi e_j}}
| \Delta_{\xi e_j}^{l_j} f(x)|^p dx d\xi\biggr)^{1/p} dt \le \\
c_9 2^{k(\alpha^{-1}, \lambda -q^{-1} \e)}
\int_0^{2^{-k}} t^{-(\alpha^{-1}, p^{-1} \e) -1} \\
\times \sum_{j =1}^d t^{-p^{-1} /\alpha_j} \biggl(\int_{c_8 t^{1/\alpha_j} B^1}
\int_{(D \cap D_{\kappa,n}^{\prime d,m,D,\alpha})_{l_j \xi e_j}}
| \Delta_{\xi e_j}^{l_j} f(x)|^p dx d\xi\biggr)^{1/p} dt \le \\
c_9 2^{k(\alpha^{-1}, \lambda -q^{-1} \e)}
\int_0^{2^{-k}} t^{-(\alpha^{-1}, p^{-1} \e) -1} \\
\times \sum_{j =1}^d t^{-p^{-1} /\alpha_j} \biggl(\int_{c_8 t^{1/\alpha_j} B^1}
\int_{D_{l_j \xi e_j} \cap D_{\kappa,n}^{\prime d,m,D,\alpha}}
| \Delta_{\xi e_j}^{l_j} f(x)|^p dx d\xi\biggr)^{1/p} dt,\\
 n \in \Z^d: Q_{\kappa,n}^d \cap G_\kappa^{d,m,D} \ne \emptyset,
\nu \in N_\kappa^{d,m,D}: \supp g_{\kappa, \nu}^{d,m} \cap Q_{\kappa,n}^d
\ne \emptyset, \mu \in \Z_+^d(\lambda).
\end{multline*}

Выбирая $ \epsilon >0 $ так, чтобы соблюдалось неравенство
$ 1 -(\alpha^{-1}, p^{-1} \e) -\epsilon >0, $ из (3.1.10), пользуясь
неравенством Г\"eльдера, выводим оценку
\begin{multline*}
\| \D^\mu (S_{\kappa, \nu_\kappa^D(\nu)}^{d,l -\e} f -
R_{\kappa, \nu_\kappa^D(\nu)}^{d,l -\e} \phi(f))
\D^{\lambda -\mu} g_{\kappa, \nu}^{d,m} \|_{L_q(Q_{\kappa,n}^d)} \le \\
c_9 2^{k(\alpha^{-1}, \lambda -q^{-1} \e)}
\int_0^{2^{-k}} t^{\epsilon -1 /q^\prime} t^{-(\alpha^{-1}, p^{-1} \e) -\epsilon -1 /q} \\
\times \sum_{j =1}^d t^{-p^{-1} /\alpha_j} \biggl(\int_{c_8 t^{1/\alpha_j} B^1}
\int_{D_{l_j \xi e_j} \cap D_{\kappa,n}^{\prime d,m,D,\alpha}}
| \Delta_{\xi e_j}^{l_j} f(x)|^p dx d\xi\biggr)^{1/p} dt \le \\
c_{10} 2^{k(\alpha^{-1}, \lambda -q^{-1} \e) -k \epsilon}
\biggl(\int_0^{2^{-k}} t^{-q(\alpha^{-1}, p^{-1} \e) -q \epsilon -1}
\sum_{j =1}^d t^{-(q p^{-1}) /\alpha_j} \\
\times \biggl(\int_{ c_8 t^{1/\alpha_j} B^1}
\int_{D_{l_j \xi e_j} \cap D_{\kappa,n}^{\prime d,m,D,\alpha}}
| \Delta_{\xi e_j}^{l_j} f(x)|^p dx d\xi \biggr)^{q /p} dt \biggr)^{1/q},\\
n \in \Z^d: Q_{\kappa,n}^d \cap G_\kappa^{d,m,D} \ne \emptyset,
\nu \in N_\kappa^{d,m,D}: \supp g_{\kappa, \nu}^{d,m} \cap Q_{\kappa,n}^d
\ne \emptyset, \mu \in \Z_+^d(\lambda).
\end{multline*}

Подставляя эту оценку в (3.1.7) и  проводя соответствующие выкладки,
в частности, используя (1.3.17), применяя неравенство Г\"eльдера
с показателем $ p /q \le 1 $ и (1.3.52), для $ f \in \mathcal K $ при $ p \le q,
\mu \in \Z_+^d(\lambda) $ получаем
\begin{multline*} \tag{3.1.11}
\biggl\| \sum_{\nu \in N_\kappa^{d,m,D}}
\D^\mu ((S_{\kappa, \nu_\kappa^D(\nu)}^{d,l -\e} f) -
(R_{\kappa, \nu_\kappa^D(\nu)}^{d,l -\e} \phi(f)))
\D^{\lambda -\mu} g_{\kappa, \nu}^{d,m} \biggr\|_{L_q(\R^d)}^q \le \\
\sum_{n \in \Z^d: Q_{\kappa,n}^d \cap G_\kappa^{d,m,D} \ne \emptyset}
\biggl(\sum_{\nu \in N_\kappa^{d,m,D}: Q_{\kappa,n}^d \cap
\supp g_{\kappa, \nu}^{d,m} \ne \emptyset}
c_{10} 2^{k(\alpha^{-1}, \lambda -q^{-1} \e) -k \epsilon}\\
\times\biggl(\int_0^{2^{-k}} t^{-q(\alpha^{-1}, p^{-1} \e) -q \epsilon -1}
\sum_{j =1}^d t^{-(q p^{-1}) /\alpha_j} \\
\times \biggl(\int_{ c_8 t^{1/\alpha_j} B^1}
\int_{D_{l_j \xi e_j} \cap D_{\kappa,n}^{\prime d,m,D,\alpha}}
| \Delta_{\xi e_j}^{l_j} f(x)|^p dx d\xi \biggr)^{q /p} dt \biggr)^{1/q}\biggr)^q \le \\
\sum_{n \in \Z^d: Q_{\kappa,n}^d \cap G_\kappa^{d,m,D} \ne \emptyset}
\biggl(c_{11} 2^{k(\alpha^{-1}, \lambda -q^{-1} \e) -k \epsilon}
\biggl(\int_0^{2^{-k}} t^{-q(\alpha^{-1}, p^{-1} \e) -q \epsilon -1}
\sum_{j =1}^d t^{-(q p^{-1}) /\alpha_j} \\
\times \biggl(\int_{ c_8 t^{1/\alpha_j} B^1}
\int_{D_{l_j \xi e_j} \cap D_{\kappa,n}^{\prime d,m,D,\alpha}}
| \Delta_{\xi e_j}^{l_j} f(x)|^p dx d\xi \biggr)^{q /p} dt \biggr)^{1/q}\biggr)^q = \\
(c_{11} 2^{k(\alpha^{-1}, \lambda -q^{-1} \e) -k \epsilon})^q
\sum_{n \in \Z^d: Q_{\kappa,n}^d \cap G_\kappa^{d,m,D} \ne \emptyset}
\int_0^{2^{-k}} t^{-q(\alpha^{-1}, p^{-1} \e) -q \epsilon -1}
\sum_{j =1}^d t^{-(q p^{-1}) /\alpha_j} \\
\times \biggl(\int_{ c_8 t^{1/\alpha_j} B^1}
\int_{D_{l_j \xi e_j} \cap D_{\kappa,n}^{\prime d,m,D,\alpha}}
| \Delta_{\xi e_j}^{l_j} f(x)|^p dx d\xi \biggr)^{q /p} dt = \\
(c_{11} 2^{k(\alpha^{-1}, \lambda -q^{-1} \e) -k \epsilon})^q
\int_0^{2^{-k}} t^{-q(\alpha^{-1}, p^{-1} \e) -q \epsilon -1}
\sum_{n \in \Z^d: Q_{\kappa,n}^d \cap G_\kappa^{d,m,D} \ne \emptyset}
\sum_{j =1}^d t^{-(q p^{-1}) /\alpha_j} \\
\times \biggl(\int_{ c_8 t^{1/\alpha_j} B^1}
\int_{D_{l_j \xi e_j} \cap D_{\kappa,n}^{\prime d,m,D,\alpha}}
| \Delta_{\xi e_j}^{l_j} f(x)|^p dx d\xi \biggr)^{q /p} dt = \\
(c_{11} 2^{k(\alpha^{-1}, \lambda -q^{-1} \e) -k \epsilon})^q
\int_0^{2^{-k}} t^{-q(\alpha^{-1}, p^{-1} \e) -q \epsilon -1}
\sum_{j =1}^d t^{-(q p^{-1}) /\alpha_j}  \times \\
\sum_{n \in \Z^d: Q_{\kappa,n}^d \cap G_\kappa^{d,m,D} \ne \emptyset}
\biggl(\int_{ c_8 t^{1/\alpha_j} B^1}
\int_{D_{l_j \xi e_j} \cap D_{\kappa,n}^{\prime d,m,D,\alpha}}
| \Delta_{\xi e_j}^{l_j} f(x)|^p dx d\xi \biggr)^{q /p} dt \le \\
(c_{11} 2^{k(\alpha^{-1}, \lambda -q^{-1} \e) -k \epsilon})^q
\int_0^{2^{-k}} t^{-q(\alpha^{-1}, p^{-1} \e) -q \epsilon -1}
\sum_{j =1}^d t^{-(q p^{-1}) /\alpha_j}  \\
\times\biggl(\sum_{n \in \Z^d: Q_{\kappa,n}^d \cap G_\kappa^{d,m,D} \ne \emptyset}
\int_{ c_8 t^{1/\alpha_j} B^1}
\int_{D_{l_j \xi e_j} \cap D_{\kappa,n}^{\prime d,m,D,\alpha}}
| \Delta_{\xi e_j}^{l_j} f(x)|^p dx d\xi \biggr)^{q /p} dt = \\
(c_{11} 2^{k(\alpha^{-1}, \lambda -q^{-1} \e) -k \epsilon})^q
\int_0^{2^{-k}} t^{-q(\alpha^{-1}, p^{-1} \e) -q \epsilon -1}
\sum_{j =1}^d t^{-(q p^{-1}) /\alpha_j} \\
\times\biggl(\int_{ c_8 t^{1/\alpha_j} B^1}
\sum_{n \in \Z^d: Q_{\kappa,n}^d \cap G_\kappa^{d,m,D} \ne \emptyset}
\int_{D_{l_j \xi e_j}} \chi_{ D_{\kappa,n}^{\prime d,m,D,\alpha}}(x)
| \Delta_{\xi e_j}^{l_j} f(x)|^p dx d\xi \biggr)^{q /p} dt = \\
(c_{11} 2^{k(\alpha^{-1}, \lambda -q^{-1} \e) -k \epsilon})^q
\int_0^{2^{-k}} t^{-q(\alpha^{-1}, p^{-1} \e) -q \epsilon -1}
\sum_{j =1}^d t^{-(q p^{-1}) /\alpha_j}  \\
\times\biggl(\int_{ c_8 t^{1/\alpha_j} B^1}
\int_{D_{l_j \xi e_j}} (\sum_{n \in \Z^d: Q_{\kappa,n}^d \cap
G_\kappa^{d,m,D} \ne \emptyset}
\chi_{ D_{\kappa,n}^{\prime d,m,D,\alpha}}(x))
| \Delta_{\xi e_j}^{l_j} f(x)|^p dx d\xi \biggr)^{q /p} dt \le \\
(c_{11} 2^{k(\alpha^{-1}, \lambda -q^{-1} \e) -k \epsilon})^q
\int_0^{2^{-k}} t^{-q(\alpha^{-1}, p^{-1} \e) -q \epsilon -1}
\sum_{j =1}^d t^{-(q p^{-1}) /\alpha_j}  \\
\times\biggl(\int_{ c_8 t^{1/\alpha_j} B^1} \int_{D_{l_j \xi e_j}} c_{12}
| \Delta_{\xi e_j}^{l_j} f(x)|^p dx d\xi \biggr)^{q /p} dt \le \\
(c_{13} 2^{k(\alpha^{-1}, \lambda -q^{-1} \e) -k \epsilon})^q
\int_0^{2^{-k}} t^{-q(\alpha^{-1}, p^{-1} \e) -q \epsilon -1}
\sum_{j =1}^d \biggl((2 c_8 t^{1/\alpha_j})^{-p^{-1}} \\
\times\biggl(\int_{ c_8 t^{1/\alpha_j} B^1} \int_{D_{l_j \xi e_j}}
| \Delta_{\xi e_j}^{l_j} f(x)|^p dx d\xi\biggr)^{1 /p} \biggr)^q dt = \\
(c_{13} 2^{k(\alpha^{-1}, \lambda -q^{-1} \e) -k \epsilon})^q
\int_0^{2^{-k}} t^{-q(\alpha^{-1}, p^{-1} \e) -q \epsilon -1}
\sum_{j =1}^d \biggl(\Omega_j^{\prime l_j}(f, c_8 t^{1/\alpha_j})_{L_p(D)}  \biggr)^q dt \le \\
(c_{13} 2^{k(\alpha^{-1}, \lambda -q^{-1} \e) -k \epsilon})^q
\int_0^{2^{-k}} t^{-q(\alpha^{-1}, p^{-1} \e) -q \epsilon -1} c_{14} t^q dt = \\
(c_{15} 2^{k(\alpha^{-1}, \lambda -q^{-1} \e) -k \epsilon})^q
\int_0^{2^{-k}} t^{q(1 -(\alpha^{-1}, p^{-1} \e) -\epsilon) -1} dt = \\
(c_{15} 2^{k(\alpha^{-1}, \lambda -q^{-1} \e) -k \epsilon})^q
c_{16} 2^{-k q(1 -(\alpha^{-1}, p^{-1} \e) -\epsilon)} = \\
(c_{17} 2^{-k (1 -(\alpha^{-1}, \lambda +p^{-1} \e -q^{-1} \e))})^q.
\end{multline*}

Соединяя (3.1.6) и (3.1.11), приходим к оценке
\begin{equation*} \tag{3.1.12}
\|\D^\lambda E_\kappa^{d,l-\e,m,D} f -A \circ \phi(f)\|_{L_q(D)} \le
c_{18} 2^{-k (1 -(\alpha^{-1}, \lambda +p^{-1} \e -q^{-1} \e))}.
\end{equation*}
Подстановка оценок (1.3.63) и (3.1.12) в (3.1.5) приводит к неравенству
\begin{equation*} \tag{3.1.13}
\| \D^\lambda f -A \circ \phi(f) \|_{L_q(D)} \le
c_{19} 2^{-k (1 -(\alpha^{-1}, \lambda +(p^{-1} -q^{-1})_+ \e))},
f \in \mathcal K, p \le q.
\end{equation*}
В условиях теоремы 3.1.2 выполнение неравенства (3.1.13) при $ q < p $
вытекает из соблюдения неравенства (3.1.13) при $ q = p $ и того факта, что
при $ q < p $ для $ h \in L_p(D) $ справедливо неравенство
$ \|h\|_{L_q(D)} \le c_{20} \|h\|_{L_p(D)}. $

Из (3.1.13) с уч\"eтом (3.1.3) следует, что
\begin{equation*}
\overline \sigma_n(U,\mathcal K,X) \le \overline \sigma_{n(k,\alpha,D)}(U,\mathcal K,X) \le
c_{19} 2^{-k(1 -(\alpha^{-1}, \lambda +(p^{-1} -q^{-1})_+ \e))}.
\end{equation*}
Отсюда, пользуясь тем, что для $ n, $ удовлетворяющего (3.1.3),
справедливо соотношение $ n \le c_{21}(\alpha,D) 2^{k(\alpha^{-1},\e)}, $
приходим к (3.1.2). $ \square $
\bigskip

3.2. В этом пункте выводится оценка  снизу величины наилучшей
точности восстановления в $ L_q(D) $ производной $ \D^\lambda f $ по
значениям в $ n $ точках функций $ f $ из $ (\mathcal H_p^\alpha)^\prime(D),
(\mathcal B_{p,\theta}^\alpha)^\prime(D). $

Теорема 3.2.1

В условиях теоремы 3.1.2 существует константа $ c_1(d,\alpha,D,p,\theta,q,\lambda) >0 $
такая, что при $ n \in \N $ выполняется неравенство
\begin{equation*} \tag{3.2.1}
\sigma_n(U,\mathcal K,X) \ge c_1 n^{-(1 -(\alpha^{-1}, \lambda +(p^{-1}
-q^{-1})_+ \e)) /(\alpha^{-1},\e)}.
\end{equation*}

В основе доказательства теоремы 3.2.1 лежит лемма 3.2.2, установленная в
[6].

Лемма 3.2.2

Пусть $ d \in \N, \alpha \in \R_+^d, \theta \in \R: \theta \ge 1,
1 \le p < \infty, 1 \le q \le \infty, \lambda \in \Z_+^d $ и
$ \l = \l(\alpha) \in \Z_+^d. $
Тогда существует константа $ c_2(d,\alpha,p,\theta,q,\lambda) >0 $ такая,
что для любого $ n \in \N $ для любого набора точек $ x^1,\ldots,x^n \in I^d $
можно построить функцию
$ f \in C_0^\infty(I^d) \cap (\mathcal B_{p,\theta}^\alpha)^{\l}(\R^d), $
для которой
\begin{equation*} \tag{3.2.2}
f(x^i) =0, i =1,\ldots,n,
\end{equation*}
и
\begin{equation*} \tag{3.2.3}
\| \D^\lambda f\|_{L_q(I^d)} \ge
c_2 n^{-(1 -(\alpha^{-1}, \lambda +(p^{-1} -q^{-1})_+ \e)) /(\alpha^{-1}, \e)}.
\end{equation*}

Отметим, что функция $ f $ в формулировке леммы 3.2.2, приведенной в
[6], принадлежит $ C_0^\infty(\R^d) \cap
(\mathcal B_{p,\theta}^\alpha)^{\l}(I^d). $ Однако, как видно из
доказательства леммы, на самом деле $ f \in C_0^\infty(I^d) \cap
(\mathcal B_{p,\theta}^\alpha)^{\l}(\R^d). $

Доказательство теоремы 3.2.1.

Принимая во внимание включение (1.1.7), доказательство проводим в
случае $ \mathcal K = (\mathcal B_{p,\theta}^\alpha)^\prime(D). $

Фиксируем точку $ x^0 \in \R^d $ и вектор $ \delta \in \R_+^d $ такие, что
$ Q = (x^0 +\delta I^d) \subset D, $ и пусть $ n \in \N, A \in \mathcal A^n(X)  $
и $ \phi \in \Phi_n(C(D)) $ определяется системой точек $ y^1,\ldots,y^n \in D. $

Для $ \lambda \in \Z_+^d, $ удовлетворяющего условиям теоремы 3.2.1,
возьм\"eм функцию $ f \in (\mathcal B_{p,\theta}^\alpha)^{\l}(\R^d) \cap
C_0^\infty(I^d), $ обладающую свойствами (3.2.2) и (3.2.3), где
$ x^i = \delta^{-1}(y^i -x^0), i =1,\ldots,n. $ Рассмотрим функцию
$ F = h_{\delta,x^0}^{-1} f $ (см. п. 2.2.). Тогда, поскольку ввиду (1.1.8)
функция $ f \in (\mathcal B_{p,\theta}^\alpha)^\prime(\R^d) \cap
C_0^\infty(I^d), $ то в силу (2.2.1), (2.2.7) функция
$ F \in (c_3 (\mathcal B_{p,\theta}^\alpha)^\prime(\R^d)) \cap C_0^\infty(Q). $
Полагая $ \mathcal F = (1 / c_3) F \mid_D, $ получаем, что
$ \supp \mathcal F \subset Q, \mathcal F \in (\mathcal B_{p,\theta}^\alpha)^\prime(D),
\mathcal F(y^i) = (1 / c_3) F(y^i) = (1 / c_3) f(\delta^{-1}(y^i -x^0)) =
(1 / c_3) f(x^i) =0, i =1,\ldots,n. $

При этом ввиду (3.2.3) соблюдается неравенство
\begin{multline*}
c_2 n^{-(1 -(\alpha^{-1}, \lambda +(p^{-1} -q^{-1})_+ \e)) /(\alpha^{-1}, \e)}
\le \|\D^\lambda f\|_{L_q(I^d)} \\
= \|\D^\lambda (h_{\delta,x^0} F)\|_{L_q(I^d)} =
c_4(\delta,\lambda,q) \| \D^\lambda F\|_{L_q(Q)} =
c_5 \| \D^\lambda \mathcal F\|_{L_q(Q)} = c_5 \| \D^\lambda \mathcal F\|_{L_q(D)} = \\
c_5 \|\D^\lambda \mathcal F -A \circ \phi(\mathcal F) +A \circ \phi(0) -\D^\lambda 0\|_{L_q(D)}
\le 2 c_5 \sup_{g \in \mathcal K} \|\D^\lambda g -A \circ \phi(g)\|_X.
\end{multline*}
Отсюда, в силу произвольности $ A \in \mathcal A^n(X) $ и $ \phi \in
\Phi_n(C(D)), $ заключаем, что имеет место (3.2.1). $\square$
\bigskip

\centerline{\S 4. Задача С.Б. Стечкина}
\centerline{для оператора частного дифференцирования}
\centerline{на неизотропных классах функций Никольского и Бесова}

\bigskip

4.1. В этом пункте будет получена оценка сверху наилучшей точности
приближения в $ L_q(D) $ оператора $ \D^\lambda $ ограниченными
операторами, действующими из $ L_s(D) $ в $ L_q(D), $ на
классах $ (\mathcal H_p^\alpha)^\prime(D) $ и
$ (\mathcal B_{p,\theta}^\alpha)^\prime(D), $ где $ D $ -- ограниченная
область $ \alpha $-типа в $ \R^d. $

Напомним постановку общей задачи,  частным случаем которой
является задача, рассматриваемая в этом параграфе.

Пусть $ X,Y $ -- банаховы пространства, $ \mathcal B(X, Y) $ --
банахово пространство непрерывных линейных операторов $ V: X \mapsto Y $
с обычной нормой, $ A: D(A) \mapsto Y $ -- линейный оператор с областью
определения $ D(A) \subset X. $ Пусть ещ\"e множество $ \mathcal K \subset D(A) $ и
$ \rho >0. $

Требуется описать поведение в зависимости от $ \rho $ величины
\begin{equation*}
E(A,X,Y,\mathcal K,\rho) = \inf_{\{ V \in \mathcal B(X,Y): \|V\|_{\mathcal B(X,Y)}
\le \rho\}} \sup_{x \in \mathcal K} \|Ax -Vx\|_Y.
\end{equation*}

Теорема 4.1.1

Пусть $ d \in \N, \alpha \in \R_+^d, 1 \le p < \infty, 1 \le q,s \le \infty $ и
$ \lambda \in \Z_+^d $ удовлетворяют  условию
(1.3.62), а также при $ p < s $ соблюдается неравенство
\begin{equation*}
1 -(\alpha^{-1}, (p^{-1} -s^{-1}) \e) >0,
\end{equation*}
и, кроме того, выполняется неравенство
$$
(\alpha^{-1}, \lambda +(s^{-1} -q^{-1})_+ \e) >0.
$$
Пусть ещ\"e $ D \subset \R^d $-- ограниченная область $ \alpha $-типа,
$ A = \D^\lambda, D(A) = \{f \in L_s(D): \D^\lambda f \in L_q(D)\}, X =
L_s(D), Y = L_q(D), \mathcal K =
(\mathcal H_p^\alpha)^\prime(D),
(\mathcal B_{p,\theta}^\alpha)^\prime(D), $ где $ \theta \in \R:
\theta \ge 1, \gamma = 1-(\alpha^{-1}, \lambda +(p^{-1} -q^{-1})_+ \e),
\tau = (\alpha^{-1}, \lambda +(s^{-1} -q^{-1})_+ \e). $ Тогда существуют
константы $ \rho_1(d,\alpha,D,q,s,\lambda) >0 $ и
$ c_1(d,\alpha,D,p,q,s,\lambda) >0 $ такие, что для $ \rho \ge \rho_1 $
справедливо неравенство
\begin{equation*} \tag{4.1.1}
E(A,X,Y,\mathcal K,\rho) \le c_1 \rho^{-\gamma /\tau}.
\end{equation*}

Для доказательства теоремы 4.1.1 достаточно заметить, что в условиях
теоремы множество $ \mathcal K \subset D(A), $ и если взять $ l = l(\alpha),
m \in \N^d: \lambda \in \Z_+^d(m), $ то для каждого достаточно большого $ k \in \Z_+ $
для оператора $ V: L_s(D) \mapsto L_q(D), $ определяемого равенством
$$
V f = (\D^\lambda (E_k^{d,l -\e,m,D,\alpha} f)) \mid_D,
$$
в силу (1.3.63), (1.3.24) (см. также (1.1.7)) выполняются соотношения
\begin{equation*}
\sup_{f \in \mathcal K} \| \D^\lambda f -V f \|_{L_q(D)} \le
c_2 2^{-k(1 -(\alpha^{-1}, \lambda +(p^{-1} -q^{-1})_+ \e))},
\end{equation*}
и
\begin{equation*}
\| V \|_{\mathcal B(L_s(D), L_q(D))} \le
c_3 2^{k (\alpha^{-1}, \lambda +(s^{-1} -q^{-1})_+ \e)},
\end{equation*}
а, следовательно,
\begin{multline*}
\inf_{\{\mathcal V \in \mathcal B(X,Y): \| \mathcal V \|_{\mathcal B(X,Y)} \le
c_3 2^{k (\alpha^{-1}, \lambda +(s^{-1} -q^{-1})_+ \e)}\}}
\sup_{f \in \mathcal K} \| \D^\lambda f -\mathcal V f \|_{Y} \le\\
c_2 2^{-k(1 -(\alpha^{-1}, \lambda +(p^{-1} -q^{-1})_+ \e))}. \square
\end{multline*}
\bigskip

4.2. В этом пункте проводится оценка снизу величины, которая
рассматривалась в п. 4.1. Для доказательства этой оценки используем
утверждение, справедливость которого по существу установлена в [6]
при доказательстве теоремы 3.2.4, опирающейся на леммы 3.2.2, 3.2.3.
При этом будем пользоваться следующим обозначением. Для множества $ S, $
состоящего из функций $ f, $ область определения которых содержит множество
$ D \subset \R^d, $ через $ S \mid_D $ обозначим множество $ S \mid_D =
\{ f \mid_D: f \in S\}. $

Теорема 4.2.1

Пусть $ d, \alpha, p, q,s, \lambda $  удовлетворяют условиям теоремы 4.1.1 и
$ \l = \l(\alpha) $ (см. п. 1.1.). Пусть ещ\"e $ A = \D^\lambda,
D(A) = \{f  \in L_s(I^d): \D^\lambda f \in L_q(I^d)\},
X = L_s(I^d), Y = L_q(I^d), \mathcal K =
((\mathcal B_{p,\theta}^\alpha)^{\l}(\R^d) \cap C_0^\infty(I^d)) \mid_{I^d}, $
где $ \theta \in \R: \theta \ge 1, \gamma = 1 -(\alpha^{-1}, \lambda +
(p^{-1} -q^{-1})_+ \e), \tau = (\alpha^{-1}, \lambda +(s^{-1} -q^{-1})_+ \e). $
Тогда существуют константы $ c_1(d,\alpha,p,\theta,q,s,\lambda) >0 $ и
$ \rho_1(d,\alpha,q,s,\lambda) >0 $ такие, что для $ \rho > \rho_1 $
выполняется неравенство
\begin{equation*} \tag{4.2.1}
E(A,X,Y,\mathcal K,\rho) \ge c_1 \rho^{-\gamma /\tau}.
\end{equation*}

Теорема 4.2.2

В условиях и обозначениях теоремы 4.1.1 существуют константы
$ \rho_2(d,\alpha,D,q,s,\lambda) >0, c_2(d,\alpha,D,p,\theta,q,s,\lambda) >0 $
такие, что при $ \rho > \rho_2 $ соблюдается неравенство (4.2.1) с
константой $ c_2 $ вместо $ c_1. $

Доказательство.

Ввиду включения (1.1.7) рассмотрим лишь случай, когда $ \mathcal K =
(\mathcal B_{p,\theta}^\alpha)^\prime(D). $
Фиксируем $ x^0 \in \R^d $ и $ \delta \in \R_+^d $ такие, что $ Q =
(x^0 +\delta I^d) \subset D. $
Прежде всего, отметим, что в силу (1.1.8), (4.2.1) при $ \rho > \rho_1$
верно неравенство
\begin{equation*} \tag{4.2.2}
c_1 \rho^{-\gamma /\tau} \le E(\D^\lambda,L_s(I^d),L_Q(I^d),
((\mathcal B_{p,\theta}^\alpha)^\prime(\R^d) \cap C_0^\infty(I^d)) \mid_{I^d},
\rho).
\end{equation*}

Далее, заметим, что для $ f \in L_s(I^d): \D^\lambda f \in L_q(I^d),
V \in \mathcal B(L_s(Q),L_q(Q)) $ имеет место соотношение
\begin{multline*}
\| \D^\lambda f -(h_{\delta,x^0} V h_{\delta,x^0}^{-1}) f\|_{L_q(I^d)} =
\biggl(\int_{I^d}| \D^\lambda f -(h_{\delta,x^0} V h_{\delta,x^0}^{-1}) f|^q dx\biggr)^{1 /q} = \\
\biggl(\int_{Q}| h_{\delta,x^0}^{-1}(\D^\lambda f) -V (h_{\delta,x^0}^{-1} f)|^q \delta^{-\e} dy\biggr)^{1 /q} = \\
\delta^{-q^{-1} \e} \biggl(\int_{Q}| \delta^\lambda \D^\lambda (h_{\delta,x^0}^{-1} f) -
V (h_{\delta,x^0}^{-1} f)|^q dy\biggr)^{1 /q} = \\
\delta^{\lambda -q^{-1} \e} \| \D^\lambda (h_{\delta,x^0}^{-1} f) -
\delta^{-\lambda} V (h_{\delta,x^0}^{-1} f)\|_{L_q(Q)}.
\end{multline*}
Отметим еще, что для $ V \in \mathcal B(L_s(Q),L_q(Q)), $ вследствие
(2.2.2), (2.2.3), норма
\begin{multline*}
\|h_{\delta,x^0} V h_{\delta,x^0}^{-1}\|_{\mathcal B(L_s(I^d),L_q(I^d))} \le
\|h_{\delta,x^0} \|_{\mathcal B(L_q(Q),L_q(I^d))} \\
\| V\|_{\mathcal B(L_s(Q),L_q(Q))} \| h_{\delta,x^0}^{-1}\|_{\mathcal B(L_s(I^d),L_s(Q))} \le \\
\delta^{-q^{-1} \e} \| V\|_{\mathcal B(L_s(Q),L_q(Q))} \delta^{s^{-1} \e} =
\delta^{s^{-1} \e -q^{-1} \e} \| V\|_{\mathcal B(L_s(Q),L_q(Q))}.
\end{multline*}
Учитывая эти обстоятельства, а также (2.2.1), (2.2.7), при $ \rho \in \R_+ $
для $ V \in \mathcal B(L_s(Q),L_q(Q)): \|V\|_{\mathcal B(L_s(Q),L_q(Q))} \le
\delta^\lambda \rho, $ получаем, что
\begin{multline*}
\inf_{\substack{\{ U \in \mathcal B(L_s(I^d),L_q(I^d)):\\
\|U\|_{\mathcal B(L_s(I^d),L_q(I^d))} \le \delta^{\lambda +s^{-1} \e -q^{-1} \e} \rho\}}}
\sup_{f \in ((\mathcal B_{p,\theta}^\alpha)^\prime(\R^d) \cap C_0^\infty(I^d)) \mid_{I^d}}
\| \D^\lambda f -U f\|_{L_q(I^d)} \le \\
\sup_{f \in ((\mathcal B_{p,\theta}^\alpha)^\prime(\R^d) \cap C_0^\infty(I^d)) \mid_{I^d}}
\| \D^\lambda f -(h_{\delta,x^0} V h_{\delta,x^0}^{-1}) f\|_{L_q(I^d)} = \\
\sup_{f \in ((\mathcal B_{p,\theta}^\alpha)^\prime(\R^d) \cap C_0^\infty(I^d)) \mid_{I^d}}
\delta^{\lambda -q^{-1} \e} \| \D^\lambda (h_{\delta,x^0}^{-1} f) -
\delta^{-\lambda} V (h_{\delta,x^0}^{-1} f)\|_{L_q(Q)} = \\
\delta^{\lambda -q^{-1} \e}
\sup_{f \in ((\mathcal B_{p,\theta}^\alpha)^\prime(\R^d) \cap C_0^\infty(I^d)) \mid_{I^d}}
\| \D^\lambda (h_{\delta,x^0}^{-1} f) -
\delta^{-\lambda} V (h_{\delta,x^0}^{-1} f)\|_{L_q(Q)} \le \\
\delta^{\lambda -q^{-1} \e}
\sup_{F \in (c_3 (\mathcal B_{p,\theta}^\alpha)^\prime(\R^d)) \cap C_0^\infty(Q)}
\| \D^\lambda (F \mid_Q) -\delta^{-\lambda} V (F \mid_Q)\|_{L_q(Q)},
\end{multline*}
и, значит,
\begin{multline*} \tag{4.2.3}
c_3 \inf_{\substack{ \mathcal V \in \mathcal B(L_s(Q),L_q(Q)):\\
\|\mathcal V\|_{\mathcal B(L_s(Q),L_q(Q))} \le \rho}}
\sup_{F \in (\mathcal B_{p,\theta}^\alpha)^\prime(\R^d) \cap C_0^\infty(Q)}
\| \D^\lambda (F \mid_Q) -\mathcal V (F \mid_Q)\|_{L_q(Q)} \ge\\
\delta^{-\lambda +q^{-1} \e}
\inf_{\substack{\{ U \in \mathcal B(L_s(I^d),L_q(I^d)):\\
\|U\|_{\mathcal B(L_s(I^d),L_q(I^d))} \le \delta^{\lambda +s^{-1} \e -q^{-1} \e} \rho\}}}
\sup_{f \in ((\mathcal B_{p,\theta}^\alpha)^\prime(\R^d) \cap C_0^\infty(I^d)) \mid_{I^d}}
\| \D^\lambda f -U f\|_{L_q(I^d)}.
\end{multline*}

Пользуясь тем, что для $ F \in C_0^\infty(Q) $ справедливо равенство
$ F \mid_D = (\mathcal I^Q(F \mid_Q)) \mid_D (\mathcal I^Q $ см. в п. 1.3.),
а для $ \mathcal V \in \mathcal B(L_s(D),L_q(D)) $ имеет место неравенство
\begin{multline*}
\| (\mathcal V((\mathcal I^Q f) \mid_D)) \mid_Q \|_{L_q(Q)} \le
\| \mathcal V((\mathcal I^Q f) \mid_D) \|_{L_q(D)} \le\\
\| \mathcal V\|_{\mathcal B(L_s(D),L_q(D))}
\| (\mathcal I^Q f) \mid_D \|_{L_s(D)} = \\
\| \mathcal V\|_{\mathcal B(L_s(D),L_q(D))}
\| ((\mathcal I^Q f) \mid_D) \mid_Q \|_{L_s(Q)} = \\
\| \mathcal V\|_{\mathcal B(L_s(D),L_q(D))} \| f\|_{L_s(Q)},
f \in L_s(Q),
\end{multline*}
для $ \mathcal V \in \mathcal B(L_s(D),L_q(D)):
\| \mathcal V\|_{\mathcal B(L_s(D),L_q(D))} \le \rho, $ выводим
\begin{multline*}
\sup_{f \in (\mathcal B_{p,\theta}^\alpha)^\prime(D)}
\| \D^\lambda f -\mathcal V f\|_{L_q(D)} \ge \\
\sup_{F \in (\mathcal B_{p,\theta}^\alpha)^\prime(\R^d) \cap C_0^\infty(Q)}
\| \D^\lambda (F \mid_D) -\mathcal V (F \mid_D)\|_{L_q(D)} \ge \\
\sup_{F \in (\mathcal B_{p,\theta}^\alpha)^\prime(\R^d) \cap C_0^\infty(Q)}
\| (\D^\lambda (F \mid_D)) \mid_Q -(\mathcal V (F \mid_D)) \mid_Q\|_{L_q(Q)} = \\
\sup_{F \in (\mathcal B_{p,\theta}^\alpha)^\prime(\R^d) \cap C_0^\infty(Q)}
\| \D^\lambda (F \mid_Q) -(\mathcal V ((\mathcal I^Q(F \mid_Q)) \mid_D)) \mid_Q\|_{L_q(Q)} \ge \\
\inf_{ V \in \mathcal B(L_s(Q),L_q(Q)): \|V\|_{\mathcal B(L_s(Q),L_q(Q))} \le \rho}
\sup_{F \in (\mathcal B_{p,\theta}^\alpha)^\prime(\R^d) \cap C_0^\infty(Q)}
\| \D^\lambda (F \mid_Q) -V (F \mid_Q)\|_{L_q(Q)},
\end{multline*}
а, следовательно,
\begin{multline*} \tag{4.2.4}
\inf_{ \mathcal V \in \mathcal B(L_s(D),L_q(D)):
\| \mathcal V\|_{\mathcal B(L_s(D),L_q(D))} \le \rho}
\sup_{f \in (\mathcal B_{p,\theta}^\alpha)^\prime(D)}
\| \D^\lambda f -\mathcal V f\|_{L_q(D)} \ge \\
\inf_{ V \in \mathcal B(L_s(Q),L_q(Q)): \|V\|_{\mathcal B(L_s(Q),L_q(Q))} \le \rho}
\sup_{F \in (\mathcal B_{p,\theta}^\alpha)^\prime(\R^d) \cap C_0^\infty(Q)}
\| \D^\lambda (F \mid_Q) -V (F \mid_Q)\|_{L_q(Q)}.
\end{multline*}

Соединяя (4.2.4), (4.2.3), (4.2.2), заключаем, что при $ \rho > \rho_2 =
\delta^{-\lambda -s^{-1} \e +q^{-1} \e} \rho_1 $ соблюдается неравенство
(4.2.1), где $ A,X,Y,\mathcal K $ имеют тот же смысл, что в формулировке
теоремы 4.1.1. $ \square $
\bigskip

\centerline{\S 5. Поперечники классов $ \D^\lambda(B((H_p^\alpha)^\prime(D))) $ и
$\D^\lambda(B((B_{p,\theta}^\alpha)^\prime(D)))$}
\centerline{в пространстве $ L_q(D)$}
\bigskip

5.1. В этом пункте введ\"eм в рассмотрение подходящие нам пространства 
кусочно-полиномиальных функций и установим некоторые вспомогательные 
утверждения, касающиеся самих пространств и линейных отображений в эти пространства,
с помощью которых проводятся верхние оценки изучаемых поперечников. Но
сначала приведем взятую из [6] лемму 5.1.1.

Лемма 5.1.1

Пусть $ d \in \N, l \in \Z_+^d, 1 \le p \le \infty. $ Тогда существует
константа $ c_1(d,l) >0 $ такая, что для
любого $ \delta \in \R_+^d $ и любого $ x^0 \in \R^d $ для $  Q  =
x^0  +\delta I^d $ можно построить систему линейных операторов
$ P_{\delta,x^0}^{d,l,\lambda}: L_1(Q) \mapsto \mathcal  P^{d,l}, \lambda
\in \Z_+^d(l), $ обладающих следующими свойствами:

1) для $ f \in L_p(Q) $ при $ \lambda \in \Z_+^d(l) $ справедливо неравенство
\begin{equation*} \tag{5.1.1}
\| P_{\delta,x^0}^{d,l,\lambda} f \|_{L_p(Q)} \le c_1 \|f\|_{L_p(Q)};
\end{equation*}

2) для $ f \in \mathcal P^{d,l} $ при  $ \lambda \in \Z_+^d(l) $
имеет место равенство
\begin{equation*} \tag{5.1.2}
P_{\delta,x^0}^{d,l,\lambda}(f \mid_Q) = f;
\end{equation*}

3) если при $ \lambda \in \Z_+^d(l) $  для $ f \in L_1(Q) $ е\"e обобщ\"eнная
производная $ \D^\lambda f \in L_1(Q), $ то соблюдается равенство
\begin{equation*} \tag{5.1.3}
\D^\lambda P_{\delta,x^0}^{d,l,\lambda} f = P_{\delta,x^0}^{d,l -\lambda,0} \D^\lambda f.
\end{equation*}

Введем в рассмотрение следующие пространства кусочно-полиномиальных
функций.

При $ d \in \N, l \in \Z_+^d, \kappa \in \Z_+^d $ и области
$ D \subset \R^d $ положим
\begin{equation*}
N_\kappa^{d,D} = \{\nu \in \Z^d: Q_{\kappa, \nu}^d \cap D \ne \emptyset\}
\end{equation*}
и через $ \mathcal P_\kappa^{d,l,D} $ обозначим линейное пространство,
состоящее из функций $ f: \R^d \mapsto \R, $ для каждой из которых существует
набор полиномов
$ \{f_\nu \in \mathcal P^{d,l}, \nu \in N_\kappa^{d,D}\} $ такой, что
для $ x \in \R^d $ выполняется равенство
\begin{equation*} \tag{5.1.4}
f(x) = \sum_{\nu \in N_\kappa^{d,D}} f_\nu(x) \chi_{\kappa,\nu}^d(x),
\end{equation*}
где $ \chi_{\kappa,\nu}^d $ --- характеристическая функция множества
$ Q_{\kappa,\nu}^d = 2^{-\kappa} \nu +2^{-\kappa} I^d. $

Нетрудно видеть, что при $ d \in \N, l \in \Z_+^d, \kappa \in \Z_+^d $ и
области $ D \subset \R^d $ отображение,
которое каждому набору полиномов
$ \{f_\nu \in \mathcal P^{d, l}, \nu \in N_\kappa^{d,D} \} $ ставит в
соответствие функцию $ f, $ задаваемую равенством (5.1.4), является
изоморфизмом прямого произведения $ \card N_\kappa^{d,D} $ экземпляров
пространства $ \mathcal P^{d,l} $ на пространство
$ \mathcal P_\kappa^{d,l,D}. $
Обозначая через $ R_\kappa^{d,l,D} = \dim \mathcal P_\kappa^{d,l,D} =
(l +\e)^{\e} \card N_\kappa^{d,D}, $
отметим, что для $ d \in \N, l \in \Z_+^d, $ и ограниченной области
$ D \subset \R^d $ существуют константы $ c_2(d,l,D) >0 $ и $ c_3(d,l,D) >0 $
такие, что при $ \kappa \in \Z_+^d $ выполняется неравенство
$$
c_2 2^{(\kappa, \e)} < R_\kappa^{d,l,D} < c_3 2^{(\kappa, \e)}.
$$
При тех же $ d,l,D $ и $ \alpha \in \R_+^d, k \in \Z_+ $ будем обозначать
$ \mathcal P_k^{d,l,D,\alpha} = \mathcal P_\kappa^{d,l,D}, $
$ R_k^{d,l,D,\alpha} = R_\kappa^{d,l,D} $ при $ \kappa = \kappa(k,\alpha). $

Нетрудно убедиться в том, что справедлива

Лемма 5.1.2

Пусть $ d \in \N, l \in \Z_+^d, D $ -- область в $ \R^d, \kappa \in \Z_+^d,
\kappa^\prime \in \Z_+^d: \kappa^\prime \le \kappa. $
Тогда для $ f \in \mathcal P_{\kappa^\prime}^{d,l,D} $ имеет место включение
\begin{equation*} \tag{5.1.5}
\biggl(\sum_{\nu \in N_\kappa^{d,D}} \chi_{\kappa,\nu}^d\biggr) f \in
\mathcal P_\kappa^{d,l,D}.
\end{equation*}

Доказательство.

Прежде всего, заметим, что поскольку при $ \kappa, \kappa^\prime \in \Z_+^d $
для каждого $ \nu \in \Z^d $ соблюдается соотношние
$$
Q_{\kappa,\nu}^d = Q_{\kappa,\nu}^d \cap \R^d = Q_{\kappa,\nu}^d \cap
(\cup_{\nu^\prime \in \Z^d} \overline Q_{\kappa^\prime,\nu^\prime}^d) =
\cup_{\nu^\prime \in \Z^d} (Q_{\kappa,\nu}^d \cap \overline Q_{\kappa^\prime,\nu^\prime}^d)
\ne \emptyset,
$$
то существует $ \nu^\prime \in \Z^d: Q_{\kappa,\nu}^d \cap
\overline Q_{\kappa^\prime,\nu^\prime}^d \ne \emptyset, $ а для такого
$ \nu^\prime $ ввиду открытости $ Q_{\kappa,\nu}^d $ множество
$ Q_{\kappa,\nu}^d \cap Q_{\kappa^\prime,\nu^\prime}^d \ne \emptyset. $

Кроме того, при $ \kappa, \kappa^\prime \in \Z_+^d: \kappa^\prime \le \kappa, $
если $ \nu, \nu^\prime \in \Z^d $ таковы, что $ Q_{\kappa,\nu}^d \cap
Q_{\kappa^\prime,\nu^\prime}^d \ne \emptyset, $ то, выбирая $ x \in
Q_{\kappa,\nu}^d \cap Q_{\kappa^\prime,\nu^\prime}^d, $ имеем
$$
2^{-\kappa} \nu < x < 2^{-\kappa^\prime} \nu^\prime +2^{-\kappa^\prime};
$$
и
$$
2^{-\kappa^\prime} \nu^\prime < x < 2^{-\kappa} \nu +2^{-\kappa},
$$
или
$$
\nu < 2^{\kappa -\kappa^\prime} \nu^\prime +2^{\kappa -\kappa^\prime};
2^{\kappa -\kappa^\prime} \nu^\prime < \nu +\e,
$$
а это значит, что
$$
\nu \le 2^{\kappa -\kappa^\prime} \nu^\prime +2^{\kappa -\kappa^\prime} -\e;
2^{\kappa -\kappa^\prime} \nu^\prime \le \nu,
$$
откуда

$$
2^{-\kappa^\prime} \nu^\prime \le 2^{-\kappa} \nu,
2^{-\kappa} \nu +2^{-\kappa} \le 2^{-\kappa^\prime} \nu^\prime +2^{-\kappa^\prime},
$$
т.е. $ Q_{\kappa,\nu}^d \subset Q_{\kappa^\prime,\nu^\prime}^d. $

Учитывая это обстоятельство, отметим еще, что при $ \kappa, \kappa^\prime
\in \Z_+^d: \kappa^\prime \le \kappa, $ для $ \nu \in N_\kappa^{d,D}, $ если
$ \nu^\prime \in \Z^d: Q_{\kappa,\nu}^d \cap Q_{\kappa^\prime,\nu^\prime}^d
\ne \emptyset, $ то $ \nu^\prime \in N_{\kappa^\prime}^{d,D}, $ ибо в этом
случае
$$
D \cap Q_{\kappa^\prime,\nu^\prime}^d \supset D \cap Q_{\kappa,\nu}^d
\ne \emptyset.
$$
И обратно, при $ \kappa, \kappa^\prime \in \Z_+^d: \kappa^\prime \le \kappa, $
для каждого $ \nu^\prime \in N_{\kappa^\prime}^{d,D} $ существует
$ \nu \in N_\kappa^{d,D}, $ для которого $ Q_{\kappa,\nu}^d \cap
Q_{\kappa^\prime,\nu^\prime}^d \ne \emptyset, $ поскольку, выбирая
$ x \in D \cap Q_{\kappa^\prime,\nu^\prime}^d, $ можно найти $ \nu \in \Z^d, $
такой, что $ x \in \overline Q_{\kappa,\nu}^d, $ что в силу открытости
$ D \cap Q_{\kappa^\prime,\nu^\prime}^d $ влечет непустоту множества
$ D \cap Q_{\kappa^\prime,\nu^\prime}^d \cap Q_{\kappa,\nu}^d. $

На основании сказанного заключаем, что в условиях леммы для
каждого $ \nu \in \Z^d $ существует единственный мультииндекс
$ \nu^\prime(\nu) \in \Z^d $ такой, что
$ Q_{\kappa,\nu}^d \cap Q_{\kappa^\prime,\nu^\prime(\nu)}^d \ne \emptyset, $
причем
\begin{equation*} \tag{5.1.6}
Q_{\kappa,\nu}^d \subset Q_{\kappa^\prime,\nu^\prime(\nu)}^d,
\end{equation*}
и отображение $ \Z^d \ni \nu \mapsto \nu^\prime(\nu) \in \Z^d $ отображает
$ N_\kappa^{d,D} $ на $ N_{\kappa^\prime}^{d,D}. $
Принимая во внимание эти обстоятельства, для $ f = \sum_{\nu^\prime \in
N_{\kappa^\prime}^{d,D}} f_{\nu^\prime} \chi_{\kappa^\prime,\nu^\prime}^d
\in \mathcal P_{\kappa^\prime}^{d,l,D} $ с учетом (5.1.6) имеем
\begin{multline*} \tag{5.1.7}
(\sum_{\nu \in N_\kappa^{d,D}} \chi_{\kappa,\nu}^d) f =
(\sum_{\nu \in N_\kappa^{d,D}} \chi_{\kappa,\nu}^d) (\sum_{\nu^\prime \in
N_{\kappa^\prime}^{d,D}} f_{\nu^\prime} \chi_{\kappa^\prime,\nu^\prime}^d) = \\
\sum_{\nu \in N_\kappa^{d,D}} \chi_{\kappa,\nu}^d (\sum_{\nu^\prime \in
N_{\kappa^\prime}^{d,D}} f_{\nu^\prime} \chi_{\kappa^\prime,\nu^\prime}^d) = 
\sum_{\nu \in N_\kappa^{d,D}} (\sum_{\nu^\prime \in N_{\kappa^\prime}^{d,D}}
f_{\nu^\prime} \chi_{\kappa,\nu}^d \chi_{\kappa^\prime,\nu^\prime}^d) = \\
\sum_{\nu \in N_\kappa^{d,D}} (\sum_{\nu^\prime \in N_{\kappa^\prime}^{d,D}:
Q_{\kappa,\nu}^d \cap Q_{\kappa^\prime,\nu^\prime}^d \ne \emptyset}
f_{\nu^\prime} \chi_{\kappa,\nu}^d \chi_{\kappa^\prime,\nu^\prime}^d) = \\
\sum_{\nu \in N_\kappa^{d,D}} f_{\nu^\prime(\nu)} \chi_{\kappa,\nu}^d
\chi_{\kappa^\prime,\nu^\prime(\nu)}^d = 
\sum_{\nu \in N_\kappa^{d,D}} f_{\nu^\prime(\nu)} \chi_{\kappa,\nu}^d \in
\mathcal P_\kappa^{d,l,D},
\end{multline*}
что завершает доказательство (5.1.5). $ \square $

Для формулировки леммы 5.1.3 введ\"eм следующие обозначения.

Для $ d \in \N, l \in \Z_+^d, \lambda \in \Z_+^d(l), \kappa \in \Z_+^d,
\nu \in \Z^d $ обозначим через $ S_{\kappa,\nu}^{d,l,\lambda}:
L_1(Q_{\kappa,\nu}^d) \mapsto \mathcal  P^{d,l} $ линейный оператор,
определяемый равенством
$$
S_{\kappa,\nu}^{d,l,\lambda} = P_{\delta,x^0}^{d,l,\lambda}
\text{ при} \delta = 2^{-\kappa}, x^0 = 2^{-\kappa} \nu
$$ (см. лемму 5.1.1).

Для $ d \in \N, l \in \Z_+^d, \lambda,  \mu \in \Z_+^d(l), \alpha \in \R_+^d, $
ограниченной области $ D $ в $ \R^d $ при $ k \in \Z_+, $ для котоого при
$ \kappa = \kappa(k,\alpha) $ существует $ \nu \in \Z^d: Q_{\kappa,\nu}^d
\subset D, $ определим непрерывный линейный оператор
$ \mathcal U_{k,\lambda}^{d,l,\mu,D,\alpha}: L_1(D) \mapsto
\mathcal P_k^{d,l -\lambda,D,\alpha} \cap L_\infty(\R^d), $ полагая для
$ f \in L_1(D) $ значение
$$
\mathcal U_{k,\lambda}^{d,l,\mu,D,\alpha} f = \sum_{\nu \in N_\kappa^{d,D}}
(\D^\lambda S_{\kappa,\nu_\kappa^D(\nu)}^{d,l,\mu}
(f \mid_{Q_{\kappa,\nu_\kappa^D(\nu)}^d})) \chi_{\kappa,\nu}^d
$$
при $ \kappa = \kappa(k,\alpha). $

Лемма 5.1.3

Пусть $ d \in \N, \alpha \in \R_+^d, l = l(\alpha), \lambda \in \Z_+^d,
D \subset \R^d $ -- ограниченная область $ \alpha $-типа. Пусть ещ\"e
$ 1 \le p < \infty, 1 \le q \le \infty $ и соблюдается условие (1.3.62).
Тогда существует константа $ c_4(d,\alpha,D,\lambda,p,q) >0 $ такая, что для
любой функции $ f \in (\mathcal H_p^\alpha)^\prime(D) $ при $ k \in \Z_+:
k \ge K^0 $ ($K^0$ см. в п. 1.3.) соблюдается неравенство
\begin{equation*} \tag{5.1.8}
\| \D^\lambda f -\mathcal U_{k,0}^{d,l-\e-\lambda,0,D,\alpha}
\D^\lambda f \|_{L_q(D)} \le c_4 2^{-k(1-(\alpha^{-1}, \lambda
+(p^{-1} -q^{-1})_+ \e))}.
\end{equation*}

Доказательство.

В условиях леммы фиксируем $ m \in \N^d: \lambda \in \Z_+^d(m). $
Тогда для $ f \in (\mathcal H_p^\alpha)^\prime(D) $
при $ k \in \Z_+: k \ge K^0, $ и $ \kappa = \kappa(k,\alpha) $ имеем
\begin{multline*} \tag{5.1.9}
\| \D^\lambda f -\mathcal U_{k,0}^{d,l-\e-\lambda,0,D,\alpha} \D^\lambda f \|_{L_q(D)} \le
\| \D^\lambda f -\D^\lambda E_k^{d,l -\e,m,D,\alpha} f\|_{L_q(D)} \\
+\| \D^\lambda E_k^{d,l -\e,m,D,\alpha} f -
\mathcal U_{k,0}^{d,l-\e-\lambda,0,D,\alpha} \D^\lambda f\|_{L_q(D)}.
\end{multline*}

Оценку первого слагаемого в правой части (5.1.9) да\"eт (1.3.63).
А для оценки второго слагаемого в правой части (5.1.9), учитывая (1.3.23) и тот
факт, что ввиду (1.2.4) выполняется равенство
\begin{equation*}
\sum_{\nu \in N_\kappa^{d,m,D}} (\D^\mu g_{\kappa, \nu}^{d,m}) \mid_D =
\begin{cases} 1, \text{ при } \mu =0; \\ 0, \text{ при } \mu \in \Z_+^d(m)
\setminus \{0\},
\end{cases}
\end{equation*}
имеем
\begin{multline*} \tag{5.1.10}
\| \D^\lambda E_k^{d,l -\e,m,D,\alpha} f -
\mathcal U_{k,0}^{d,l-\e-\lambda,0,D,\alpha} \D^\lambda f\|_{L_q(D)} = \\
\| \D^\lambda (E_\kappa^{d,l -\e,m,D} f) -
\mathcal U_{k,0}^{d,l-\e-\lambda,0,D,\alpha} \D^\lambda f\|_{L_q(D)} = \\
\| \D^\lambda (\sum_{\nu \in N_\kappa^{d,m,D}}
(S_{\kappa, \nu_\kappa^D(\nu)}^{d,l -\e} f) g_{\kappa, \nu}^{d,m}) -
\mathcal U_{k,0}^{d,l-\e-\lambda,0,D,\alpha} \D^\lambda f \|_{L_q(D)} = \\
\| \sum_{\nu \in N_\kappa^{d,m,D}}
\D^\lambda ((S_{\kappa, \nu_\kappa^D(\nu)}^{d,l -\e} f) g_{\kappa, \nu}^{d,m}) -
\mathcal U_{k,0}^{d,l-\e-\lambda,0,D,\alpha} \D^\lambda f \|_{L_q(D)} = \\
\| \sum_{\nu \in N_\kappa^{d,m,D}} \sum_{\mu \in \Z_+^d(\lambda)} C_\lambda^\mu
(\D^\mu (S_{\kappa, \nu_\kappa^D(\nu)}^{d,l -\e} f))
\D^{\lambda -\mu} g_{\kappa, \nu}^{d,m} -
\mathcal U_{k,0}^{d,l-\e-\lambda,0,D,\alpha} \D^\lambda f \|_{L_q(D)} = \\
\| \sum_{\mu \in \Z_+^d(\lambda)} C_\lambda^\mu
\sum_{\nu \in N_\kappa^{d,m,D}}
(\D^\mu (S_{\kappa, \nu_\kappa^D(\nu)}^{d,l -\e} f))
\D^{\lambda -\mu} g_{\kappa, \nu}^{d,m} -
\mathcal U_{k,0}^{d,l-\e-\lambda,0,D,\alpha} \D^\lambda f \|_{L_q(D)} = \\
\biggl\| \sum_{\mu \in \Z_+^d(\lambda)} C_\lambda^\mu
\sum_{\nu \in N_\kappa^{d,m,D}}
(\D^\mu (S_{\kappa, \nu_\kappa^D(\nu)}^{d,l -\e} f))
\D^{\lambda -\mu} g_{\kappa, \nu}^{d,m} -\\
\sum_{\mu \in \Z_+^d(\lambda)} C_\lambda^\mu
(\mathcal U_{k,0}^{d,l-\e-\mu,0,D,\alpha} \D^\mu f)
(\sum_{\nu \in N_\kappa^{d,m,D}}
\D^{\lambda -\mu} g_{\kappa, \nu}^{d,m}) \biggr\|_{L_q(D)} = \\
\biggl\| \sum_{\mu \in \Z_+^d(\lambda)} C_\lambda^\mu
\sum_{\nu \in N_\kappa^{d,m,D}}
(\D^\mu (S_{\kappa, \nu_\kappa^D(\nu)}^{d,l -\e} f))
\D^{\lambda -\mu} g_{\kappa, \nu}^{d,m} -\\
\sum_{\mu \in \Z_+^d(\lambda)} C_\lambda^\mu
\sum_{\nu \in N_\kappa^{d,m,D}}
(\mathcal U_{k,0}^{d,l-\e-\mu,0,D,\alpha} \D^\mu f)
\D^{\lambda -\mu} g_{\kappa, \nu}^{d,m} \biggr\|_{L_q(D)} = \\
\| \sum_{\mu \in \Z_+^d(\lambda)} C_\lambda^\mu
\sum_{\nu \in N_\kappa^{d,m,D}}
(\D^\mu (S_{\kappa, \nu_\kappa^D(\nu)}^{d,l -\e} f) -
(\mathcal U_{k,0}^{d,l-\e-\mu,0,D,\alpha} \D^\mu f))
\D^{\lambda -\mu} g_{\kappa, \nu}^{d,m} \|_{L_q(D)} \le \\
\sum_{\mu \in \Z_+^d(\lambda)} C_\lambda^\mu
\| \sum_{\nu \in N_\kappa^{d,m,D}}
(\D^\mu (S_{\kappa, \nu_\kappa^D(\nu)}^{d,l -\e} f) -
(\mathcal U_{k,0}^{d,l-\e-\mu,0,D,\alpha} \D^\mu f))
\D^{\lambda -\mu} g_{\kappa, \nu}^{d,m} \|_{L_q(D)}.
\end{multline*}

Для проведения оценки правой части (5.1.10) заметим, что в силу (5.1.3)
справедливо равенство
\begin{multline*} \tag{5.1.11}
\mathcal U_{k,0}^{d,l-\e-\mu,0,D,\alpha} \D^\mu f =
\sum_{\nu \in N_\kappa^{d,D}}
(S_{\kappa,\nu_\kappa^D(\nu)}^{d,l-\e-\mu,0}
((\D^\mu f) \mid_{Q_{\kappa,\nu_\kappa^D(\nu)}^d})) \chi_{\kappa,\nu}^d = \\
\sum_{\nu \in N_\kappa^{d,D}}
(\D^\mu S_{\kappa,\nu_\kappa^D(\nu)}^{d,l-\e,\mu}
(f \mid_{Q_{\kappa,\nu_\kappa^D(\nu)}^d})) \chi_{\kappa,\nu}^d = 
\mathcal U_{k,\mu}^{d,l-\e,\mu,D,\alpha} f, \mu \in \Z_+^d(\lambda),
\end{multline*}
и имеет место соотношение
\begin{equation*} \tag{5.1.12}
D = \cup_{\nu \in N_\kappa^{d,D}} (D \cap \overline Q_{\kappa,\nu}^d) =
(\cup_{\nu \in N_\kappa^{d,D}} (D \cap Q_{\kappa,\nu}^d)) \cup A_\kappa^D
\subset (\cup_{\nu \in N_\kappa^{d,D}} Q_{\kappa,\nu}^d) \cup A_\kappa^D,
\end{equation*}
причем множества в правой части (5.1.12) попарно не пересекаются и
$ \mes A_\kappa^D =0. $
Принимая во внимание (5.1.11), (5.1.12),  при $ \mu \in \Z_+^d(\lambda) $
выводим
\begin{multline*} \tag{5.1.13}
\biggl\| \sum_{\nu \in N_\kappa^{d,m,D}}
(\D^\mu (S_{\kappa, \nu_\kappa^D(\nu)}^{d,l -\e} f) -
(\mathcal U_{k,0}^{d,l-\e-\mu,0,D,\alpha} \D^\mu f))
\D^{\lambda -\mu} g_{\kappa, \nu}^{d,m} \biggr\|_{L_q(D)}^q = \\
\int_D \biggl| \sum_{\nu \in N_\kappa^{d,m,D}}
(\D^\mu (S_{\kappa, \nu_\kappa^D(\nu)}^{d,l -\e} f) -
(\mathcal U_{k,\mu}^{d,l-\e,\mu,D,\alpha} f))
\D^{\lambda -\mu} g_{\kappa, \nu}^{d,m} \biggr|^q dx \le \\
\int_{\cup_{n \in N_\kappa^{d,D}} Q_{\kappa,n}^d}
\biggl| \sum_{\nu \in N_\kappa^{d,m,D}}
(\D^\mu (S_{\kappa, \nu_\kappa^D(\nu)}^{d,l -\e} f) -
(\mathcal U_{k,\mu}^{d,l-\e,\mu,D,\alpha} f))
\D^{\lambda -\mu} g_{\kappa, \nu}^{d,m} \biggr|^q dx = \\
\sum_{n \in N_\kappa^{d,D}} \int_{Q_{\kappa,n}^d}
\biggl| \sum_{\nu \in N_\kappa^{d,m,D}}
(\D^\mu (S_{\kappa, \nu_\kappa^D(\nu)}^{d,l -\e} f) -
(\mathcal U_{k,\mu}^{d,l-\e,\mu,D,\alpha} f))
\D^{\lambda -\mu} g_{\kappa, \nu}^{d,m} \biggr|^q dx = \\
\sum_{n \in N_\kappa^{d,D}} \int_{Q_{\kappa,n}^d}
\biggl| \sum_{\nu \in N_\kappa^{d,m,D}}
(\D^\mu (S_{\kappa, \nu_\kappa^D(\nu)}^{d,l -\e} f) -
\D^\mu (S_{\kappa,\nu_\kappa^D(n)}^{d,l-\e,\mu} f))
\D^{\lambda -\mu} g_{\kappa, \nu}^{d,m} \biggr|^q dx = \\
\sum_{n \in N_\kappa^{d,D}} \int_{Q_{\kappa,n}^d}
\biggl| \sum_{\nu \in N_\kappa^{d,m,D}: Q_{\kappa,n}^d \cap
\supp g_{\kappa,\nu}^{d,m} \ne \emptyset}
(\D^\mu (S_{\kappa, \nu_\kappa^D(\nu)}^{d,l -\e} f) -
\D^\mu (S_{\kappa,\nu_\kappa^D(n)}^{d,l-\e,\mu} f))
\D^{\lambda -\mu} g_{\kappa, \nu}^{d,m} \biggr|^q dx = \\
\sum_{n \in N_\kappa^{d,D}} \biggl\| \sum_{\nu \in N_\kappa^{d,m,D}: Q_{\kappa,n}^d
\cap \supp g_{\kappa,\nu}^{d,m} \ne \emptyset}
(\D^\mu (S_{\kappa, \nu_\kappa^D(\nu)}^{d,l -\e} f) -
\D^\mu (S_{\kappa,\nu_\kappa^D(n)}^{d,l-\e,\mu} f))
\D^{\lambda -\mu} g_{\kappa, \nu}^{d,m} \biggr\|_{L_q(Q_{\kappa,n}^d)}^q \le \\
\sum_{n \in N_\kappa^{d,D}} \biggl(\sum_{\nu \in N_\kappa^{d,m,D}: Q_{\kappa,n}^d
\cap \supp g_{\kappa,\nu}^{d,m} \ne \emptyset}
\biggl\|(\D^\mu (S_{\kappa, \nu_\kappa^D(\nu)}^{d,l -\e} f) -
\D^\mu (S_{\kappa,\nu_\kappa^D(n)}^{d,l-\e,\mu} f))
\D^{\lambda -\mu} g_{\kappa, \nu}^{d,m} \biggr\|_{L_q(Q_{\kappa,n}^d)}\biggr)^q.
\end{multline*}

Используя (1.2.5), (1.1.1), для $ \mu \in \Z_+^d(\lambda), n \in
N_\kappa^{d,D}, \nu \in N_\kappa^{d,m,D}: Q_{\kappa,n}^d \cap
\supp g_{\kappa,\nu}^{d,m} \ne \emptyset, $ находим, что
\begin{multline*} \tag{5.1.14}
\|(\D^\mu (S_{\kappa, \nu_\kappa^D(\nu)}^{d,l -\e} f) -
\D^\mu (S_{\kappa,\nu_\kappa^D(n)}^{d,l-\e,\mu} f))
\D^{\lambda -\mu} g_{\kappa, \nu}^{d,m} \|_{L_q(Q_{\kappa,n}^d)} \le \\
\| \D^{\lambda -\mu} g_{\kappa, \nu}^{d,m} \|_{L_\infty(\R^d)}
\| \D^\mu ((S_{\kappa, \nu_\kappa^D(\nu)}^{d,l -\e} f) -
(S_{\kappa,\nu_\kappa^D(n)}^{d,l-\e,\mu} f))
\|_{L_q(Q_{\kappa,n}^d)} \le \\
c_5 2^{(\kappa, \lambda -\mu)} c_6 2^{(\kappa, \mu +p^{-1} \e -q^{-1} \e)}
\| S_{\kappa, \nu_\kappa^D(\nu)}^{d,l -\e} f -
S_{\kappa,\nu_\kappa^D(n)}^{d,l-\e,\mu} f \|_{L_p(Q_{\kappa,n}^d)} = \\
c_7 2^{(\kappa, \lambda +p^{-1} \e -q^{-1} \e)}
\| S_{\kappa, \nu_\kappa^D(\nu)}^{d,l -\e} f -
S_{\kappa,\nu_\kappa^D(n)}^{d,l-\e,\mu} f \|_{L_p(Q_{\kappa,n}^d)}.
\end{multline*}

Оценим норму в правой части (5.1.14). Фиксировав $ n \in N_\kappa^{d,D},
\nu \in N_\kappa^{d,m,D}: Q_{\kappa,n}^d \cap \supp g_{\kappa, \nu}^{d,m} \ne
\emptyset, $ для $ Q_{\kappa, \nu_\kappa^D(n)}^d \subset D $
и $ Q_{\kappa, \nu_\kappa^D(\nu)}^d \subset D $ выберем
последовательности $ \nu^\iota \in \Z^d, j^\iota \in \Nu_{1,d}^1, \epsilon^\iota
\in \{-1, 1\}, \iota =0,\ldots,\Iota, $ для которых $ \nu^0 = \nu_\kappa^D(\nu),
\nu^{\Iota} = \nu_\kappa^D(n) $ и соблюдаются соотношения (1.3.31),
(1.3.32), (1.3.33). Тогда
\begin{multline*} \tag{5.1.15}
\| S_{\kappa, \nu_\kappa^D(\nu)}^{d,l -\e} f -
S_{\kappa, \nu_\kappa^D(n)}^{d,l -\e,\mu} f \|_{L_p(Q_{\kappa,n}^d)} = 
\| S_{\kappa, \nu^0}^{d,l -\e} f -
S_{\kappa, \nu_\kappa^D(n)}^{d,l -\e,\mu} f\|_{L_p(Q_{\kappa,n}^d)} = \\
\| S_{\kappa, \nu^0}^{d,l -\e} f -\sum_{\iota =1}^{\Iota}
(S_{\kappa, \nu^\iota}^{d,l -\e} f -S_{\kappa, \nu^\iota}^{d,l -\e} f) -
S_{\kappa, \nu_\kappa^D(n)}^{d,l -\e,\mu} f\|_{L_p(Q_{\kappa,n}^d)} = \\
\| \sum_{\iota =0}^{\Iota -1} (S_{\kappa, \nu^\iota}^{d,l -\e} f -
S_{\kappa, \nu^{\iota +1}}^{d,l -\e} f) +
S_{\kappa, \nu^{\Iota}}^{d,l -\e} f -
S_{\kappa, \nu_\kappa^D(n)}^{d,l -\e,\mu} f\|_{L_p(Q_{\kappa,n}^d)} \le \\
\sum_{\iota =0}^{\Iota -1} \| S_{\kappa, \nu^\iota}^{d,l -\e} f -
S_{\kappa, \nu^{\iota +1}}^{d,l -\e} f\|_{L_p(Q_{\kappa,n}^d)} +
\| S_{\kappa, \nu^{\Iota}}^{d,l -\e} f -
S_{\kappa, \nu_\kappa^D(n)}^{d,l -\e,\mu} f\|_{L_p(Q_{\kappa,n}^d)}.
\end{multline*}

Для проведения оценки слагаемых в правой части (5.1.15) приведем некоторые
полезные для нас факты. Прежде всего, заметим, что $ N_\kappa^{d,D} \subset
N_\kappa^{d,m,D} \cap \{ \nu^\prime \in \Z^d: Q_{\kappa, \nu^\prime}^d
\cap G_\kappa^{d,m,D} \ne \emptyset\} $ (см. п. 1.3.).
С учетом этого замечания отметим, что вследствие (1.3.31), (1.3.13), (1.3.12)
имеет место оценка
\begin{multline*} \tag{5.1.16}
\Iota \le c_{8} \| \nu^0 -\nu^{\Iota} \|_{l_\infty^d} = c_{8} \| \nu_\kappa^D(\nu) -
\nu_\kappa^D(n) \|_{l_\infty^d}  \le \\
 c_8 (\| \nu_\kappa^D(\nu) -n \|_{l_\infty^d} +
\| n -\nu_\kappa^D(n) \|_{l_\infty^d}) \le \\
c_8(\| \gamma^1\|_{l_\infty^d} +c_9) =
c_{10}(d,m,D,\alpha).
\end{multline*}

Далее, при $ \iota =0,\ldots,\Iota, j =1,\ldots,d $ для $ x \in
\overline Q_{\kappa,\nu^\iota}^d, $ используя (5.1.16), (1.3.32) и (1.3.13),
заключаем, что справедливо неравенство
\begin{multline*}
| x_j -2^{-\kappa_j} n_j | \le | x_j -2^{-\kappa_j} (\nu^\iota)_j | +
| 2^{-\kappa_j} (\nu^\iota)_j -2^{-\kappa_j} n_j | \le \\
 2^{-\kappa_j} +
| 2^{-\kappa_j} (\nu^\iota)_j -2^{-\kappa_j} n_j | = 2^{-\kappa_j} +
| 2^{-\kappa_j} (\nu^\iota)_j -\sum_{i =0}^{\iota -1} (2^{-\kappa_j} (\nu^i)_j
-2^{-\kappa_j} (\nu^i)_j) -2^{-\kappa_j} n_j | = \\
2^{-\kappa_j} +| \sum_{i =0}^{\iota -1} (2^{-\kappa_j} (\nu^{i +1})_j
-2^{-\kappa_j} (\nu^i)_j) +2^{-\kappa_j} (\nu^0)_j -2^{-\kappa_j} n_j | \le \\
2^{-\kappa_j} +\sum_{i =0}^{\iota -1} | 2^{-\kappa_j} (\nu^{i +1})_j
-2^{-\kappa_j} (\nu^i)_j| +| 2^{-\kappa_j} (\nu^0)_j -2^{-\kappa_j} n_j | \le \\
2^{-\kappa_j} +\sum_{i =0}^{\iota -1} 2^{-\kappa_j} \| \nu^{i +1} -\nu^i \|_{l_\infty^d}
+| 2^{-\kappa_j} (\nu^0)_j -2^{-\kappa_j} n_j | = \\
2^{-\kappa_j} +\sum_{i =0}^{\iota -1} 2^{-\kappa_j} \| \epsilon^i e_{j^i} \|_{l_\infty^d}
+| 2^{-\kappa_j} (\nu_\kappa^D(\nu))_j  -2^{-\kappa_j} n_j | = \\
2^{-\kappa_j} +\sum_{i =0}^{\iota -1} 2^{-\kappa_j}
+| 2^{-\kappa_j} (\nu_\kappa^D(\nu))_j  -2^{-\kappa_j} n_j | = \\
2^{-\kappa_j} (\iota +1) +| 2^{-\kappa_j} (\nu_\kappa^D(\nu))_j -
2^{-\kappa_j} n_j | \le \\
2^{-\kappa_j} (\Iota +1) +\gamma_j^1 2^{-\kappa_j} \le 
2^{-\kappa_j} (c_{10} +1) +\gamma_j^1 2^{-\kappa_j} = \gamma_j^2 2^{-\kappa_j}.
\end{multline*}
Отсюда, учитывая (1.3.33), видим, что существует вектор
$ \gamma^2(d,m,D,\alpha) \in \R_+^d $ такой, что для
$$
\mathfrak D_{\kappa,n}^{d,m,D,\alpha} = 2^{-\kappa} n +\gamma^2 2^{-\kappa} B^d,
n \in N_\kappa^{d,D},
$$
соблюдается включение
\begin{equation*} \tag{5.1.17}
Q_{\kappa,\nu^\iota}^d \cup Q_{\kappa,\nu^{\iota +1}}^d \subset Q^\iota \subset
D \cap \mathfrak D_{\kappa,n}^{d,m,D,\alpha}, \iota =0,\ldots,\Iota -1.
\end{equation*}

Из приведенного определения с учетом того, что $ \gamma^2 > \e, $ следует,
что при $ n \in N_\kappa^{d,D} $ верно включение
\begin{equation*} \tag{5.1.18}
Q_{\kappa, n}^d \subset \mathfrak D_{\kappa, n}^{d,m,D,\alpha},
\end{equation*}
а из (5.1.17), (5.1.18) вытекает, что
\begin{equation*} \tag{5.1.19}
Q_{\kappa,n}^d \subset (2^{-\kappa} \nu^{\iota} +(\gamma^2 +\e) 2^{-\kappa} B^d),
\iota =0,\ldots,\Iota.
\end{equation*}

Учитывая (5.1.18), нетрудно видеть, что существует константа
$ c_{11}(d,m,D,\alpha) >0 $ такая, что для каждого $ x \in \R^d $ число
\begin{equation*} \tag{5.1.20}
\card \{ n \in N_\kappa^{d,D}:
x \in \mathfrak D_{\kappa,n}^{d,m,D,\alpha} \} \le c_{11}.
\end{equation*}

Принимая во внимание, что в силу (1.3.32), (1.3.33) при
$ \iota =0,\ldots,\Iota -1 $ выполняется равенство
\begin{equation*}
Q^\iota = \begin{cases} 2^{-\kappa} \nu^\iota +2^{-(\kappa -e_{j^\iota})} I^d,
\text{ при } \epsilon^\iota =1; \\
2^{-\kappa} \nu^{\iota +1} +2^{-(\kappa -e_{j^\iota})} I^d,
\text{ при } \epsilon^\iota =-1,
\end{cases}
\end{equation*}
определим линейный оператор $ S^\iota: L_1(Q^\iota) \mapsto \mathcal P^{d,l-\e}, $
полагая
$$
S^\iota = P_{\delta,x^0}^{d,l -\e}
$$
при $ \delta = 2^{-(\kappa -e_{j^\iota})}, $
\begin{equation*}
x^0 = \begin{cases} 2^{-\kappa} \nu^\iota, \text{ при } \epsilon^\iota =1; \\
2^{-\kappa} \nu^{\iota +1}, \text{ при } \epsilon^\iota =-1,
\end{cases}
\iota =0,\ldots,\Iota -1.
\end{equation*}

Теперь проведем оценку слагаемых в правой части (5.1.15).
При $ \iota =0,\ldots,\Iota -1, $ благодаря (5.1.19) применяя (1.1.1),
а затем используя (1.1.2), (1.1.3), (5.1.17), и, наконец, пользуясь (1.1.4)
и учитывая (5.1.17), приходим к неравенству
\begin{multline*} \tag{5.1.21}
\| S_{\kappa, \nu^\iota}^{d,l -\e} f -
S_{\kappa, \nu^{\iota +1}}^{d,l -\e} f\|_{L_p(Q_{\kappa,n}^d)} = \\
\| S_{\kappa, \nu^\iota}^{d,l -\e} f -S^\iota f +S^\iota f -
S_{\kappa, \nu^{\iota +1}}^{d,l -\e} f\|_{L_p(Q_{\kappa,n}^d)} \le \\
\| S_{\kappa, \nu^\iota}^{d,l -\e} f -S^\iota f \|_{L_p(Q_{\kappa,n}^d)}
+\| S^\iota f -S_{\kappa, \nu^{\iota +1}}^{d,l -\e} f\|_{L_p(Q_{\kappa,n}^d)} \le \\
c_{12} \| S_{\kappa, \nu^\iota}^{d,l -\e} f -S^\iota f \|_{L_p(Q_{\kappa,\nu^\iota}^d)}
+c_{12} \| S^\iota f -S_{\kappa, \nu^{\iota +1}}^{d,l -\e} f
\|_{L_p(Q_{\kappa,\nu^{\iota +1}}^d)} = \\
c_{12} \| S_{\kappa, \nu^\iota}^{d,l -\e} (f -S^\iota f)
\|_{L_p(Q_{\kappa,\nu^\iota}^d)}
+c_{12} \| S_{\kappa, \nu^{\iota +1}}^{d,l -\e} (f -S^\iota f)
\|_{L_p(Q_{\kappa,\nu^{\iota +1}}^d)} \le \\
c_{13} \| f -S^\iota f \|_{L_p(Q_{\kappa,\nu^\iota}^d)}
+c_{13} \| f -S^\iota f\|_{L_p(Q_{\kappa,\nu^{\iota +1}}^d)} \le 
c_{14} \| f -S^\iota f \|_{L_p(Q^\iota)} \le \\
 c_{15} \sum_{j =1}^d
((2^{-(\kappa -e_{j^\iota})})_j)^{-1/p} \biggl(\int_{(2^{-(\kappa -e_{j^\iota})})_j B^1}
\int_{Q_{l_j \xi e_j}^\iota}| \Delta_{\xi e_j}^{l_j} f(x)|^p dx d\xi\biggr)^{1/p} \le \\
c_{16} \sum_{j =1}^d 2^{\kappa_j /p} \biggl(\int_{2^{-\kappa_j +1} B^1}
\int_{Q_{l_j \xi e_j}^\iota}| \Delta_{\xi e_j}^{l_j} f(x)|^p dx d\xi\biggr)^{1/p} \le \\
c_{16} \sum_{j =1}^d 2^{\kappa_j /p} \biggl(\int_{2^{-\kappa_j +1} B^1}
\int_{(D \cap \mathfrak D_{\kappa,n}^{d,m,D,\alpha})_{l_j \xi e_j}}
| \Delta_{\xi e_j}^{l_j} f(x)|^p dx d\xi\biggr)^{1/p} \le \\
c_{16} \sum_{j =1}^d 2^{\kappa_j /p} \biggl(\int_{2^{-\kappa_j +1} B^1}
\int_{D_{l_j \xi e_j} \cap \mathfrak D_{\kappa,n}^{d,m,D,\alpha}}
| \Delta_{\xi e_j}^{l_j} f(x)|^p dx d\xi\biggr)^{1/p}.
\end{multline*}

Для оценки последнего слагаемого в правой части (5.1.15) при $ \mu \in
\Z_+^d(\lambda), $ ввиду (5.1.19)
используя (1.1.1), затем пользуясь (5.1.2), (5.1.1) и применяя (1.1.4),
с учетом (5.1.17) получаем
\begin{multline*} \tag{5.1.22}
\| S_{\kappa, \nu^{\Iota}}^{d,l -\e} f -
S_{\kappa, \nu_\kappa^D(n)}^{d,l -\e,\mu} f\|_{L_p(Q_{\kappa,n}^d)} = 
\| S_{\kappa, \nu^{\Iota}}^{d,l -\e} f -
S_{\kappa, \nu^{\Iota}}^{d,l -\e,\mu} f\|_{L_p(Q_{\kappa,n}^d)} \le \\
c_{17} \| S_{\kappa, \nu^{\Iota}}^{d,l -\e} f -
S_{\kappa, \nu^{\Iota}}^{d,l -\e,\mu} f\|_{L_p(Q_{\kappa,\nu^{\Iota}}^d)} = 
c_{17} \| S_{\kappa, \nu^{\Iota}}^{d,l -\e,\mu} (f -
S_{\kappa, \nu^{\Iota}}^{d,l -\e} f)\|_{L_p(Q_{\kappa,\nu^{\Iota}}^d)} \le \\
c_{18} \| f -S_{\kappa, \nu^{\Iota}}^{d,l -\e} f\|_{L_p(Q_{\kappa,\nu^{\Iota}}^d)} \le 
c_{19} \sum_{j =1}^d 2^{\kappa_j /p} \biggl(\int_{ 2^{-\kappa_j} B^1}
\int_{(Q_{\kappa,\nu^{\Iota}}^d)_{l_j \xi e_j}}
| \Delta_{\xi e_j}^{l_j} f(x)|^p dx d\xi\biggr)^{1/p} \le \\
c_{19} \sum_{j =1}^d 2^{\kappa_j /p} \biggl(\int_{ 2^{-\kappa_j} B^1}
\int_{(D \cap \mathfrak D_{\kappa,n}^{d,m,D,\alpha})_{l_j \xi e_j}}
| \Delta_{\xi e_j}^{l_j} f(x)|^p dx d\xi\biggr)^{1/p} \le \\
c_{19} \sum_{j =1}^d 2^{\kappa_j /p} \biggl(\int_{ 2^{-\kappa_j} B^1}
\int_{D_{l_j \xi e_j} \cap \mathfrak D_{\kappa,n}^{d,m,D,\alpha}}
| \Delta_{\xi e_j}^{l_j} f(x)|^p dx d\xi\biggr)^{1/p}.
\end{multline*}
Объединяя (5.1.15), (5.1.21), (5.1.22) и учитывая (5.1.16), при $ \mu \in
\Z_+^d(\lambda), n \in N_\kappa^{d,D}, \nu \in N_\kappa^{d,m,D}: Q_{\kappa,n}^d
\cap \supp g_{\kappa, \nu}^{d,m} \ne \emptyset, $ имеем
\begin{multline*}
\| S_{\kappa, \nu_\kappa^D(\nu)}^{d,l -\e} f -
S_{\kappa, \nu_\kappa^D(n)}^{d,l -\e,\mu} f \|_{L_p(Q_{\kappa,n}^d)} \le \\
\Iota c_{16} \sum_{j =1}^d 2^{\kappa_j /p} \biggl(\int_{2^{-\kappa_j +1} B^1}
\int_{D_{l_j \xi e_j} \cap \mathfrak D_{\kappa,n}^{d,m,D,\alpha}}
| \Delta_{\xi e_j}^{l_j} f(x)|^p dx d\xi\biggr)^{1/p} +\\
c_{19} \sum_{j =1}^d 2^{\kappa_j /p} \biggl(\int_{2^{-\kappa_j} B^1}
\int_{D_{l_j \xi e_j} \cap \mathfrak D_{\kappa,n}^{d,m,D,\alpha}}
| \Delta_{\xi e_j}^{l_j} f(x)|^p dx d\xi\biggr)^{1/p} \le \\
c_{10} c_{16} \sum_{j =1}^d 2^{\kappa_j /p} \biggl(\int_{ 2^{-\kappa_j +1} B^1}
\int_{D_{l_j \xi e_j} \cap \mathfrak D_{\kappa,n}^{d,m,D,\alpha}}
| \Delta_{\xi e_j}^{l_j} f(x)|^p dx d\xi\biggr)^{1/p} +\\
c_{19} \sum_{j =1}^d 2^{\kappa_j /p} \biggl(\int_{ 2^{-\kappa_j +1} B^1}
\int_{D_{l_j \xi e_j} \cap \mathfrak D_{\kappa,n}^{d,m,D,\alpha}}
| \Delta_{\xi e_j}^{l_j} f(x)|^p dx d\xi\biggr)^{1/p} = \\
c_{20} \sum_{j =1}^d 2^{\kappa_j /p} \biggl(\int_{ 2^{-\kappa_j +1} B^1}
\int_{D_{l_j \xi e_j} \cap \mathfrak D_{\kappa,n}^{d,m,D,\alpha}}
| \Delta_{\xi e_j}^{l_j} f(x)|^p dx d\xi\biggr)^{1/p}.
\end{multline*}

Подставляя эту оценку в (5.1.14), для $ \mu \in \Z_+^d(\lambda),
n \in N_\kappa^{d,D}, \nu \in N_\kappa^{d,m,D}: Q_{\kappa,n}^d \cap
\supp g_{\kappa,\nu}^{d,m} \ne \emptyset, $ получаем неравенство
\begin{multline*} \tag{5.1.23}
\|(\D^\mu (S_{\kappa, \nu_\kappa^D(\nu)}^{d,l -\e} f) -
\D^\mu (S_{\kappa,\nu_\kappa^D(n)}^{d,l-\e,\mu} f))
\D^{\lambda -\mu} g_{\kappa, \nu}^{d,m} \|_{L_q(Q_{\kappa,n}^d)} \le \\
c_{21} 2^{(\kappa, \lambda +p^{-1} \e -q^{-1} \e)}
\sum_{j =1}^d 2^{\kappa_j /p} \biggl(\int_{ 2^{-\kappa_j +1} B^1}
\int_{D_{l_j \xi e_j} \cap \mathfrak D_{\kappa,n}^{d,m,D,\alpha}}
| \Delta_{\xi e_j}^{l_j} f(x)|^p dx d\xi\biggr)^{1/p}.
\end{multline*}

Соединяя (5.1.13) с (5.1.23), и учитывая (1.3.17), а затем применяя
неравенство Г\"eльдера с показателем $ q $ и неравенство Г\"eльдера с показателем
$ p/q \le 1, $ и, наконец, принимая во внимание (5.1.20) и используя
неравенство Г\"eльдера с показателем $ 1/q, $ при $ \mu \in \Z_+^d(\lambda) $
выводим
\begin{multline*} \tag{5.1.24}
\| \sum_{\nu \in N_\kappa^{d,m,D}}
(\D^\mu (S_{\kappa, \nu_\kappa^D(\nu)}^{d,l -\e} f) -
(\mathcal U_{k,0}^{d,l-\e-\mu,0,D,\alpha} \D^\mu f))
\D^{\lambda -\mu} g_{\kappa, \nu}^{d,m} \|_{L_q(D)}^q \le \\
\sum_{n \in N_\kappa^{d,D}} \biggl(\sum_{\substack{\nu \in N_\kappa^{d,m,D}:\\
 Q_{\kappa,n}^d \cap \supp g_{\kappa, \nu}^{d,m} \ne \emptyset}}
c_{21} 2^{(\kappa, \lambda +p^{-1} \e -q^{-1} \e)}\times\\
\sum_{j =1}^d 2^{\kappa_j /p} \biggl(\int_{2^{-\kappa_j +1} B^1}
\int_{D_{l_j \xi e_j} \cap \mathfrak D_{\kappa,n}^{d,m,D,\alpha}}
| \Delta_{\xi e_j}^{l_j} f(x)|^p dx d\xi\biggr)^{1/p}\biggr)^q \le \\
\sum_{n \in N_\kappa^{d,D}}
(c_{22} c_{21} 2^{(\kappa, \lambda +p^{-1} \e -q^{-1} \e)}
\sum_{j =1}^d 2^{\kappa_j /p} (\int_{2^{-\kappa_j +1} B^1}
\int_{D_{l_j \xi e_j} \cap \mathfrak D_{\kappa,n}^{d,m,D,\alpha}}
| \Delta_{\xi e_j}^{l_j} f(x)|^p dx d\xi)^{1/p})^q = \\
(c_{23} 2^{(\kappa, \lambda +p^{-1} \e -q^{-1} \e)})^q
\sum_{n \in N_\kappa^{d,D}}
\biggl(\sum_{j =1}^d 2^{\kappa_j /p} \biggl(\int_{2^{-\kappa_j +1} B^1}
\int_{D_{l_j \xi e_j} \cap \mathfrak D_{\kappa,n}^{d,m,D,\alpha}}
| \Delta_{\xi e_j}^{l_j} f(x)|^p dx d\xi\biggr)^{1/p}\biggr)^q \le \\
(c_{24} 2^{(\kappa, \lambda +p^{-1} \e -q^{-1} \e)})^q
\sum_{n \in N_\kappa^{d,D}}
\sum_{j =1}^d 2^{\kappa_j q /p} \biggl(\int_{2^{-\kappa_j +1} B^1}
\int_{D_{l_j \xi e_j} \cap \mathfrak D_{\kappa,n}^{d,m,D,\alpha}}
| \Delta_{\xi e_j}^{l_j} f(x)|^p dx d\xi\biggr)^{q/p} = \\
(c_{24} 2^{(\kappa, \lambda +p^{-1} \e -q^{-1} \e)})^q
\sum_{j =1}^d 2^{\kappa_j q /p}
\sum_{n \in N_\kappa^{d,D}}
\biggl(\int_{2^{-\kappa_j +1} B^1}
\int_{D_{l_j \xi e_j} \cap \mathfrak D_{\kappa,n}^{d,m,D,\alpha}}
| \Delta_{\xi e_j}^{l_j} f(x)|^p dx d\xi\biggr)^{q/p} \le \\
(c_{24} 2^{(\kappa, \lambda +p^{-1} \e -q^{-1} \e)})^q
\sum_{j =1}^d 2^{\kappa_j q /p}
\biggl(\sum_{n \in N_\kappa^{d,D}}
\int_{2^{-\kappa_j +1} B^1}
\int_{D_{l_j \xi e_j} \cap \mathfrak D_{\kappa,n}^{d,m,D,\alpha}}
| \Delta_{\xi e_j}^{l_j} f(x)|^p dx d\xi\biggr)^{q/p} = \\
(c_{24} 2^{(\kappa, \lambda +p^{-1} \e -q^{-1} \e)})^q
\sum_{j =1}^d 2^{\kappa_j q /p}
\biggl(\int_{2^{-\kappa_j +1} B^1}
\sum_{n \in N_\kappa^{d,D}}
\int_{D_{l_j \xi e_j}} \chi_{\mathfrak D_{\kappa,n}^{d,m,D,\alpha}}(x)
| \Delta_{\xi e_j}^{l_j} f(x)|^p dx d\xi\biggr)^{q/p} = \\
(c_{24} 2^{(\kappa, \lambda +p^{-1} \e -q^{-1} \e)})^q
\sum_{j =1}^d 2^{\kappa_j q /p}
\biggl(\int_{2^{-\kappa_j +1} B^1}
\int_{D_{l_j \xi e_j}}
\biggl(\sum_{n \in N_\kappa^{d,D}}
\chi_{\mathfrak D_{\kappa,n}^{d,m,D,\alpha}}(x)\biggr)
| \Delta_{\xi e_j}^{l_j} f(x)|^p dx d\xi\biggr)^{q/p} \le \\
(c_{24} 2^{(\kappa, \lambda +p^{-1} \e -q^{-1} \e)})^q
\sum_{j =1}^d 2^{\kappa_j q /p}
\biggl(\int_{2^{-\kappa_j +1} B^1}
\int_{D_{l_j \xi e_j}} c_{11}
| \Delta_{\xi e_j}^{l_j} f(x)|^p dx d\xi\biggr)^{q/p} \le \\
(c_{25} 2^{(\kappa, \lambda +p^{-1} \e -q^{-1} \e)})^q
\biggl(\sum_{j =1}^d 2^{\kappa_j /p}
\biggl(\int_{2^{-\kappa_j +1} B^1}
\int_{D_{l_j \xi e_j}} | \Delta_{\xi e_j}^{l_j} f(x)|^p dx d\xi\biggr)^{1/p}\biggr)^q = \\
(c_{25} 2^{(\kappa, \lambda +p^{-1} \e -q^{-1} \e)}
\sum_{j =1}^d 2^{\kappa_j /p}
\biggl(\int_{2^{-\kappa_j +1} B^1}
\int_{D_{l_j \xi e_j}} | \Delta_{\xi e_j}^{l_j} f(x)|^p dx d\xi\biggr)^{1/p})^q \le \\
\biggl( c_{25} 2^{k (\alpha^{-1}, \lambda  +(p^{-1} -q^{-1})_+ \e)}
\sum_{j =1}^d 2^{k /(\alpha_j p)}
\biggl(\int_{4 2^{-k /\alpha_j} B^1}\int_{D_{l_j \xi e_j}} 
| \Delta_{\xi e_j}^{l_j} f(x)|^p dx d\xi\biggr)^{1/p}\biggr)^q \le \\
(c_{26} 2^{k (\alpha^{-1}, \lambda  +(p^{-1} -q^{-1})_+ \e)}
\sum_{j =1}^d \Omega_j^{\prime l_j}(f, 4 2^{-k /\alpha_j})_{L_p(D)})^q.
\end{multline*}

Используя (5.1.24) при оценке правой части (5.1.10), приходим к неравенству
\begin{multline*} \tag{5.1.25}
\| \D^\lambda E_k^{d,l -\e,m,D,\alpha} f -
\mathcal U_{k,0}^{d,l-\e-\lambda,0,D,\alpha} \D^\lambda f\|_{L_q(D)} \le \\
\sum_{\mu \in \Z_+^d(\lambda)} C_\lambda^\mu
c_{26} 2^{k (\alpha^{-1}, \lambda +(p^{-1} -q^{-1})_+ \e)}
\sum_{j =1}^d \Omega_j^{\prime l_j}(f, 4 2^{-k /\alpha_j})_{L_p(D)} \le \\
c_{27} 2^{-k(1 -(\alpha^{-1}, \lambda  +(p^{-1} -q^{-1})_+ \e))}.
\end{multline*}
Соединяя (5.1.9) с (1.3.63) и (5.1.25), получаем (5.1.8) при $ p \le q. $
Неравенство (5.1.8) при $ q < p $ вытекает из неравенства (5.1.8) при
$ q = p $ и того факта, что при $ q < p $ для $ h \in L_p(D) $ и ограниченной
области $ D $ справедливо неравенство
$ \|h\|_{L_q(D)} \le c_{28} \|h\|_{L_p(D)}. \square$

Непосредственным следствием леммы 5.1.3 является одно из утверждений леммы
5.1.4, для формулировки которой потребуются следующие обозначения.

Для $ d \in \N, l \in \Z_+^d, \alpha \in \R_+^d, $ ограниченной области $ D $
в $ \R^d $ при $ k \in \Z_+, $ для которого при
$ \kappa = \kappa(k,\alpha), \kappa^\prime = \kappa(k -1,\alpha) $
существуют $ \nu \in \Z^d: Q_{\kappa,\nu}^d \subset D, $ и
$ \nu^\prime \in \Z^d: Q_{\kappa^\prime,\nu^\prime}^d \subset D, $ с учетом (5.1.5) определим непрерывный линейный оператор
$ V_k^{d,l,D,\alpha}: L_1(D) \mapsto \mathcal P_k^{d,l,D,\alpha} \cap
L_\infty(\R^d) $ равенством
\begin{equation*}
V_k^{d,l,D,\alpha} = \mathcal U_{k,0}^{d,l,0,D,\alpha} -
(\sum_{\nu \in N_\kappa^{d,D}} \chi_{\kappa,\nu}^d)
\mathcal U_{k-1,0}^{d,l,0,D,\alpha}.
\end{equation*}

Лемма 5.1.4

В условиях леммы 5.1.3 существует константа
$ c_{29}(d,\alpha,D,\lambda,p,q) >0 $ такая, что
для $ f \in (\mathcal H_p^\alpha)^\prime(D) $:

1) в $ L_q(D) $ имеет место равенство
\begin{equation*} \tag{5.1.26}
\D^\lambda f = (\mathcal U_{K^0,0}^{d,l-\e-\lambda,0,D,\alpha} \D^\lambda f) \mid_D
+\sum_{k = K^0 +1}^\infty (V_k^{d,l-\e-\lambda,D,\alpha}(\D^\lambda f)) \mid_D,
\end{equation*}
и
2) при $ k \in \Z_+: k > K^0, $ справедливо неравенство
\begin{equation*} \tag{5.1.27}
\| V_k^{d,l-\e-\lambda,D,\alpha} (\D^\lambda f)\|_{L_q(\R^d)} \le
c_{29} 2^{-k(1 -(\alpha^{-1}, \lambda +(p^{-1} -q^{-1})_+ \e))}.
\end{equation*}

Доказательство.

Начнем с доказательства (5.1.26). При соблюдении условий леммы для
$ f \in (\mathcal H_p^\alpha)^\prime(D), \k \in \Z_+: \k > K^0 (K^0 $ см. в п. 1.3.),
при $ \kappa = \kappa(k,\alpha), k \ge K^0, $ с учетом того, что
$ \sum_{\nu \in N_\kappa^{d,D}} \chi_{\kappa,\nu}^d =1 $ почти всюду на $ D, $

в силу (5.1.8) выполняется соотношение
\begin{multline*}
\| \D^\lambda f -((\mathcal U_{K^0,0}^{d,l-\e-\lambda,0,D,\alpha} \D^\lambda f) \mid_D
+\sum_{k = K^0 +1}^{\k} (V_k^{d,l-\e-\lambda,D,\alpha}(\D^\lambda f)) \mid_D)
\|_{L_q(D)} = \\
\| \D^\lambda f -((\mathcal U_{K^0,0}^{d,l-\e-\lambda,0,D,\alpha} \D^\lambda f) \mid_D
+\sum_{k = K^0 +1}^{\k} ((\mathcal U_{k,0}^{d,l-\e-\lambda,0,D,\alpha}
\D^\lambda f) \mid_D -\\
(\sum_{\nu \in N_\kappa^{d,D}} \chi_{\kappa,\nu}^d) \mid_D
(\mathcal U_{k-1,0}^{d,l-\e-\lambda,0,D,\alpha} \D^\lambda f) \mid_D))
\|_{L_q(D)} = \\
\| \D^\lambda f -((\mathcal U_{K^0,0}^{d,l-\e-\lambda,0,D,\alpha} \D^\lambda f) \mid_D
+\sum_{k = K^0 +1}^{\k} ((\mathcal U_{k,0}^{d,l-\e-\lambda,0,D,\alpha}
\D^\lambda f) \mid_D -\\
(\mathcal U_{k-1,0}^{d,l-\e-\lambda,0,D,\alpha} \D^\lambda f) \mid_D))
\|_{L_q(D)} = \\
\| \D^\lambda f -(\mathcal U_{\k,0}^{d,l-\e-\lambda,0,D,\alpha} \D^\lambda f) \mid_D
\|_{L_q(D)} \le 
c_4 2^{-\k(1-(\alpha^{-1}, \lambda +(p^{-1} -q^{-1})_+ \e))} \to 0 \text{ при }
\k \to \infty,
\end{multline*}
т.е. справедливо (5.1.26).

Перейдем к выводу неравенства (5.1.27). В условиях леммы при $ k \in \Z_+: k > K^0,
\kappa = \kappa(k,\alpha), \kappa^\prime = \kappa(k -1,\alpha),
\k = \kappa -\kappa^\prime $ для $ f \in (\mathcal H_p^\alpha)^\prime(D), $
учитывая (5.1.7), (5.1.3), имеем
\begin{multline*} \tag{5.1.28}
\| V_k^{d,l-\e-\lambda,D,\alpha} (\D^\lambda f)\|_{L_q(\R^d)}^q = \\
\biggl\| (\mathcal U_{k,0}^{d,l-\e-\lambda,0,D,\alpha} \D^\lambda f) -
(\sum_{\nu \in N_\kappa^{d,D}} \chi_{\kappa,\nu}^d)
(\mathcal U_{k-1,0}^{d,l-\e-\lambda,0,D,\alpha} \D^\lambda f) \biggr\|_{L_q(\R^d)}^q = \\
\biggl\| (\sum_{\nu \in N_\kappa^{d,D}}
(S_{\kappa,\nu_\kappa^D(\nu)}^{d,l-\e-\lambda,0} (\D^\lambda f))
\chi_{\kappa,\nu}^d) -
(\sum_{\nu \in N_\kappa^{d,D}} \chi_{\kappa,\nu}^d)
(\sum_{\nu^\prime \in N_{\kappa^\prime}^{d,D}}
(S_{\kappa^\prime,\nu_{\kappa^\prime}^D(\nu^\prime)}^{d,l-\e-\lambda,0}
(\D^\lambda f)) \chi_{\kappa^\prime,\nu^\prime}^d) \biggr\|_{L_q(\R^d)}^q = \\
\biggl\| (\sum_{\nu \in N_\kappa^{d,D}}
(S_{\kappa,\nu_\kappa^D(\nu)}^{d,l-\e-\lambda,0} (\D^\lambda f))
\chi_{\kappa,\nu}^d) -
(\sum_{\nu \in N_\kappa^{d,D}}
(S_{\kappa^\prime,\nu_{\kappa^\prime}^D(\nu^\prime(\nu))}^{d,l-\e-\lambda,0}
(\D^\lambda f)) \chi_{\kappa,\nu}^d) \biggr\|_{L_q(\R^d)}^q = \\
\biggl\| (\sum_{\nu \in N_\kappa^{d,D}}
(\D^\lambda S_{\kappa,\nu_\kappa^D(\nu)}^{d,l-\e,\lambda} f)
\chi_{\kappa,\nu}^d) -
(\sum_{\nu \in N_\kappa^{d,D}}
(\D^\lambda S_{\kappa^\prime,
\nu_{\kappa^\prime}^D(\nu^\prime(\nu))}^{d,l-\e,\lambda} f)
\chi_{\kappa,\nu}^d) \biggr\|_{L_q(\R^d)}^q = \\
\biggl\| \sum_{\nu \in N_\kappa^{d,D}}
((\D^\lambda S_{\kappa,\nu_\kappa^D(\nu)}^{d,l-\e,\lambda} f) -
(\D^\lambda S_{\kappa^\prime,
\nu_{\kappa^\prime}^D(\nu^\prime(\nu))}^{d,l-\e,\lambda} f))
\chi_{\kappa,\nu}^d \biggr\|_{L_q(\R^d)}^q = \\
\int_{ \cup_{n \in N_\kappa^{d,D}} Q_{\kappa,n}^d }
\biggl| \sum_{\nu \in N_\kappa^{d,D}}
((\D^\lambda S_{\kappa,\nu_\kappa^D(\nu)}^{d,l-\e,\lambda} f) -
(\D^\lambda S_{\kappa^\prime,
\nu_{\kappa^\prime}^D(\nu^\prime(\nu))}^{d,l-\e,\lambda} f))
\chi_{\kappa,\nu}^d \biggr|^q dx  = \\
\sum_{n \in N_\kappa^{d,D}} \int_{ Q_{\kappa,n}^d}
\biggl| \sum_{\nu \in N_\kappa^{d,D}}
((\D^\lambda S_{\kappa,\nu_\kappa^D(\nu)}^{d,l-\e,\lambda} f) -
(\D^\lambda S_{\kappa^\prime,
\nu_{\kappa^\prime}^D(\nu^\prime(\nu))}^{d,l-\e,\lambda} f))
\chi_{\kappa,\nu}^d \biggr|^q dx = \\
\sum_{n \in N_\kappa^{d,D}} \int_{ Q_{\kappa,n}^d}
\biggl| ((\D^\lambda S_{\kappa,\nu_\kappa^D(n)}^{d,l-\e,\lambda} f) -
(\D^\lambda S_{\kappa^\prime,
\nu_{\kappa^\prime}^D(\nu^\prime(n))}^{d,l-\e,\lambda} f))
\chi_{\kappa,n}^d \biggr|^q dx = \\
\sum_{n \in N_\kappa^{d,D}}
\| (\D^\lambda S_{\kappa,\nu_\kappa^D(n)}^{d,l-\e,\lambda} f) -
(\D^\lambda S_{\kappa^\prime,
\nu_{\kappa^\prime}^D(\nu^\prime(n))}^{d,l-\e,\lambda} f)
\|_{L_q(Q_{\kappa,n}^d)}^q.
\end{multline*}

Применение неравенства (1.1.1) для оценки слагаемых в правой части (5.1.28)
дает
\begin{multline*} \tag{5.1.29}
\| (\D^\lambda S_{\kappa,\nu_\kappa^D(n)}^{d,l-\e,\lambda} f) -
(\D^\lambda S_{\kappa^\prime,
\nu_{\kappa^\prime}^D(\nu^\prime(n))}^{d,l-\e,\lambda} f)
\|_{L_q(Q_{\kappa,n}^d)} = \\
\| \D^\lambda (S_{\kappa,\nu_\kappa^D(n)}^{d,l-\e,\lambda} f -
S_{\kappa^\prime,
\nu_{\kappa^\prime}^D(\nu^\prime(n))}^{d,l-\e,\lambda} f)
\|_{L_q(Q_{\kappa,n}^d)} \le \\
c_{30} 2^{(\kappa, \lambda +p^{-1} \e -q^{-1} \e)}
\| S_{\kappa,\nu_\kappa^D(n)}^{d,l-\e,\lambda} f -
S_{\kappa^\prime,
\nu_{\kappa^\prime}^D(\nu^\prime(n))}^{d,l-\e,\lambda} f
\|_{L_p(Q_{\kappa,n}^d)}, n \in N_\kappa^{d,D}.
\end{multline*}

Для оценки правой части (5.1.29), замечая, что при $ n \in N_\kappa^{d,D} $
справедливы включения
\begin{multline*} \tag{5.1.30}
Q_{\kappa,\nu_\kappa^D(n)}^d \subset D \supset Q_{\kappa^\prime,
\nu_{\kappa^\prime}^D(\nu^\prime(n))}^d = 2^{-\kappa^\prime}
\nu_{\kappa^\prime}^D(\nu^\prime(n)) +2^{-\kappa^\prime} I^d =\\
 (2^{-\kappa}
2^{\k} \nu_{\kappa^\prime}^D(\nu^\prime(n)) +2^{-\kappa^\prime} I^d)
\supset (2^{-\kappa} 2^{\k} \nu_{\kappa^\prime}^D(\nu^\prime(n)) +2^{-\kappa} I^d) =
Q_{\kappa, 2^{\k} \nu_{\kappa^\prime}^D(\nu^\prime(n))}^d,
\end{multline*}
фиксируя произвольное $ n \in N_\kappa^{d,D}, $ выберем последовательности
$ \nu^\iota \in \Z^d, j^\iota \in \Nu_{1,d}^1, \epsilon^\iota \in \{-1, 1\},
\iota =0,\ldots,\Iota, $ для которых
$ \nu^0 = \nu_\kappa^D(n), \nu^{\Iota} = 2^{\k} \nu_{\kappa^\prime}^D(\nu^\prime(n)) $
и соблюдаются соотношения (1.3.31), (1.3.32), (1.3.33). Тогда
\begin{multline*} \tag{5.1.31}
\| S_{\kappa,\nu_\kappa^D(n)}^{d,l-\e,\lambda} f -
S_{\kappa^\prime,
\nu_{\kappa^\prime}^D(\nu^\prime(n))}^{d,l-\e,\lambda} f
\|_{L_p(Q_{\kappa,n}^d)} = 
\| S_{\kappa, \nu^0}^{d,l -\e,\lambda} f -
S_{\kappa^\prime,
\nu_{\kappa^\prime}^D(\nu^\prime(n))}^{d,l-\e,\lambda} f
\|_{L_p(Q_{\kappa,n}^d)} = \\
\| S_{\kappa, \nu^0}^{d,l -\e,\lambda} f -\sum_{\iota =0}^{\Iota}
(S_{\kappa, \nu^\iota}^{d,l -\e} f -S_{\kappa, \nu^\iota}^{d,l -\e} f) -
S_{\kappa^\prime,
\nu_{\kappa^\prime}^D(\nu^\prime(n))}^{d,l-\e,\lambda} f
\|_{L_p(Q_{\kappa,n}^d)} = \\
\| S_{\kappa, \nu^0}^{d,l -\e,\lambda} f -
S_{\kappa, \nu^0}^{d,l -\e} f +\sum_{\iota =0}^{\Iota -1}
(S_{\kappa, \nu^\iota}^{d,l -\e} f -S_{\kappa, \nu^{\iota +1}}^{d,l -\e} f) +
S_{\kappa, \nu^{\Iota}}^{d,l -\e} f -
S_{\kappa^\prime,
\nu_{\kappa^\prime}^D(\nu^\prime(n))}^{d,l-\e,\lambda} f
\|_{L_p(Q_{\kappa,n}^d)} \le \\
\| S_{\kappa, \nu^0}^{d,l -\e,\lambda} f -
S_{\kappa, \nu^0}^{d,l -\e} f\|_{L_p(Q_{\kappa,n}^d)} +
\sum_{\iota =0}^{\Iota -1} \| S_{\kappa, \nu^\iota}^{d,l -\e} f -
S_{\kappa, \nu^{\iota +1}}^{d,l -\e} f\|_{L_p(Q_{\kappa,n}^d)} +\\
\| S_{\kappa, \nu^{\Iota}}^{d,l -\e} f -
S_{\kappa^\prime,
\nu_{\kappa^\prime}^D(\nu^\prime(n))}^{d,l-\e,\lambda} f
\|_{L_p(Q_{\kappa,n}^d)}, n \in N_\kappa^{d,D}.
\end{multline*}

Для проведения оценки слагаемых в правой части (5.1.31) отметим некоторые
полезные для этого факты.

Сначала заметим, что при $ n \in N_\kappa^{d,D} $
вследствие (1.3.12) (при $ \nu^\prime(n) $ вместо $ \nu, \kappa^\prime $ вместо
$ \kappa $ ), а также благодаря (5.1.6),
(1.3.42), для $ j =1,\ldots,d $ соблюдается неравенство
\begin{multline*} \tag{5.1.32}
| 2^{-\kappa_j} 2^{\k_j} (\nu_{\kappa^\prime}^D(\nu^\prime(n)))_j -
2^{-\kappa_j} n_j | = 
| 2^{-\kappa^\prime_j} (\nu_{\kappa^\prime}^D(\nu^\prime(n)))_j -
2^{-\kappa_j} n_j | = \\
| 2^{-\kappa^\prime_j} (\nu_{\kappa^\prime}^D(\nu^\prime(n)))_j -
2^{-\kappa^\prime_j} (\nu^\prime(n))_j
+2^{-\kappa^\prime_j} (\nu^\prime(n))_j -2^{-\kappa_j} n_j | \le \\
| 2^{-\kappa^\prime_j} (\nu_{\kappa^\prime}^D(\nu^\prime(n)))_j -
2^{-\kappa^\prime_j} (\nu^\prime(n))_j | +| 2^{-\kappa^\prime_j}
(\nu^\prime(n))_j -2^{-\kappa_j} n_j | \le \\
2^{-\kappa^\prime_j} \| \nu_{\kappa^\prime}^D(\nu^\prime(n)) -
\nu^\prime(n) \|_{l_\infty^d} +| 2^{-\kappa^\prime_j} (\nu^\prime(n))_j -
2^{-\kappa_j} n_j | \le \\
c_{31} 2^{-\kappa^\prime_j} +2^{-\kappa^\prime_j} = 
(c_{31} +1) 2^{\k_j} 2^{-\kappa_j} \le (c_{31} +1) 2^{1 +1 /\alpha_j} 2^{-\kappa_j}=\\
\rho_j 2^{-\kappa_j} \le \| \rho \|_{l_\infty^d} 2^{-\kappa_j},
\end{multline*}
откуда для $ x \in \overline Q_{\kappa^\prime, \nu_{\kappa^\prime}^D(\nu^\prime(n))}^d $ --
\begin{multline*}
| x_j -2^{-\kappa_j} n_j | = | x_j -2^{-\kappa^\prime_j}
(\nu_{\kappa^\prime}^D(\nu^\prime(n)))_j +
2^{-\kappa^\prime_j} (\nu_{\kappa^\prime}^D(\nu^\prime(n)))_j
-2^{-\kappa_j} n_j | \le \\
| x_j -2^{-\kappa^\prime_j} (\nu_{\kappa^\prime}^D(\nu^\prime(n)))_j |
+| 2^{-\kappa^\prime_j} (\nu_{\kappa^\prime}^D(\nu^\prime(n)))_j
-2^{-\kappa_j} n_j | \le \\
2^{-\kappa^\prime_j} +\rho_j 2^{-\kappa_j} = (2^{\k_j} +\rho_j) 2^{-\kappa_j} \le
(2^{1 +1 /\alpha_j} +\rho_j) 2^{-\kappa_j},
\end{multline*}
т.е.
\begin{equation*} \tag{5.1.33}
Q_{\kappa^\prime, \nu_{\kappa^\prime}^D(\nu^\prime(n))}^d \subset
\overline Q_{\kappa^\prime, \nu_{\kappa^\prime}^D(\nu^\prime(n))}^d
\subset (2^{-\kappa} n +(2^{\e +\alpha^{-1}} +\rho) 2^{-\kappa} B^d).
\end{equation*}

Отметим еще, что в силу (1.3.31), (1.3.32), (1.3.12) (при $ \nu = n), $
а также ввиду (5.1.32) при $ n \in N_\kappa^{d,D} $ справедлива оценка
\begin{multline*} \tag{5.1.34}
\Iota \le c_{32} \| \nu^0 -\nu^{\Iota} \|_{l_\infty^d} = c_{32} \| \nu_\kappa^D(n) -
2^{\k} \nu_{\kappa^\prime}^D(\nu^\prime(n)) \|_{l_\infty^d} \le\\
c_{32} (\| \nu_\kappa^D(n) -n \|_{l_\infty^d}
+\| n -2^{\k} \nu_{\kappa^\prime}^D(\nu^\prime(n)) \|_{l_\infty^d}) \le 
c_{32} (c_{33} +\| \rho\|_{l_\infty^d}) = c_{34}.
\end{multline*}

Далее, при $ n \in N_\kappa^{d,D}, \iota =0,\ldots,\Iota, j =1,\ldots,d $ для
$ x \in \overline Q_{\kappa,\nu^\iota}^d $ ввиду (1.3.32), (5.1.34),
(1.3.12) (при $ \nu = n $) выполняется неравенство
\begin{multline*}
| x_j -2^{-\kappa_j} n_j | \le | x_j -2^{-\kappa_j} (\nu^\iota)_j | +
| 2^{-\kappa_j} (\nu^\iota)_j -2^{-\kappa_j} n_j | \le 2^{-\kappa_j} +
| 2^{-\kappa_j} (\nu^\iota)_j -2^{-\kappa_j} n_j | = \\
2^{-\kappa_j} +| 2^{-\kappa_j} (\nu^\iota)_j -\sum_{i =0}^{\iota -1}
(2^{-\kappa_j} (\nu^i)_j -2^{-\kappa_j} (\nu^i)_j) -2^{-\kappa_j} n_j | = \\
2^{-\kappa_j} +| \sum_{i =0}^{\iota -1} (2^{-\kappa_j} (\nu^{i +1})_j
-2^{-\kappa_j} (\nu^i)_j) +2^{-\kappa_j} (\nu^0)_j -2^{-\kappa_j} n_j | \le \\
2^{-\kappa_j} +\sum_{i =0}^{\iota -1} | 2^{-\kappa_j} (\nu^{i +1})_j
-2^{-\kappa_j} (\nu^i)_j| +| 2^{-\kappa_j} (\nu^0)_j -2^{-\kappa_j} n_j | \le \\
2^{-\kappa_j} +\sum_{i =0}^{\iota -1} 2^{-\kappa_j} \| \nu^{i +1} -\nu^i \|_{l_\infty^d}
+| 2^{-\kappa_j} (\nu^0)_j -2^{-\kappa_j} n_j | = \\
2^{-\kappa_j} +\sum_{i =0}^{\iota -1} 2^{-\kappa_j} \| \epsilon^i e_{j^i} \|_{l_\infty^d}
+| 2^{-\kappa_j} (\nu_\kappa^D(n))_j  -2^{-\kappa_j} n_j | = \\
2^{-\kappa_j} +\sum_{i =0}^{\iota -1} 2^{-\kappa_j}
+| 2^{-\kappa_j} (\nu_\kappa^D(n))_j  -2^{-\kappa_j} n_j | = \\
2^{-\kappa_j} (\iota +1) +| 2^{-\kappa_j} (\nu_\kappa^D(n))_j -
2^{-\kappa_j} n_j | \le \\
2^{-\kappa_j} (\Iota +1) +2^{-\kappa_j} \| \nu_\kappa^D(n) -n \|_{l_\infty^d}   \le 
2^{-\kappa_j} (c_{34} +1) +2^{-\kappa_j} c_{35} = c_{36} 2^{-\kappa_j},
\end{multline*}
т.е. имеет место включение
\begin{equation*}
\overline Q_{\kappa,\nu^\iota}^d \subset (2^{-\kappa} n +c_{36} 2^{-\kappa} B^d),
\end{equation*}
из которого следует (см. (1.3.33)), что
\begin{equation*} \tag{5.1.35}
Q_{\kappa,\nu^\iota}^d \cup Q_{\kappa,\nu^{\iota +1}}^d \subset Q^\iota \subset
D \cap (2^{-\kappa} n +c_{36} 2^{-\kappa} B^d), \iota =0,\ldots,\Iota -1.
\end{equation*}

Итак, задавая координаты вектора $ \sigma \in \R_+^d $ соотношением
$ \sigma_j = \max(2^{1 +1 /\alpha_j} +\rho_j, c_{36}), j =1,\ldots,d, $ и обозначая
$$
\mathfrak D_{\kappa,n}^{\prime d,D,\alpha} = 2^{-\kappa} n +\sigma 2^{-\kappa} B^d,
n \in N_\kappa^{d,D},
$$
в силу (5.1.35), (5.1.30), (5.1.33) имеем
\begin{equation*} \tag{5.1.36}
Q_{\kappa,\nu^\iota}^d \cup Q_{\kappa,\nu^{\iota +1}}^d \subset Q^\iota \subset
D \cap \mathfrak D_{\kappa,n}^{\prime d,D,\alpha}, \iota =0,\ldots,\Iota -1.
\end{equation*}
\begin{equation*} \tag{5.1.37}
Q_{\kappa^\prime, \nu_{\kappa^\prime}^D(\nu^\prime(n))}^d \subset
D \cap \mathfrak D_{\kappa,n}^{\prime d,D,\alpha},
n \in N_\kappa^{d,D}.
\end{equation*}

Из приведенных определений с учетом того, что $ \sigma > \e, $ видно, что
при $ n \in N_\kappa^{d,D} $ справедливо включение
\begin{equation*} \tag{5.1.38}
Q_{\kappa, n}^d \subset \mathfrak D_{\kappa, n}^{\prime d,D,\alpha},
\end{equation*}
а из (5.1.35) следует, что
\begin{equation*} \tag{5.1.39}
Q_{\kappa,n}^d \subset (2^{-\kappa} \nu^{\iota} +(c_{36} +1) 2^{-\kappa} B^d),
\iota =0,\ldots,\Iota.
\end{equation*}

Учитывая (5.1.38), легко видеть, что существует константа
$ c_{37}(d,D,\alpha) >0 $ такая, что для каждого $ x \in \R^d $ число
\begin{equation*} \tag{5.1.40}
\card \{ n \in N_\kappa^{d,D}:
x \in \mathfrak D_{\kappa,n}^{\prime d,D,\alpha} \} \le c_{37}.
\end{equation*}

Принимая во внимание, что при $ \iota =0,\ldots,\Iota -1 $ в силу (1.3.33),
(1.3.32) имеет место равенство
\begin{equation*}
Q^\iota = \begin{cases} 2^{-\kappa} \nu^\iota +2^{-(\kappa -e_{j^\iota})} I^d,
\text{ при } \epsilon^\iota =1; \\
2^{-\kappa} \nu^{\iota +1} +2^{-(\kappa -e_{j^\iota})} I^d,

\text{ при } \epsilon^\iota =-1,
\end{cases}
\end{equation*}
определим линейный оператор $ S^\iota, $ полагая
$$
S^\iota = P_{\delta,x^0}^{d,l -\e}
$$
при $ \delta = 2^{-(\kappa -e_{j^\iota})}, $
\begin{equation*}
x^0 = \begin{cases} 2^{-\kappa} \nu^\iota, \text{ при } \epsilon^\iota =1; \\
2^{-\kappa} \nu^{\iota +1}, \text{ при } \epsilon^\iota =-1,
\end{cases}
\iota =0,\ldots,\Iota -1.
\end{equation*}

Теперь проведем оценку слагаемых в правой части (5.1.31). Оценивая первое
слагаемое, согласно (1.1.1) с учетом (5.1.39), в силу (5.1.2), (5.1.1) и на
основании (1.1.4), (5.1.36) выводим
\begin{multline*} \tag{5.1.41}
\| S_{\kappa, \nu^0}^{d,l -\e,\lambda} f -
S_{\kappa, \nu^0}^{d,l -\e} f\|_{L_p(Q_{\kappa,n}^d)} \le 
c_{38} \| S_{\kappa, \nu^0}^{d,l -\e,\lambda} f -
S_{\kappa, \nu^0}^{d,l -\e} f\|_{L_p(Q_{\kappa,\nu^0}^d)} = \\
c_{38} \| S_{\kappa, \nu^0}^{d,l -\e,\lambda} (f -
S_{\kappa, \nu^0}^{d,l -\e} f)\|_{L_p(Q_{\kappa,\nu^0}^d)} \le 
c_{39} \| f -S_{\kappa, \nu^0}^{d,l -\e} f\|_{L_p(Q_{\kappa,\nu^0}^d)} \le \\
c_{40} \sum_{j =1}^d 2^{\kappa_j /p} \biggl(\int_{2^{-\kappa_j} B^1}
\int_{(Q_{\kappa,\nu^0}^d)_{l_j \xi e_j}} |\Delta_{\xi e_j}^{l_j} f(x)|^p dx d\xi\biggr)^{1/p} \le \\
c_{40} \sum_{j =1}^d 2^{\kappa_j /p} \biggl(\int_{2^{-\kappa_j} B^1}
\int_{(D \cap \mathfrak D_{\kappa,n}^{\prime d,D,\alpha})_{l_j \xi e_j}}
| \Delta_{\xi e_j}^{l_j} f(x)|^p dx d\xi\biggr)^{1/p} \le \\
c_{40} \sum_{j =1}^d 2^{\kappa_j /p} \biggl(\int_{2^{-\kappa_j} B^1}
\int_{D_{l_j \xi e_j} \cap \mathfrak D_{\kappa,n}^{\prime d,D,\alpha}}
| \Delta_{\xi e_j}^{l_j} f(x)|^p dx d\xi\biggr)^{1/p}.
\end{multline*}

Далее, при $ \iota =0,\ldots,\Iota -1, $ благодаря (5.1.39) применяя (1.1.1),
а затем используя (1.1.2), (1.1.3), (5.1.36), и, наконец, пользуясь (1.1.4)
и учитывая (5.1.36), приходим к неравенству
\begin{multline*} \tag{5.1.42}
\| S_{\kappa, \nu^\iota}^{d,l -\e} f -
S_{\kappa, \nu^{\iota +1}}^{d,l -\e} f\|_{L_p(Q_{\kappa,n}^d)} = 
\| S_{\kappa, \nu^\iota}^{d,l -\e} f -S^\iota f +S^\iota f -
S_{\kappa, \nu^{\iota +1}}^{d,l -\e} f\|_{L_p(Q_{\kappa,n}^d)} \le \\
\| S_{\kappa, \nu^\iota}^{d,l -\e} f -S^\iota f \|_{L_p(Q_{\kappa,n}^d)}
+\| S^\iota f -S_{\kappa, \nu^{\iota +1}}^{d,l -\e} f\|_{L_p(Q_{\kappa,n}^d)} \le \\
c_{38} \| S_{\kappa, \nu^\iota}^{d,l -\e} f -S^\iota f \|_{L_p(Q_{\kappa,\nu^\iota}^d)}
+c_{38} \| S^\iota f -S_{\kappa, \nu^{\iota +1}}^{d,l -\e} f
\|_{L_p(Q_{\kappa,\nu^{\iota +1}}^d)} = \\
c_{38} \| S_{\kappa, \nu^\iota}^{d,l -\e} (f -S^\iota f)
\|_{L_p(Q_{\kappa,\nu^\iota}^d)}
+c_{38} \| S_{\kappa, \nu^{\iota +1}}^{d,l -\e} (f -S^\iota f)
\|_{L_p(Q_{\kappa,\nu^{\iota +1}}^d)} \le \\
c_{41} \| f -S^\iota f \|_{L_p(Q_{\kappa,\nu^\iota}^d)}
+c_{41} \| f -S^\iota f\|_{L_p(Q_{\kappa,\nu^{\iota +1}}^d)} \le 
c_{42} \| f -S^\iota f \|_{L_p(Q^\iota)} \\
\le c_{43} \sum_{j =1}^d
((2^{-(\kappa -e_{j^\iota})})_j)^{-1/p} \biggl(\int_{(2^{-(\kappa -e_{j^\iota})})_j B^1}
\int_{Q_{l_j \xi e_j}^\iota}| \Delta_{\xi e_j}^{l_j} f(x)|^p dx d\xi\biggr)^{1/p} \le \\
c_{43} \sum_{j =1}^d 2^{\kappa_j /p} \biggl(\int_{2^{-\kappa_j +1} B^1}
\int_{Q_{l_j \xi e_j}^\iota}| \Delta_{\xi e_j}^{l_j} f(x)|^p dx d\xi\biggr)^{1/p} \le \\
c_{43} \sum_{j =1}^d 2^{\kappa_j /p} \biggl(\int_{2^{-\kappa_j +1} B^1}
\int_{(D \cap \mathfrak D_{\kappa,n}^{\prime d,D,\alpha})_{l_j \xi e_j}}
| \Delta_{\xi e_j}^{l_j} f(x)|^p dx d\xi\biggr)^{1/p} \le \\
c_{43} \sum_{j =1}^d 2^{\kappa_j /p} \biggl(\int_{2^{-\kappa_j +1} B^1}
\int_{D_{l_j \xi e_j} \cap \mathfrak D_{\kappa,n}^{\prime d,D,\alpha}}
| \Delta_{\xi e_j}^{l_j} f(x)|^p dx d\xi\biggr)^{1/p}.
\end{multline*}

Для оценки последнего слагаемого в правой части (5.1.31), ввиду (5.1.39)
используя (1.1.1), затем принимая во внимание (5.1.30), пользуясь (1.1.2),
(5.1.2), (1.1.3), (5.1.1), применяя (1.1.4), с учетом (5.1.37) получаем
\begin{multline*} \tag{5.1.43}
\| S_{\kappa, \nu^{\Iota}}^{d,l -\e} f -
S_{\kappa^\prime, \nu_{\kappa^\prime}^D(\nu^\prime(n))}^{d,l-\e,\lambda} f
\|_{L_p(Q_{\kappa,n}^d)} \le 
c_{38} \| S_{\kappa, \nu^{\Iota}}^{d,l -\e} f -
S_{\kappa^\prime, \nu_{\kappa^\prime}^D(\nu^\prime(n))}^{d,l-\e,\lambda} f
\|_{L_p(Q_{\kappa,\nu^{\Iota}}^d)} \le \\
c_{38} \| S_{\kappa, \nu^{\Iota}}^{d,l -\e} f -
S_{\kappa^\prime, \nu_{\kappa^\prime}^D(\nu^\prime(n))}^{d,l-\e} f
\|_{L_p(Q_{\kappa,\nu^{\Iota}}^d)} +
c_{38} \| S_{\kappa^\prime, \nu_{\kappa^\prime}^D(\nu^\prime(n))}^{d,l-\e} f -
S_{\kappa^\prime, \nu_{\kappa^\prime}^D(\nu^\prime(n))}^{d,l-\e,\lambda} f
\|_{L_p(Q_{\kappa,\nu^{\Iota}}^d)} \le \\
c_{38} \| S_{\kappa, \nu^{\Iota}}^{d,l -\e} f -
S_{\kappa^\prime, \nu_{\kappa^\prime}^D(\nu^\prime(n))}^{d,l-\e} f
\|_{L_p(Q_{\kappa,\nu^{\Iota}}^d)} +
c_{38} \| S_{\kappa^\prime, \nu_{\kappa^\prime}^D(\nu^\prime(n))}^{d,l-\e} f -
S_{\kappa^\prime, \nu_{\kappa^\prime}^D(\nu^\prime(n))}^{d,l-\e,\lambda} f
\|_{L_p(Q_{\kappa^\prime, \nu_{\kappa^\prime}^D(\nu^\prime(n))}^d)} = \\
c_{38} \| S_{\kappa, \nu^{\Iota}}^{d,l -\e} (f -
S_{\kappa^\prime, \nu_{\kappa^\prime}^D(\nu^\prime(n))}^{d,l-\e} f)
\|_{L_p(Q_{\kappa,\nu^{\Iota}}^d)} +
c_{38} \| S_{\kappa^\prime, \nu_{\kappa^\prime}^D(\nu^\prime(n))}^{d,l-\e,\lambda}
(S_{\kappa^\prime, \nu_{\kappa^\prime}^D(\nu^\prime(n))}^{d,l-\e} f -f)
\|_{L_p(Q_{\kappa^\prime, \nu_{\kappa^\prime}^D(\nu^\prime(n))}^d)} \le \\
c_{44} \| f -S_{\kappa^\prime, \nu_{\kappa^\prime}^D(\nu^\prime(n))}^{d,l-\e} f
\|_{L_p(Q_{\kappa,\nu^{\Iota}}^d)} +
c_{45} \| S_{\kappa^\prime, \nu_{\kappa^\prime}^D(\nu^\prime(n))}^{d,l-\e} f -f
\|_{L_p(Q_{\kappa^\prime, \nu_{\kappa^\prime}^D(\nu^\prime(n))}^d)} \le \\
c_{46} \| f -S_{\kappa^\prime, \nu_{\kappa^\prime}^D(\nu^\prime(n))}^{d,l-\e} f
\|_{L_p(Q_{\kappa^\prime, \nu_{\kappa^\prime}^D(\nu^\prime(n))}^d)} \le \\
c_{47} \sum_{j =1}^d 2^{\kappa^\prime_j /p} \biggl(\int_{2^{-\kappa^\prime_j} B^1}
\int_{(Q_{\kappa^\prime, \nu_{\kappa^\prime}^D(\nu^\prime(n))}^d)_{l_j \xi e_j}}
| \Delta_{\xi e_j}^{l_j} f(x)|^p dx d\xi\biggr)^{1/p} \le \\
c_{47} \sum_{j =1}^d 2^{\kappa_j /p} \biggl(\int_{c_{48} 2^{-\kappa_j} B^1}
\int_{(D \cap \mathfrak D_{\kappa,n}^{\prime d,D,\alpha})_{l_j \xi e_j}}
| \Delta_{\xi e_j}^{l_j} f(x)|^p dx d\xi\biggr)^{1/p} \le \\
c_{47} \sum_{j =1}^d 2^{\kappa_j /p} \biggl(\int_{c_{48} 2^{-\kappa_j} B^1}
\int_{D_{l_j \xi e_j} \cap \mathfrak D_{\kappa,n}^{\prime d,D,\alpha}}
| \Delta_{\xi e_j}^{l_j} f(x)|^p dx d\xi\biggr)^{1/p},
\end{multline*}
где $ c_{48} = \max_{j =1,\ldots,d} 2^{1 +1 /\alpha_j} > 2. $
Объединяя (5.1.31), (5.1.41), (5.1.42), (5.1.43) и учитывая (5.1.34), при
$ n \in N_\kappa^{d,D} $ имеем
\begin{multline*} \tag{5.1.44}
\| S_{\kappa,\nu_\kappa^D(n)}^{d,l-\e,\lambda} f -
S_{\kappa^\prime,
\nu_{\kappa^\prime}^D(\nu^\prime(n))}^{d,l-\e,\lambda} f
\|_{L_p(Q_{\kappa,n}^d)} \le \\
c_{40} \sum_{j =1}^d 2^{\kappa_j /p} \biggl(\int_{2^{-\kappa_j} B^1}
\int_{D_{l_j \xi e_j} \cap \mathfrak D_{\kappa,n}^{\prime d,D,\alpha}}
| \Delta_{\xi e_j}^{l_j} f(x)|^p dx d\xi\biggr)^{1/p} + \\
\Iota c_{43} \sum_{j =1}^d 2^{\kappa_j /p} \biggl(\int_{2^{-\kappa_j +1} B^1}
\int_{D_{l_j \xi e_j} \cap \mathfrak D_{\kappa,n}^{\prime d,D,\alpha}}
| \Delta_{\xi e_j}^{l_j} f(x)|^p dx d\xi\biggr)^{1/p} +\\
c_{47} \sum_{j =1}^d 2^{\kappa_j /p} \biggl(\int_{c_{48} 2^{-\kappa_j} B^1}
\int_{D_{l_j \xi e_j} \cap \mathfrak D_{\kappa,n}^{\prime d,D,\alpha}}
| \Delta_{\xi e_j}^{l_j} f(x)|^p dx d\xi\biggr)^{1/p} \le \\
c_{40} \sum_{j =1}^d 2^{\kappa_j /p} \biggl(\int_{c_{48} 2^{-\kappa_j} B^1}
\int_{D_{l_j \xi e_j} \cap \mathfrak D_{\kappa,n}^{\prime d,D,\alpha}}
| \Delta_{\xi e_j}^{l_j} f(x)|^p dx d\xi\biggr)^{1/p} +\\
c_{34} c_{43} \sum_{j =1}^d 2^{\kappa_j /p} \biggl(\int_{c_{48} 2^{-\kappa_j} B^1}
\int_{D_{l_j \xi e_j} \cap \mathfrak D_{\kappa,n}^{\prime d,D,\alpha}}
| \Delta_{\xi e_j}^{l_j} f(x)|^p dx d\xi\biggr)^{1/p} +\\
c_{47} \sum_{j =1}^d 2^{\kappa_j /p} \biggl(\int_{c_{48} 2^{-\kappa_j} B^1}
\int_{D_{l_j \xi e_j} \cap \mathfrak D_{\kappa,n}^{\prime d,D,\alpha}}
| \Delta_{\xi e_j}^{l_j} f(x)|^p dx d\xi\biggr)^{1/p} = \\
c_{49} \sum_{j =1}^d 2^{\kappa_j /p} \biggl(\int_{c_{48} 2^{-\kappa_j} B^1}
\int_{D_{l_j \xi e_j} \cap \mathfrak D_{\kappa,n}^{\prime d,D,\alpha}}
| \Delta_{\xi e_j}^{l_j} f(x)|^p dx d\xi\biggr)^{1/p}.
\end{multline*}
Обращаясь к (5.1.29) и (5.1.44), заключаем, что при
$ n \in N_\kappa^{d,D} $ выполняется неравенство
\begin{multline*}
\| (\D^\lambda S_{\kappa,\nu_\kappa^D(n)}^{d,l-\e,\lambda} f) -
(\D^\lambda S_{\kappa^\prime,
\nu_{\kappa^\prime}^D(\nu^\prime(n))}^{d,l-\e,\lambda} f)
\|_{L_q(Q_{\kappa,n}^d)} \le \\
c_{50} 2^{(\kappa, \lambda +p^{-1} \e -q^{-1} \e)}
\sum_{j =1}^d 2^{\kappa_j /p} \biggl(\int_{c_{48} 2^{-\kappa_j} B^1}
\int_{D_{l_j \xi e_j} \cap \mathfrak D_{\kappa,n}^{\prime d,D,\alpha}}
| \Delta_{\xi e_j}^{l_j} f(x)|^p dx d\xi\biggr)^{1/p}.
\end{multline*}

Подставляя эту оценку в (5.1.28) и  применяя неравенство Г\"eльдера с
показателем $ q $ и неравенство Г\"eльдера с показателем $ p/q \le 1, $ и,
наконец, принимая во внимание (5.1.40) и используя неравенство Г\"eльдера с
показателем $ 1 /q, $ получаем
\begin{multline*} \tag{5.1.45}
\| V_k^{d,l-\e-\lambda,D,\alpha} (\D^\lambda f)\|_{L_q(\R^d)}^q \le \\
\sum_{n \in N_\kappa^{d,D}}
\biggl(c_{50} 2^{(\kappa, \lambda +p^{-1} \e -q^{-1} \e)}
\sum_{j =1}^d 2^{\kappa_j /p} \biggl(\int_{c_{48} 2^{-\kappa_j} B^1}
\int_{D_{l_j \xi e_j} \cap \mathfrak D_{\kappa,n}^{\prime d,D,\alpha}}
| \Delta_{\xi e_j}^{l_j} f(x)|^p dx d\xi\biggr)^{1/p}\biggr)^q = \\
(c_{50} 2^{(\kappa, \lambda +p^{-1} \e -q^{-1} \e)})^q
\sum_{n \in N_\kappa^{d,D}}
\biggl(\sum_{j =1}^d 2^{\kappa_j /p} \biggl(\int_{c_{48} 2^{-\kappa_j} B^1}
\int_{D_{l_j \xi e_j} \cap \mathfrak D_{\kappa,n}^{\prime d,D,\alpha}}
| \Delta_{\xi e_j}^{l_j} f(x)|^p dx d\xi\biggr)^{1/p}\biggr)^q \le \\
(c_{51} 2^{(\kappa, \lambda +p^{-1} \e -q^{-1} \e)})^q
\sum_{n \in N_\kappa^{d,D}}
\sum_{j =1}^d 2^{\kappa_j q /p} \biggl(\int_{c_{48} 2^{-\kappa_j} B^1}
\int_{D_{l_j \xi e_j} \cap \mathfrak D_{\kappa,n}^{\prime d,D,\alpha}}
| \Delta_{\xi e_j}^{l_j} f(x)|^p dx d\xi\biggr)^{q/p} = \\
(c_{51} 2^{(\kappa, \lambda +p^{-1} \e -q^{-1} \e)})^q
\sum_{j =1}^d 2^{\kappa_j q /p}
\sum_{n \in N_\kappa^{d,D}}
\biggl(\int_{c_{48} 2^{-\kappa_j} B^1}
\int_{D_{l_j \xi e_j} \cap \mathfrak D_{\kappa,n}^{\prime d,D,\alpha}}
| \Delta_{\xi e_j}^{l_j} f(x)|^p dx d\xi\biggr)^{q/p} \le \\
(c_{51} 2^{(\kappa, \lambda +p^{-1} \e -q^{-1} \e)})^q
\sum_{j =1}^d 2^{\kappa_j q /p}
\biggl(\sum_{n \in N_\kappa^{d,D}}
\int_{c_{48} 2^{-\kappa_j} B^1}
\int_{D_{l_j \xi e_j} \cap \mathfrak D_{\kappa,n}^{\prime d,D,\alpha}}
| \Delta_{\xi e_j}^{l_j} f(x)|^p dx d\xi\biggr)^{q/p} = \\
(c_{51} 2^{(\kappa, \lambda +p^{-1} \e -q^{-1} \e)})^q
\sum_{j =1}^d 2^{\kappa_j q /p}
\biggl(\int_{c_{48} 2^{-\kappa_j} B^1}
\sum_{n \in N_\kappa^{d,D}}
\int_{D_{l_j \xi e_j}} \chi_{\mathfrak D_{\kappa,n}^{\prime d,D,\alpha}}(x)
| \Delta_{\xi e_j}^{l_j} f(x)|^p dx d\xi\biggr)^{q/p} = \\
(c_{51} 2^{(\kappa, \lambda +p^{-1} \e -q^{-1} \e)})^q
\sum_{j =1}^d 2^{\kappa_j q /p}
\biggl(\int_{c_{48} 2^{-\kappa_j} B^1}
\int_{D_{l_j \xi e_j}}
(\sum_{n \in N_\kappa^{d,D}}
\chi_{\mathfrak D_{\kappa,n}^{\prime d,D,\alpha}}(x))
| \Delta_{\xi e_j}^{l_j} f(x)|^p dx d\xi\biggr)^{q/p} \le \\
(c_{51} 2^{(\kappa, \lambda +p^{-1} \e -q^{-1} \e)})^q
\sum_{j =1}^d 2^{\kappa_j q /p}
\biggl(\int_{c_{48} 2^{-\kappa_j} B^1}
\int_{D_{l_j \xi e_j}} c_{37} | \Delta_{\xi e_j}^{l_j} f(x)|^p dx d\xi\biggr)^{q/p} \le \\
(c_{52} 2^{(\kappa, \lambda +p^{-1} \e -q^{-1} \e)})^q
\biggl(\sum_{j =1}^d 2^{\kappa_j /p}
(\int_{c_{48} 2^{-\kappa_j} B^1}
\int_{D_{l_j \xi e_j}} | \Delta_{\xi e_j}^{l_j} f(x)|^p dx d\xi)^{1/p}\biggr)^q = \\
(c_{52} 2^{(\kappa, \lambda +p^{-1} \e -q^{-1} \e)}
\sum_{j =1}^d 2^{\kappa_j /p}
\biggl(\int_{c_{48} 2^{-\kappa_j} B^1}
\int_{D_{l_j \xi e_j}} | \Delta_{\xi e_j}^{l_j} f(x)|^p dx d\xi\biggr)^{1/p})^q \le \\
\biggl( c_{52} 2^{k (\alpha^{-1}, \lambda  +(p^{-1} -q^{-1})_+ \e)}
\sum_{j =1}^d 2^{k /(\alpha_j p)}
\biggl(\int_{c_{53} 2^{-k /\alpha_j} B^1}
\int_{D_{l_j \xi e_j}} | \Delta_{\xi e_j}^{l_j} f(x)|^p dx d\xi\biggr)^{1/p}\biggr)^q \le \\
(c_{54} 2^{k (\alpha^{-1}, \lambda  +(p^{-1} -q^{-1})_+ \e)}
\sum_{j =1}^d \Omega_j^{\prime l_j}(f, c_{53} 2^{-k /\alpha_j})_{L_p(D)})^q.
\end{multline*}

Из (5.1.45) приходим к (5.1.27) при $ p \le q. $
Неравенство (5.1.27) при $ q < p $ вытекает из неравенства (5.1.27) при
$ q = p $ и того факта, что при $ q < p, l \in \Z_+^d, \kappa \in \Z_+^d, $
ограниченной области $ D \subset \R^d $ для $ h \in \mathcal P_\kappa^{d,l,D} $
справедливо неравенство
$$
\|h\|_{L_q(\R^d)} = \|h\|_{L_q(D^\prime)} \le c_{55} \|h\|_{L_p(D^\prime)}
= c_{55} \|h\|_{L_p(\R^d)},
$$
где $ D^\prime $ -- ограниченная область в $ \R^d $ такая, что
$$
\cup_{n \in N_\kappa^{d,D}} Q_{\kappa,n}^d \subset D^\prime
$$
при $ \kappa \in \Z_+^d. \square $

С применением леммы 5.1.3 устанавливается также

Предложение 5.1.5

В условиях леммы 5.1.3 множество $ \D^\lambda (B((H_p^\alpha)^\prime(D))) $ --
вполне ограничено (а при $ 1<p<\infty $ компактно) в $ L_q(D). $

В заключение этого пункта для $ d \in \N, l \in \Z_+^d, $ ограниченной области
$ D \subset \R^d $ и $ \kappa \in \Z_+^d $ обозначим через $ \I_\kappa^{d,l,D}:
\mathcal P_\kappa^{d,l,D} \mapsto \R^{R_\kappa^{d,l,D}} $ линейный изоморфизм,
определяемый для $ f \in \mathcal P_\kappa^{d,l,D} $ равенством
$$
\I_\kappa^{d,l,D} f = \{ f_\nu(2^{-\kappa} \nu +2^{-\kappa} \lambda):
\nu \in N_\kappa^{d,D}, \lambda \in \Z_+^d(l) \},
$$
где $ \{f_\nu \in \mathcal P^{d,l}, \nu \in N_\kappa^{d,D}\} $ -- набор
полиномов, удовлетворяющих (5.1.4).

Несложно проверить, что имеет место

Лемма 5.1.6

Для $ d \in \N, l \in \Z_+^d, $ ограниченной области $ D \subset \R^d $
и $ 1 \le p \le \infty $ существуют константы $ c_{56}(d,l,p) >0 $ и
$ c_{57}(d,l,p) >0 $ такие, что при $ \kappa \in \Z_+^d $ для $ f \in
\mathcal P_\kappa^{d,l,D} $ соблюдается неравенство
\begin{equation*} \tag{5.1.46}
c_{56} \|f\|_{L_p(\R^d)} \le 2^{-(\kappa, \e) p^{-1}}
\| \I_\kappa^{d,l,D} f \|_{l_p^{R_\kappa^{d,l,D}}} \le c_{57} \|f\|_{L_p(\R^d)}.
\end{equation*}

Лемма 5.1.7

В условиях леммы 5.1.6 при $ \kappa \in \Z_+^d $ отображение
\begin{equation*}
U = U_\kappa^{d,l,D,p}: \mathcal P_\kappa^{d,l,D} \cap L_p(\R^d) \mapsto
\{f \mid_D: f \in \mathcal P_\kappa^{d,l,D}\} \cap L_p(D),
\end{equation*}
задаваемое равенством $ U f = f \mid_D, $ является линейным гомеоморфизмом.

Доказательство.

Линейность и сюръективность отображения $ U $ очевидны. Проверим его инъективность.
Если для $ f \in \mathcal P_\kappa^{d,l,D} $ значение $ U f =0, $
то ввиду представления (5.1.4) для $ x \in D $ выполняется равенство
$$
\sum_{\nu \in N_\kappa^{d,D}} f_\nu(x) \chi_{\kappa,\nu}^d(x) =0,
$$
а, следовательно, при $ \nu \in N_\kappa^{d,D} $ для $ x \in D \cap
Q_{\kappa,\nu}^d $ соблюдается равенство $ f_\nu(x) = f_\nu(x) \chi_{\kappa,\nu}^d(x) =0. $
Отсюда, учитывая включение $ f_\nu \in \mathcal P^{d,l} $ и открытость непустого
множества $ D \cap Q_{\kappa,\nu}^d, $ заключаем, что $ f_\nu(x) =0 $ для
$ x \in \R^d, \nu \in N_\kappa^{d,D}, $ т.е. $ f =0. $
Принимая во внимание сказанное, видим, что функционал $ \| U f \|_{L_p(D)} $
является нормой на конечномерном линейном пространстве
$ \mathcal P_\kappa^{d,l,D}. $ А поскольку любые нормы на конечномерном
линейном пространстве эквивалентны, то существует константа $ c_{58}(d,l,D,p,\kappa) >0 $
такая, что для $ f \in \mathcal P_\kappa^{d,l,D} $ справедливо неравенство
\begin{equation*} \tag{5.1.47}
c_{58} \| f \|_{L_p(\R^d)} \le \| f \mid_D\|_{L_p(D)} \le \| f \|_{L_p(\R^d)},
\end{equation*}
что завершает доказательство леммы. $ \square $
\bigskip

5.2. В этом пункте да\"eтся описание слабой асимптотики изучаемых в
настоящем параграфе поперечников, но сначала напомним их
определения (см. [10]).

Пусть $ C $ -- подмножество банахова пространства $X$ и
$ n \in \Z_+.$ Тогда $n$-поперечни\-ком по Колмогорову множества $C$
в пространстве $X$ называется величина
$$
d_n(C,X) = \inf_{M \in \mathcal M_n(X)} \sup_{x \in C} \inf_{y \in M} \|x -y\|_X,
$$
где $ \mathcal M_n(X) $ -- совокупность всех плоскостей $ M $ в $ X, $ у которых
$ \dim M \le n. $

$ n$-поперечником по Гельфанду множества $ C $ в пространстве $ X $ называется
величина
$$
d^n(C,X) = \inf_{M \in \mathcal M^n(X)} \sup_{x \in C \cap M} \|x\|_X,
$$
где $ \mathcal M^n(X)$ -- множество всех замкнутых линейных подпространств $ M $
в $ X, $ у которых $ \codim M \le n. $

Линейным $ n$-поперечником множества $ C $ в пространстве $ X $
называется величина
$$
\lambda_n(C,X) = \inf_{A \in \mathcal L_n(X)} \sup_{x \in C} \|x -Ax\|_X,
$$
где $ \mathcal L_n(X)$ -- множество всех непрерывных аффинных
отображений $ A: X \mapsto X, $ у которых $ \Rank A \le n. $

$ n$-поперечником по Александрову  множества $ C $ в пространстве $ X $
называется величина
$$
a_n(C,X) = \inf_{\phi \in \mathcal A_n(C,X)} \sup_{x \in C} \|x -\phi(x)\|_X,
$$
где $ \mathcal A_n(C,X)$ -- множество  всех непрерывных отображений
$ \phi: C \mapsto X, $ для которых существует компакт
$ K_{\phi} \subset X $ такой, что $ \Im \phi \subset K_{\phi} $ и
$ \dim K_{\phi} \le n. $

Введ\"eм ещ\"e следующую величину
$$
a^n(C,X) = \inf_{\phi \in \mathcal A^n(C,X)} \sup_{y \in \Im \phi}
\diam(\phi^{-1}(y))_X,
$$
где $ \mathcal A^n(C,X)$ -- совокупность всех непрерывных относительно метрики
пространства $ X $ отображений $ \phi:  C \mapsto K_{\phi}, $ для которых
$ K_{\phi}$ -- метрический компакт размерности не больше $ n, $ а
$ \diam(E)_X$ -- диаметр множества $ E $ в метрике пространства $ X. $

Наконец, энтропийным $n$-поперечником множества $ C $ в пространстве $ X $
будем называть величину
$$
\epsilon_n(C,X) = \inf\{\epsilon >0: \exists x^1,\ldots,x^{2^n} \in
X: C \subset \cup_{j=1}^{2^n}(x^j +\epsilon B(X))\}.
$$

С помощью утверждений 5.2.1 и 5.2.2, взятых из [11], доказывается лемма 5.2.3.

Предложение 5.2.1

Пусть $ U: X \mapsto Y$ -- непрерывное линейное отображение
банахова пространства $ X $ в банахово пространство $ Y $ и $ L \subset X$ --
замкнутое линейное подпространство, для которого сужение $ U \mid_L $ является
гомеоморфизмом $ L $ на $ U(L), $ и пусть $ C \subset L$ -- некоторое
множество. Тогда для $ p_n(\cdot,\cdot) = d_n(\cdot,\cdot), d^n(\cdot,\cdot),
\lambda_n(\cdot,\cdot),  a^n(\cdot,\cdot), \epsilon_n(\cdot,\cdot) $ при
$ n \in \Z_+ $ имеет место неравенство
\begin{equation*}  \tag{5.2.1}
p_n(U(C), Y) \le \| U \|_{\mathcal B(X,Y)} p_n(C,X).
\end{equation*}

Предложение 5.2.2

Пусть $ C $ -- подмножество банахова  пространства $ X, $ для которого
существует последовательность непрерывных линейных отображений
$ V_j: X \mapsto X_j $ пространства $ X $ в замкнутые линейные подпространства
$ X_j, j \in \Z_+, $ такая, что для каждого $ x \in C $ имеет место
представление $ x = \sum_{j =0}^{\infty} V_j x. $ Тогда для любых $ n,j_0 \in
\Z_+ $ и любого набора $ \{n_j \in \Z_+, j =0,\ldots,j_0\} $ таких, что
$ n \ge \sum_{j =0}^{j_0} n_j, $ для $ p_n(\cdot,\cdot) = d_n(\cdot,\cdot),
d^n(\cdot,\cdot), \lambda_n(\cdot,\cdot), a^n(\cdot,\cdot), \epsilon_n(\cdot,\cdot) $
выполняется неравенство
\begin{equation*} \tag{5.2.2}
p_n(C, X) \le \sum_{j =0}^{j_0} p_{n_j}(V_j(C), X_j)
+2\sum_{j =j_0 +1}^{\infty} d^0(V_j(C), X_j).
\end{equation*}

Лемма 5.2.3

Пусть выполнены условия леммы 5.1.3.  Тогда для
$ C = \D^\lambda (B((H_p^\alpha)^\prime(D))) $ и $ X = L_q(D) $
справедливы следующие утверждения:

1) существует константа $ c_1(d,\alpha,D,p,q,\lambda) >0 $ такая, что для
любых $ n,j_0,k \in \Z_+: k \ge K^0, $ и любого набора чисел $ \{n_j \in \Z_+,
j =1,\ldots,j_0\} $ таких, что
$ n \ge R_k^{d,l-\e-\lambda,D,\alpha} +\sum_{j =1}^{j_0} n_j, $ для
$ p_n(\cdot,\cdot) = d_n(\cdot,\cdot), d^n(\cdot,\cdot), \lambda_n(\cdot,\cdot),
a^n(\cdot,\cdot) $ имеет место неравенство
\begin{multline*} \tag{5.2.3}
p_n(C, X) \le c_1 \biggl(2^{-k(1 -(\alpha^{-1}, \lambda +(p^{-1} -q^{-1}) \e))}
\sum_{j =1}^{j_0} 2^{-j(1 -(\alpha^{-1}, \lambda +(p^{-1} -q^{-1}) \e))}
p_{n_j} \\
\times \bigl(B(l_p^{R_{k+j}^{d,l-\e-\lambda,D,\alpha}}),
l_q^{R_{k+j}^{d,l-\e-\lambda,D,\alpha}}\bigr) +2^{-(k+j_0)(1 -(\alpha^{-1},
\lambda +(p^{-1} -q^{-1})_+ \e))}\biggr);
\end{multline*}

2) существует константа $ c_2(d,\alpha,D,p,q,\lambda) >0 $ такая,
что для любых $ n,j_0 \in \Z_+ $ и любого набора $ \{n_j \in \Z_+,
j =0,\ldots,j_0\} $ таких, что $ n \ge \sum_{j =0}^{j_0} n_j, $ соблюдается
неравенство
\begin{multline*} \tag{5.2.4}
\epsilon_n(C,X) \le c_2 \biggl( \sum_{j =0}^{j_0}
2^{-j(1 -(\alpha^{-1}, \lambda +(p^{-1} -q^{-1}) \e))}
\epsilon_{n_j}\bigl(B(l_p^{R_{K^0 +j}^{d,l-\e-\lambda,D,\alpha}}\bigr),
l_q^{R_{K^0 +j}^{d,l-\e-\lambda,D,\alpha}}\bigr) \\
+2^{-j_0(1 -(\alpha^{-1}, \lambda +(p^{-1} -q^{-1})_+ \e))}\biggr).
\end{multline*}

Доказательство.

Прежде всего, в силу (5.1.26) при $ k \in \Z_+: k \ge K^0,,j_0 \in \Z_+ $ для
$ F \in C $ в $ X $ имеет место равенство
\begin{equation*}
F = (\mathcal U_{k,0}^{d,l-\e-\lambda,0,D,\alpha} F) \mid_D
+\sum_{j =1}^{j_0} (V_{k+j}^{d,l-\e-\lambda,D,\alpha} F) \mid_D
+\sum_{j =j_0 +1}^\infty (V_{k+j}^{d,l-\e-\lambda,D,\alpha} F) \mid_D.
\end{equation*}

Поэтому для $ n \in \Z_+ $ и $ \{n_j \in \Z_+, j =1,\ldots,j_0\} $
таких, что $ n \ge R_k^{d,l-\e-\lambda,D,\alpha}
+\sum_{j =1}^{j_0} n_j, $ ввиду (5.2.2) справедлива оценка
\begin{multline*} \tag{5.2.5}
p_n(C,X) \le p_{R_k^{d,l-\e-\lambda,D,\alpha}}
\bigl((\mathcal U_{k,0}^{d,l-\e-\lambda,0,D,\alpha}(C)) \mid_D,
(\mathcal P_k^{d,l-\e-\lambda,D,\alpha}) \mid_D \cap X\bigr) \\
+\sum_{j =1}^{j_0} p_{n_j} \bigl((V_{k+j}^{d,l-\e-\lambda,D,\alpha}(C)) \mid_D,
(\mathcal P_{k+j}^{d,l-\e-\lambda,D,\alpha}) \mid_D \cap L_q(D)\bigr) \\
+2\sum_{j =j_0 +1}^\infty d^0\bigl((V_{k+j}^{d,l-\e-\lambda,D,\alpha}(C)) \mid_D,
(\mathcal P_{k+j}^{d,l-\e-\lambda,D,\alpha}) \mid_D \cap L_q(D)\bigr).
\end{multline*}

Ясно, что
\begin{equation*} \tag{5.2.6}
p_{R_k^{d,l-\e-\lambda,D,\alpha}}
\bigl((\mathcal U_{k,0}^{d,l-\e-\lambda,0,D,\alpha}(C)) \mid_D,
(\mathcal P_k^{d,l-\e-\lambda,D,\alpha}) \mid_D \cap X \bigr) =0.
\end{equation*}

Далее, при $ j =1,\ldots,j_0, $ применяя сначала с учетом леммы 5.1.7
неравенства (5.2.1) и (5.1.47), а затем -- (5.2.1) и (5.1.46) (при $ \kappa =
\kappa(k +j,\alpha) $), наконец, используя (5.1.46) и (5.1.27) (при $ q = p $),
получаем
\begin{equation*} \tag{5.2.7}
\begin{split}
p_{n_j} \bigl((V_{k+j}^{d,l-\e-\lambda,D,\alpha}(C)) \mid_D,
(\mathcal P_{k+j}^{d,l-\e-\lambda,D,\alpha}) \mid_D \cap L_q(D) \bigr) \le \\
p_{n_j} \bigl(V_{k+j}^{d,l-\e-\lambda,D,\alpha}(C),
\mathcal P_{k+j}^{d,l-\e-\lambda,D,\alpha} \cap L_q(\R^d) \bigr) \\
\le c_3 2^{-(k+j)(1-(\alpha^{-1}, \lambda +(p^{-1} -q^{-1}) \e))}
p_{n_j} \bigl(B(l_p^{R_{k+j}^{d,l-\e-\lambda,D,\alpha}}),
l_q^{R_{k+j}^{d,l-\e-\lambda,D,\alpha}}\bigr).
\end{split}
\end{equation*}

Наконец, учитывая (5.1.47), (5.1.27), имеем
\begin{multline*} \tag{5.2.8}
\sum_{j =j_0 +1}^\infty d^0\bigl((V_{k+j}^{d,l-\e-\lambda,D,\alpha}(C)) \mid_D,
(\mathcal P_{k+j}^{d,l-\e-\lambda,D,\alpha}) \mid_D \cap L_q(D) \bigr) \le \\
\sum_{j =J_0 +1}^\infty
d^0 \bigl(V_{k+j}^{d,l-\e-\lambda,D,\alpha}(C),
\mathcal P_{k+j}^{d,l-\e-\lambda,D,\alpha} \cap L_q(\R^d) \bigr) \le \\
c_4 2^{-(k+j_0)(1 -(\alpha^{-1}, \lambda +(p^{-1} -q^{-1})_+ \e))}.
\end{multline*}

Соединяя (5.2.5) -- (5.2.8), приходим к (5.2.3).

Вывод (5.2.4) проводится аналогично. $ \square $

Как видно из доказательства лемм 4.2.4, 3.2.2, 3.2.3 в [6], имеет место

Лемма 5.2.4

Пусть $ d \in \N, \alpha \in  \R_+^d, \lambda \in \Z_+^d, 1 \le p < \infty,
1 \le q \le \infty $ удовлетворяют условию (1.3.62) и $ \theta \in \R:
1 \le \theta < \infty, \l = \l(\alpha). $ Тогда существует
константа $ c_5(d,\alpha,p,\theta,q,\lambda) >0 $ такая, что для $ C =
\D^\lambda ((B((B_{p,\theta}^\alpha)^{\l}(\R^d)) \cap C_0^\infty(I^d)) \mid_{I^d}),
X = L_q(I^d) $ и $ p_n(\cdot,\cdot) = d_n(\cdot,\cdot), d^n(\cdot,\cdot),
\lambda_n(\cdot,\cdot), a^n(\cdot,\cdot), \epsilon_n(\cdot,\cdot) $
для любых натуральных чисел $ n< N $
справедлсво неравенство
\begin{equation*} \tag{5.2.9}
p_n(C,X) \ge c_5 N^{-(1-(\alpha^{-1},  \lambda +(p^{-1} -q^{-1}) \e)) /(\alpha^{-1},\e)}
p_n(B(l_p^N), l_q^N).
\end{equation*}

С использованием леммы 5.2.4 доказывается

Лемма 5.2.5

Пусть $ d \in \N, \alpha \in \R_+^d, D $ -- ограниченная область $ \alpha $-типа
в $ \R^d, \lambda \in \Z_+^d, 1 \le p < \infty, 1 \le q \le \infty $
удовлетворяют условию (1.3.62) и $ \theta \in \R:
1 \le \theta < \infty. $ Тогда существует константа
$ c_6(d,\alpha,D,p,\theta,q,\lambda) >0 $ такая, что для $ C =
\D^\lambda (B((B_{p,\theta}^\alpha)^\prime(D))), X = L_q(D) $ и
$ p_n(\cdot,\cdot) = d_n(\cdot,\cdot), d^n(\cdot,\cdot), \lambda_n(\cdot,\cdot),
a^n(\cdot,\cdot), \epsilon_n(\cdot,\cdot) $
для любых натуральных чисел $ n < N $
выполняется неравенство (5.2.9) с константой $ c_6 $ вместо $ c_5. $

Доказательство.

Фиксируем точку $ x^0 \in \R^d $ и вектор $ \delta \in \R_+^d $
такие, что $ Q = (x^0 +\delta I^d) \subset D. $
Сначала заметим, что в условиях леммы при $ n,N \in \N: n < N, \l = \l(\alpha), $
благодаря (5.2.9) и (1.1.9), имеет место неравенство
\begin{multline*} \tag{5.2.10}
c_5 N^{-(1 -(\alpha^{-1}, \lambda +(p^{-1} -q^{-1}) \e)) /(\alpha^{-1},\e)}
p_n(B(l_p^N), l_q^N) \le \\
p_n(\D^\lambda ((B((B_{p,\theta}^\alpha)^{\l}(\R^d)) \cap C_0^\infty(I^d)) \mid_{I^d}), L_q(I^d)) \le \\
p_n(\D^\lambda ((B((B_{p,\theta}^\alpha)^\prime(\R^d)) \cap C_0^\infty(I^d)) \mid_{I^d}), L_q(I^d)).
\end{multline*}

Далее, принимая во внимание, что в силу (2.2.6) справедливо включение
\begin{multline*}
\D^\lambda ((B((B_{p,\theta}^\alpha)^\prime(\R^d)) \cap C_0^\infty(I^d)) \mid_{I^d}) =
\{ \D^\lambda (f \mid_{I^d}): f \in
(B((B_{p,\theta}^\alpha)^\prime(\R^d)) \cap C_0^\infty(I^d)) \} = \\
h_{\delta,x^0} (\{ h_{\delta,x^0}^{-1}(\D^\lambda (f \mid_{I^d})): f \in
(B((B_{p,\theta}^\alpha)^\prime(\R^d)) \cap C_0^\infty(I^d)) \}) = \\
h_{\delta,x^0} (\{ \delta^{\lambda} \D^\lambda (h_{\delta,x^0}^{-1}(f \mid_{I^d})):
f \in (B((B_{p,\theta}^\alpha)^\prime(\R^d)) \cap C_0^\infty(I^d)) \}) = \\
\delta^\lambda h_{\delta,x^0} (\{ \D^\lambda ((h_{\delta,x^0}^{-1} f) \mid_Q):
f \in (B((B_{p,\theta}^\alpha)^\prime(\R^d)) \cap C_0^\infty(I^d)) \}) = \\
\delta^\lambda h_{\delta,x^0} (\{ \D^\lambda (F \mid_Q):
F \in h_{\delta,x^0}^{-1}(B((B_{p,\theta}^\alpha)^\prime(\R^d)) \cap C_0^\infty(I^d)) \}) = \\
\delta^\lambda h_{\delta,x^0} (\{ \D^\lambda (F \mid_Q):
F \in h_{\delta,x^0}^{-1}(B((B_{p,\theta}^\alpha)^\prime(\R^d))) \cap C_0^\infty(Q) \}) \subset \\
\delta^\lambda h_{\delta,x^0} (\{ \D^\lambda (F \mid_Q):
F \in (c_7 B((B_{p,\theta}^\alpha)^\prime(\R^d))) \cap C_0^\infty(Q) \}) = \\
c_8 h_{\delta,x^0} (\{ \D^\lambda (F \mid_Q):
F \in B((B_{p,\theta}^\alpha)^\prime(\R^d)) \cap C_0^\infty(Q) \}),
\end{multline*}
благодаря (5.2.1), (2.2.2), имеем
\begin{multline*} \tag{5.2.11}
p_n(\D^\lambda ((B((B_{p,\theta}^\alpha)^\prime(\R^d)) \cap C_0^\infty(I^d)) \mid_{I^d}), L_q(I^d)) \le \\
p_n(c_8 h_{\delta,x^0} (\{ \D^\lambda (F \mid_Q):
F \in B((B_{p,\theta}^\alpha)^\prime(\R^d)) \cap C_0^\infty(Q) \}),L_q(I^d)) \le \\
c_8 \delta^{-q^{-1} \e} p_n(\{ (\D^\lambda F) \mid_Q:
F \in B((B_{p,\theta}^\alpha)^\prime(\R^d)) \cap C_0^\infty(Q) \},L_q(Q)).
\end{multline*}

Теперь, определяя в $ L_q(D) $ замкнутое подпространство
$$
L_q^0(Q,D) = \{f \in L_q(D): f = \chi_Q f\},
$$
заметим, что
$$
\{ \D^\lambda (F \mid_D):
F \in B((B_{p,\theta}^\alpha)^\prime(\R^d)) \cap C_0^\infty(Q) \} \subset
L_q^0(Q,D),
$$
а для оператора $ U: L_q(D) \ni f \mapsto U f = f \mid_Q \in L_q(Q), $
его сужение $ U \mid_{L_q^0(Q,D)} $ является изометрическим изоморфизмом
пространства $ L_q^0(Q,D) \cap L_q(D) $ на $ L_q(Q). $
С учетом сказанного применяя (5.2.1), получаем
\begin{multline*} \tag{5.2.12}
p_n(\{ (\D^\lambda F) \mid_Q:
F \in B((B_{p,\theta}^\alpha)^\prime(\R^d)) \cap C_0^\infty(Q) \},L_q(Q)) = \\
p_n(\{ (\D^\lambda (F \mid_D)) \mid_Q:
F \in B((B_{p,\theta}^\alpha)^\prime(\R^d)) \cap C_0^\infty(Q) \},L_q(Q)) = \\
p_n(U(\{ \D^\lambda (F \mid_D):
F \in B((B_{p,\theta}^\alpha)^\prime(\R^d)) \cap C_0^\infty(Q) \}),L_q(Q)) \le \\
p_n(\{ \D^\lambda (F \mid_D):
F \in B((B_{p,\theta}^\alpha)^\prime(\R^d)) \cap C_0^\infty(Q) \},L_q(D)) \le \\
p_n(\{ \D^\lambda f: f \in B((B_{p,\theta}^\alpha)^\prime(D)),L_q(D)).
\end{multline*}

Соединяя (5.2.12), (5.2.11), (5.2.10), видим, что для $ C $ и $ X $ из
формулировки леммы соблюдается (5.2.9) с константой $ c_6 $ вместо $ c_5. \square $

Опираясь на соотношения (5.1.8), (5.2.3), (5.2.4) и (5.2.9) с константой
$ c_6 $ вместо $ c_5, $ подобно тому, как это сделано в [6], устанавливается
слабая асимптотика колмогоровского, гельфандовского, линейного,
александровского и энтропийного поперечников классов
$ \D^\lambda (B((H_p^\alpha)^\prime(D))) $ и
$ \D^\lambda (B((B_{p,\theta}^\alpha)^\prime(D))) $ в пространстве $ L_q(D), $
заданных в ограниченной области $ D \alpha $-типа.

Теорема 5.2.6

Пусть $ d \in \N, \alpha \in \R_+^d, \lambda \in \Z_+^d, 1 \le p < \infty,
1 \le q \le \infty $ удовлетворяют условию (1.3.62) и $ \theta \in \R:
1 \le \theta < \infty, $ а $ D \subset \R^d $ -- ограниченная область
$ \alpha $-типа. Тогда для $ C = \D^\lambda(B((H_p^\alpha)^\prime(D))),
\D^\lambda(B((B_{p,\theta}^\alpha)^\prime(D))) $ и $ X = L_q(D) $
справедливо соотношение
\begin{equation*}
d_n(C,X) \asymp 
\begin{cases} n^{-(1-(\alpha^{-1}, \lambda))
/(\alpha^{-1},\e)
+(p^{-1} -q^{-1})_+}, \text{ при } q \le p \text{ или } p < q \le 2; \\
n^{-(1-(\alpha^{-1}, \lambda)) /(\alpha^{-1},\e) +(p^{-1} -2^{-1})_+},\parbox[t]{.4\textwidth} { при $q > \max(2,p),\\
1 -(\alpha^{-1},\lambda) -(\alpha^{-1},\e) (1/p +(1/2 -1/p)_+)
>0$.} 
\end{cases}
\end{equation*}

Теорема 5.2.7

Пусть $ d, \alpha, \lambda,p,q, \theta, D, $ а также $ C $ и $ X $ имеют тот же
смысл, что и в теореме 5.2.6. Тогда имеет место
соотношение
\begin{equation*}
d^n(C,X) \asymp 
\begin{cases} n^{-(1-(\alpha^{-1},  \lambda))
/(\alpha^{-1},\e) +(p^{-1} -q^{-1})_+},
\parbox[t]{.4\textwidth} { при $(q \le p$ или $ 2 \le p < q)$
и соблюдении условия (1.3.62);} \\
n^{-(1-(\alpha^{-1}, \lambda)) /(\alpha^{-1},\e) +(2^{-1} -q^{-1})_+},\
\parbox[t]{.4\textwidth} { при $p < \min(2,q),
1 -(\alpha^{-1},\lambda) -(\alpha^{-1},\e) (1/q^\prime +(1/2 -1/q^\prime)_+) >0.$}
\end{cases}
\end{equation*}

Теорема 5.2.8

В условиях теоремы 5.2.6 верно соотношение
\begin{equation*}
\lambda_n(C,X) \asymp 
\begin{cases} n^{-(1-(\alpha^{-1}, \lambda))
/(\alpha^{-1},\e) +(p^{-1} -q^{-1})_+},  \parbox[t]{.4\textwidth} { при $(q  \le p$
или $p < q \le 2$ или
$2 \le p < q)$;} \\
n^{-(1-(\alpha^{-1}, \lambda)) /(\alpha^{-1},\e) +1/2 -1/q +(1/p
+1/q -1)_+}, \parbox[t]{.4\textwidth} { при $p < 2 < q,
1-(\alpha^{-1},\lambda) -(\alpha^{-1},\e) (1+(1/p^\prime -1/q)_+) >0.$}
\end{cases}
\end{equation*}

Теорема 5.2.9

В условиях теоремы 5.2.6 соблюдается соотношение
\begin{equation*}
a_n(C,X) \asymp n^{-(1-(\alpha^{-1}, \lambda)) /(\alpha^{-1},\e)}.
\end{equation*}

Теорема 5.2.10

В условиях теоремы 5.2.6 при $ 1-(\alpha^{-1},\lambda)
-(\alpha^{-1},\e) >0 $ выполняется соотношение
\begin{equation*}
\epsilon_n(C,X) \asymp n^{-(1-(\alpha^{-1}, \lambda)) /(\alpha^{-1},\e)}.
\end{equation*}


\bigskip
\newpage
\begin{thebibliography}{1}
\bibitem{1.}
Бесов О. В., Ильин В. П.
Естественное расширение класса областей в теоремах вложения //
Матем. сб., 75(117):4 (1968), 483--495.
\bibitem{2.}
Буренков В. И., Файн Б. Л.
О продолжении функций из анизотропных пространств с сохранением класса //
Тр. МИАН СССР, 150, 1979, 52--66.
\bibitem{3.}
Шварцман П. А.
Теоремы продолжения с сохранением локально-аппроксимационных свойств функций в неизотропном случае //
Зап. научн. сем. ЛОМИ, 113 (1981), 247--252
\bibitem{4.}
Файн Б. Л.
О продолжении функций из анизотропных пространств Соболева //
Тр. МИАН СССР, 170 (1984), 248--272.
\bibitem{5.}
Мамедов Г. А.
Продолжение анизотропного класса Г\"eльдера //
Изв. вузов. Матем., 1985, №6, 36--40.
\bibitem{6.}
Кудрявцев С. Н.
Приближение производных функций конечной гладкости из неизотропных классов //
Изв. РАН. Сер. матем. 68:1 (2004), 79--122.
\bibitem{7.}
Никольский С. М.
Приближение функций многих переменных и теоремы вложения // М.: Наука. 1977.
\bibitem{8.}
Кудрявцев С. Н.
Обобщенные ряды Хаара и их применение //
Analysis Mathematica, 37:2 (2011), 103--150.
\bibitem{9.}
Чуи Ч. К.
Введение в вейвлеты // М.: Мир, 2001.
\bibitem{10.}
Тихомиров В. М.
Некоторые вопросы теории приближения // М.: Изд-во МГУ. 1976.
\bibitem{11.}
Кудрявцев С. Н.
Поперечники классов гладких функций //
Изв. РАН. Сер. матем. 59:4 (1995), 81--104.
\end{thebibliography}
\end{document}